\documentclass[a4paper,reqno,oneside]{amsart}
\usepackage{amsfonts}
\usepackage{geometry}
\usepackage{tkz-graph,subcaption}

\usetikzlibrary{calc}
\usepackage{wrapfig}
\usepackage{pgf}
\usetikzlibrary{decorations.pathmorphing}
\usetikzlibrary{shapes,arrows,automata,decorations.pathreplacing,angles,quotes,calc}
\usetikzlibrary{shapes.geometric}
\tikzset{snake it/.style={decorate, decoration=snake}}
\tikzset{
	partial ellipse/.style args={#1:#2:#3}{
		insert path={+ (#1:#3) arc (#1:#2:#3)}
	}
}
\usepackage{tkz-euclide}
\usepackage{amsmath}
\usepackage{url}
\usepackage[colorlinks=true,citecolor={black!40!green},urlcolor={black!40!blue},linkcolor={blue}]{hyperref}
\usepackage{amssymb}
\usepackage{amsthm}
\usepackage{enumerate}
\usepackage{enumitem}
\usepackage{color}

\theoremstyle{plain}
\newtheorem{thm}{Theorem}
\newtheorem{lem}{Lemma}
\newtheorem{defi}{Definition}
\newtheorem{cor}{Corollary}
\newtheorem{pro}{Proposition}
\newtheorem{rem}{Remark}
\newtheorem{rem*}{Remark}

\newtheorem*{def*}{Definition}
\newtheorem*{prob*}{Direct Monodromy Problem}
\newtheorem*{probinv*}{Inverse Monodromy Problem}
\newtheorem*{probinvred*}{Reduced Inverse Monodromy Problem}
\newtheorem*{lem*}{Lemma}
\newtheorem*{cor*}{\textsc{Corollary}}
\newtheorem*{con*}{\textsc{Conjecture}}
\newtheorem*{thm*}{Theorem}
\newcommand{\e}{\varepsilon}
\newcommand{\C}{\mathbb{C}}
\newcommand{\Cb}{\overline{\mathbb{C}}}
\newcommand{\uo}{\underline{0}}
\newcommand{\Zq}{\mathbb{Z}_4}
\newcommand{\bb}[1]{\mathbb{ #1 }}
\newcommand{\mc}[1]{\mathcal{ #1 }}
\newcommand{\av}[1]{\left| #1 \right|}

\newcommand{\om}{\omega}
\newcommand{\la}{\lambda}

\newcommand{\beq}{\begin{equation}}
\newcommand{\eeq}{\end{equation}}

\newcommand{\norm}[1]{\left| #1 \right|}
\DeclareMathOperator*{\res}{Res}

\title{Roots of generalised Hermite polynomials when both parameters are large}
\author{Davide Masoero and Pieter Roffelsen}

\address{Grupo de F\'isica Matem\'atica e Departamento da Matem\'atica da Universidade de Lisboa,
Campo Grande Edif\'icio C6,
1749-016 Lisboa, Portugal.}
\email{dmasoero@gmail.com}

\address{SISSA, via Bonomea 265, 34136 Trieste, Italy.}
\email{proffels@sissa.it}

\date{}
\begin{document}
\pagestyle{plain}

\begin{abstract}
	We study the roots of the generalised Hermite polynomials $H_{m,n}$ when both $m$ and $n$ are large.
	We prove that the roots, when appropriately rescaled, densely fill a bounded quadrilateral region, called the elliptic region, and
	organise themselves on a deformed rectangular grid, as was numerically observed by Clarkson.
	We describe the elliptic region and the deformed grid in terms
	of elliptic integrals and their degenerations.
	
\end{abstract}

 \maketitle

{\hypersetup{linkcolor=blue,linktocpage=true}
	\tableofcontents}

\section{Introduction}

For $m,n\in\mathbb{N}$, the generalised Hermite polynomial $H_{m,n}$ is the polynomial
of degree $m\times n$ defined by the determinantal
formula
\begin{equation*}
H_{m,n}(z)=\gamma_{m,n}\begin{vmatrix}
H_m(z)  &  H_{m+1}(z)  &  \ldots  & H_{m+n-1}(z)\\
H_m^{(1)}(z)  &  H_{m+1}^{(1)}(z)  &  \ldots  & H_{m+n-1}^{(1)}(z)\\
\vdots	  &\vdots & \ddots   & \vdots\\
H_m^{(n-1)}(z)  &  H_{m+1}^{(n-1)}(z)  &  \ldots  & H_{m+n-1}^{(n-1)}(z)\\
\end{vmatrix}
\end{equation*}
where $H_k^{(l)}(z)$ denotes the $l$-th derivative of the $k$-th Hermite polynomial
$$H_k(z)=(-1)^ke^{z^2}\frac{d^k}{d z^k}\left[e^{-z^2}\right]$$ and $\gamma_{m,n}\in\mathbb{C}^*$ is a,
for the purpose of this paper, irrelevant constant multiplier.

The goal of the present paper is to study
the asymptotic distribution of the roots of $H_{m,n}$ as $m+n \to \infty$. We call the asymptotics unrestricted since,
in contrast to our previous paper \cite{maro18},  we do not require
any of the two parameters to remain bounded.

The most striking property of the generalised Hermite polynomials is that
they yield families of rational solutions to the fourth Painlev\'e equation
\begin{equation}\label{eq:piv}
\om''=\displaystyle\frac{1}{2\om}\left(\om'\right)^2+ \frac{3}{2} \om^3 + 4 z\om^2
+2(z^2+1-2\theta_\infty)\om-\frac{8\theta_0^2}{\om}
,\quad \theta:=(\theta_0,\theta_\infty) \in \mathbb{C}^2.
\end{equation}
For example, the functions $\omega_{m,n}^{(\text{I})}=\frac{d}{dz}\log\frac{H_{m+1,n}}{H_{m,n}}$
solve the above equation with parameters $\theta_0=\tfrac{1}{2}n,  \theta_\infty=m+\tfrac{1}{2}n+1$ \cite{noumiyamada}. This
establishes an explicit relation among poles of rational solutions of Painleve IV and
roots of generalised Hermite polynomials; hence the problem of our interest can be restated as the study of the
asymptotic distribution of singularities of the Hermite-family of rational solutions of the Painleve IV equation \cite{maro18}.

We mentioned in our previous paper \cite{maro18} that the asymptotic analysis of
rational solutions of Painlev\'e equations (and more generally special solutions)
has recently been the object of intense study, see e.g. \cite{costin12,piwkb}. This is even more the case at the time of writing since
in the meanwhile many new results have appeared in the literature \cite{bothner18,buckingham18,vanass18}.
There is a clear reason for this: most often when a problem in
applied or pure mathematics is solved by means of a solution of a Painlev\'e equation (see e.g. \cite{dubrovin2009}),
this is singled-out by some special (i.e. non generic) asymptotic expansion,
which reflects itself in a special associated Riemann-Hilbert problem  amenable to a thorough analysis
to a degree not attainable for generic solutions \footnote{On the other side, the structure of  the Painlev\'e equations can
be properly understood only when studying the general solution, see \cite{fokas,guzzetti12,iorgov15}.}.

The above is clearly the case for generalised Hermite polynomials. In fact in \cite{maro18} we showed that
in the case of the generalised Hermite polynomials
the isomonodromic deformation method for Painleve IV simplifies dramatically. To be more precise,
we proved the following theorem, which characterises the roots of $H_{m,n}$ by means of an inverse monodromy problem
for a specific class of anharmonic oscillators known as biconfluent Heun equations.
\begin{thm}[Theorem 2.2 in \cite{maro18}]\label{thm:inversemonodromy}
For $m,n\in\mathbb{N}$, the point  $a \in \bb{C}$ is a root of $H_{m,n}$ if and only if there exists $b \in \bb{C}$ such that
the anharmonic oscillator
	\begin{align}\label{eq:anharmonic}
	&\psi''(\lambda)=(\la^2+2 a \la + a^2 - (2m+n) -\frac{b}{\la}+\frac{n^2-1}{4 \lambda^2})\psi(\lambda),
	\end{align}
	satisfies the following two properties:
	\begin{enumerate}
		\item\textbf{Apparent Singularity Condition.} The resonant singularity at $\la=0$ is
		apparent or equivalently the monodromy around the singularity is scalar. In a formula,
		\begin{equation*}
		 \psi(e^{2 \pi i} \la)=(-1)^{n+1} \psi(\la) \, , \quad \forall \psi \mbox{ solution  of \eqref{eq:anharmonic}.}
		\end{equation*}
		\item\textbf{Quantisation Condition.} There exists a non-zero solution of \eqref{eq:anharmonic} which solves
		the following boundary value problem
		\begin{equation*}
		 \lim_{\la \to +\infty} \psi(\la)=\lim_{\la \to 0^+} \psi(\la)=0 \; .
		\end{equation*}
	\end{enumerate}
\end{thm}
Studying the asymptotic solution of the above inverse problem, we obtain our main result, which we name the \textit{bulk asymptotics}:
\begin{enumerate}
 \item We determine
the region of the complex plane which asymptotically gets filled densely by the roots. This is a quadrilateral domain that we call the
elliptic region following Buckingham \cite{buckingham18}, who had already described it.
\item We describe the bulk asymptotics of the roots, that is, we obtain an asymptotic description of the
roots uniformly on any fixed compact subset of the interior of the elliptic region. In particular we show that the roots
asymptotically organise themselves on a
deformed grid, thus confirming Clarkson's numerical findings \cite{clarkson2003piv,clarksonpiv}.
\end{enumerate}
The paper is organised as follows. In Section \ref{sec:results} we state our main results, and
we announce a conjectural formula concerning the critical asymptotics, that is, the asymptotics of 
roots approaching the four corners of the elliptic region.
In Section \ref{section:elliptic} we study of the elliptic region.
Section \ref{sec:WKB} is devoted to complex WKB method for equation \eqref{eq:anharmonic}.
In Section \ref{section:bulkasymptotic} we prove the main theorems
concerning the bulk asymptotics. Finally in the appendix we collect some formulae from the theory
of elliptic integrals which
are used in the main body.

\begin{figure}[ht]
	\centering
		\includegraphics[width=10cm]{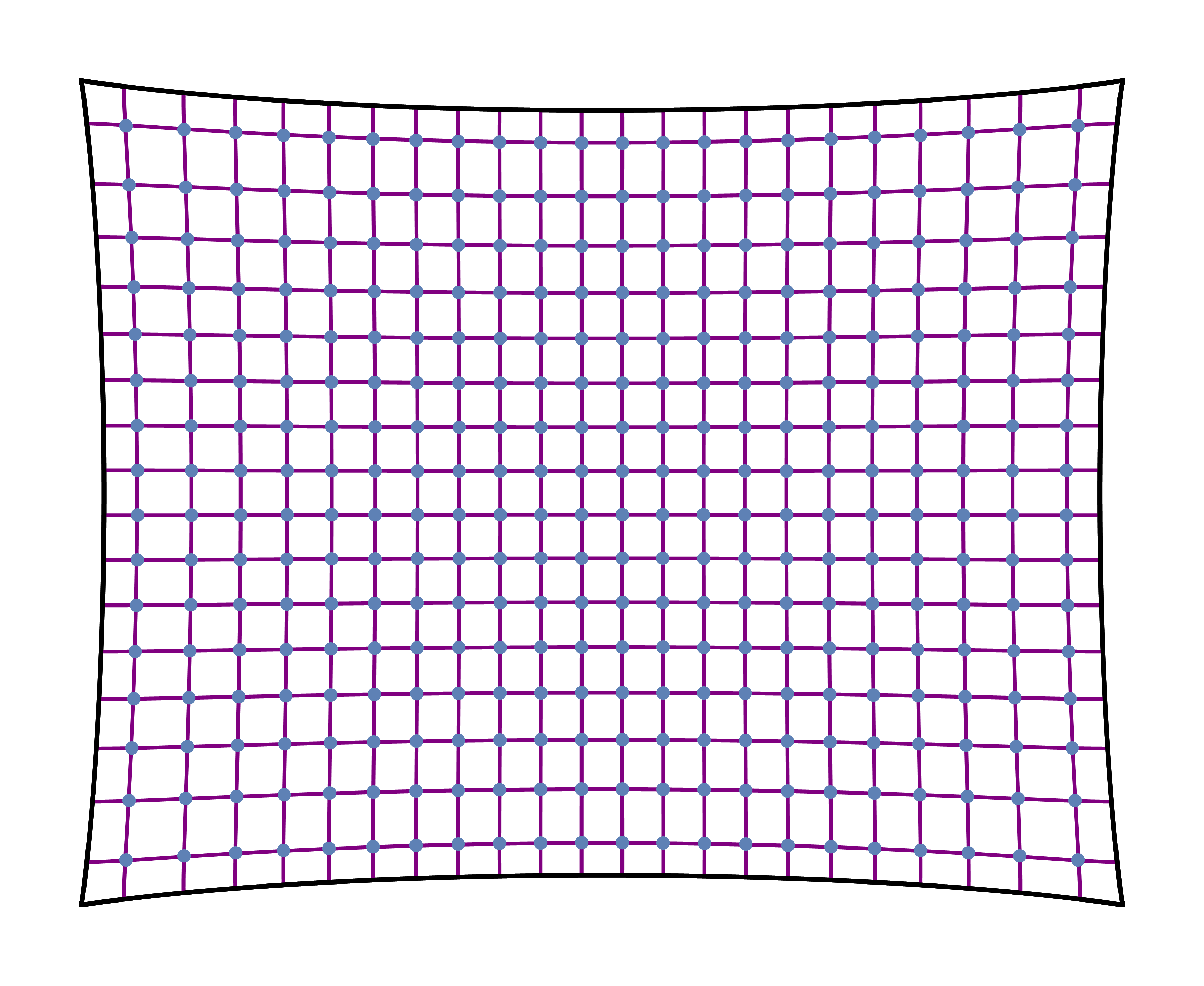} 
	\caption{Asymptotically the roots lie on the vertices of a deformed grid. The asymptotic prediction
	is stunningly precise even for moderate $m,n$. In the picture, the elliptic region with in purple the deformed grid, and in blue
	the roots of $H_{m,n}(z)$ superimposed
		with $(m,n)=(22,16)$.
		} \label{fig:intro}
\end{figure}

Before we begin our paper, in Figure \ref{fig:intro} we show the reader a pictorial description
of our results, which is worth a thousand lemmas.

\subsection*{Note on the Literature}
The problem of describing the distribution of roots of generalised Hermite polynomials when both parameters are large
were raised by Peter Clarkson in \cite{clarksonpiv}. The slightly more general problem of describing the asymptotics
of generalised Hermite polynomials was the object of much recent interest, see e.g. \cite{maro18,buckingham18,vanass18}.
In our previous paper \cite{maro18} we dealt with the asymptotic
location of roots in case when one of the two parameters stays bounded.
Buckingham studied the asymptotics of generalised Hermite polynomials when both parameters are large in the complement of
the elliptic region, that is, in the region where no roots are present \cite{buckingham18}.
We present here the first correct (a previous attempt is carefully analysed in \cite{buckingham18,vanass18}, and shown to be erroneous)
asymptotic description of roots of generalised Hermite polynomials when both parameters are large.

\subsection*{Acknowledgements}
D.M. is partially supported by the FCT Project PTDC/MAT-PUR/ 30234/2017 `Irregular
connections on algebraic curves and Quantum Field Theory' and by the
FCT Investigator grant IF/00069/2015 `A mathematical framework for the ODE/IM correspondence'. P.R. acknowledges the support of the H2020-MSCA-RISE-2017 PROJECT No. 778010
IPADEGAN.

\section{Results}\label{sec:results}
We introduce a large parameter $E=2m+n$, and the scaled
parameters $\alpha, \beta, \nu$,
\begin{equation}\label{eq:scaledparameters}
 E=2m+n,\quad \alpha=E^{-\frac12}a,\quad \beta=E^{-\frac32}b, \quad \nu=\frac{n}{E}.
\end{equation}
The asymmetric choice of $E$ is motivated by the particular asymmetric dependence of the oscillator \eqref{eq:anharmonic} on $m,n$.
Without loss in generality, since $H_{n,m}(z)=H_{m,n}(iz)$,  we assume that $m\geq n$, so $\nu\in[0,\frac13]$.
We further assume that $\nu \in (0,\frac13]$. The case $\nu=0$ requires a different approach and was dealt with in our previous paper
\cite{maro18}.

By applying the scalings \eqref{eq:scaledparameters} to the oscillator \eqref{eq:anharmonic} in Theorem \ref{thm:inversemonodromy}, and also scaling the independent variable $\la \to E^{\frac12}\la$, the problem laying before us is the rigorous study of
the no-logarithm and quantisation condition on the anharmonic oscillator
\begin{align}\label{eq:scaled}
&\psi''(\la)=\left(E^2V(\lambda;\alpha,\beta,\nu)-\frac{1}{4 \la^2} \right)\psi(\la),\\ \label{eq:scaledpotential}
&V(\lambda;\alpha,\beta,\nu)=\lambda^2+2\alpha \lambda+\alpha^2-1-\beta\lambda^{-1}+\frac{\nu^2}{4}\lambda^{-2},
\end{align}
in the $E\rightarrow +\infty$ limit.
Our approach to the inverse monodromy problem is based on
the complex WKB analysis of equation \eqref{eq:scaled}. The latter builds on the approximation of solutions by means
of the (multivalued) WKB
functions
\begin{equation}\label{eq:wkb_functions}
\psi=V^{-\frac14}e^{\pm E \int^\la \sqrt{V(\mu)}d\mu},
\end{equation}
 where $V$ is the function appearing in \eqref{eq:scaledpotential}.
Notice that in the above formula we have neglected the term $-\frac{1}{4 \la^2}$. This is called the Langer modification and it is
necessary to obtain a correct approximation when regular singularities are present; we will discuss it further when proving our results.

\subsection{The Elliptic Region}
In this subsection we introduce the compact region $K_a$ in the complex $\alpha$-plane which
asymptotically gets filled with the roots of $H_{m,n}(E^\frac{1}{2}\alpha)$ as $E\rightarrow \infty$.
This region is defined as the projection onto the $\alpha$-plane of a compact set
$K\subseteq \{(\alpha,\beta)\in\mathbb{C}^2\}$ which asymptotically gets filled with the solutions $(\alpha,\beta)$
of the no-logarithm and quantisation condition on \eqref{eq:scaled}.

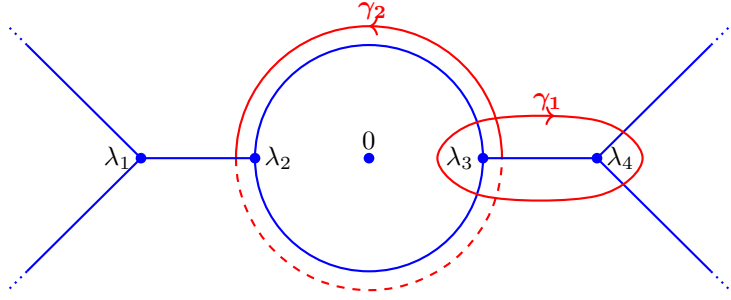
\begin{figure}[t]
	\centering
	\begin{tikzpicture}[scale=1.0]
	
	\draw[red,thick,decoration={markings, mark=at position 0.5 with {\arrow{>}}},postaction={decorate}] (1.75,0) arc (0:180:1.75cm);
	\draw[red,thick,dashed] (-1.75,0) arc (-180:0:1.75cm);

	\tikzstyle{star}  = [circle, minimum width=4pt, fill, inner sep=0pt];
	\node[blue,star]     at (-3,0){};
	\node[left]     at (-3,0) {$\lambda_1$};
	
	\node[blue,star]     at (-1.5,0) {};
	\node[right]     at (-1.5,0) {$\lambda_2$};
	
	\node[blue,star]     at (0,0) {};
	\node[above]     at (0,0) {0};
	
	\node[blue,star]     at (1.5,0) {};
	\node[left]     at (1.5,0) {$\lambda_3$};
	
	\node[blue,star]     at (3,0) {};
	\node[right]     at (3,0) {$\lambda_4$};
	
	\draw[blue,thick] (-3,0) -- (-1.5,0);
	\draw[blue,thick] (1.5,0) -- (3,0);
	\draw[blue,thick] (0,0) circle (1.5cm);
	
	\draw[blue,thick] (-3,0) -- (-4.5,1.5);
	\draw[blue,thick,dotted] (-4.5,1.5) -- (-4.75,1.75);
	
	\draw[blue,thick] (3,0) -- (4.5,1.5);
	\draw[blue,thick,dotted] (4.5,1.5) -- (4.75,1.75);
	
	\draw[blue,thick] (-3,0) -- (-4.5,-1.5);
	\draw[blue,thick,dotted] (-4.5,-1.5) -- (-4.75,-1.75);
	
	\draw[blue,thick] (3,0) -- (4.5,-1.5);
	\draw[blue,thick,dotted] (4.5,-1.5) -- (4.75,-1.75);
	
	\draw [red,thick,decoration={markings, mark=at position 0.15 with {\arrow{>}}},
	postaction={decorate}] plot [smooth cycle,tension=0.6] coordinates {(0.9,0) (1.5,0.5)  (3, 0.5) (3.6,0) (3,-0.5)  (1.5,-0.5)  };	
	
	\node[red]     at (2.35,0.75) {$\boldsymbol{\gamma_1}$};
	\node[red]     at (0.05,1.95) {$\boldsymbol{\gamma_2}$};

	\end{tikzpicture}
	\caption{Topological representation of Stokes complex in blue at $(\alpha,\beta)=(0,0)$ and more generally for $(\alpha,\beta)\in R$, where $\lambda_k$, $1\leq k\leq 4$, are the zeros of $V(\lambda;\alpha,\beta)$. Furthermore in red are depicted the two cycles $\gamma_1$ and $\gamma_2$
		used in equations \eqref{eq:defs1s2}.} \label{fig:approximate}
\end{figure}

In order to define $K$, we briefly introduce the Stokes complex associated with the potential $V(\lambda)$.
In the WKB analysis of \eqref{eq:scaled} the Stokes lines of $V(\lambda)$ play an important role. These are defined as
the levels sets $\Re\int_{\lambda^*}^\lambda \sqrt{V(\la)}d\la=0$ in $\mathbb{P}^1$, starting from any zero $\lambda^*$ of the potential $V(\lambda)$.
The Stokes complex $\mathcal{C}=\mathcal{C}(\alpha,\beta)\subseteq \mathbb{P}^1$ of $V(\lambda)$ is defined as the union of all
its Stokes lines and zeros. 

For example, if we set $(\alpha,\beta)=(0,0)$, then the potential simplifies to
\begin{equation*}
V(\lambda)=\lambda^2-1+\frac{\nu^2}{4}\lambda^{-2},
\end{equation*}
which has four real zeros
\begin{align*}
&\lambda_1=-\sqrt{\tfrac{1}{2}(1+\sqrt{1-\nu^2})},\quad \lambda_2=-\sqrt{\tfrac{1}{2}(1-\sqrt{1-\nu^2})},\\ &\lambda_3=+\sqrt{\tfrac{1}{2}(1-\sqrt{1-\nu^2})},\quad \lambda_4=+\sqrt{\tfrac{1}{2}(1+\sqrt{1-\nu^2})},
\end{align*}
and its Stokes complex is depicted topologically\footnote{For an accurate numerical plot of the Stokes complex when $(\alpha,\beta)=(0,0)$ with $\nu=\frac{1}{3}$, see Figure \ref{fig:stokesrealplot}.} in Figure \ref{fig:approximate}.

\begin{defi}\label{defi:RK}
	We define $R$ as the set of $(\alpha,\beta)\in\mathbb{C}^2$ such that the Stokes
	complex $\mathcal{C}(\alpha,\beta)$ is homeomorphic to the Stokes complex at $(0,0)$,
	 denote its closure by $K=\overline{R}$ and define $K_a$ as the projection of $K$ onto
	 the $\alpha$-plane.
	 We call $K_a$ the elliptic region.
\end{defi}
We note that $R$, $K$ and $K_a$ each depend parametrically on $\nu$.
The region $R$ can locally be described within $\{(\alpha,\beta)\in\mathbb{C}^2\}$ by the vanishing of the real parts of integrals $\int_{\gamma}{y d\la}$ along two homologically independent cycles $\gamma_{1,2}$ on the Riemann surface $y^2=V(\lambda)$. In particular $R$ is a smooth real surface. The region $K=\overline{R}$ is a smooth real surface with boundaries and corners and the same holds true for $K_a$. We proceed with describing the elliptic region $K_a\subseteq \{\alpha\in\mathbb{C}\}$ in detail. 

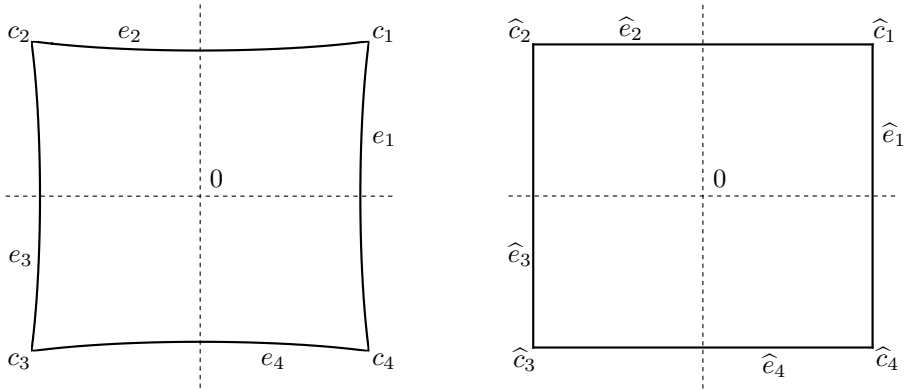
\begin{figure}
	\begin{tikzpicture}
	\node[anchor=south west,inner sep=0] (imageKa) at (0,0) {\includegraphics[width=0.35\textwidth]{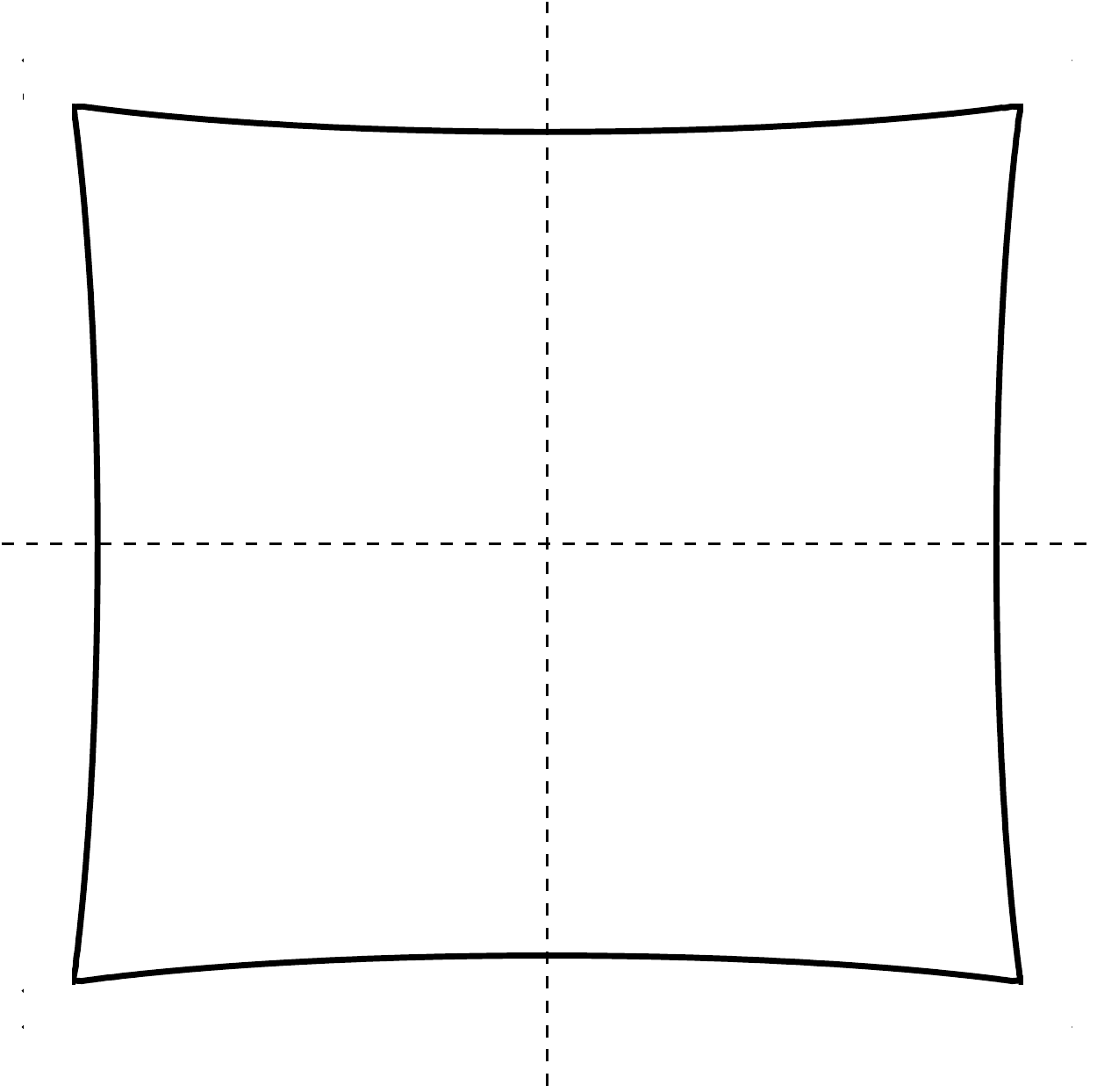}};
	\node[anchor=south west,inner sep=0] (imageQ) at (0.45\textwidth,0) {\includegraphics[width=0.35\textwidth]{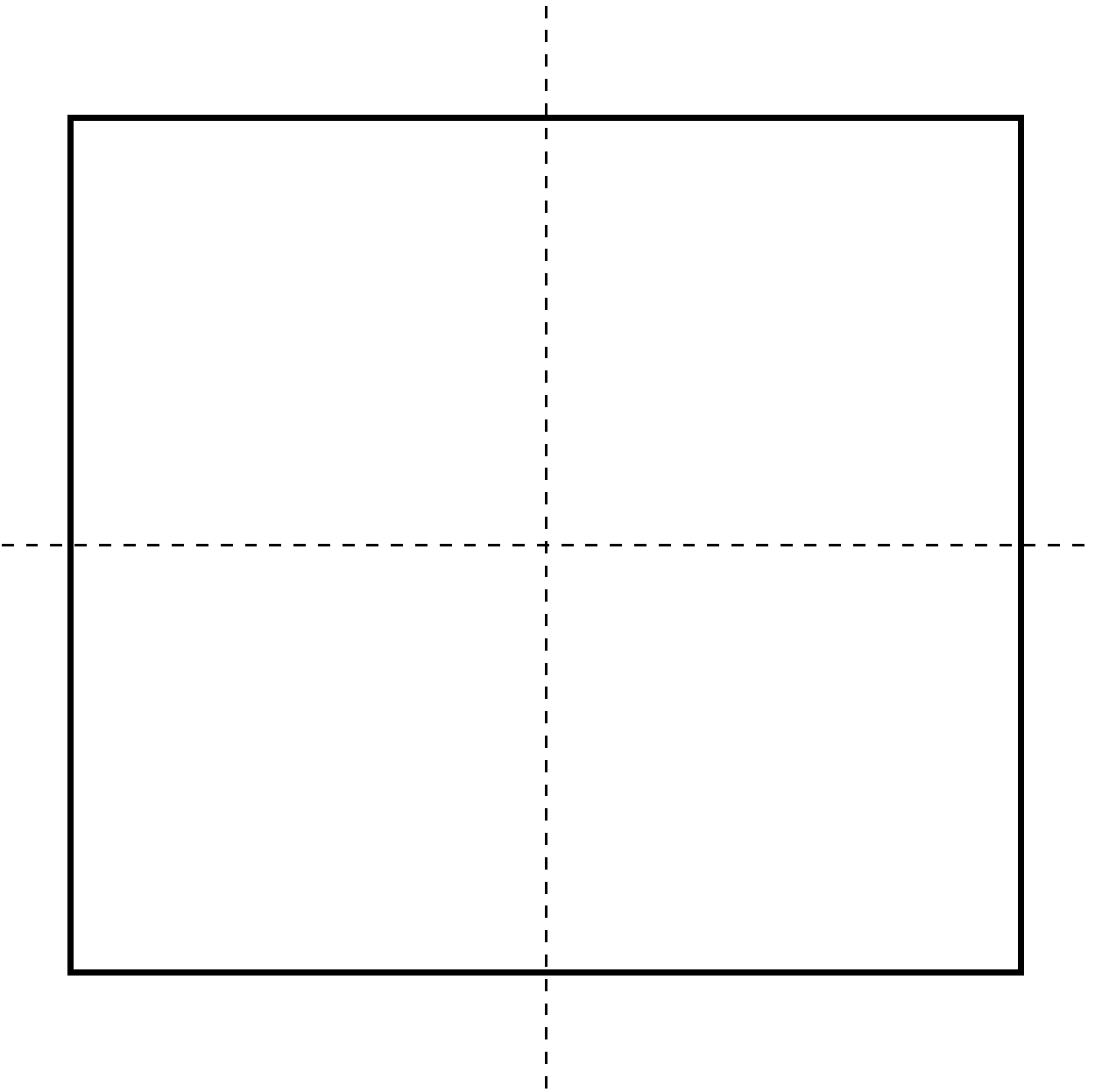}};
	\begin{scope}[x={(imageQ.south east)},y={(imageKa.north west)}]
	\node[above right] at (0.96,0.875) {$\widehat{c}_1$};
	\node[above right] at (0.40,0.875) {$c_1$};
	\node[above left] at (0.60,0.875) {$\widehat{c}_2$};
	\node[above left] at (0.04,0.875) {$c_2$};
	\node[below left] at (0.605,0.134) {$\widehat{c}_3$};
	\node[below left] at (0.04,0.12) {$c_3$};
	\node[below right] at (0.963,0.137) {$\widehat{c}_4$};
	\node[below right] at (0.40,0.12) {$c_4$};
	
	\node[right] at (0.97,0.66) {$\widehat{e}_1$};
	\node[right] at (0.40,0.65) {$e_1$};
	
	\node[above] at (0.7,0.88) {$\widehat{e}_2$};
	\node[above] at (0.14,0.875) {$e_2$};
	
	\node[left] at (0.6,0.35) {$\widehat{e}_3$};
	\node[left] at (0.042,0.34) {$e_3$};
	
	\node[below] at (0.86,0.118) {$\widehat{e}_4$};
	\node[below] at (0.3,0.12) {$e_4$};
	
	\node[above right] at (0.21875,0.5) {$0$};
	\node[above right] at (0.78125,0.5) {$0$};
	\end{scope}
	\end{tikzpicture}
	\caption{Edges and corners of the elliptic region $K_a$, defined in Definition \ref{defi:RK}, and of the quadrilateral $Q$, defined in equation \eqref{eq:defiQ}. These pictures were generated with the parameter choice $\nu=\frac{1}{3}$.}
	\label{fig:marked}
\end{figure}

$K_a$ is a quadrilateral domain, invariant under complex conjugation and reflection in the origin. It is a simply connected region whose boundary is a Jordan curve composed of four analytic pieces, which we call edges, meeting at four corners, as in Figure \ref{fig:marked}.
The interior of $K_a$ corresponds to potentials $V(\lambda)$ with four distinct zeros, whereas at the 
edges and corners coalescence of these zeros occurs: at edges a pair of zeros coalesces and
at corners three zeros coalesce.

Call 
$e_k, 1\leq k\leq 4$, the edges of $K_a$ and
$c_k, 1\leq k\leq 4$, the corners of $K_a$,
see Figure \ref{fig:marked}, so that
\begin{equation*}
\partial K_a=e_1\sqcup e_2\sqcup e_3\sqcup e_4 \sqcup\{c_1,c_2,c_3,c_4\}.
\end{equation*}

An elementary calculation shows that, upon imposing that the potential $V(\lambda;\alpha,\beta)$ has a zero with multiplicity three, $\alpha$ must be a root of the octic polynomial
\begin{equation}\label{eq:defiCintro}
C(\alpha):=\alpha^8-6(3\nu^2+1)\alpha^4+8(1-9\nu^2)\alpha^2-3(9\nu^4+6\nu^2+1).
\end{equation}
This polynomial has precisely two real roots, two purely imaginary roots, and four additional complex roots related to each other by reflections in the real and imaginary axes, see Lemma \ref{lem:zerospolynomial} for details.
 For $1\leq k\leq 4$, the corner $c_k$ of $K_a$ equals the unique root of $C(\alpha)$
in the $k$-th quadrant of the complex $\alpha$-plane.

\begin{figure}
	\begin{center}
		\includegraphics[width=0.4\textwidth]{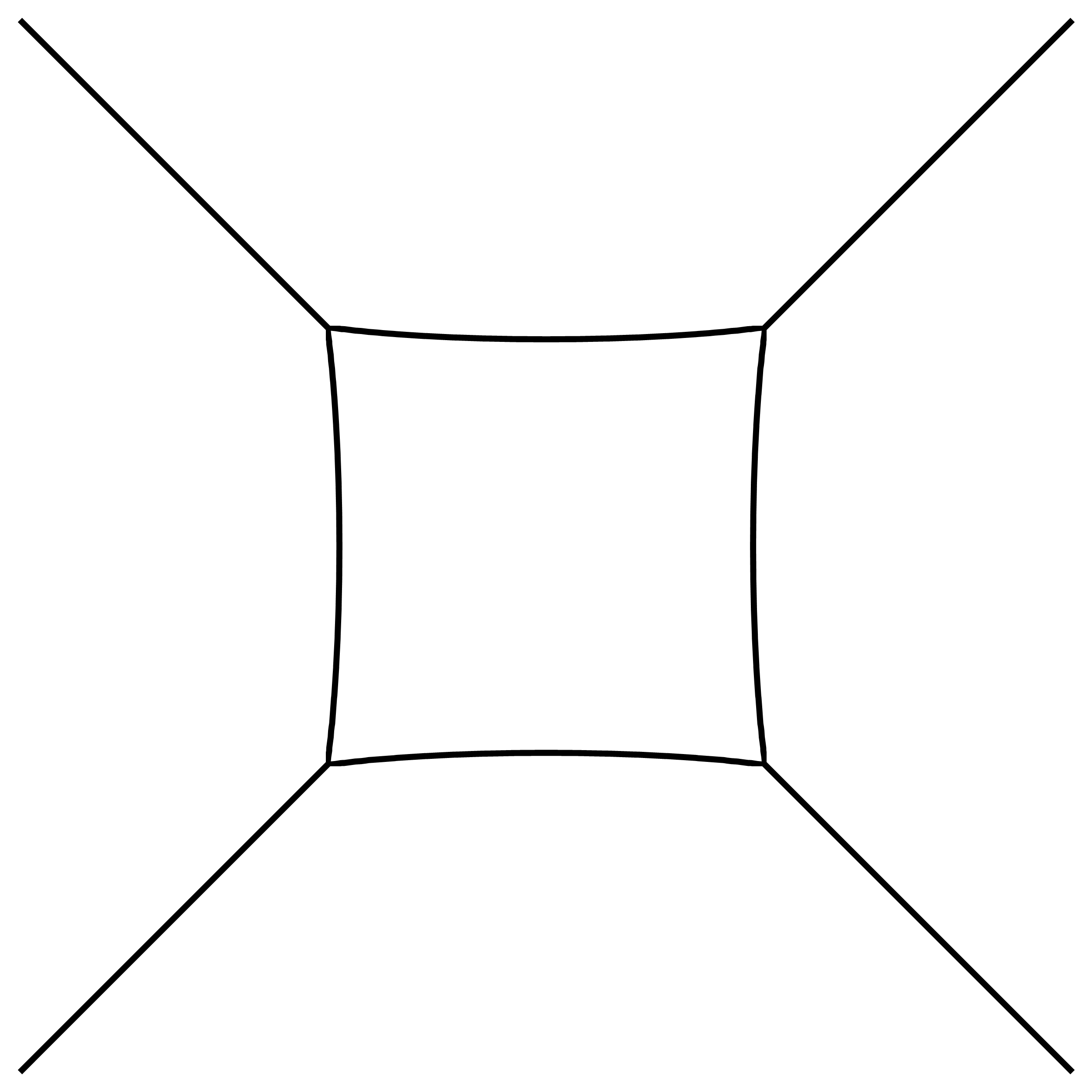} 
	\end{center}
	\caption{The zero set $\{\psi(\alpha)=0\}$, where $\psi$ is the harmonic function defined in equation \eqref{eq:parametrisation}, here depicted with $\nu=\frac{1}{3}$.} \label{fig:parametrisation}
\end{figure}

For $1\leq k\leq 4$, the edge $e_k$ is a smooth curve in the half-plane $\{\tfrac{1}{2}\pi(k-2)<\arg{\alpha}<\tfrac{1}{2}\pi k\}$
with end-points $c_{k-1}$ and $c_k$, where $c_{-1}:=c_4$. We have the following implicit parametrisation of $\partial K_a$:
let $x=x(\alpha)$ be the unique algebraic function which solves the quartic
\begin{equation*}
3 x^4+4\alpha x^3+(\alpha^2-1)x^2-\frac{\nu^2}{4}=0
\end{equation*}
analytically in the complex $\alpha$-plane with $x(\alpha)\sim\frac{\nu}{2}\alpha^{-1}$ as $\alpha\rightarrow \infty$ and branch-cuts
 the diagonals
$[c_1,c_3]$ and $[c_2,c_4]$. On the same cut plane there exists a unique algebraic function $y=y(\alpha)$ which solves
\begin{equation*}
y^2=\alpha^2+6 x\alpha+6x^2-1
\end{equation*}
with $y(\alpha)\sim \alpha$ as $\alpha\rightarrow \infty$. We set
\begin{equation}\label{eq:parametrisation}
\psi(\alpha)=\tfrac{1}{2}\Re\left[\alpha y+\tfrac{1}{2}(1-\nu)\log(p_1)-\log(p_2)+\nu \log(p_3)\right],
\end{equation}
 where
\begin{equation*}
p_1=1-2 x\alpha-2x^2,\quad p_2=2x+\alpha+y,\quad p_3=\frac{x(\alpha^2+5 x\alpha+4x^2-1) +\frac{1}{2}\nu  y}{x^2}.
\end{equation*}
Then $\psi$ is a univalued harmonic function on the cut plane, satisfying $\psi(-\alpha)=\psi(\alpha)$, and its level set $\{\psi(\alpha)=0\}$ consists of the boundary $\partial K_a$ plus four additional lines which emanate from the corners and go to infinity along the asymptotic directions $e^{\frac{\pi}{4}(2k-1)}\infty$, $1\leq k\leq 4$, see Figure \ref{fig:parametrisation}.
Buckingham \cite[Sections 1.1 and 3.2]{buckingham18} gives a slightly different but equivalent parametrisation of the boundary, which also includes the four additional lines emanating from the corners.

%

To formulate our asymptotic results, we have to introduce a particular function which maps the elliptic region $K_a$ homeomorphically onto the quadrilateral
\begin{equation}\label{eq:defiQ}
Q:=\left[-\tfrac{1}{2}(1-\nu)\pi,+\tfrac{1}{2}(1-\nu)\pi\right]\times \left[-\nu\pi,+\nu\pi\right].
\end{equation}
To this end, consider for general $(\alpha,\beta)$ the (possibly singular) affine elliptic curve 
\begin{equation}\label{eq:defigamma}
\Gamma_{\alpha,\beta} = \lbrace (\la,y)\in \C^2, \, y^2= \lambda^4+2\alpha \lambda^3+(\alpha^2-1)\la^2-
\beta\lambda+\frac{\nu^2}{4}=\la^2 V(\la)\rbrace \; .
\end{equation}
We denote by $\widehat{\Gamma}:=\widehat{\Gamma}(\alpha,\beta)$ its compactification, obtained by adding two points at infinity.
The elliptic curve is endowed with the projection $(\la,y)\mapsto \la$ onto the Riemann sphere $\Cb$.
The pull-back of the multivalued differential $\sqrt{V(\la)}d\la$ is the meromorphic differential $\omega=\frac{y}{\la}d\la$ on $\hat{\Gamma}$. 

Suppose that $(\alpha,\beta) \in R$, if we cut $\mathbb{P}^1$ by the Stokes lines connecting $\la_1$ with $\la_2$ and
$\la_3$ with $\la_4$, then
the function $\sqrt{V(\la)}$ is single-valued and
the smooth elliptic curve $\widehat{\Gamma}$ can be realised as a doubly sheeted cover of this cut Riemann sphere. We name 
the lower sheet the one fixed by the requirement $\lim_{\la \to +\infty }\Re \sqrt{V(\la)}=+\infty$
(equivalently, $\res_{\la=0}\sqrt{V(\la)}d\la=-\frac{\nu}{2}$ on the lower sheet).
To represent a curve in $\widehat{\Gamma}$ as a curve in the two-sheeted covering, we draw it dashed in the lower sheet and solidly
in the upper sheet. With this standard notation, the cycles $\gamma_1,\gamma_2$, defined in Figure \ref{fig:approximate}, form a basis
basis of $H_1(\widehat{\Gamma},\bb{Z})$.

By Definition \ref{defi:RK}, we have  $\Re \oint_{\gamma_1} \omega= \Re \oint_{\gamma_2} \omega=0$, thus yielding
a real mapping from $R$ to $\bb{R}^2$:
\begin{equation}\label{def:S}
	\mathcal{S}: R\rightarrow \mathbb{R}^2, (\alpha,\beta)\mapsto
	(-is_1(\alpha,\beta),-is_2(\alpha,\beta)),
\end{equation}
with
\begin{subequations}\label{eq:defs1s2}
	\begin{align}\label{eq:s1}
	s_1(\alpha,\beta)&=\int_{\gamma_1} \omega+\frac{i \pi (1-\nu)}{2},\\ \label{eq:s2}
	s_2(\alpha,\beta)&=\int_{\gamma_2} \omega.
	\end{align}
	\end{subequations}
In the following theorem some important properties of the mapping $\mathcal{S}$ are stated
and a homeomorphism between the elliptic region $K_a$ and the quadrilateral $Q$ is constructed which is crucial in the formulation of our asymptotic results.
\begin{thm}\label{thm:mappingS}
The mapping $\mathcal{S}$ has a unique continuous extension to $K=\overline{R}$ with range equal to the quadrilateral $Q$ defined in \eqref{eq:defiQ}, and this continuous extension
\begin{equation*}
\mathcal{S}:K\rightarrow Q
\end{equation*}
is a homeomorphism. Furthermore, the projection
\begin{equation}\label{eq:defiproj}
\Pi_a:K\rightarrow K_a, (\alpha,\beta)\mapsto \alpha
\end{equation}	
of $K$ onto $K_a$ is injective and the composition
\begin{equation}\label{eq:defiSa}
\mathcal{S}_a=\mathcal{S}\circ \Pi_a^{-1}:K_a\rightarrow Q
\end{equation}
is a homeomorphism, which
\begin{itemize}
	\item maps the interior of the elliptic region $C^\infty$-diffeomorphically onto the interior of $Q$,
	\item maps the edge $e_k$ $C^\infty$-diffeomorphically onto the (open) edge $\hat{e}_k$ of the rectangle $Q$ for $1\leq k\leq 4$,
	\item maps the corner $c_k$ to the corner $\hat{c}_k$ of the rectangle $Q$ for $1\leq k\leq 4$,
\end{itemize} 
where the (open) edges and corners of $Q$ are labelled as in Figure \ref{fig:marked}.
\end{thm}

\begin{rem}
For $1\leq k\leq 4$, the edges $e_k$ and $e_{k+1}$ meet at the 
corner $c_k$ with interior angle equal to $\frac{2}{5}\pi$. In particular $\mathcal{S}_a$ is not conformal.
\end{rem}

\subsection{Bulk Asymptotics}
\label{section:bulk}
For sake of simplicity we state our result when $m,n\to \infty$ with the ratio $\frac{m}{n}$ fixed.
We thus choose $p\geq q$, $p,q$ either equal or co-prime, and fix the ratio $\frac{m}{n}=\frac{p}{q}$.
Hence the numbers $m,n$ take values in the sequences $m=tq,n=tp, t \in \bb{N}^*$. 
Correspondingly $\nu=\frac{p}{2q+p}\in(0,\frac13]$ is fixed and the large parameter $E$ belongs to the sequence $(2q+p)t, t \in \bb{N}^*$.

\begin{defi}\label{def:Im}
Given an integer number $m \in \bb{N}^*$ we denote $I_m=\lbrace -m+1,-m+3,\dots,m-1\rbrace \subseteq \bb{Z}$.
 $\forall (j,k)  \in I_m \times I_n$
we let $(\alpha_{j,k},\beta_{j,k}) \in K$ be the unique solution of $\mathcal{S}(\alpha,\beta)=(\frac{\pi j}{E},\frac{\pi k}{E})$.
\end{defi}

\begin{defi}\label{def:Qsigma}
A filling fraction is a real number $\sigma \in (0,1)$.
We let $I_m^{\sigma}=I_m \cap [\sigma(-m+1), \sigma(m-1)]$ and define $Q^{\sigma} \subset Q$ as the closed rectangle $[-\frac{\pi \lfloor \sigma (m-1) \rfloor}{E},
\frac{\pi \lfloor \sigma (m-1) \rfloor}{E}]
\times[-\frac{\pi \lfloor \sigma (n-1) \rfloor}{E},
\frac{\pi \lfloor \sigma (n-1) \rfloor}{E}]$, which in the large $E$ limit converges to
$\sigma\cdot Q$.

Finally we define $K^{\sigma}=\mathcal{S}^{-1}(Q^{\sigma})$ and $K_a^{\sigma}=\Pi_a(K^{\sigma})$ as the projection of $K^{\sigma}$ on the $\alpha$-plane.
\end{defi}

By means of the complex WKB method, we will prove -- in Section \ref{section:bulkasymptotic} -- the following asymptotic description
of the solutions $(\alpha,\beta)$ of
the inverse monodromy problem characterising the roots of generalised Hermite polynomials.
\begin{thm}\label{thm:bulkbeta}
 Fix $\sigma \in (0,1)$. Then there exists an $R_{\sigma}>0$ such that for $E$ large enough the following hold true:
 \begin{enumerate}
  \item In each
 ball of centre $(\alpha_{j,k},\beta_{j,k})$, $(j,k) \in I_m^{\sigma} \times I_n^{\sigma}$ and radius $R_{\sigma}E^{-2}$ there exists
 a unique point $(\alpha,\beta)$ such that the anharmonic oscillator \eqref{eq:scaled} satisfies
 the inverse monodromy problem characterising the roots of generalised Hermite polynomials.
 \item In the $\epsilon$ neighbourhood of $K^{\sigma}$ with radius $R_{\sigma}E^{-2}$, there are exactly
  $\lfloor \sigma m \rfloor \times \lfloor \sigma n \rfloor$ points $(\alpha,\beta)$ such that the anharmonic oscillator
 (\ref{eq:scaled}) satisfies the inverse monodromy problem characterising the roots of generalised Hermite polynomials.
 \end{enumerate}
\end{thm}

Our main result is the
following theorem which provides an asymptotic formula for
the roots of the generalised Hermite polynomials
as well as an estimate of the error. It is a straightforward corollary of Theorem \ref{thm:bulkbeta}.
\begin{thm}\label{thm:bulkintro}
 Fix  a filling fraction $\sigma \in (0,1)$. Then there exists an $R_{\sigma}>0$ such that for $E$ large enough the following hold true:
  Each
 disc with center $E^{\frac12}\alpha_{j,k}$ and radius $R_{\sigma}E^{-\frac32}$, $(j,k) \in I_m^{\sigma} \times I_n^{\sigma} $,
 contains a root of the generalised Hermite polynomial $H_{m,n}$.
 \begin{proof}
  A point $a \in \bb{C}$ is a root of $H_{m,n}$ if and only if there exists $\beta$ such that $(E^{-\frac{1}{2}}a,\beta)$
  provides a solution to the inverse monodromy problem for the scaled oscillator (\ref{eq:scaled}). Since 
  the projection $\Pi_a:K^{\sigma} \to K^{\sigma}_a$ is a local diffeomorphism (Proposition \ref{pro:projection}), the thesis follows
  from Theorem \ref{thm:bulkbeta}.
 \end{proof}
\end{thm}

\begin{figure}[ht]
	\centering
	\begin{tabular}{ c c c }
		\includegraphics[width=4.5cm]{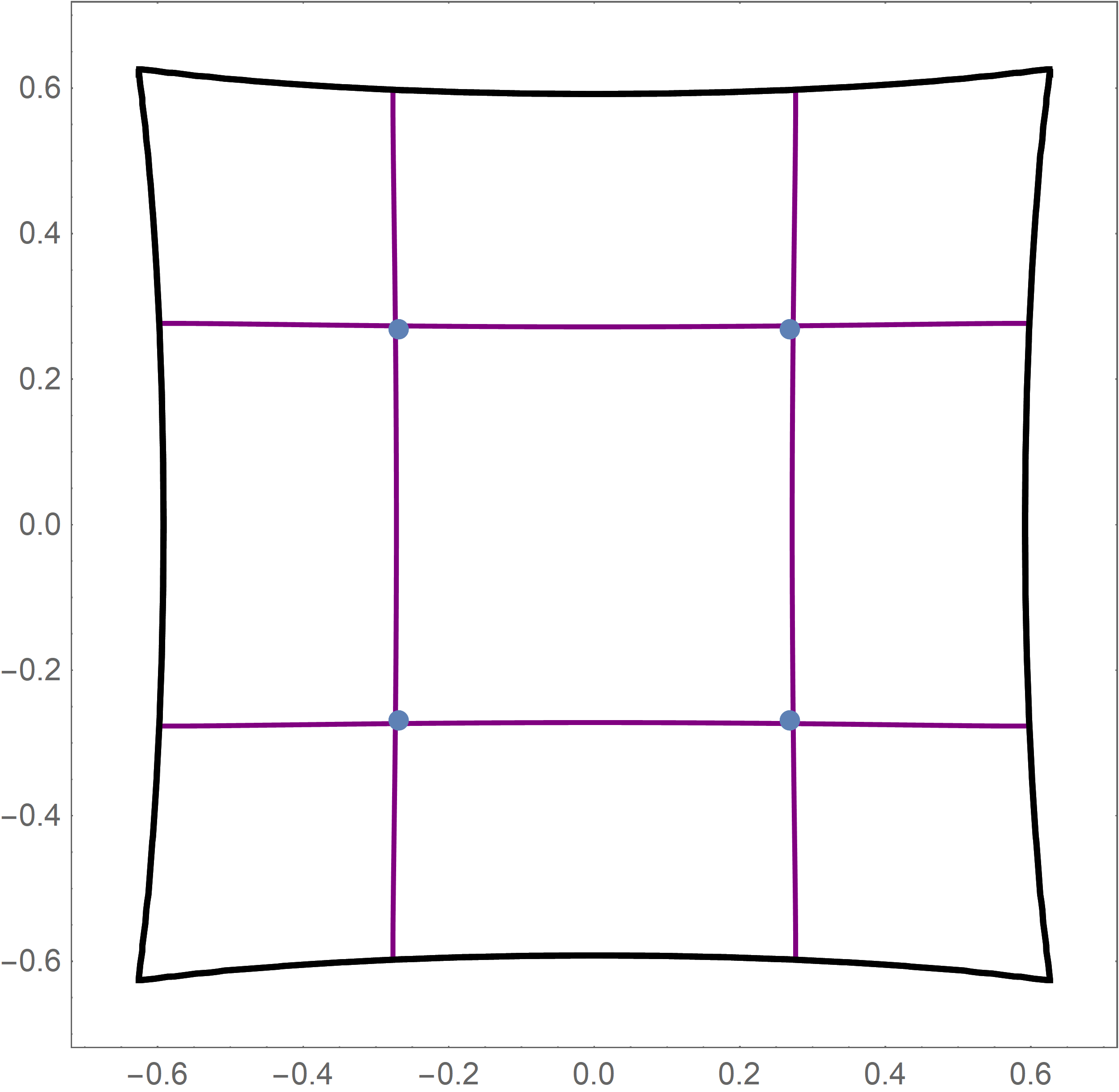} & \includegraphics[width=4.5cm]{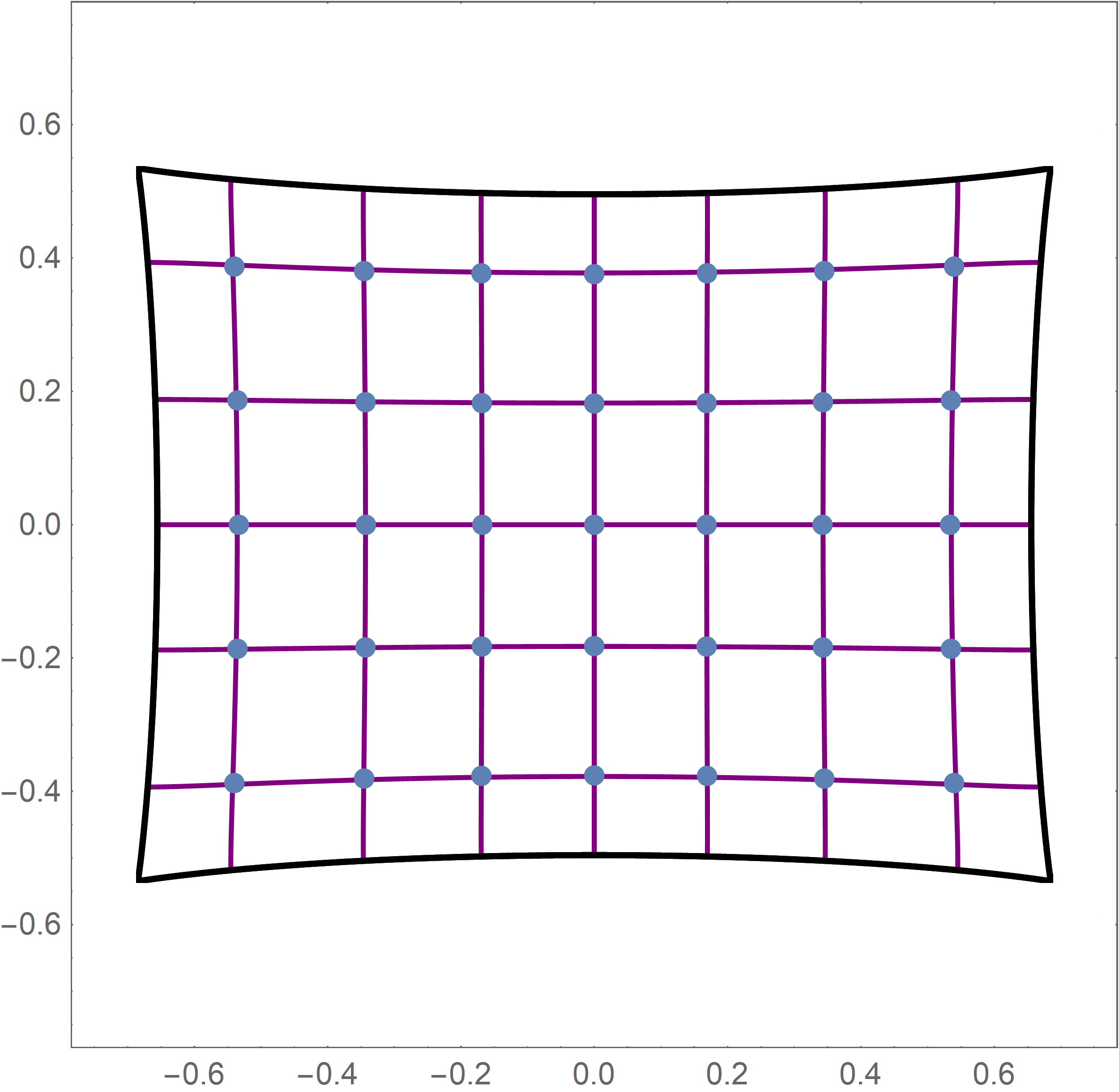} & \includegraphics[width=4.5cm]{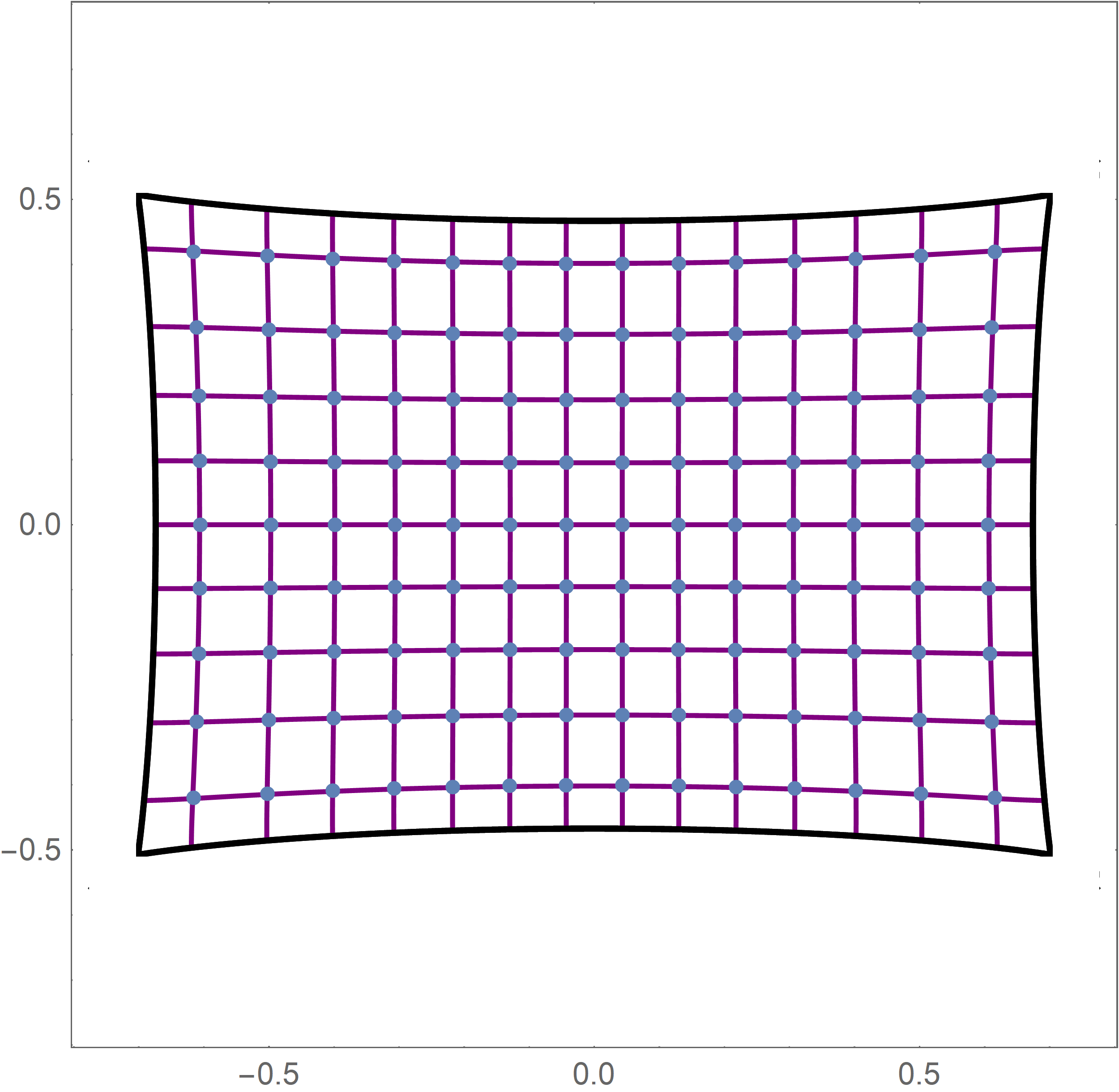}\\
		\includegraphics[width=4.5cm]{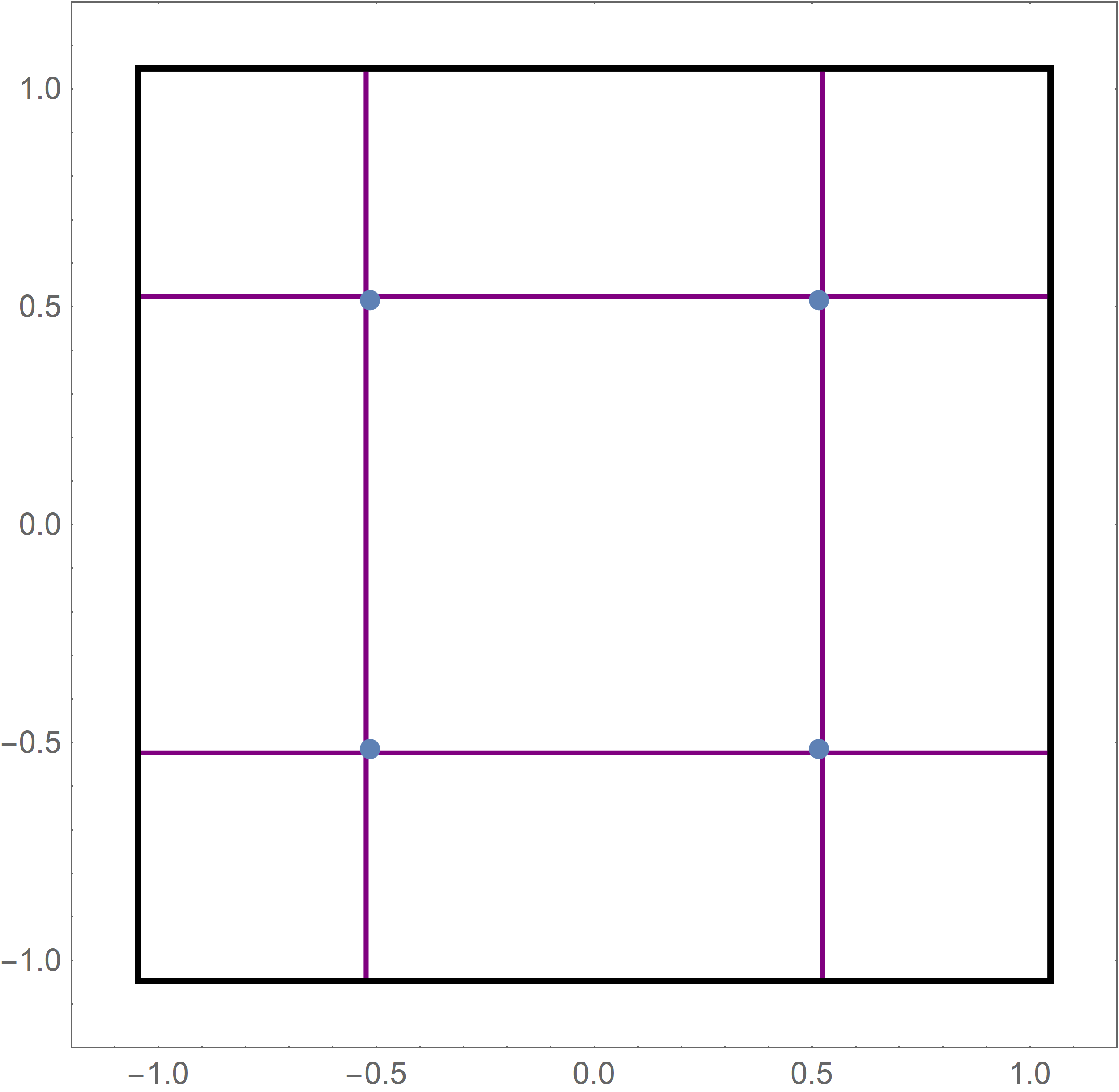} & \includegraphics[width=4.5cm]{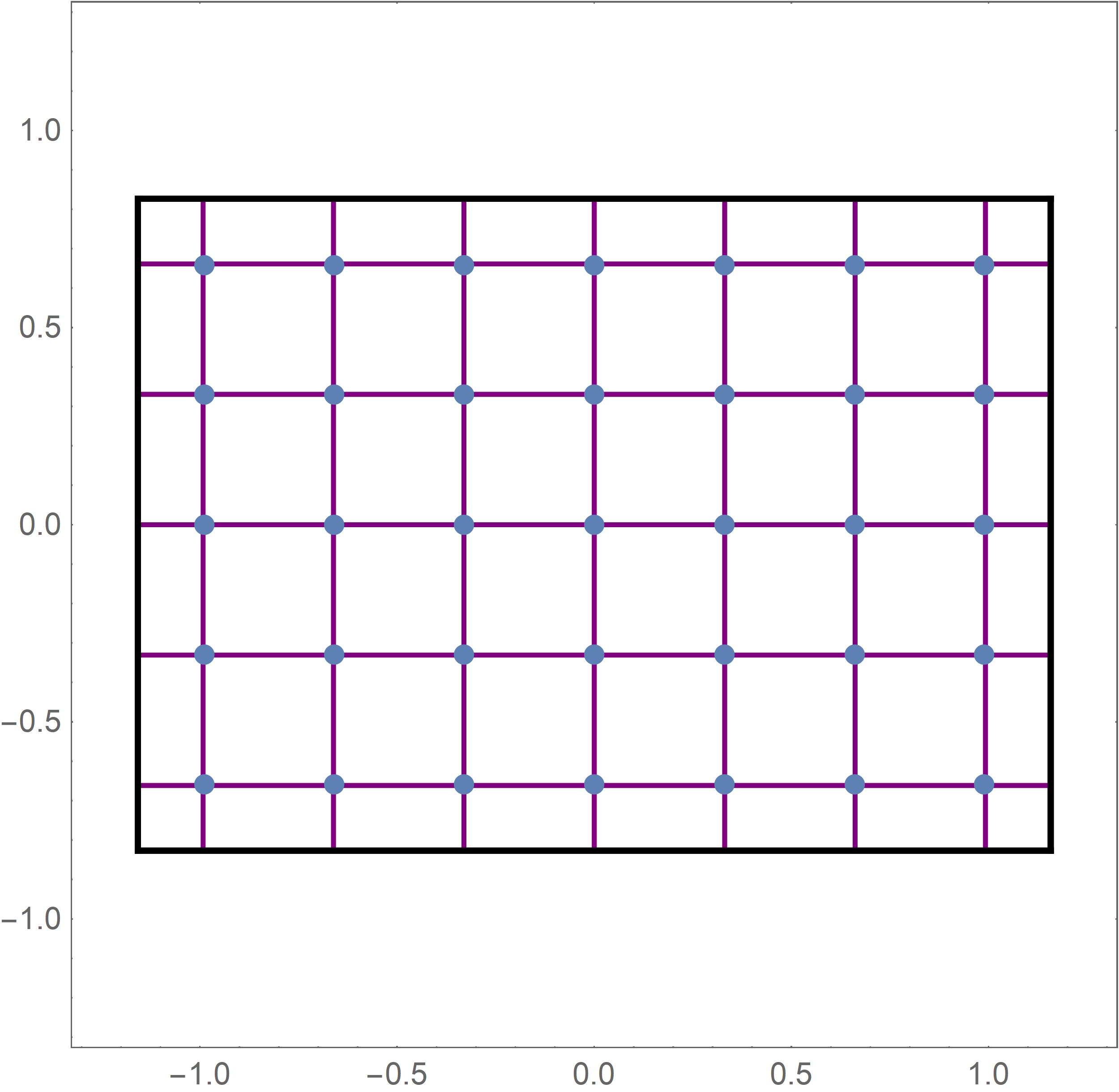} & \includegraphics[width=4.5cm]{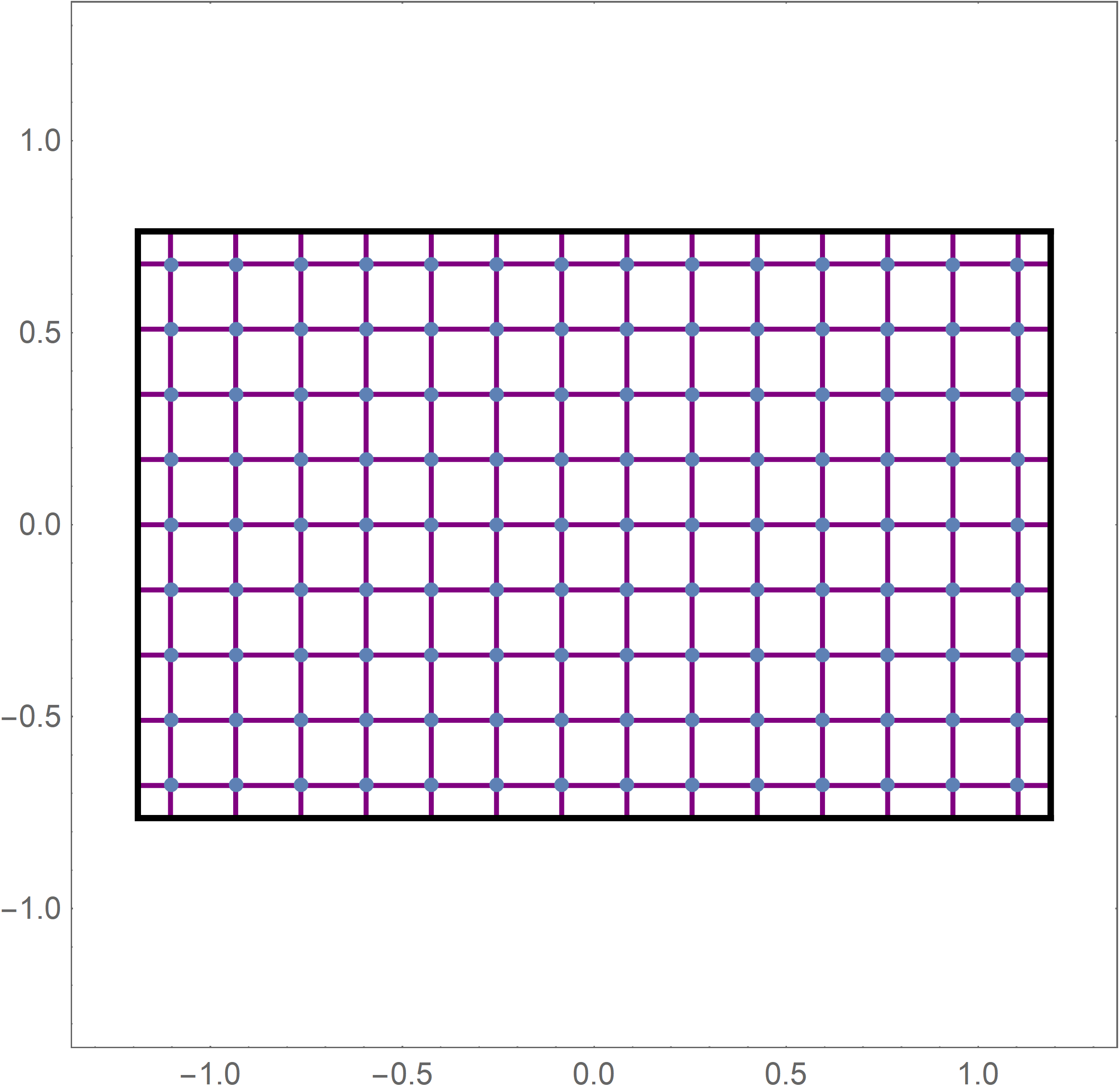}
	\end{tabular}
	\caption{On the top row the elliptic region $K_a$ with in purple the inverse images under $\mathcal{S}_a$ of the grid lines \eqref{eq:gridlines} and the roots of $H_{m,n}(E^\frac{1}{2}\alpha)$ in blue superimposed, and on the bottom row the rectangle $Q$ with in purple the grid lines \eqref{eq:gridlines} and the images of the roots of $H_{m,n}(E^\frac{1}{2}\alpha)$ under $\mathcal{S}_a$ in blue superimposed
		for the respective values $(m,n)=(2,2)$, $(7,5)$, $(14,9)$.} \label{fig:grid}
\end{figure}

Clarkson \cite{clarkson2003piv,clarksonspecialpol} observed numerically that the zeros
of generalised Hermite polynomial $H_{m,n}$ seem to organise themselves on the intersection points of
a deformed grid of $m$ vertical and $n$ horizontal lines, see the top pictures in Figure \ref{fig:grid}. The mapping $\mathcal{S}_a$ rectifies this deformed grid.
Namely, after Theorem \ref{thm:bulkintro}, the images under $\mathcal{S}_a$ of the rescaled roots of $H_{m,n}$ asymptotically organise
themselves along the intersection points of the true rectangular grid made of the respective equally spaced $m$
vertical and $n$ horizontal lines
\begin{equation}\label{eq:gridlines}
l_v^{(j)}= \{(x,y) \in Q, x=\frac{\pi j}{E}\}\quad  (j \in I_m), \quad l_h^{(k)}=\lbrace (x,y) \in Q, y=\frac{\pi k}{E}\}\quad (k \in I_n),
\end{equation} 
Consequently, the rescaled roots of $H_{m,n}$ asymptotically organise themselves along the intersection points of the deformed  grid made out of the deformed grid lines $\mathcal{S}_a^{-1}(l_v^{(j)}), j\in I_m$ and $\mathcal{S}_a^{-1}(l_h^{(k)}), k\in I_n$, see Figure \ref{fig:grid}.

The deformed grid is actually a regular grid for small $\alpha$'s. Indeed we have the following immediate corollary
of Theorem \ref{thm:bulkintro}.
\begin{cor}\label{cor:smallalphaintro}
Fix $N_0$ and suppose $|j|,|k| \leq N_0 $.
The approximate roots $\alpha_{j,k}$ lie on the regular grid generated by the vectors
$\frac{\mathcal{K}\left(\frac{1-\nu}{1+\nu}\right)}{E \sqrt{1+\nu}}$
and $i \frac{\mathcal{K}\left(\frac{2\nu}{1+\nu}\right)}{E \sqrt{1+\nu}}$,  where $\mathcal{K}(m)$ denotes
the standard complete elliptic integral of the first kind with respect to the parameter $m=k^2$.
More precisely, the following asymptotic formula holds
\begin{equation}\label{eq:smalljkintro}
 \alpha_{j,k}= \frac{1}{\sqrt{1+\nu}}\left(\mathcal{K}\left(\frac{1-\nu}{1+\nu}\right)\frac{j}{E}+i
  \mathcal{K}\left(\frac{2\nu}{1+\nu}\right)\frac{k}{E}\right)+ O\big(E^{-2}\big)
\end{equation}
Note that the above equation describes a square grid when $\nu=\frac{1}{3}$.
\begin{proof}
 By Theorem \ref{thm:bulkbeta}, we need to solve $S(\alpha_{j,k},\beta_{j,k})=(\frac{\pi j}{E},\frac{\pi k}{E})$ for $j,k$ bounded.
 Since $j,k$ are bounded, up to a uniform
(and immaterial) $O(E^{-2})$ error we can solve the above equation using the first order Taylor expansion of $\mc{S}$ at
$(\alpha,\beta)=(0,0)$ (recall that $\mc{S}(0,0)=(0,0)$).
More precisely, if we let $J$ be the Jacobian of $\mc{S}$ at $(\alpha,\beta)=0$, then we have that
$\alpha_{j,k}$ is, up to the $O(E^{-2})$ error, the first component of the solution of the linear equation
$J\cdot(\alpha,\beta)=(\frac{\pi j}{E},\frac{\pi k}{E})$. That is,
\begin{align*}
 \alpha_{j,k}= (\det J)^{-1} \frac{\pi i}{E} \big( J_{22} j- J_{12} k \big)+O\big(E^{-2}\big) \; .
 \end{align*}
 The thesis follows from the above formula upon substituting the actual values of the entries of $J$
 given in equation \eqref{eq:jacobianexplicit} in the appendix.
 \end{proof}
\end{cor}

\begin{figure}[htpb]
		\centering
		\includegraphics[width=8cm]{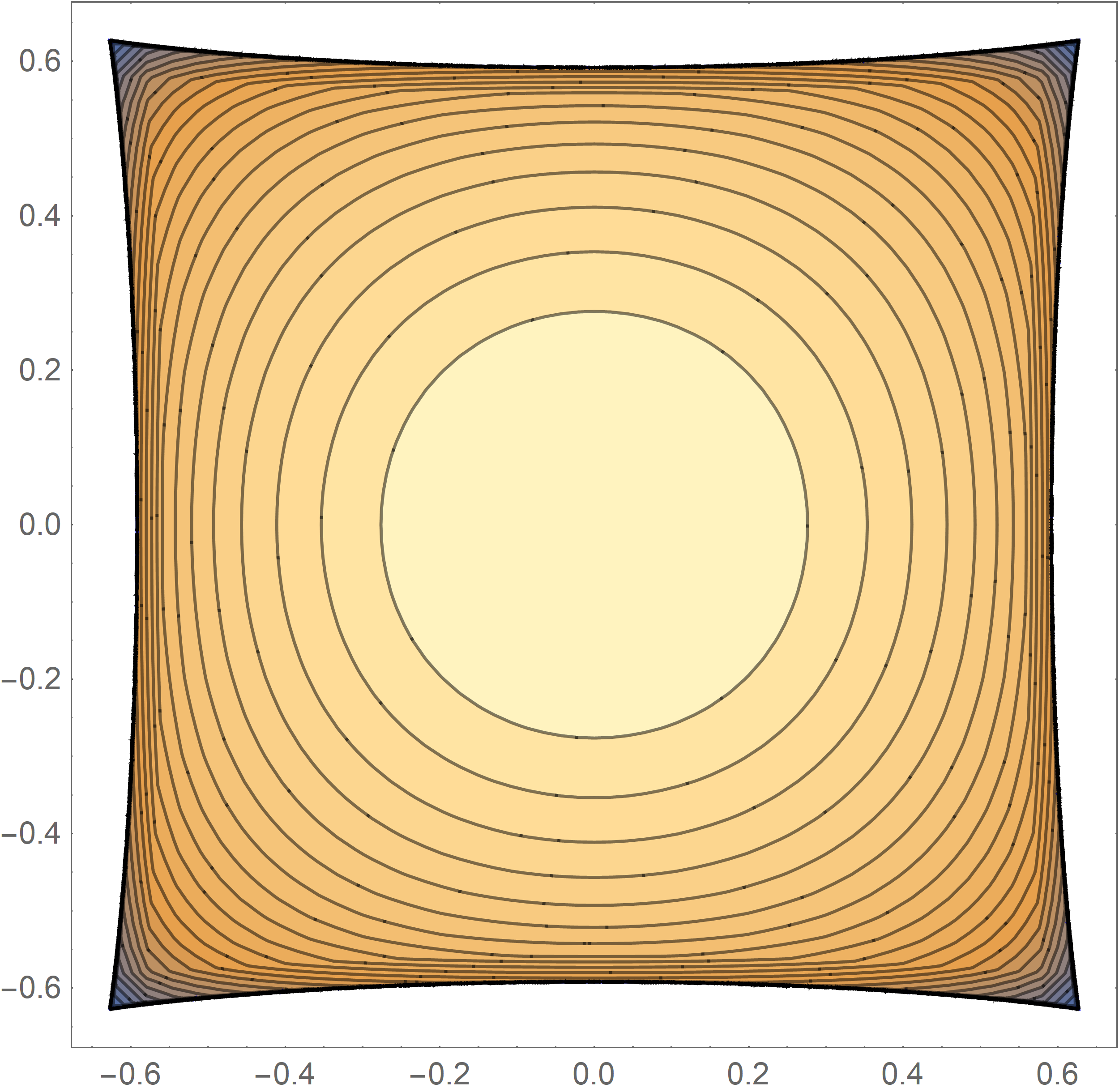}\hspace{1cm}
		\includegraphics[width=1cm]{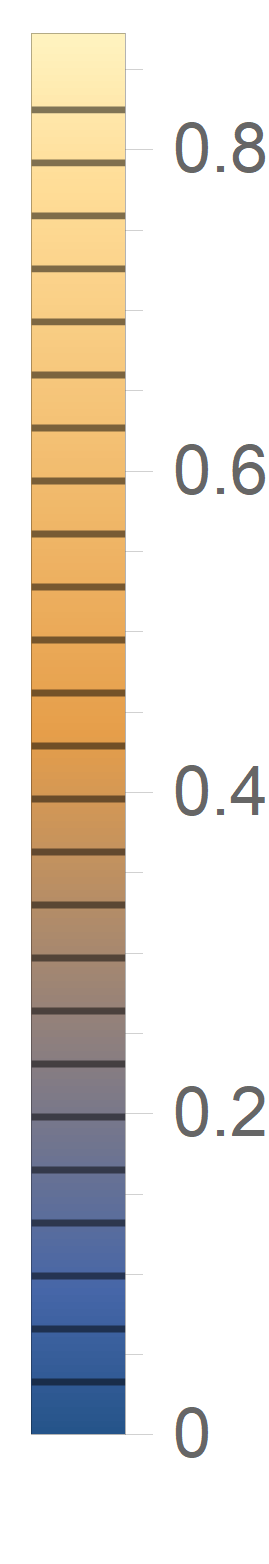} 
		
		\caption{Asymptotic density of roots of $H_{m,n}(E^{\frac{1}{2}}\alpha)$, described by the density function \eqref{eq:density_function}, as $E\rightarrow \infty$ with $\nu=\frac{1}{3}$ within the elliptic region $K_a$.}\label{fig:density}
	\end{figure}

\begin{rem}
	After Theorem \ref{thm:bulkintro}, we know that, asymptotically, almost all roots of $H_{m,n}(E^{\frac{1}{2}}\alpha)$ lie in the elliptic region $K_a$ and that the images of those under $S_a$ are uniformly distributed within the quadrilateral $Q$. Therefore, the asymptotic distribution of the roots in the $\alpha$-plane is described by the density function	
	\begin{equation}\label{eq:density_function}
	\Phi_\nu(\alpha)=\begin{cases}
	\frac{1}{2\nu(1-\nu)\pi^2}\begin{vmatrix}
	\frac{\partial S_1}{\partial \alpha_R} & \frac{\partial S_1}{\partial \alpha_I}\\
	\frac{\partial S_2}{\partial \alpha_R} & \frac{\partial S_2}{\partial \alpha_I}
	\end{vmatrix} &\text{if $\alpha\in R_a$,}\\
	0 & \text{if $\alpha\in \mathbb{C}\setminus R_a$,}
	\end{cases} 
	\end{equation}
	where $\alpha=\alpha_R+\alpha_I i$, with $\alpha_R,\alpha_I\in\mathbb{R}$, which is continuous on $\mathbb{C}$. See Figure \ref{fig:density} for a plot of this density function with $\nu=\frac{1}{3}$.
\end{rem}

\paragraph*{\textbf{A conjectural formula describing the critical behaviour}}
We announce here a conjectural formula describing the critical asymptotics of roots of generalised Hermite polynomials.

Theorem \ref{thm:bulkintro} does not cover the asymptotics of roots approaching the corner of the elliptic region. These are called
critical asymptotics, and, according to our computations, are described by means of the Tritronqu\'ee solution $y_T$
of the Painlev\'e I equation
\begin{equation}\label{eq:PI}
 y''(z)=6y^2(z)-z.
\end{equation}
Let us be more precise.
The Tritronqu\'ee solution, which was discovered by Boutroux \cite{boutroux}, can be defined as the unique
solution of equation \eqref{eq:PI} which does not have poles in the closed sector $|\arg{z}|\leq \frac{4 \pi}{5}$ \cite{costin12}.
It has however  an infinite number of poles \cite{piwkb}, which a fortiori lie in the sector $|\arg{z}|> \frac{4 \pi}{5}$.
Let $c_1$ be the upper right corner of the elliptic region, and $p$ be a pole of the Tritronqu\'ee solution $y_T$.
Our computations led us to the following conjectural formula (whose proof will be provided in a forthcoming publication \cite{maroprep}):
for $E$ large enough, there is a unique root $\alpha_p$ of $H_{m,n}(E^{\frac{1}{2}}\alpha)$
with the following asymptotic behaviour,
 \begin{equation}\label{eq:tritronqueeas}
 \alpha_p=c_1-\kappa p E^{-\frac{2}{5}}+\mathcal{O}(E^{-\frac{3}{5}}) \mbox{ as } E\rightarrow  \infty \; .
 \end{equation}
 Here $\kappa$ is a constant, independent of $p$, defined by the equation
 \begin{equation*}
 \kappa^5=\frac{2 c_1^3(2+c_1^2)^3}{(c_1^4-3\nu^2-1)(c_1^4+4 c_1^2+3\nu^2+1)},
 \quad -\frac{3}{4}\pi-\frac{\pi}{10}\leq \arg \kappa<-\frac{3}{4}\pi+\frac{\pi}{10}.
 \end{equation*}
The critical asymptotics for the other corners are given by the above formulas,
upon substituting $c_1$ for the corresponding $c_k,k\neq 1$, and by choosing
the appropriate solution of the quintic equation.

The above conjectural formula is not totally unexpected,
as the poles of the Tritronqu\'ee solution describe simultaneous quantisation of the cubic oscillator 
\cite{piwkb} and near the corners of the elliptic region, three of the zeros of the potential coalesce, so that,
in appropriate local coordinates, the anharmonic oscillator \eqref{eq:anharmonic} becomes a cubic oscillator.

\section{The Elliptic Region}
\label{section:elliptic}
This section is devoted to the study of the regions $R$, $K$ and $K_a$, defined in Definition \ref{defi:RK},
as well as the mapping $\mathcal{S}:R\rightarrow \mathbb{R}^2$ and the projection $\Pi_a:K\rightarrow K_a$ introduced in equations \eqref{def:S} and \eqref{eq:defiproj} respectively. The main objectives of this section are proving Theorem \ref{thm:mappingS} and proving that the implicit parametrisation of the elliptic region given in Section \ref{sec:results} is correct. The proofs revolve around the Stokes geometry of the potential $V$ and involve the following main steps.
\begin{enumerate}[label={S.\arabic*}]
	\item\label{step1} We show that $R$ is a smooth real two-dimensional regular submanifold of $\{(\alpha,\beta)\in\mathbb{C}^2\}$ and that the mapping $\mathcal{S}$ is a local diffeomorphism. 
	\item\label{step2} 	We prove that $\mathcal{S}$ has a unique continuous extension to $K=\overline{R}$ with corresponding range a subset of the quadrilateral $Q$ defined in equation \eqref{eq:defiQ}.
	\item\label{step3} We show that the continuous extension
	\begin{equation*}
	\mathcal{S}:K\rightarrow Q
	\end{equation*}
	 is a homeomorphism which maps the interior $R$ of $K$ diffeomorphically onto the interior of $Q$.	
	 \item\label{step4} We derive an implicit parametrisation of the boundary $\delta R:=K\setminus R$ of $R$, as well as of $\Pi_a(\delta R)$.
	 \item\label{step5} We show that $\partial K_a=\Pi_a(\delta R)$, which in conjunction with the parametrisation of $\Pi_a(\delta R)$ in step \ref{step4}, proves the correctness of the description of the elliptic region $K_a$ given in Section \ref{sec:results}.
	 \item\label{step6} We show that the projection $\Pi_a:K\rightarrow K_a$ of $K$ onto the elliptic region is a homeomorphism, which maps the interior of its domain diffeomorphically onto the interior of its co-domain. 
	  \item\label{step7} After steps \ref{step3} and \ref{step6}, we know that $\mathcal{S}_a=\mathcal{S}\circ \Pi_a^{-1}:K_a\rightarrow Q$ maps the elliptic region $K_a$ homeomorphically onto the quadrilateral $Q$, and the interior of the former diffeomorphically onto the interior of the latter. In the final step, we show that $\mathcal{S}_a$ maps edges and corners of the elliptic region to the edges and corners of $Q$, as indicated in Theorem \ref{thm:mappingS}, which finishes the proof of Theorem \ref{thm:mappingS}.
\end{enumerate}
We carry out the various steps in the order they are presented above, due to their logical dependencies in our approach.

The present section is organised as follows.
In Subsection \ref{subsec:stokes} we recall some basic notions of the theory of Stokes complexes tailored to the potential
\begin{equation}\label{eq:potential}
V(\lambda;\alpha,\beta)=\lambda^2+2\alpha \lambda+\alpha^2-1-\beta\lambda^{-1}+\frac{\nu^2}{4}\lambda^{-2}.
\end{equation}
We introduce the notion of Boutroux curves \cite{bertola2011,bertola2009,boutroux,piwkb}
and collect some basic results relevant to this paper.
In Subsection \ref{subsec:results}, we carry out the Steps \ref{step1}-\ref{step7} of the proof of Theorem \ref{thm:mappingS} and
of the implicit parametrisation of the elliptic region. Finally, Subsection \ref{subsec:anti} is devoted to
to the study of the anti-Stokes complex of the potential
$V(\lambda;\alpha,\beta)$ for $(\alpha,\beta) \in R$. The anti-Stokes complex does not play any role
in the proof of Theorem \ref{thm:mappingS} but it is an important technical tool for the complex WKB method.

To the best of our knowledge, such a meticulous study of the elliptic region, for this or similar problems, has not been carried out before.
Our approach is unfortunately complex and lengthy. Ultimately the difficulty in the proofs lies
in the fact that we want to prove global results of mappings, like $\mathcal{S}$ and $\Pi_a$,
which are not conformal, over which we only have strong control locally and only implicitly, and whose
domains, $K$ and $K_a$, have a parametrically complicated and merely piece-wise smooth boundary.

We hope that our approach will be simplified with the help of future research.

\subsection{Stokes Complexes}\label{subsec:stokes}
In this subsection we review some basic notions of WKB theory, tailored to the potential \eqref{eq:potential}, and we refer to \cite{fedoryuk,strebel,jenkins} for the more general setting.

Firstly we recall the notion of turning points.
\begin{defi}\label{def:turn}
	A zero of multiplicity $n$ of \eqref{eq:potential} is called a turning point of degree $n$ of the potential,
	for $n\geq 1$. The turning points and $\lambda=0$ are called the critical points of the potential. Other points are called generic.
\end{defi}
Let $\lambda_0$ be a generic point and consider, for any choice of sign, the so called action integral
\begin{equation}\label{eq:actionintegral}
S(\lambda_0,\lambda)=\int_{\lambda_0}^{\lambda}\sqrt{V(\mu;\alpha,\beta)}d\mu
\end{equation}
defined on the universal covering of the $\lambda$-plane minus critical points.  Away from critical points we have
\begin{equation}\label{eq:actionregular}
\frac{\partial}{\partial \lambda}S(\lambda_0,\lambda)\neq 0,
\end{equation}
on the universal covering. Let $i_{\lambda_0}$ be the level curve $\Re S(\lambda_0,\lambda)=0$ through $\lambda_0$ on the universal covering, then it follows from equation \eqref{eq:actionregular} that $i_{\lambda_0}$ is a non-self intersecting line in aforementioned universal covering. Now consider its projection $\widetilde{i}_{\lambda_0}$ onto the $\lambda$ plane minus critical points. Since the action integral can only differ by a minus sign on different sheets, it follows that
 $\widetilde{i}_{\lambda_0}$ is everywhere locally diffeomorphic to a line (and in particular does not self-intersect), and since it is connected, it is either diffeomorphic to a circle or a line.
We call $\widetilde{i}_{\lambda_0}$ the pojected level curve through $\lambda_0$.

Note that different projected level curves cannot intersect and thus the set of projected level curves forms a complete foliation
of the $\lambda$-plane minus critical points. We have the following important dichotomy.
\begin{lem}\label{lem:stokes_curves}
	Let $\lambda_0$ be a generic point of the potential \eqref{eq:potential} and $\widetilde{i}_{\lambda_0}$
	its corresponding level curve. If $\widetilde{i}_{\lambda_0}$ is diffeomorphic to a circle, then it is
	homotopic to a simple encircling of $\lambda=0$ in the $\lambda$-plane minus critical points.
	Otherwise $\widetilde{i}_{\lambda_0}$ is diffeomorphic to a line. Let $x\mapsto \gamma_{\lambda_0}(x)$
	be a diffeomorphism of $\mathbb{R}$ onto $\widetilde{i}_{\lambda_0}$, then, for $\epsilon\in \{\pm 1\}$,
	we have the following dichotomy: either
	\begin{enumerate}[label=(\roman*)]
		\item $\lim_{x\rightarrow \epsilon \infty}=\infty$ and the curve is asymptotic to one of the
		four rays $e^{\frac{1}{4}(2k-1)\pi i}\mathbb{R}_+$, $k\in\mathbb{Z}_4$, in the $\lambda$-plane.
		\item or $\lim_{x\rightarrow \epsilon \infty}=\lambda_*$ with $\lambda_*$ a turning point of the potential.
	\end{enumerate}
\end{lem}
\begin{proof}
	See Strebel \cite{strebel}.
\end{proof}
In alignment with the above lemma, we make the following definition.
\begin{defi}\label{defi:dichotomy}
	Considering the dichotomy in Lemma \ref{lem:stokes_curves}, when $\widetilde{i}_{\lambda_0}$ is
	diffeomorphic to a line, we call, in case (i) or (ii) respectively the asymptotic direction $\infty_k:=e^{\frac{1}{4}(2k-1)\pi i}\infty$ or the corresponding turning point $\lambda_*$ an endpoint of $\widetilde{i}_{\lambda_0}$.
	We call $\widetilde{i}_{\lambda_0}$ a Stokes line if at least one endpoint is a turning point.
\end{defi}
We define $\mathbb{C}_\infty$ as the complex plane with the addition of four marked points at
infinity $\infty_k:=e^{\frac{1}{4}(2k-1)\pi i}\infty$, $1\leq k \leq 4$, with the unique topology
making the following map a homeomorphism,
\begin{align}
&L:\mathbb{C}_\infty\rightarrow \mathbb{D}\cup\{e^{\frac{1}{4}(2k-1)\pi i}:1\leq k\leq 4\},\label{eq:defiL}\\
&L(\rho e^{i\phi})=\frac{2}{\pi} \arctan(\rho)e^{i\phi}\quad (\rho\in \mathbb{R}_{\geq 0}, \phi \in \mathbb{R}),\nonumber\\
&L(\infty_k)=e^{\frac{1}{4}(2k-1)\pi i}\quad (1\leq k\leq 4),\nonumber
\end{align}
where $\mathbb{D}$ denotes the open unit disc. We denote $\mathbb{C}_\infty^*=\mathbb{C}_\infty\setminus\{0\}$.

\begin{defi}\label{def:stokescomplex}
	The Stokes complex $\mc{C}=\mc{C}(\alpha,\beta)\subseteq \mathbb{C}_\infty^*$ of the potential
	\eqref{eq:potential} is the union of the marked points at infinity and the Stokes lines and turning points
	of the potential. The corresponding internal Stokes complex is defined as the union of the turning points and
	those Stokes lines with only turning points as endpoints.\\
	 The Stokes complex $\mc{C}$ is an embedded marked graph into $\mathbb{C}_\infty^*$ with
vertices equal to the turning points and the four points at infinity, with edges given by the Stokes lines. We call
two Stokes complexes isomorphic if they are isomorphic as embedded marked graphs.
\end{defi}

In the following proposition we summarise some basic facts concerning the Stokes complexes under consideration.
\begin{pro}\label{pro:basic}
	Let $\mathcal{C}$ be the Stokes complex of the potential \eqref{eq:potential}, then
	\begin{itemize}
		\item For any turning point $\lambda^*$, say of degree $n$, there are precisely $n+2$ Stokes
		lines emanating from it, counting Stokes lines with all end points equal to $\lambda_*$ double.
		\item $\mathcal{C}$ is connected;
		\item The internal Stokes complex contains precisely one Jordan curve, the interior of which
		contains $\lambda=0$;
		\item If two different Stokes lines have the same endpoints, then they are homotopically
		inequivalent in $\mathbb{C}_\infty^*$.
	\end{itemize}
\end{pro}
\begin{proof}
	See \cite{fedoryuk,strebel}.
\end{proof}

\begin{figure}[htpb]
	\centering
	\begin{tikzpicture}[scale=0.7]
	\tikzstyle{star}  = [circle, minimum width=4pt, fill, inner sep=0pt];
	
	\draw[black] (0,0) ellipse (6.3cm and 2.15cm);
	\node[blue,star]     at (-3,0){};
	\node[left]     at (-3,0) {$\lambda_1$};
	
	\node[blue,star]     at (-1.5,0) {};
	\node[right]     at (-1.5,0) {$\lambda_2$};
	
	\node[blue,star]     at (0,0) {};
	\node[above]     at (0,0) {0};
	
	\node[blue,star]     at (1.5,0) {};
	\node[left]     at (1.5,0) {$\lambda_3$};
	
	\node[blue,star]     at (3,0) {};
	\node[right]     at (3,0) {$\lambda_4$};
	
	\draw[blue,thick] (-3,0) -- (-1.5,0);
	\draw[blue,thick] (1.5,0) -- (3,0);
	\draw[blue,thick] (0,0) circle (1.5cm);
	
	\draw[blue,thick] (-3,0) -- (-4.5,1.5);
	\node[blue,star]     at (-4.5,1.5) {};
	\node[above left]     at (-4.5,1.5) {$\infty_2$};
	
	\draw[blue,thick] (3,0) -- (4.5,1.5);
	\node[blue,star]     at (4.5,1.5) {};
	\node[above right]     at (4.5,1.5) {$\infty_1$};
	
	\draw[blue,thick] (-3,0) -- (-4.5,-1.5);
	\node[blue,star]     at (-4.5,-1.5) {};
	\node[below left]     at (-4.5,-1.5) {$\infty_3$};
	
	\draw[blue,thick] (3,0) -- (4.5,-1.5);
	\node[blue,star]     at (4.5,-1.5) {};
	\node[below right]     at (4.5,-1.5) {$\infty_4$};

	\end{tikzpicture} 
	\caption{The Stokes complex $\mathcal{C}(0,0)$ as an embedded graph in $\mathbb{C}_\infty^*$.} \label{fig:complex00}
\end{figure}
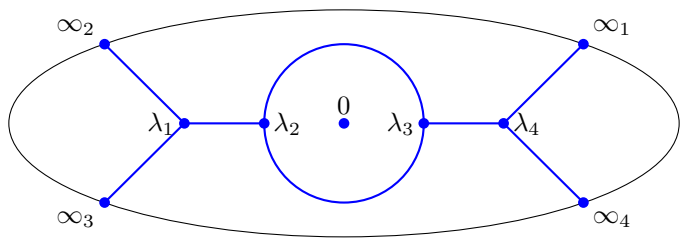

As an example, let us consider the Stokes complex corresponding to the potential \eqref{eq:potential} with
$(\alpha,\beta)=(0,0)$. It has real turning points
\begin{align}
&\lambda_1=-\sqrt{\tfrac{1}{2}(1+\sqrt{1-\nu^2})},\quad \lambda_2=-\sqrt{\tfrac{1}{2}(1-\sqrt{1-\nu^2})},
\label{eq:turning}\\ &\lambda_3=+\sqrt{\tfrac{1}{2}(1-\sqrt{1-\nu^2})},\quad \lambda_4=+\sqrt{\tfrac{1}{2}(1+\sqrt{1-\nu^2})},\nonumber
\end{align}
and a topological representation in $\mathbb{C}$ of the Stokes complex is given in Figure \ref{fig:approximate}, a topological representation as embedded graph in $\mathbb{C}_\infty^*$ is given in Figure \ref{fig:complex00} and a numerical plot is given in Figure \ref{fig:complex00anti}.

\subsubsection{Boutroux Curves}
\label{append:boutroux}
We introduce the notion of Boutroux curves \cite{boutroux,piwkb,bertola2011,bertola2009} tailored to our setting.
Recall from Section \ref{sec:results} the definition of the elliptic curve $\widehat{\Gamma}$, as the compactification
of the affine curve $\Gamma=\lbrace y^2=\la^2 V(\la) \rbrace $, and of the meromorphic differential $\omega=\sqrt{V(\la)}d\la$.
\begin{defi}\label{def:boutroux}
	We call the elliptic curve $\widehat{\Gamma}$ a
	Boutroux curve if $\Re\oint_\gamma\omega=0$ for
	any closed cycle $\gamma$ in $\widehat{\Gamma}$. We denote by $\Omega$ the set of all
	$(\alpha,\beta)\in\mathbb{C}^2$ such that $\widehat{\Gamma}(\alpha,\beta)$ is a Boutroux curve.
\end{defi}
Boutroux curves are highly relevant in the subsequent analysis in the section, as $R\subseteq \Omega$, by Definition 
\ref{defi:RK} of the region $R$.
We collect below some basic results on Boutroux curves and of their deformations which are
needed in the process of pursuing Steps 1-7.

We have the following equivalent characterisations of Boutroux curves.
\begin{lem}\label{lem:boutrouxclas}
	Let $(\alpha,\beta)\in \mathbb{C}^2$, then the following are equivalent:
	\begin{enumerate}[label=(\roman*)]
		\item $\widehat{\Gamma}(\alpha,\beta)$ is a Boutroux curve,
		\item each external vertex of the Stokes complex $\mathcal{C}(\alpha,\beta)$ has valency one,
		\item the inner Stokes complex corresponding to $V(\lambda;\alpha,\beta)$, defined in Definition \ref{def:stokescomplex}, is connected.
	\end{enumerate}
\end{lem}
\begin{proof}
	Let us note that (iii) trivially implies (i) and it is easy to see that (ii) implies (iii) using the fact that the Stokes complex is connected, see Proposition \ref{pro:basic}. It remains to be shown that (i) implies (ii).
	
	Suppose $\widehat{\Gamma}(\alpha,\beta)$ is a Boutroux curve. By Proposition \ref{pro:basic}, each external vertex has valency at least one. Suppose the valency of any particular external vertex, say $e^{\frac{1}{4}(2k-1)\pi i}$, is greater than one. Then there exist distinct Stokes lines $l_1$ and $l_2$, both asymptotic to $e^{\frac{1}{4}(2k-1)\pi i}\mathbb{R}_+$, such that
	\begin{equation}\label{eq:realint}
	\Re\int_{\Lambda_1}^{\Lambda_2}\sqrt{V(\lambda)}d\lambda=0,
	\end{equation}
	for any choice of $\Lambda_1\in l_1$, $\Lambda_2\in l_2$ and choice of connecting contour. Let $C,R>0$ be such that
	\begin{equation*}
	|\sqrt{V(\lambda)}-\lambda|\leq C,
	\end{equation*}
	for $|\lambda|>R$. Then, for any $\Lambda_1\in l_1$ and $\Lambda_2\in l_2$ with $|\Lambda_1|,|\Lambda_2|>R$, equation \eqref{eq:realint} with a choice of contour lying in $\{|\lambda|>R\}$, implies
	\begin{equation}\label{eq:inequality}
	|\Re(\Lambda_2^2-\Lambda_1^2)|\leq 2C|\Lambda_2-\Lambda_1|.
	\end{equation}
	We choose $\Lambda_j=re^{i\theta_i(r)}\in l_j$ for $r\gg0$, with $\theta_j(r)\sim \tfrac{1}{4}(2k-1)\pi$ as $r\rightarrow +\infty$, for $j\in 1,2$. Note that $\theta_1(r)-\theta_2(r)\notin 2\pi \mathbb{Z}$ for any $r$.
	Now inequality \eqref{eq:inequality} translates to
	\begin{equation*}
	r\leq 2 C \frac{|e^{i(\theta_2(r)-\theta_1(r))}-1|}{|\cos{(2\theta_2(r))}-\cos{(2\theta_1(r))}|},
	\end{equation*}
	but the right-hand side converges to $C$ as $r\rightarrow +\infty$ and we have arrived at a contradiction. We infer that each external vertex must has valency one, which completes the proof of the lemma.
\end{proof}

We proceed with discussing deformations of Stokes complexes, and in particular Boutroux curves.
Firstly we consider generic deformations of Stokes lines in the following lemma.

\begin{lem}[deformations of Stokes lines]\label{lem:deformationlines}
	Let $(\alpha^*,\beta^*)\in \mathbb{C}^2$ and consider the Stokes complex of $V(\lambda;\alpha^*,\beta^*)$. Let $\mu^*$ be a turning point of $V(\lambda;\alpha^*,\beta^*)$ and let $l^*$ be a Stokes line with $\mu^*$ as one of its endpoints. Let $\gamma_0:\mathbb{R}_{\geq 0}\rightarrow \mathbb{C}$ be a homeomorphism onto $l^*\cup\{\mu^*\}$ such that its restriction to $\mathbb{R}_+$ is a diffeomorphism onto $l^*$. Then, for any $0<T<\infty$, there exists
	\begin{itemize}
		\item  a simply connected open neighbourhood $W \subseteq \mathbb{C}^2$ of $(\alpha^*,\beta^*)$,
		\item an analytic function $\mu(\alpha,\beta)$ on $W$ with $\mu(\alpha^*,\beta^*)=\mu^*$, and
		\item a smooth mapping $\gamma:W\times [0,T)\rightarrow \mathbb{C},(\alpha,\beta,t)\mapsto \gamma_{(\alpha,\beta)}(t)$,
	\end{itemize}
	such that $\gamma_{(\alpha^*,\beta^*)}([0,T))=\gamma_0([0,T))$, and, for all $(\alpha,\beta)\in W$,
	\begin{itemize}
		\item $\mu(\alpha,\beta)$ is a turning point of $V(\lambda;\alpha,\beta)$,
		\item $\gamma_{(\alpha,\beta)}(0)=\mu(\alpha,\beta)$, and
		\item $\gamma_{(\alpha,\beta)}:[0,T)\rightarrow \mathbb{C}$ is, when restricted to $(0,T)$, a diffeomorphism onto part of a unique Stokes line $l=l(\alpha,\beta)$ of $V(\lambda;\alpha,\beta)$.
	\end{itemize}	
	Furthermore, if $l^*$ is asymptotic to $e^{\frac{1}{4}(2k-1)\pi i}\mathbb{R}_+$ for some $k\in \mathbb{Z}_4$, then the above also holds for $T=+\infty$, so that $\gamma_{(\alpha,\beta)}:[0,\infty)\rightarrow \mathbb{C}$ is, when restricted to $(0,\infty)$, a  diffeomorphism onto a unique Stokes line $l=l(\alpha,\beta)$ of $V(\lambda;\alpha,\beta)$ asymptotic to $e^{\frac{1}{4}(2k-1)\pi i}\mathbb{R}_+$, for all $(\alpha,\beta)\in W$.
\end{lem}
\begin{proof}
The proof is a straightforward but technical exercise in analysis. We leave it to the interested reader.
\end{proof}

\begin{lem}[generic deformations of Boutroux curves]\label{lem:deformboutroux}
	Let $T$ be a both locally and globally simply connected metric space, together with a continuous mapping
	\begin{equation*}
	(\alpha,\beta):T\rightarrow \mathbb{C}^2,t\mapsto (\alpha(t),\beta(t))
	\end{equation*}
	such that $(\alpha(t),\beta(t))\in \Omega$ for all $t\in T$. Further suppose that turning points of $V(\lambda;\alpha(t),\beta(t))$ do not merge or split on $T$, i.e. there exists an $m\in \{2,3,4\}$ and continuous functions $\mu_j:T\rightarrow \mathbb{C}$, for $1\leq j\leq m$, such that $\{\mu_j(t): 1\leq j\leq m\}$ are the turning points of $V(\lambda;\alpha(t),\beta(t))$ and $|\{\mu_j(t): 1\leq j\leq m\}|=m$ for all $t\in T$.\\	
	Then the isomorphism class of the Stokes complex $\mathcal{C}(\alpha(t),\beta(t))$ is constant on $T$.
\end{lem}
\begin{proof}
	Note that it is enough to show that, for every $t^*\in T$, there exists an open neighbourhood $B\subseteq T$ of $t^*$ such that the Stokes complex $\mathcal{C}(\alpha(t),\beta(t))$ is isomorphic to $\mathcal{C}(\alpha(t^*),\beta(t^*))$ for $t\in B$.
	
	Let $t^*\in T$ and write $(\alpha^*,\beta^*)=T(t^*)$. We denote $V^*(\lambda)=V(\lambda;\alpha^*,\beta^*)$ and $V_t(\lambda)=V(\lambda;\alpha(t),\beta(t))$ for $t\in T$.
	Let us pick one of the turning points, say $\mu_j(t)$, $1\leq j\leq m$. It is convenient to number the Stokes lines of $V_t(\lambda)$ emanating from $\mu_j(t)$ uniquely. Let $r_j\geq 1$ be the degree of the turning point and
	$l_1^*,\ldots,l_{r_j+2}^*$ be the Stokes lines of $V^*(\lambda)$ emanating from $\mu_j^*$\footnote{Here we allow for duplicity of a Stokes line if all its endpoints equal $\mu_j(t)$.}.
	We may choose a unique $\theta_1^*\in [0,\frac{1}{r_j+2}2\pi)$, write $\theta_s^*=\theta_1^*+\frac{s-1}{r_j+2}2\pi$, and order the Stokes lines such that $l_s^*$ emanates from $\mu_j^*$ with angle $\theta_s^*$ for $1\leq s\leq r_j+2$. It is now easy to see that there exists a unique continuous mapping $\theta_1: T\rightarrow \mathbb{R}$ with $\theta_1(t^*)=\theta_1^*$, so that, for every $t\in T$ and $1\leq s\leq r_j+2$, there exists a unique Stokes line $l_s(t)$ of $V_t(\lambda)$ which emanates from $\mu_j(t)$ with angle $\theta_s=\theta_1(t)+\frac{s-1}{r_j+2}2\pi$.
	We call $l_s(t)$ the continuous extension to $T$ of $l_s^*$ with respect to the turning point $\mu_j^*$ and angle $\theta_j^*$. Note that every Stokes line of $V^*(\lambda)$ has a unique continuous extension to $T$ with respect to every turning point from which it emanates for every of the angles by which it does so. In particular, a priori it is not excluded that an internal Stoke line might have two different continuous extensions to $T$.
	
	Let $1\leq s\leq r_j+2$. Suppose $l_s^*$ is asymptotic to $e^{\frac{1}{4}(2k-1)\pi i}\mathbb{R}_+$ for a $k\in\mathbb{Z}_4$, then we may apply Lemma \ref{lem:deformationlines} to find an open neighbourhood $B_{j,s}\subseteq T$ of $t^*$ such that the Stokes line $l_s(t)$ is asymptotic to $e^{\frac{1}{4}(2k-1)\pi i}\mathbb{R}_+$ for all $t\in B_{j,s}$.
	
	Suppose instead that $l_s^*$ is an internal Stokes line and has endpoints $\{\mu_j^*,\mu_k^*\}$. It might be that $\mu_j^*=\mu_k^*$, in which case there is an $1\leq s'\leq r_j+2$ with $s'\neq s$ such that $l_s^*=l_{s'}^*$. Regardless of this, let $\theta_\diamond^*\in [0,2\pi)$ by the angle by which $l_s^*$ emanates from $\mu_k^*$ and denote by $l_{\diamond}(t)$ the unique continuous extension of $l_s^*$ to $T$ with respect to $\mu_k^*$ and angle $\theta_\diamond^*$. We wish to show that there exists an open neighbourhood $B_{j,s}\subseteq T$ of $t^*$ such that $l_s(t)=l_\diamond(t)$, and thus $l_s(t)$ is an internal Stokes line with endpoints $\{\mu_j(t),\mu_k(t)\}$, for all $t\in B_{j,s}$.
	
	Assume, for the sake of contradiction, that such a set does not exist. Then there exists a sequence $(t_n)_{n\geq 1}$ in $T$ with $t_n\rightarrow t^*$ as $n\rightarrow \infty$ such that $l_s(t_n)\neq l_\diamond(t_n)$ for all $n\geq 1$. Fix a point $\Lambda\in l_s(t^*)=l_\diamond(t^*)=l_s^*$. Using Lemma \ref{lem:deformationlines}, it is easy to see that there exists a sequence $(\Lambda_n)_{n\geq 1}$ in $\mathbb{C}$ with $\Lambda_n\rightarrow \Lambda$ as $n\rightarrow \infty$ such that $\Lambda_n\in l_s(t_n)$ for $n\geq 1$. Similarly there exists a sequence $(\Lambda_n^\diamond)_{n\geq 1}$ in $\mathbb{C}$ with $\Lambda_n^\diamond\rightarrow \Lambda$ as $n\rightarrow \infty$ such that $\Lambda_n^\diamond\in l_\diamond(t_n)$ for $n\geq 1$.
	In fact, it is not difficult to see that we may choose these sequences such that $\Lambda_n^\diamond-\Lambda_n$ is approximately orthogonal to the Stokes line $l_s^*$ through $\Lambda$ for large $n$, namely
	\begin{equation}\label{eq:orthogonal}
	\Lambda_n^\diamond-\Lambda_n=\epsilon_n |\Lambda_n^\diamond-\Lambda_n| \overline{\sqrt{V^*(\Lambda)}}+o(\Lambda_n^\diamond-\Lambda_n)
	\end{equation}
	as $n\rightarrow \infty$, for some irrelevant choice of signs $\epsilon_n\in\{\pm 1\}$, $n\geq 1$.

	Since $(\alpha(t),\beta(t))\in \Omega$ for $t\in T$, we know, by Lemma \ref{lem:boutrouxclas}, that the internal Stokes complex of $V_{t_n}(\lambda)$ is connected  and in particular
	\begin{equation}\label{eq:vanish}
	\Re\int_{\Lambda_n}^{\Lambda_n^\diamond}{\sqrt{V_{t_n}(\lambda)}d\lambda}=0,
	\end{equation} 
	for every choice of contour avoiding critical points and $n\geq 1$. For large $n$, say $n\geq N\gg0$, we may choose as contour the line segment between $\Lambda_n$ and $\Lambda_n^\diamond$ and find a $c$, independent of $n$, such that
	\begin{equation*}
	\left|\int_{\Lambda_n}^{\Lambda_n^\diamond}{\sqrt{V_{t_n}(\lambda)}d\lambda}-(\Lambda_n^\diamond-\Lambda_n)\sqrt{V_{t_n}(\Lambda)}\right|\leq c|\Lambda_n^\diamond-\Lambda_n|^2
	\end{equation*}
	for $n\geq N$, so that equation \eqref{eq:vanish} yields
	\begin{equation*}
	|\Re(\Lambda_n^\diamond-\Lambda_n)\sqrt{V_{t_n}(\Lambda)}|\leq c|\Lambda_n^\diamond-\Lambda_n|^2
	\end{equation*}
	for $n\geq N$. Therefore, using \eqref{eq:orthogonal}, we obtain
	\begin{equation*}
	|\Re \overline{\sqrt{V_{0}(\Lambda)}}\sqrt{V_{t_n}(\Lambda)}|=o(1)
	\end{equation*}
	as $n\rightarrow \infty$. However, the left-hand side converges to $|V_{0}(\Lambda)|$ as $n\rightarrow \infty$ and therefore $V_{0}(\Lambda)=0$. Clearly $\Lambda$ is not a turning point of $V^*(\lambda)$ and we have arrived at a contradiction.
	We conclude that there exists an open neighbourhood $B_{j,s}\subseteq T$ of $t^*$ such that $l_s(t)=l_\diamond(t)$, and thus $l_s(t)$ is an internal Stokes line with endpoints $\{\mu_j(t),\mu_k(t)\}$, for all $t\in B_{j,s}$.
	
	Now
	\begin{equation*}
	B=\bigcap_{1\leq j\leq m,1\leq s\leq r_j+2} {B_{j,s}}
	\end{equation*}
	is an open neighbourhood of $t^*$ and the Stokes complex $\mathcal{C}(\alpha(t),\beta(t))$ is isomorphic to $\mathcal{C}(\alpha^*,\beta^*)$ for $t\in B$.	
\end{proof}

In the following two lemmas perturbations of singular Boutroux curves are discussed. In their formulation we make use of the discrimiant $\Delta$ of the polynomial
\begin{equation*}
\lambda^2V(\lambda;\alpha,\beta)=\lambda^4+2\alpha\lambda^3+(\alpha^2-1)\lambda^2-\beta \lambda+\tfrac{1}{4}\nu^2,
\end{equation*}
with respect to $\lambda$, explicitly
\begin{align}\label{eq:discriminant}
\Delta(\alpha,\beta)=&-27 \beta^4+4\alpha(9-\alpha^2)\beta^3+2\left(3\nu^2(5\alpha^2-6)+2(\alpha^2-1)^2\right)\beta^2\\
&+4\nu^2 \alpha\left(6\nu^2-(1-\alpha^2)(10-\alpha^2)\right)\beta+\nu^2\left(4\nu^4+\nu^2(\alpha^4-20\alpha^2-8)+4(1-\alpha^2)^3\right).\nonumber
\end{align}

\begin{lem}\label{lem:double}
	Let $(\alpha^*,\beta^*)\in \Omega$ be such that $V(\lambda;\alpha^*,\beta^*)$ has one double turning point and two simple turning points, say $\mu_1^*,\mu_2^*\in\mathbb{C}^*$. Choose a simply connected open neighbourhood $W\subseteq \mathbb{C}^2$ of $(\alpha^*,\beta^*)$ such that, for $j\in \{1,2\}$, there exists a unique analytic function $\mu_j:W\rightarrow \mathbb{C}^*$ with $\mu_j(\alpha^*,\beta^*)=\mu_j^*$, such that $\mu_j(\alpha,\beta)$ is a simple turning point of $V(\lambda;\alpha,\beta)$ for $(\alpha,\beta)\in W$.\\
	Then there exists an open neighbourhood $W_0\subseteq W$ of $(\alpha^*,\beta^*)$, such that, for all $(\alpha,\beta)\in W_0\cap \Omega$, if $\Delta(\alpha,\beta)\neq 0$, then the two turning points of $V(\lambda;\alpha,\beta)$, not equal to $\mu_1$ or $\mu_2$, are connected by a unique Stokes line which is homotopic to the straight line segment between the two in $\mathbb{C}^*$ minus critical points.
\end{lem}
\begin{proof}
	Note that, for $(\alpha,\beta)\in \Omega$, considering the Stokes $\mathcal{C}(\alpha,\beta)$, the sum of the valencies of the vertices equals twice the number of edges minus $4$. Using this identity in conjunction with Lemma \ref{lem:deformationlines}, the lemma follows quite directly. We leave the details to the reader.
\end{proof}

The following lemma is proven similarly to the above.
\begin{lem}\label{lem:triple}
	Let $(\alpha^*,\beta^*)\in \Omega$ be such that $V(\lambda;\alpha^*,\beta^*)$ has one triple turning point and one simple turning point, say $\mu_1^*\in\mathbb{C}^*$. Choose a simply connected open neighbourhood $W\subseteq \mathbb{C}^2$ of $(\alpha^*,\beta^*)$ such that there exists a unique analytic function $\mu_1:W\rightarrow \mathbb{C}^*$ with $\mu_1(\alpha^*,\beta^*)=\mu_1^*$, such that $\mu_1(\alpha,\beta)$ is a simple turning point of $V(\lambda;\alpha,\beta)$ for $(\alpha,\beta)\in W$.\\
	Then there exists an open neighbourhood $W_0\subseteq W$ of $(\alpha^*,\beta^*)$, such that, for all $(\alpha,\beta)\in W_0\cap \Omega$ not equal to $(\alpha^*,\beta^*)$, 
	\begin{itemize}
		\item if $\Delta(\alpha,\beta)\neq 0$, then the three simple turning points $\{\mu_2,\mu_3,\mu_4\}$ of $V(\lambda;\alpha,\beta)$, not equal to $\mu_1(\alpha,\beta)$, can be labelled such that $\mu_2$ and $\mu_3$
		are connected by a Stokes line which is homotopic to the straight line segment between the two in $\mathbb{C}^*$ minus critical points, and the same for the pair $\{\mu_3,\mu_4\}$;
		\item if $\Delta(\alpha,\beta)=0$, then the two turning points $\{\mu_2,\mu_3\}$, one simple and one double, of $V(\lambda;\alpha,\beta)$, not equal to $\mu_1(\alpha,\beta)$, 
		are connected by a Stokes line which is homotopic to the straight line segment between the two in $\mathbb{C}^*$ minus critical points,
	\end{itemize}
	where $\Delta(\alpha,\beta)$ denotes the discriminant \eqref{eq:discriminant}.
\end{lem}

\subsection{The Seven Steps}\label{subsec:results}
We start this subsection recalling the main objects related to the elliptic region.
Firstly, in the following definition, we introduce a new region $R_a$ and recall the definition of the sets $R$, $K$
and $K_a$ as well as the projection $\Pi_a$ defined in equation \eqref{eq:defiproj}.
\begin{def*}
By $R\subseteq \mathbb{C}^2$ we denote the set of all $(\alpha,\beta)\in \bb{C}^2$ for which the
	corresponding Stokes complex $\mc{C}(\alpha,\beta)$ is isomorphic to $\mc{C}(0,0)$. Furthermore,
	\begin{itemize}
		\item $K$ denotes the closure of $R$,
		\item $R_a$ denotes the projection of $R$ onto the $\alpha$-plane,
		\item $K_a$ denotes the projection of $K$ onto the
		$\alpha$-plane and $\Pi_a$ denotes the corresponding projection.
		\item $K_a$ is the elliptic region.
	\end{itemize}
\end{def*}

Let us recall the mapping $\mathcal{S}$ introduced in \eqref{def:S}, which is defined by means of the complete
elliptic integrals $s_{1,2}$
\footnote{For convenience of the reader we have collected some explicit formulas (that we will need in the sequel)
for the complete elliptic integrals $s_{1,2}$ and their
Jacobian in Appendix \ref{append:computation}.}
defined in equations \eqref{eq:defs1s2}. For $(\alpha,\beta)\in R$, the following Boutroux conditions hold,
\begin{equation}\label{eq:res1s2}
\Re s_1=0,\quad \Re s_2=0,
\end{equation}
due to the Stokes geometry of the potential on $R$ and thus $\mathcal{S}$ is indeed a real-valued mapping from $R$ into $\mathbb{R}^2$.

Next, we briefly discuss some symmetries of the region $R$ and mapping $\mathcal{S}$ in the following
\begin{lem}\label{lem:sym}
	The set $R$ is invariant under complex conjugation $(\alpha,\beta)\mapsto (\overline{\alpha},\overline{\beta})$
	and reflection $(\alpha,\beta)\mapsto (-\alpha,-\beta)$. Furthermore, for $(\alpha,\beta)\in R$,
	\begin{align*}
	\mathcal{S}(\overline{\alpha},\overline{\beta})&=\left(\mathcal{S}_1(\alpha,\beta),-\mathcal{S}_2(\alpha,\beta)\right),\\
	\mathcal{S}(-\alpha,-\beta)&=\left(-\mathcal{S}_1(\alpha,\beta),-\mathcal{S}_2(\alpha,\beta)\right),
	\end{align*}
	and in particular $\mathcal{S}(0,0)=(0,0)$.
\end{lem}
\begin{proof}
	Note that we have the following symmetries
	\begin{equation}\label{eq:csym}
	V(\lambda;\overline{\alpha},\overline{\beta})=\overline{V(\overline{\lambda};\alpha,\beta)},\quad
	V(\lambda;-\alpha,-\beta)=V(-\lambda;\alpha,\beta),
	\end{equation} from which the lemma easily follows.\end{proof}
	
It will be helpful to consider the following factorisation of the potential,
\begin{equation*}
\lambda^2V(\lambda;\alpha,\beta)=(\lambda-\lambda_1)(\lambda-\lambda_2)(\lambda-\lambda_3)(\lambda-\lambda_4),
\end{equation*}
where the turning points $\Lambda:=(\lambda_1,\lambda_2,\lambda_3,\lambda_4)$ are defined, up to permutations, by the following algebraic constraints
\begin{subequations}
	\label{eq:constraint}
	\begin{align}
	&\lambda_1\lambda_2\lambda_3\lambda_4=\tfrac{1}{4}\nu^2,\label{eq:constraint1}\\
	&\lambda_1\lambda_2+\lambda_1\lambda_3+\lambda_1\lambda_4+\lambda_2\lambda_3+
	\lambda_2\lambda_4+\lambda_3\lambda_4=\tfrac{1}{2}(\lambda_1^2+\lambda_2^2+\lambda_3^2+\lambda_4^2)-2,\label{eq:constraint2}\\
	&
	-\tfrac{1}{2}(\lambda_1+\lambda_2+\lambda_3+\lambda_4)=\alpha, \label{eq:alphaoflambdas}\\ \label{eq:betaoflambdas}
	&\tfrac{1}{4}\nu^2 (\lambda_1^{-1}+\lambda_2^{-1}+\lambda_3^{-1}+\lambda_4^{-1})= \beta.
	\end{align}
\end{subequations}
When $(\alpha,\beta)\in R$, we may unambiguously, i.e. not just up to permutation, define the turning points
via the Stokes complex and we will always do this in alignment with Figure \ref{fig:complex00}.

We proceed with discussing steps \ref{step1}-\ref{step7} in more detail.
\subsubsection*{Step \ref{step1}}
Let us take a point $(\alpha^*,\beta^*)\in R$. To understand what $R$ looks like in a
small neighbourhood of $(\alpha^*,\beta^*)$, we consider local deformations in $(\alpha,\beta)$ around this point. To this end, let us take a small open neighbourhood $U\times V\ni (\alpha^*,\beta^*)$ such that the turning points $\Lambda=\Lambda(\alpha,\beta)$ do not coalesce as we move $(\alpha,\beta)$ within in. Note that the complete elliptic integrals $s_{1,2}$ have a unique analytic extension to $U\times V$. The main insight on which much of the subsequent analysis relies, is that, locally within $U\times V$, the region $R$ is described exactly by the Boutroux conditions \eqref{eq:res1s2}.

Namely, as $\gamma_{1,2}$ form a basis of cycles, it is clear that, as we deform $(\alpha,\beta)$ in such a way that the Boutroux conditions \eqref{eq:res1s2} remain valid, the underlying elliptic curve remains a Boutroux curve, i.e.
\begin{equation*}
(U\times V)\cap \Omega = \{(\alpha,\beta)\in U\times V: \Re s_1=\Re s_2=0\},
\end{equation*}
where the set $\Omega$ of Boutroux curves is defined in Definition \ref{def:boutroux}. Furthermore, due to Lemma \ref{lem:deformboutroux},
we know that deformations of Boutroux curves preserve the underlying Stokes complex, as long as none of the turning points coalesce or split. Therefore, we have, by shrinking $U$ and $V$ if necessary,
\begin{equation}\label{eq:Rlocal}
(U\times V)\cap R= (U\times V)\cap \Omega = \{(\alpha,\beta)\in U\times V: \Re s_1=\Re s_2=0\}.
\end{equation}
This local characterisation of $R$ has two important consequences. Firstly, since $R$ is locally prescribed by
the vanishing of the real parts of two analytic functions, we can endow it with a geometric structure via an application
of the implicit function theorem. Secondly, as one approaches the border $\delta R=K\setminus R$ from within $R$,
local application of the implicit function theorem fails, which means that some of the turning points must coalesce there.
The first consequence is made precise in the following lemma.
\begin{lem}\label{lem:deformation}
	The set $R$ is a smooth $2$-dimensional regular submanifold of $\{(\Re{\alpha},\Im{\alpha},\Re{\beta},\Im{\beta})\in\mathbb{R}^4\}$
	and the mapping $\mathcal{S}$
	is locally diffeomorphic\footnote{In this paper, diffeomorphic always means $C^\infty$-diffeomorphic.} on $R$.
	Furthermore, for any given $(\alpha^*,\beta^*)\in R$, there exist simply connected open sets
	$U\subseteq\mathbb{C}$ and $V\subseteq\mathbb{C}$, containing respectively $\alpha^*$ and $\beta^*$,
	and a  diffeomorphism $B:U\rightarrow V$, such that
	\begin{equation*}
	R\cap(U\times V)=\{(\alpha,B(\alpha)):\alpha\in U\}.
	\end{equation*}
\begin{proof}
	As $(\alpha^*,\beta^*)\in R$, we know that the turning points of $V(\lambda;\alpha^*,\beta^*)$
	are all simple. Hence we can take small open discs $U_0\subseteq \mathbb{C}$ and $V_0\subseteq \mathbb{C}$,
	centred respectively at $\alpha^*$ and $\beta^*$, such that on $U_0\times V_0$ the turning points of
	$V(\lambda;\alpha,\beta)$ do not coalesce. On this set the turning points $\lambda_i=\lambda_i(\alpha,\beta)$
	are analytic in $(\alpha,\beta)$, $1\leq i \leq 4$, and the defining equations \eqref{eq:defs1s2} of
	$s_1(\alpha,\beta)$ and $s_2(\alpha,\beta)$ have unique analytic continuation to $U_0\times V_0$,
	the result of which we denote by $s_1(\alpha,\beta)$ and $s_2(\alpha,\beta)$ as well.
	The idea of the proof is to make use of the fact that $R$ is locally defined by equations \eqref{eq:res1s2}.
	
	Using the standard decomposition into real and imaginary part,
	\begin{equation*}
	\alpha=\alpha_R+i\alpha_I,\quad \beta=\beta_R+i\beta_I,\quad s_i=s_i^R+is_i^I\quad (i=1,2),
	\end{equation*}
	the analytic mapping
	\begin{equation*}
	(s_1,s_2):U_0\times V_0\rightarrow \mathbb{C}^2,(\alpha,\beta)\mapsto (s_1(\alpha,\beta),s_2(\alpha,\beta))
	\end{equation*}
	can be rewritten as a real smooth mapping 
	\begin{equation*}
	H:(\alpha_R,\alpha_I,\beta_R,\beta_I)\mapsto (s_1^R,s_1^I,s_2^R,s_2^I)
	\end{equation*}
	on an open neighbourhood of $(\alpha_R^*,\alpha_I^*,\beta_R^*,\beta_I^*)$.
	
	Due to equation \eqref{eq:jacobiandet}, we know that $(s_1,s_2)$ is locally biholomorphic near
	$(\alpha^*,\beta^*)$ and thus $H$ is locally diffeomorphic near $(\alpha_R^*,\alpha_I^*,\beta_R^*,\beta_I^*)$.
	To prove the thesis, we proceed in showing that the mapping $(\beta_R,\beta_I)\mapsto (s_1^R,s_2^R)$ is a
	local  diffeomorphism near $(\beta_R^*,\beta_I^*)$, for  $(\alpha_R,\alpha_I)=(\alpha_R^*,\alpha_I^*)$ fixed.

	To this end we apply the implicit function theorem to show that $\{\Re{\alpha},\Im{\alpha}\}$
	are good local coordinates on the zero set
	\begin{equation*}
	\Re{s_1(\alpha,\beta)}=0,\quad \Re{s_2(\alpha,\beta)}=0.
	\end{equation*}
	This requires that the Jacobian determinant
	\begin{equation*}
	\delta=\delta(\alpha,\beta)=\begin{vmatrix}
	\frac{\partial s_1^R}{\partial \beta_R} & \frac{\partial s_1^R}{\partial \beta_I}\\
	\frac{\partial s_2^R}{\partial \beta_R} & \frac{\partial s_2^R}{\partial \beta_I}
	\end{vmatrix}
	\end{equation*}
	does not vanish at $(\alpha,\beta)=(\alpha^*,\beta^*)$.
	By the Cauchy-Riemann equations, we have
	$\delta=\Im{\left[\frac{\partial s_1}{\partial \beta}\cdot \overline{\frac{\partial s_2}{\partial \beta}}\right]}$.
	Now, note that
	\begin{equation*}
	\frac{\partial s_1}{\partial \beta}=-\frac{1}{2}\int_{\gamma_1}\widetilde{\omega},\quad
	\frac{\partial s_2}{\partial \beta}=-\frac{1}{2}\int_{\gamma_2}\widetilde{\omega},
	\end{equation*}
	where $\widetilde{\omega}$ is the holomorphic differential form $\widetilde{\omega}=\frac{d\lambda}{y}$. As
	$\gamma_1$ and $\gamma_2$ are homologically independent, 	
	$(p_1,p_2):=(-2 \frac{\partial s_1}{\partial \beta},-2\frac{\partial s_2}{\partial \beta})$ form a pair of
	$\mathbb{R}$-linearly independent periods of the elliptic curve $\widehat{\Gamma}(\alpha,\beta)$. Therefore
	$p_1\overline{p}_2\notin \mathbb{R}$ and thus $\delta(\alpha,\beta)\neq 0$ on $U_0\times V_0$ and in particular at
	$(\alpha,\beta)=(\alpha^*,\beta^*)$\footnote{We note that this can also be proven using the explicit formulae in
	Appendix \ref{append:computation}.}.
	
	By the implicit function theorem, and the fact that $(s_1,s_2)$ is a local biholomorphism, there exist simply
	connected open sets $U\subseteq U_0$ and  $V\subseteq V_0$ with $\alpha^*\in U$ and $\beta^*\in V$, and a diffeomorphism
	$B:U\rightarrow V$ such that
	\begin{equation*}
	\{(\alpha,B(\alpha)):\alpha\in U\}=\{(\alpha,\beta)\in U\times V: \Re{s_1(\alpha,\beta)}=\Re{s_2(\alpha,\beta)}=0\}
	\end{equation*}
	and hence in particular $B(\alpha^*)=\beta^*$. 	
	Furthermore, note that
	\begin{equation*}
	\{(\alpha,B(\alpha)):\alpha\in U\}=\Omega\cap(U\times V).
	\end{equation*}
	
	Applying Lemma \ref{lem:deformboutroux}, it follows that the Stokes complex $\mathcal{C}(\alpha,B(\alpha))$
	is homeomorphic to $\mathcal{C}(\alpha^*,\beta^*)$ and hence $(\alpha,B(\alpha))\in R$ for all $\alpha\in U$.
	Therefore
	\begin{equation}\label{eq:ROmega}
	R\cap(U\times V)=\Omega\cap(U\times V)=\{(\alpha,B(\alpha)):\alpha\in U\}.
	\end{equation}
	
	Since $(s_1,s_2)$ is a local diffeomorphism, we may choose a simply connected open neighbourhood $W\subseteq \{(\Re{\alpha},\Im{\alpha},\Re{\beta},
	\Im{\beta})\in\mathbb{R}^4\}$ of $(\alpha^*,\beta^*)$, with $W\subseteq U\times V$, such that
	$H|_W:W\rightarrow \mathbb{R}^4$ is a  diffeomorphism onto its open image $H(W)\subseteq \mathbb{R}^4$.
	So $(W,H|_W)$ is a local chart of $\{(\Re{\alpha},\Im{\alpha},\Re{\beta},\Im{\beta})\in\mathbb{R}^4\}$ and
	\begin{equation*}
	H|_W(R\cap W)=H(W)\cap \{(\sigma_1^R,\sigma_1^I,\sigma_2^R,\sigma_2^I)\in\mathbb{R}^4:\sigma_1^R=\sigma_2^R=0\}.
	\end{equation*}
	We conclude that $R$ is a smooth 2-dimensional regular submanifold of
	$\{(\Re{\alpha},\Im{\alpha},\Re{\beta},\Im{\beta})\in\mathbb{R}^4\}$.
	Furthermore, $(R\cap W,(s_1^I,s_2^I))$ is a local chart of $R$, and as $\mathcal{S}=(s_1^I,s_2^I)$,
	we find that $\mathcal{S}$ is locally diffeomorphic on $R$.
\end{proof}
\end{lem}
As a corollary of this lemma we obtain that $R_a=\Pi_a(R)\subseteq \mathbb{C}$ is open and that the projection $\Pi_a$ 
is locally diffeomorphic in $R$.
Furthermore, the turning points $\lambda_k=\lambda_k(\alpha,\beta)$, $1\leq k\leq 4$,  are smooth functions on $R$.
\subsubsection*{Step \ref{step2}}
 To carry out step \ref{step2}, we have to continuously
extend the mapping $\mathcal{S}$ to $K$ and show that the range of the extension is a subset of the quadrilateral $Q$.
Now, a direct application of the residue theorem to the meromorphic differential $\omega$,
yields the following lemma.
\begin{lem}\label{lem:range}
	The image of $R$ under $\mathcal{S}$ is contained in the open rectangle
	\begin{equation*}
	Q^{\circ}:=\left(-\tfrac{1}{2}(1-\nu)\pi,+\tfrac{1}{2}(1-\nu)\pi\right)\times \left(-\nu\pi,+\nu\pi\right).
	\end{equation*}
	
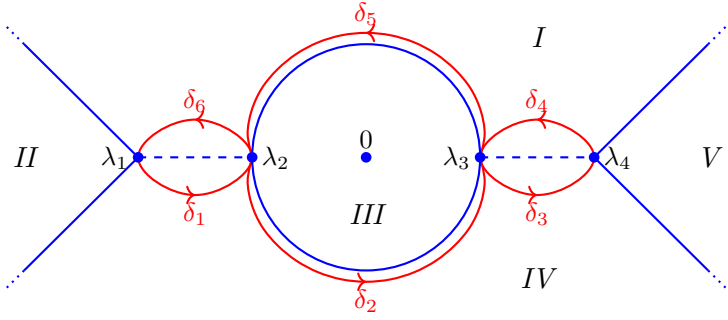
\begin{figure}[htpb]
	\centering
	\begin{tikzpicture}[scale=1.0]

	\draw [red,thick,decoration={markings, mark=at position 0 with {\arrow{>}}},decoration={markings, mark=at position 0.5 with {\arrow{>}}},
	postaction={decorate}] plot [smooth cycle,tension=0.8] coordinates {(1.5,0) (2.25,-0.5)  (3, 0) (2.25,0.5) };	
	
	\draw [red,thick,decoration={markings, mark=at position 0 with {\arrow{>}}},decoration={markings, mark=at position 0.5 with {\arrow{>}}},
	postaction={decorate}] plot [smooth cycle,tension=0.8] coordinates {(-3,0) (-2.25,-0.5)  (-1.5, 0) (-2.25,0.5) };
	
	\draw [red,thick,decoration={markings, mark=at position 0.5 with {\arrow{>}}},
	postaction={decorate}] plot [smooth,tension=0.6] coordinates 
	{(-1.5,0) 
		( {1.59*cos(-168.75)} , {1.59*sin(-168.75)} )
		( {1.63*cos(-157.5)} , {1.63*sin(-157.5)} ) 
		( {1.65*cos(-146.25)} , {1.65*sin(-146.25)} ) 
		( {1.65*cos(-135)} , {1.65*sin(-135)} ) 
		( {1.65*cos(-123.75)} , {1.65*sin(-123.75)} ) 
		( {1.65*cos(-112.5)} , {1.65*sin(-112.5)} )
		( {1.65*cos(-101.25)} , {1.65*sin(-101.25)} )
		(0,-1.65)
		( {1.65*cos(-78.75)} , {1.65*sin(-78.75)} )
		( {1.65*cos(-67.5)} , {1.65*sin(-67.5)} )
		( {1.65*cos(-56.25)} , {1.65*sin(-56.25)} )
		( {1.65*cos(-45)} , {1.65*sin(-45)} )
		( {1.65*cos(-33.75)} , {1.65*sin(-33.75)} )
		( {1.63*cos(-22.5)} , {1.63*sin(-22.5)} )
		( {1.59*cos(-11.25)} , {1.59*sin(-11.25)} )
		(1.5,0) };
	
	\draw [red,thick,decoration={markings, mark=at position 0.5 with {\arrow{>}}},
	postaction={decorate}] plot [smooth,tension=0.6] coordinates 
	{(1.5,0) 
		( {1.59*cos(-168.75+180)} , {1.59*sin(-168.75+180)} )
		( {1.63*cos(-157.5+180)} , {1.63*sin(-157.5+180)} ) 
		( {1.65*cos(-146.25+180)} , {1.65*sin(-146.25+180)} ) 
		( {1.65*cos(-135+180)} , {1.65*sin(-135+180)} ) 
		( {1.65*cos(-123.75+180)} , {1.65*sin(-123.75+180)} ) 
		( {1.65*cos(-112.5+180)} , {1.65*sin(-112.5+180)} )
		( {1.65*cos(-101.25+180)} , {1.65*sin(-101.25+180)} )
		(0,1.65)
		( {1.65*cos(-78.75+180)} , {1.65*sin(-78.75+180)} )
		( {1.65*cos(-67.5+180)} , {1.65*sin(-67.5+180)} )
		( {1.65*cos(-56.25+180)} , {1.65*sin(-56.25+180)} )
		( {1.65*cos(-45+180)} , {1.65*sin(-45+180)} )
		( {1.65*cos(-33.75+180)} , {1.65*sin(-33.75+180)} )
		( {1.63*cos(-22.5+180)} , {1.63*sin(-22.5+180)} )
		( {1.59*cos(-11.25+180)} , {1.59*sin(-11.25+180)} )
		(-1.5,0) };
	
	\node[red]     at (-2.25,-0.75) {$\delta_1$};
	\node[red]     at (0,-1.9) {$\delta_2$};
	\node[red]     at (2.25,-0.75) {$\delta_3$};
	\node[red]     at (2.25,0.75) {$\delta_4$};
	\node[red]     at (0,1.9) {$\delta_5$};
	\node[red]     at (-2.25,0.75) {$\delta_6$};

	\tikzstyle{star}  = [circle, minimum width=4pt, fill, inner sep=0pt];
	\node[blue,star]     at (-3,0){};
	\node[left]     at (-3,0) {$\lambda_1$};
	
	\node[blue,star]     at (-1.5,0) {};
	\node[right]     at (-1.5,0) {$\lambda_2$};
	
	\node[blue,star]     at (0,0) {};
	\node[above]     at (0,0) {$0$};
	
	\node[blue,star]     at (1.5,0) {};
	\node[left]     at (1.5,0) {$\lambda_3$};
	
	\node[blue,star]     at (3,0) {};
	\node[right]     at (3,0) {$\lambda_4$};
	
	\draw[dashed,blue,thick] (-3,0) -- (-1.5,0);
	\draw[dashed,blue,thick] (1.5,0) -- (3,0);
	\draw[blue,thick] (0,0) circle (1.5cm);
	
	\draw[blue,thick] (-3,0) -- (-4.5,1.5);
	\draw[blue,thick,dotted] (-4.5,1.5) -- (-4.75,1.75);
	
	\draw[blue,thick] (3,0) -- (4.5,1.5);
	\draw[blue,thick,dotted] (4.5,1.5) -- (4.75,1.75);
	
	\draw[blue,thick] (-3,0) -- (-4.5,-1.5);
	\draw[blue,thick,dotted] (-4.5,-1.5) -- (-4.75,-1.75);
	
	\draw[blue,thick] (3,0) -- (4.5,-1.5);
	\draw[blue,thick,dotted] (4.5,-1.5) -- (4.75,-1.75);
	
	\node     at (-4.5,0) {{\large II}};
	\node     at (4.5,0) {{\large V}};
	\node     at (2.25,1.6) {{\large I}}; 
	\node     at (2.25,-1.6) {{\large IV}}; 
	\node     at (0,-0.75) {{\large III}};

	\end{tikzpicture}
	\caption{Schematic representation of contours used in the proof of Lemma \ref{lem:range} in red,
		with the Stokes complex in blue, and the regions $I,\ldots,V$ used in the proof of Proposition \ref{prop:Sextension}.}
	\label{fig:contours}
\end{figure}

\begin{proof}
	The meromorphic differential $\omega$ has four poles,
	$0_{\pm}:= (\la=0,y=\pm \frac{\nu}{2})$ and $\infty_{\pm}$, where $\infty_{\pm}$ is characterised by $y\sim \pm \lambda^2$ as $(\lambda,y)\rightarrow \infty_{\pm}$, so that $\infty_+$ lies in the lower sheet of the elliptic curve $\widehat{\Gamma}$ and $\infty_-$ lies in the upper sheet. Their corresponding residues are easily computed,
	\begin{equation}
	\res_{(\lambda,y)=0_\pm}\omega=\pm \frac{\nu}{2}, \quad \res_{(\lambda,y)=\infty_{\pm}}\omega=\mp \frac12 \; .
	\end{equation}
	Consider the oriented contours $\delta_k$, $1\leq k\leq 6$,
	defined as in Figure \ref{fig:contours}, in which the Stokes complex is depicted in blue and the dashed Stokes lines between $\lambda_1$ and
	$\lambda_2$ and $\lambda_3$ and $\lambda_4$ act as branch cuts, such that all
	the contours $\delta_k$ lie in the upper sheet of the elliptic curve $\widehat{\Gamma}$, characterised by
	$\omega\sim \tfrac{1}{2}\nu \lambda^{-1}d\lambda$ as $\lambda\rightarrow 0$.	
	We define
	\begin{equation}\label{eq:defir}
	r_i=-i\int_{\delta_i}{\omega}\quad (1\leq i\leq 6),
	\end{equation}
	all of which are real and positive due to the Stokes geometry of the potential. Furthermore, $r_1=r_6$ and $r_3=r_4$.
	
	Since $0_+$ and $\infty_-$ lie in the same sheet as the above contours, the residue theorem yields 
	\begin{align*}
	r_2+r_5&=2\pi\res_{(\lambda,y)=0_+}{\omega}=\pi\nu,\\
	r_1+r_2+r_3+r_4+r_5+r_6&=2\pi \res_{(\lambda,y)=\infty_-}{\omega} =\pi.
	\end{align*}
	It follows in particular that
	\begin{equation}\label{eq:r4r6}
	r_4+r_6=\tfrac{1}{2}(1-\nu)\pi.
	\end{equation}
	Note that $s_2(\alpha,\beta)=(r_5-r_2)i$
	and since $r_2+r_5=\pi\nu$ with $r_2,r_5>0$, we have
	\begin{equation*}
	-\pi\nu<-is_2(\alpha,\beta)< +\pi\nu.
	\end{equation*}
	Similarly, $s_1(\alpha,\beta)=-2r_4i+\tfrac{1}{2}(1-\nu)\pi i$
	and it follows from equation \eqref{eq:r4r6} that $0<r_4<\tfrac{1}{2}(1-\nu)\pi$, so
	\begin{equation*}
	-\tfrac{1}{2}(1-\nu)\pi<-is_1(\alpha,\beta)< +\tfrac{1}{2}(1-\nu)\pi,
	\end{equation*}
	which finishes the proof of the lemma.
\end{proof}
\end{lem}
Due to the above lemma, the range of any continuous extension of $\mathcal{S}$ to $K=\overline{R}$ must therefore be a subset of $Q$. This means we merely have to prove that such an extension exists to finalise step \ref{step2}, since uniqueness is guaranteed.

 To extend the mapping $\mathcal{S}$ to $K$, we have to analyse the border $\delta R=K\setminus R$ of $R$ and understand what happens to the turning points $\Lambda$, the cycles $\gamma_{1,2}$ and ultimately the complete elliptic integrals $s_{1,2}$ upon approaching it. Let us stipulate that $\delta R$ is not the topological boundary $\partial R$ of $R$, and that, by definition, $K=R\sqcup \delta R$.

As we have already alleged to below equation \eqref{eq:Rlocal}, the border $\delta R$ is characterised by the merging of turning points. Since the set of Boutroux curves $\Omega$ is topologically closed, the border $\delta R$ thus corresponds to singular Boutroux curves. Now, the possible mergers of turning points as one approaches the border are heavily constrained, since two turning points can only merge if they are connected by a Stokes line, due to Lemma \ref{lem:double}, and Lemma \ref{lem:triple} gives a similar constraint for the merging of three turning points. 

\begin{figure}[htpb]
	\centering
	\begin{tabular}{ c | c | c }
		\begin{tikzpicture}[scale=0.5]
		\node     at (2.5,1.25) {$C_2$};

		\tikzstyle{star}  = [circle, minimum width=4pt, fill, inner sep=0pt];
		\node[blue,star]     at (3,0){};
		\node[right]     at (3,0) {$\eta_1$};
		
		\node[blue,star]     at (0,0) {};
		\node[above]     at (0,0) {0};
		
		\node[blue,star]     at (-0.75,-0.75) {};
		\node[left]     at (-0.78,-0.75) {$\eta_2$};

		\draw[blue,thick] (3,0) -- (4.75,1.75);

		\draw[blue,thick] (3,0) -- (4.75,-1.75);
		
		\draw [blue,thick] plot [smooth,tension=0.6] coordinates {(3,0) (2,-0.1) (0.6,-1.4) (-0.75,-0.75)  };
		
		\draw [blue,thick] plot [smooth,tension=0.9] coordinates {(-0.75,-0.75) (-0.5,0.6) (0.85,0.85) (0.5,-0.5) (-0.75,-0.75) };
		
		\draw [blue,thick] plot [smooth,tension=0.8] coordinates {(-0.75,-0.75) (-1,-0.5) (-1,0.8) (-1.75,1.75) };
		\draw [blue,thick] plot [smooth,tension=0.8] coordinates {(-0.75,-0.75) (-1.75,-1.75) };

		\end{tikzpicture} 
		& 
		\begin{tikzpicture}[scale=0.5]
		\node     at (-2.45,1.125) {$E_2$};

		\tikzstyle{star}  = [circle, minimum width=4pt, fill, inner sep=0pt];
		\node[blue,star]     at (-3,0){};
		\node[left]     at (-3,0) {$\mu_1$};
		
		\node[blue,star]     at (0,-1) {};
		\node[below]     at (0,-1) {$\mu_3$};
		
		\node[blue,star]     at (0,0) {};
		\node[above]     at (0,0) {0};
		
		\node[blue,star]     at (3,0) {};
		\node[right]     at (3,0) {$\mu_2$};
		
		\draw [blue,thick] plot [smooth,tension=0.6] coordinates {(-3,0) (-2.5,-0.2) (-1,-1.3) (0,-1)  };
		\draw [blue,thick] plot [smooth,tension=0.6] coordinates {(0,-1) (1,0) (0,1) (-1,0) (0,-1)  };
		
		\draw [blue,thick] plot [smooth,tension=0.6] coordinates {(3,0) (2.5,-0.2) (1,-1.3) (0,-1)  };
		
		\draw[blue,thick] (-3,0) -- (-4.75,1.75);
		
		\draw[blue,thick] (3,0) -- (4.75,1.75);
		
		\draw[blue,thick] (-3,0) -- (-4.75,-1.75);
		
		\draw[blue,thick] (3,0) -- (4.75,-1.75);
		
		\draw [red,thick,decoration={markings, mark=at position 0.65 with {\arrow{>}}},
		postaction={decorate}] plot [smooth cycle,tension=0.7] coordinates {(4.3,0) (3,-0.7) 
			(2,-1.5) (1,-1.5) (0,-1) (1,0) (2,0.3) (3,0.7)  };	
		
		\node[red]     at (2.35,1) {$\boldsymbol{\gamma_1}$};
		
		\end{tikzpicture}
		
		& 
		\begin{tikzpicture}[scale=0.5]
		\node     at (-2.5,1.25) {$C_1$};

		\tikzstyle{star}  = [circle, minimum width=4pt, fill, inner sep=0pt];
		\node[blue,star]     at (-3,0){};
		\node[left]     at (-3,0) {$\eta_1$};
		
		\node[blue,star]     at (0,0) {};
		\node[above]     at (0,0) {0};
		
		\node[blue,star]     at (0.75,-0.75) {};
		\node[right]     at (0.78,-0.75) {$\eta_2$};

		\draw[blue,thick] (-3,0) -- (-4.75,1.75);

		\draw[blue,thick] (-3,0) -- (-4.75,-1.75);
		
		\draw [blue,thick] plot [smooth,tension=0.6] coordinates {(-3,0) (-2,-0.1) (-0.6,-1.4) (0.75,-0.75)  };
		
		\draw [blue,thick] plot [smooth,tension=0.9] coordinates {(0.75,-0.75) (0.5,0.6) (-0.85,0.85) (-0.5,-0.5) (0.75,-0.75) };
		
		\draw [blue,thick] plot [smooth,tension=0.8] coordinates {(0.75,-0.75) (1,-0.5) (1,0.8) (1.75,1.75) };
		\draw [blue,thick] plot [smooth,tension=0.8] coordinates {(0.75,-0.75) (1.75,-1.75) };

		\end{tikzpicture}  \\ \hline
		\begin{tikzpicture}[scale=0.5]
		\node     at (2.5,1.25) {$E_3$};
		\draw[red,thick,decoration={markings, mark=at position 0.5 with {\arrow{>}}},postaction={decorate}]
		(1.75,0) arc (0:180:1.625cm and 1.75cm);
		\draw[red,thick,dashed] (-1.5,0) arc (-180:0:1.625cm and 1.75cm);

		\tikzstyle{star}  = [circle, minimum width=4pt, fill, inner sep=0pt];

		\node[blue,star]     at (-1.5,0) {};
		\node[left]     at (-1.5,0) {$\mu_3$};
		
		\node[blue,star]     at (0,0) {};
		\node[above]     at (0,0) {0};
		
		\node[blue,star]     at (1.5,0) {};
		\node[left]     at (1.5,0) {$\mu_2$};
		
		\node[blue,star]     at (3,0) {};
		\node[right]     at (3,0) {$\mu_1$};
		
		\draw[blue,thick] (1.5,0) -- (3,0);
		\draw[blue,thick] (0,0) circle (1.5cm);
		
		\draw[blue,thick] (-1.5,0) -- (-3.25,1.75);
		
		\draw[blue,thick] (3,0) -- (4.75,1.75);
		
		\draw[blue,thick] (-1.5,0) -- (-3.25,-1.75);
		
		\draw[blue,thick] (3,0) -- (4.75,-1.75);

		\node[red]     at (0.05,2.1) {$\boldsymbol{\gamma_2}$}; 
		\end{tikzpicture}
		& 
		\begin{tikzpicture}[scale=0.5]
		\node     at (-2.5,1.25) {$G_0$};
		\draw[red,thick,decoration={markings, mark=at position 0.5 with {\arrow{>}}},postaction={decorate}] (1.75,0) arc (0:180:1.75cm);
		\draw[red,thick,dashed] (-1.75,0) arc (-180:0:1.75cm);

		\tikzstyle{star}  = [circle, minimum width=4pt, fill, inner sep=0pt];
		\node[blue,star]     at (-3,0){};
		\node[left]     at (-3,0) {$\lambda_1$};
		
		\node[blue,star]     at (-1.5,0) {};
		\node[right]     at (-1.5,0) {$\lambda_2$};
		
		\node[blue,star]     at (0,0) {};
		\node[above]     at (0,0) {0};
		
		\node[blue,star]     at (1.5,0) {};
		\node[left]     at (1.5,0) {$\lambda_3$};
		
		\node[blue,star]     at (3,0) {};
		\node[right]     at (3,0) {$\lambda_4$};
		
		\draw[blue,thick] (-3,0) -- (-1.5,0);
		\draw[blue,thick] (1.5,0) -- (3,0);
		\draw[blue,thick] (0,0) circle (1.5cm);
		
		\draw[blue,thick] (-3,0) -- (-4.75,1.75);
		
		\draw[blue,thick] (3,0) -- (4.75,1.75);
		
		\draw[blue,thick] (-3,0) -- (-4.75,-1.75);
		
		\draw[blue,thick] (3,0) -- (4.75,-1.75);
		
		\draw [red,thick,decoration={markings, mark=at position 0.25 with {\arrow{<}}},
		postaction={decorate}] (2.25,0) ellipse (1.85cm and 0.7cm)	;
		
		\node[red]     at (2.35,1.0) {$\boldsymbol{\gamma_1}$};
		\node[red]     at (0.05,2.1) {$\boldsymbol{\gamma_2}$}; 
		\end{tikzpicture}
		& 
		\begin{tikzpicture}[scale=0.5]
		\node     at (-2.5,1.25) {$E_1$};
		\draw[red,thick,decoration={markings, mark=at position 0.5 with {\arrow{>}}},postaction={decorate}] (1.5,0) arc (0:180:1.625cm and 1.75cm);
		\draw[red,thick,dashed] (-1.75,0) arc (-180:0:1.625cm and 1.75cm);

		\tikzstyle{star}  = [circle, minimum width=4pt, fill, inner sep=0pt];
		\node[blue,star]     at (-3,0){};
		\node[left]     at (-3,0) {$\mu_1$};
		
		\node[blue,star]     at (-1.5,0) {};
		\node[right]     at (-1.5,0) {$\mu_2$};
		
		\node[blue,star]     at (0,0) {};
		\node[above]     at (0,0) {0};
		
		\node[blue,star]     at (1.5,0) {};
		\node[right]     at (1.5,0) {$\mu_3$};

		\draw[blue,thick] (-3,0) -- (-1.5,0);
		\draw[blue,thick] (0,0) circle (1.5cm);
		
		\draw[blue,thick] (-3,0) -- (-4.75,1.75);
		
		\draw[blue,thick] (1.5,0) -- (3.25,1.75);
		
		\draw[blue,thick] (-3,0) -- (-4.75,-1.75);
		
		\draw[blue,thick] (1.5,0) -- (3.25,-1.75);

		\node[red]     at (0.05,2.1) {$\boldsymbol{\gamma_2}$}; 
		\end{tikzpicture}\\ \hline
		\begin{tikzpicture}[scale=0.5]
		\node     at (2.5,1.25) {$C_3$};

		\tikzstyle{star}  = [circle, minimum width=4pt, fill, inner sep=0pt];
		\node[blue,star]     at (3,0){};
		\node[right]     at (3,0) {$\eta_1$};
		
		\node[blue,star]     at (0,0) {};
		\node[above]     at (0,0) {0};
		
		\node[blue,star]     at (-0.75,0.75) {};
		\node[left]     at (-0.78,0.75) {$\eta_2$};

		\draw[blue,thick] (3,0) -- (4.75,1.75);

		\draw[blue,thick] (3,0) -- (4.75,-1.75);
		
		\draw [blue,thick] plot [smooth,tension=0.6] coordinates {(3,0) (2,0.1) (0.6,1.4) (-0.75,0.75)  };
		
		\draw [blue,thick] plot [smooth,tension=0.9] coordinates {(-0.75,0.75) (-0.5,-0.5) (0.85,-0.85) (0.5,0.6) (-0.75,0.75) };
		
		\draw [blue,thick] plot [smooth,tension=0.8] coordinates {(-0.75,0.75) (-1,0.5) (-1,-0.8) (-1.75,-1.75) };
		\draw [blue,thick] plot [smooth,tension=0.8] coordinates {(-0.75,0.75) (-1.75,1.75) };

		\end{tikzpicture} 
		& 
		\begin{tikzpicture}[scale=0.5]
		\node     at (-2.5,1.2) {$E_4$};

		\tikzstyle{star}  = [circle, minimum width=4pt, fill, inner sep=0pt];
		\node[blue,star]     at (-3,0){};
		\node[left]     at (-3,0) {$\mu_1$};
		
		\node[blue,star]     at (0,1) {};
		\node[above]     at (0,1) {$\mu_3$};
		
		\node[blue,star]     at (0,0) {};
		\node[above]     at (0,0) {0};
		
		\node[blue,star]     at (3,0) {};
		\node[right]     at (3,0) {$\mu_2$};
		
		\draw [blue,thick] plot [smooth,tension=0.6] coordinates {(-3,0) (-2.5,0.2) (-1,1.3) (0,1)  };
		\draw [blue,thick] plot [smooth,tension=0.6] coordinates {(0,1) (1,0) (0,-1) (-1,0) (0,1)  };
		
		\draw [blue,thick] plot [smooth,tension=0.6] coordinates {(3,0) (2.5,0.2) (1,1.3) (0,1)  };
		
		\draw[blue,thick] (-3,0) -- (-4.75,1.75);
		
		\draw[blue,thick] (3,0) -- (4.75,1.75);
		
		\draw[blue,thick] (-3,0) -- (-4.75,-1.75);
		
		\draw[blue,thick] (3,0) -- (4.75,-1.75);
		
		\draw [red,thick,decoration={markings, mark=at position 0.1 with {\arrow{<}}},
		postaction={decorate}] plot [smooth cycle,tension=0.7] coordinates {(4.3,0) (3,0.7)  (2,1.5) (1,1.5) (0,1) (1,0) (2,-0.3) (3,-0.7)  };	
		
		\node[red]     at (2.35,1.75) {$\boldsymbol{\gamma_1}$};
		
		\end{tikzpicture}		 
		&
		\begin{tikzpicture}[scale=0.5]
		\node     at (-2.5,1.25) {$C_4$};

		\tikzstyle{star}  = [circle, minimum width=4pt, fill, inner sep=0pt];
		\node[blue,star]     at (-3,0){};
		\node[left]     at (-3,0) {$\eta_1$};
		
		\node[blue,star]     at (0,0) {};
		\node[above]     at (0,0) {0};
		
		\node[blue,star]     at (0.75,0.75) {};
		\node[right]     at (0.78,0.75) {$\eta_2$};

		\draw[blue,thick] (-3,0) -- (-4.75,1.75);

		\draw[blue,thick] (-3,0) -- (-4.75,-1.75);
		
		\draw [blue,thick] plot [smooth,tension=0.6] coordinates {(-3,0) (-2,0.1) (-0.6,1.4) (0.75,0.75)  };
		
		\draw [blue,thick] plot [smooth,tension=0.9] coordinates {(0.75,0.75) (0.5,-0.5) (-0.85,-0.85) (-0.5,0.6) (0.75,0.75) };
		
		\draw [blue,thick] plot [smooth,tension=0.8] coordinates {(0.75,0.75) (1,0.5) (1,-0.8) (1.75,-1.75) };
		\draw [blue,thick] plot [smooth,tension=0.8] coordinates {(0.75,0.75) (1.75,1.75) };

		\end{tikzpicture} 
	\end{tabular}
	\caption{Topological representation of isomorphism classes $E_1,\ldots, E_4$ and $C_1,\ldots, C_4$ and
		$G_0$ of Stokes complexes in blue with corresponding contours used in the proof of Lemma \ref{lem:extension}
		in red. Note that $G_0$ is the isomorphism class of $\mathcal{C}(0,0)$.} \label{fig:stokes_complexes}
\end{figure}
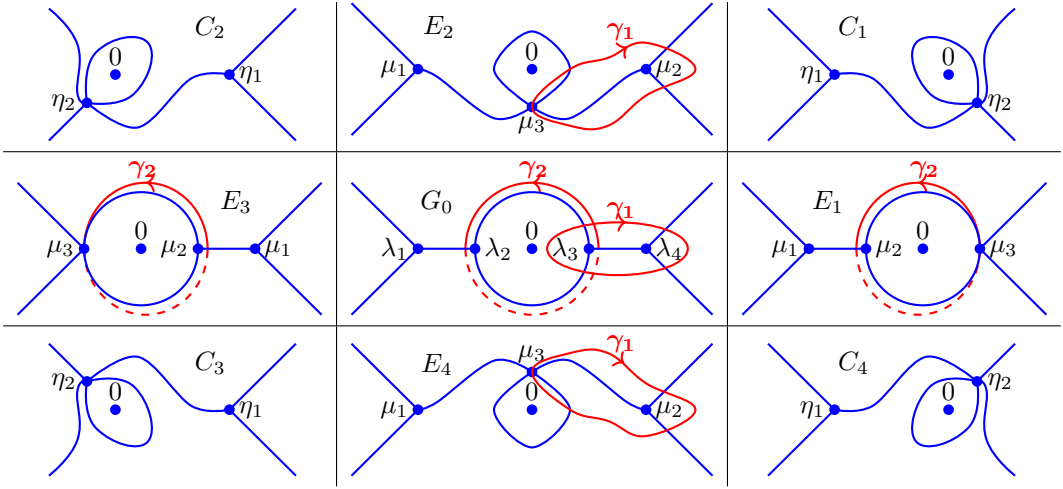

Consider Figure \ref{fig:stokes_complexes}, which defines 9 isomorphism classes of Stokes complexes, where we recognise the class $G_0$ as the one corresponding to points in $R$. An example of a legal merger is given by the merging of the turning points $\lambda_3$ and $\lambda_4$, leading to a Stokes complex in the class $E_1$. Similarly, there are two topologically inequivalent ways that the turning points $\lambda_2$ and $\lambda_3$ can merge, leading to Stokes complexes in the classes $E_2$ and $E_4$. Note that in the classes $E_1,\ldots E_4$, the Stokes complexes have one double turning point and two simple turning points. Similarly, in the classes $C_1,\ldots,C_4$, the Stokes complexes have one triple turning point and one remaining simple turning point.
The Stokes complexes in the classes $C_1,\ldots, C_4$ can be obtained through the merging of three turning points.

The following lemma asserts that, for each point on the border $\delta R$, the corresponding Stokes complex falls
in (precisely) one of these eight classes $E_1,\ldots, E_4,C_1,\ldots,C_4$. 
\begin{lem}\label{lem:partition}
	For any point $(\alpha,\beta)\in \delta R$, the corresponding elliptic curve $\widehat{\Gamma}(\alpha,\beta)$
	is a singular Boutroux curve and the associated Stokes complex $\mathcal{C}(\alpha,\beta)$ falls within (precisely)
	one of the classes $E_k,1\leq k\leq 4$ or $C_k,1\leq k\leq 4$, defined in Figure \ref{fig:stokes_complexes}.
\end{lem}
	\begin{proof}
		Firstly we show that $\delta R\subseteq \Omega$. Let $(\alpha^*,\beta^*)\in \delta R$, then there exists a sequence
		$(\alpha_n,\beta_n)_{n\geq 1}$ in $R$ such that $(\alpha_n,\beta_n)\rightarrow (\alpha^*,\beta^*)$ as
		$n\rightarrow \infty$. 	
		Take any Jordan curve $\gamma$ in $\mathbb{C}^*$ minus turning points of $V(\lambda;\alpha^*,\beta^*)$, which encircles an even number of turning points counting multiplicity, then $\sqrt{V(\lambda)}$ can be chosen single-valued along $\gamma$ and
		\begin{equation*}
		\Re\oint_{\gamma}\sqrt{V(\lambda;\alpha^*,\beta^*)}d\lambda=\lim_{n\rightarrow \infty} 
		\Re\oint_{\gamma}\sqrt{V(\lambda;\alpha_n,\beta_n)}d\lambda=0.
		\end{equation*}
		It follows that $(\alpha^*,\beta^*)\in \Omega$ and hence $\delta R\subseteq \Omega$.
		
		Next, we show that $\delta R\subseteq \{\Delta(\alpha,\beta)=0\}$.  Let $(\alpha^*,\beta^*)\in \delta R$ and suppose that $\Delta(\alpha^*,\beta^*)\neq 0$. The proof of
		equation \eqref{eq:ROmega} can be used line for line to show that there exists an open neighbourhood
		$W\subseteq \{(\Re{\alpha},\Im{\alpha},\Re{\beta},\Im{\beta})\in\mathbb{R}^4\}$ of $(\alpha^*,\beta^*)$ such that
		\begin{equation*}
		W\cap \Omega=\{(\alpha,\beta)\in W:\mathcal{C}(\alpha,\beta)\sim \mathcal{C}(\alpha^*,\beta^*)\}.
		\end{equation*}
		As $W\cap \Omega\cap R=W\cap R\neq \emptyset$, it follows that $\mathcal{C}(\alpha^*,\beta^*)\sim \mathcal{C}(0,0)$
		and hence $(\alpha^*,\beta^*)\in R$. But $R\cap \delta R=\emptyset$ and we have arrived at a contradiction.
		It follows that $\Delta(\alpha^*,\beta^*)=0$.	
		
		We conclude that
		\begin{equation*}
		\delta R\subseteq \Omega\cap \{\Delta(\alpha,\beta)=0\}.
		\end{equation*}
		Namely, take any $(\alpha^*,\beta^*)\in \delta R$, then $\widehat{\Gamma}(\alpha^*,\beta^*)$ is a singular Boutroux curve. Next, we want to prove that the Stokes complex $\mathcal{C}(\alpha^*,\beta^*)$ must fall in one of the eight isomorphism classes introduced in Figure \ref{fig:stokes_complexes}. To this end, take a sequence
		$(\alpha_n,\beta_n)_{n\geq 1}$ in $R$ such that $(\alpha_n,\beta_n)\rightarrow (\alpha^*,\beta^*)$ as
		$n\rightarrow \infty$. 
		
		Let $\{\mu_1,\ldots,\mu_m\}$ be the turning points of $V(\lambda;\alpha^*,\beta^*)$, discarding multiplicity, so that $1\leq m\leq 3$.
		It is easy to see by direct computation, that either $V(\lambda;\alpha^*,\beta^*)$ has one double turning point and two simple turning points, or it has one triple turning point and one simple turning point.
		Indeed, if $V(\lambda;\alpha^*,\beta^*)$ would have two double turning points or one quadruple turning point, then it is easy to see that $\nu^2=1$, which contradicts our assumption $0<\nu \leq \tfrac{1}{3}$. In particular, it follows that either $m=2$ or $m=3$.
		
		Let us denote the turning points of $V(\lambda;\alpha_n,\beta_n)$ by $\lambda_j^n=\lambda_j(\alpha_n,\beta_n)$
		for $1\leq j\leq 4$ and $n\geq 1$, then we may assume, by replacing $(\alpha_n,\beta_n)_{n\geq 1}$ by an appropriate
		subsequence if necessary, that there exists a surjective mapping $\sigma:\{1,2,3,4\}\rightarrow \{1,2,\ldots,m\}$
		such that $\lambda_j^n\rightarrow \mu_{\sigma(j)}$ as $n\rightarrow \infty$ for $1\leq j\leq 4$.
		
		Let us first consider the case $m=3$ and label the turning points of $V(\lambda;\alpha^*,\beta^*)$
		such that $\mu_1$ and $\mu_2$ are simple and $\mu_3$ is the double turning point. Applying Lemma \ref{lem:double}, either
		\begin{enumerate}[label=(\roman*)]
			\item $\lambda_1^n,\lambda_2^n\rightarrow \mu_3$ as $n\rightarrow \infty$,
			\item $\lambda_2^n,\lambda_3^n\rightarrow \mu_3$ as $n\rightarrow \infty$, or
			\item $\lambda_3^n,\lambda_4^n\rightarrow \mu_3$ as $n\rightarrow \infty$.
		\end{enumerate}
		In case (i), it follows from Lemma \ref{lem:deformationlines} that
		\begin{itemize}
			\item  $\mu_3$ has one emanating Stokes line asymptotic to $e^{+\frac{3}{4}\pi i}\mathbb{R}_+$ and
			another asymptotic to $e^{-\frac{3}{4}\pi i}\mathbb{R}_+$,
			\item $\mu_1$ or $\mu_2$ has one emanating Stokes line asymptotic to $e^{+\frac{1}{4}\pi i}\mathbb{R}_+$
			and another asymptotic to $e^{-\frac{1}{4}\pi i}\mathbb{R}_+$.
		\end{itemize}
		It is now easy to see that $\mathcal{C}(\alpha^*,\beta^*)\in E_3$, using Proposition \ref{pro:basic} and
		Lemma \ref{lem:boutrouxclas}.
		
		In case (iii) it follows analogously that $\mathcal{C}(\alpha^*,\beta^*)\in E_1$. 
		
		In case (ii), it follows from Lemma \ref{lem:deformationlines} that we may, if necessary,
		renumber $\{\mu_1,\mu_2\}$ such that
		\begin{itemize}
			\item  $\mu_1$ has one emanating Stokes line asymptotic to $e^{+\frac{3}{4}\pi i}\mathbb{R}_+$
			and another asymptotic to $e^{-\frac{3}{4}\pi i}\mathbb{R}_+$,
			\item $\mu_2$ has one emanating Stokes line asymptotic to $e^{+\frac{1}{4}\pi i}\mathbb{R}_+$
			and another asymptotic to $e^{-\frac{1}{4}\pi i}\mathbb{R}_+$.
		\end{itemize}
		It is now easy to see that either $\mathcal{C}(\alpha^*,\beta^*)\in E_2$ or $\mathcal{C}(\alpha^*,\beta^*)\in E_4$, using Proposition \ref{pro:basic} and Lemma \ref{lem:boutrouxclas}.
		
		Next, let us consider $m=2$ and label the turning points of $V(\lambda;\alpha^*,\beta^*)$ so
		that $\mu_1$ is simple and $\mu_2$ is the triple turning point. By Lemma \ref{lem:triple}, either
		\begin{enumerate}
			\item $\lambda_1^n,\lambda_2^n,\lambda_3^n\rightarrow \mu_2$ as $n\rightarrow \infty$, or
			\item $\lambda_2^n,\lambda_3^n,\lambda_4^n\rightarrow \mu_2$ as $n\rightarrow \infty$.
		\end{enumerate}
		Similarly as above, it follows in case (1) that
		either $\mathcal{C}(\alpha^*,\beta^*)\in C_2$ or $\mathcal{C}(\alpha^*,\beta^*)\in C_3$,
		and in case (2) that either $\mathcal{C}(\alpha^*,\beta^*)\in C_1$ or $\mathcal{C}(\alpha^*,\beta^*)\in C_4$.	
		The lemma follows.
	\end{proof}

Due to Lemma \ref{lem:partition}, the border $\delta R$ can be partitioned according to the isomorphism class of the corresponding Stokes complex. Let us accordingly define
\begin{equation*}
\widetilde{e}_k=\{(\alpha,\beta)\in \delta R:\mathcal{C}(\alpha,\beta)\in E_k\},
\quad \widetilde{c}_k=\{(\alpha,\beta)\in \delta R:\mathcal{C}(\alpha,\beta)\in C_k\}\quad (1\leq k\leq 4),
\end{equation*}
so that we obtain the partition
\begin{equation}\label{eq:disjoint}
\delta R=\widetilde{e}_1\sqcup \widetilde{e}_2\sqcup \widetilde{e}_3\sqcup 
\widetilde{e}_4\sqcup \widetilde{c}_1\sqcup \widetilde{c}_2\sqcup \widetilde{c}_3\sqcup \widetilde{c}_4.
\end{equation}
We remark that we have chosen the notation such that, for example, part $\widetilde{e}_1$ of the border $\delta R$  a posteriori corresponds to part $e_1$ of the boundary of the elliptic region $K_a$ via the projection $\Pi_a$, and such that it corresponds to the open edge $\widehat{e}_1$ of the boundary of the quadrilateral $Q$ via the mapping $\mathcal{S}$, see Figure \ref{fig:marked}.

Now, take a point $(\alpha^*,\beta^*)\in\delta R$ on the border, then it falls in one of the eight disjoint subsets on the right-hand side of \eqref{eq:disjoint}. This gives us the Stokes complex at this point and, more importantly, tells us what happens to the turning points $\Lambda$, cycles $\gamma_{1,2}$ and thus complete elliptic integrals $s_{1,2}$ corresponding to a point $(\alpha,\beta)\in R$ as it approaches $(\alpha^*,\beta^*)$. For example, it follows that, as one approaches $\widetilde{e}_1$ from within $R$, the turning points $\lambda_3$ and $\lambda_4$ merge. This means that $\gamma_2$ becomes a simple encircling of $\lambda_{3}^*=\lambda_4^*$ at $\widetilde{e}_1$ and thus $\int_{\gamma_2}\omega=0$, i.e. $s_1=\tfrac{1}{2}i\pi(1-\nu)$ on $\widetilde{e}_1$.
In this way, we obtain the unique continuous extension of $\mathcal{S}$ to the border $\delta R$, as further detailed
in the following lemma, which completes step \ref{step2}.
\begin{lem}\label{lem:extension}
	The mapping $\mathcal{S}$ on $R$ has a unique continuous extension to $K$,
	\begin{equation*}
	\mathcal{S}:K\rightarrow Q,\quad Q=\left[-\tfrac{1}{2}(1-\nu)\pi,+\tfrac{1}{2}(1-\nu)\pi\right]\times \left[-\nu\pi,+\nu\pi\right],
	\end{equation*} 
	with
	\begin{align*}
	\mathcal{S}(\widetilde{e}_k)\subseteq \widehat{e}_k,\quad \mathcal{S}(\widetilde{c}_k)\subseteq \{\widehat{c}_k\}\quad (1\leq k \leq 4),
	\end{align*}
	where the open edges and corners $\widehat{e}_k, \widehat{c}_k$, $1\leq k\leq 4$ of $Q$ are
	defined as in Figure \ref{fig:marked}.\\
	Furthermore, for $1\leq k\leq 4$, the `edge' $\widetilde{e}_k$ is a smooth $1$-dimensional
	regular submanifold of $\{(\Re{\alpha},\Im{\alpha},\Re{\beta},\Im{\beta})\in\mathbb{R}^4\}$ and the restricted mapping
	\begin{equation}\label{eq:edgediffeo}
	\mathcal{S}|_{\widetilde{e}_k}:\widetilde{e}_k\rightarrow \widehat{e}_k
	\end{equation} is locally diffeomorphic.
\begin{proof}
	Let us denote by $l_+$ and $l_-$ respectively the upper and lower Stokes line with respect to $0$,
	as depicted in Figure \ref{fig:complex00}, with endpoints $\lambda_2$ and $\lambda_3$. Then we can characterise the separate parts of the border in the disjoint union \eqref{eq:disjoint} as follows.

		Near the edges $\widetilde{e}_k$, $1\leq k\leq 4$, two of the turning points of $V(\lambda)$ merge,
		\begin{itemize}
			\item[$\widetilde{e}_1$:] here $\lambda_3$ and $\lambda_4$ merge,
			\item[$\widetilde{e}_2$:] here $\lambda_2$ and $\lambda_3$ merge via the vanishing of $l_-$,
			\item[$\widetilde{e}_3$:] here $\lambda_2$ and $\lambda_3$ merge via the vanishing of $l_+$,
			\item[$\widetilde{e}_4$:] here $\lambda_1$ and $\lambda_2$ merge.
		\end{itemize} 
		Similarly, the corners $\widetilde{c}_k$ are characterised by the merging of three turning points,
		\begin{itemize}
			\item[$\widetilde{c}_1$:] here $\lambda_2$, $\lambda_3$ and $\lambda_4$ merge via the vanishing of $l_-$,
			\item[$\widetilde{c}_2$:]  here $\lambda_1$, $\lambda_2$ and $\lambda_3$ merge via the vanishing of $l_-$,
			\item[$\widetilde{c}_3$:]  here $\lambda_1$, $\lambda_2$ and $\lambda_3$ merge via the vanishing of $l_+$,
			\item[$\widetilde{c}_4$:]  here $\lambda_2$, $\lambda_3$ and $\lambda_4$ merge via the vanishing of $l_+$.
		\end{itemize}
		Indeed, let us for example prove the characterisation of $\widetilde{e}_1$.
		Let $(\alpha^*,\beta^*)\in \widetilde{e}_1$ and let $(\alpha_n,\beta_n)_{n\geq 1}$ be a sequence in $R$ which converges to $(\alpha^*,\beta^*)$. Let $\lambda_k^n=\lambda_k(\alpha_n,\beta_n)$ denote the turning points of $V(\lambda;\alpha_n,\beta_n)$
		for $1\leq k\leq 4$ and $n\geq 1$, and $\mu_1,\mu_2,\mu_3$ denote the turning points of $V(\lambda;\alpha^*,\beta^*)$ as in Figure \ref{fig:stokes_complexes} for the class $E_1$.
		
		We have to show that $\lambda_k^n\rightarrow \mu_{\sigma_0(k)}$ for $1\leq k\leq 4$, where $\sigma_0(k)=k$ for $1\leq k\leq 3$ and $\sigma_0(4)=3$. Suppose this is not the case, then there exists a subsequence $(\alpha_{n_m},\beta_{n_m})_{m\geq 1}$
		together with a surjective mapping $\sigma:\{1,2,3,4\}\rightarrow \{1,2,3\}$ not equal to $\sigma_0$,
		such that $\lambda_k^{m_n}\rightarrow \mu_{\sigma(k)}$ as $n\rightarrow \infty$ for $1\leq k\leq 4$. But, by applying the argument in Lemma \ref{lem:partition} to the subsequence $(\alpha_{n_m},\beta_{n_m})_{m\geq 1}$, it then follows that the Stokes complex $\mathcal{C}(\alpha^*,\beta^*)$ must fall in one of the classes $E_2,E_3$ or $E_4$, which contradicts that $(\alpha^*,\beta^*)\in \widetilde{e}_1$. We conclude that $\lambda_k^n\rightarrow \mu_{\sigma_0(k)}$ for $1\leq k\leq 4$ so that $\widetilde{e}_1$ is indeed characterised by the merging of the turning points $\lambda_3$ and $\lambda_4$.
		The other characterisations are shown analogously.

	Using aforementioned characterisation of the different parts of the border $\delta R$, it is straightforward to extend $\mathcal{S}$ to $K$. We first discuss the extension of $\mathcal{S}$ to $R\cup \widetilde{e}_1$. Let
	$(\alpha,\beta)\in \widetilde{e}_1$ and take any sequence $(\alpha_n,\beta_n)_{n\geq 1}$ in
	$R$ with $(\alpha_n,\beta_n)\rightarrow (\alpha^*,\beta^*)$ as $n\rightarrow \infty$. We
	label the turning points $\{\mu_1,\mu_2,\mu_3\}$ of $V(\lambda;\alpha,\beta)$ such that $\mu_3$
	is the double turning point, $\mu_2$ is the unique simple turning point adjacent to $\mu_3$ and
	$\mu_1$ is the remaining simple turning point, in accordance with Figure \ref{fig:stokes_complexes}, as before.
	
	Let us denote the turning points of $V(\lambda;\alpha_n,\beta_n)$ by $\lambda_j^n=\lambda_j(\alpha_n,\beta_n)$
	for $1\leq j\leq 4$ and $n\geq 1$. Then we know that
	\begin{equation*}
	\lambda_1^n\rightarrow \mu_1,\quad \lambda_2^n\rightarrow \mu_2,\quad \lambda_3^n\rightarrow \mu_3,\quad \lambda_4^n\rightarrow \mu_3
	\end{equation*}
	as $n\rightarrow \infty$ and thus
	\begin{equation*}
	s_1(\alpha_n,\beta_n)\rightarrow +\tfrac{1}{2}(1-\nu)\pi i,\qquad s_2(\alpha_n,\beta_n)\rightarrow
	s_2^{e_1}(\alpha,\beta)\quad (n\rightarrow \infty),
	\end{equation*}
	where
	\begin{equation}\label{eq:defiS2}
	s_2^{e_1}(\alpha,\beta)=\int_{\gamma_2}\omega
	\end{equation}
	with $\gamma_2$ as defined in Figure \ref{fig:stokes_complexes} on the class $E_1$.
	In summary,
	\begin{equation}\label{eq:limitSn}
	\mathcal{S}(\alpha_n,\beta_n)\rightarrow (+\tfrac{1}{2}(1-\nu)\pi, -is_2^{e_1}(\alpha,\beta))\quad (n\rightarrow \infty)
	\end{equation}
	and we therefore set
	\begin{equation*}
	\mathcal{S}(\alpha,\beta)=(+\tfrac{1}{2}(1-\nu)\pi, -is_2^{e_1}(\alpha,\beta)).
	\end{equation*}	
	It is easy to see that $\mathcal{S}(\alpha,\beta)\in \widehat{e}_1$, following the lines in the proof of Lemma \ref{lem:range}.

	To prove continuity of $\mathcal{S}$ on $R\cup \widetilde{e}_1$, it remains to be shown that $s_2^{e_1}$ is continuous on $\widetilde{e}_1$. We proceed in showing that $\widetilde{e}_1$ is a regular smooth submanifold of $\mathbb{R}^4$ and that $s_2^{e_1}$ is locally diffeomorphic on $\widetilde{e}_1$, which in particular implies the required continuity.
	
	Let $(\alpha^*,\beta^*)\in \widetilde{e}_1$, then $\Delta(\alpha^*,\beta^*)=0$. Take a simply connected
	open neighbourhood $U$ of $\alpha^*$, a simply connected open neighbourhood $V$ of $\beta^*$ and a
	biholomorphism $B_\Delta:U\rightarrow V$ such that
	\begin{equation*}
	W:=\{\Delta(\alpha,\beta)=0\}\cap (U\times V) = \{(\alpha,B_\Delta(\alpha)):\alpha\in U\}.
	\end{equation*}
	Note that $W$ is a regular complex submanifold of $\mathbb{C}^2$ and $\widetilde{e}_1\cap (U\times V) \subseteq W$,
	which allows us to study $\widetilde{e}_1$ locally within $W$. To this end, we extend $s_2^{e_1}$ to an analytic mapping on $W$ and show that it is biholomorphic at $(\alpha^*,\beta^*)$.
	
	Let $\{\mu_1^*,\mu_2^*,\mu_3^*\}$ be the turning points of $V(\lambda;\alpha^*,\beta^*)$, labelled
	in accordance with Figure \ref{fig:stokes_complexes} as above.	
	Let $D_j\subseteq \mathbb{C}^*$ with $\mu_j^*\in D_j$, $1\leq j\leq 3$, be mutually disjoint open discs.
	We define the following analytic function on $D_1\times D_2\times D_3$,
	\begin{equation*}
	\widetilde{S}_2(\mu_1,\mu_2,\mu_3)=\oint_{\gamma_2}{\frac{\sqrt{(\mu-\mu_1)(\mu-\mu_2)}(\mu-\mu_3)}{\mu}d\mu},
	\end{equation*}
	with branch and contour chosen consistently with $s_2^{e_1}$ in equation \eqref{eq:defiS2}.\\	
	For $1\leq j\leq 3$, there exists a unique analytic function $\mu_j: W\rightarrow \mathbb{C}^*$
	with $\mu_j(\alpha^*,\beta^*)=\mu_j^*$, $1\leq j\leq 3$, such that $\mu_j(\alpha,\beta)$ is a
	turning point of $V(\lambda;\alpha,\beta)$ for $(\alpha,\beta)\in W$.	 
	We extend $s_2^{e_1}$, defined in \eqref{eq:defiS2}, analytically to $W$, by setting
	\begin{equation*}
	s_2^{e_1}(\alpha,\beta)=s_2^{e_1}(\alpha,B_\Delta(\alpha))=\widetilde{S}_2(\mu_1(\alpha,B_\Delta(\alpha)),\mu_2(\alpha,B_\Delta(\alpha)),\mu_3(\alpha,B_\Delta(\alpha))),
	\end{equation*}
	for $(\alpha,\beta)\in W$.
	
	Now, note that
	\begin{equation*}
	\Omega\cap W=\{(\alpha,\beta)\in W:\Re{s_2^{e_1}(\alpha,\beta)}=0\}.
	\end{equation*}
	We wish to show that $s_2^{e_1}:W\rightarrow \mathbb{C}$ is locally biholomorphic at $(\alpha^*,\beta^*)$, i.e.
	\begin{equation}\label{eq:derivative}
	\frac{\partial}{\partial \alpha}s_2^{e_1}(\alpha,B_\Delta(\alpha))
	\end{equation}
	does not vanish at $\alpha=\alpha^*$.
	
	In order to compute \eqref{eq:derivative}, we introduce,
	inspired by \eqref{eq:constraint1} and \eqref{eq:constraint2},
	\begin{align*}
	\widetilde{S}_3(\mu_1,\mu_2,\mu_3)&=\mu_1\mu_2\mu_3^2-\tfrac{1}{4}\nu^2,\\
	\widetilde{S}_4(\mu_1,\mu_2,\mu_3)&=\mu_1\mu_2+2\mu_1\mu_3+2\mu_2\mu_3-\tfrac{1}{2}(\mu_1^2+\mu_2^2)+2,
	\end{align*}
	and compute the Jacobian
	\begin{equation*}
	|J_{(\widetilde{S}_2,\widetilde{S}_3,\widetilde{S}_4)}(\mu)|=
	\pm 2  \mu_3(\mu_2-\mu_1)(\mu_3-\mu_1)^{\frac{3}{2}}(\mu_3-\mu_2)^{\frac{3}{2}},
	\quad J_{(\widetilde{S}_2,\widetilde{S}_3,\widetilde{S}_4)}(\mu):=
	\left(\frac{\partial \widetilde{S}_{m+1}}{\partial\mu_n}\right)_{1\leq m,n\leq 3}.
	\end{equation*}
	It follows that $M:=J_{(\widetilde{S}_2,\widetilde{S}_3,\widetilde{S}_4)}(\mu^*)$ is invertible.
	By the chain rule, we have
	\begin{equation*}
	\begin{pmatrix}
	\frac{\partial \widetilde{S}_2}{\partial\alpha} \\
	\frac{\partial \widetilde{S}_3}{\partial\alpha} \\
	\frac{\partial \widetilde{S}_4}{\partial\alpha} 
	\end{pmatrix}=
	J_{(\widetilde{S}_2,\widetilde{S}_3,\widetilde{S}_4)}(\mu)\cdot \begin{pmatrix}
	\frac{\partial \mu_1}{\partial\alpha}\\
	\frac{\partial \mu_2}{\partial\alpha}\\
	\frac{\partial \mu_3}{\partial\alpha}
	\end{pmatrix}
	\end{equation*}
	which, when specialised to $(\alpha,\beta)=(\alpha^*,\beta^*)$, gives
	\begin{equation*}
	\begin{pmatrix}
	\frac{\partial}{\partial \alpha}s_2^{e_1}(\alpha,B(\alpha))|_{\alpha=\alpha^*} \\
	0\\
	0
	\end{pmatrix}=
	M\cdot \begin{pmatrix}
	\frac{\mu_1^*}{\mu_3^*-\mu_1^*}\\
	\frac{\mu_2^*}{\mu_3^*-\mu_2^*}\\
	\frac{\mu_3^*(\mu_1^*+\mu_2^*-2\mu_3^*)}{2(\mu_3^*-\mu_1^*)(\mu_3^*-\mu_2^*)}
	\end{pmatrix}
	\end{equation*}	
	Since $M$ is invertible, it thus follows that \eqref{eq:derivative} does not vanishes at $\alpha=\alpha^*$
	and hence $s_2^{e_1}(\alpha,\beta)$ is a local biholomorphism at $(\alpha,\beta)=(\alpha^*,\beta^*)$.
	Therefore, there exists a simply connected open neighbourhood $W_0\subseteq W$ of $(\alpha^*,\beta^*)$ such that 
	\begin{equation*}
	(F_1,F_2):W_0\rightarrow \mathbb{R}^2,(\alpha,\beta)\mapsto (\Re \mathcal{S}_2(\alpha,\beta),\Im \mathcal{S}_2(\alpha,\beta))
	\end{equation*}
	is a diffeomorphism onto its open image $F(W_0)\subseteq \mathbb{R}^2$ with
	\begin{equation*}
	F(W_0)\cap(\{0\}\times \mathbb{R})=\{0\}\times I,
	\end{equation*}
	where $I\subseteq \mathbb{R}$ is an open interval. We apply Lemma \ref{lem:deformboutroux} with 
	\begin{equation*}
	T:=\{(\alpha,\beta)\in W_0:\Re \mathcal{S}_2(\alpha,\beta)=0\},
	\end{equation*}
	giving
	\begin{equation*}
	T=W_0\cap \Omega =W_0\cap \widetilde{e}_1.
	\end{equation*}
	Let us in particular note that, by choosing an open neighbourhood
	$Z\subseteq \{(\Re{\alpha},\Im{\alpha},\Re{\beta},\Im{\beta})\in\mathbb{R}^4\}$ of
	$(\alpha^*,\beta^*)$ with $Z\cap \{\Delta(\alpha,\beta)=0\}\subseteq W_0$, we have
	\begin{equation}\label{eq:edge_isolation}
	Z\cap \{\Delta(\alpha,\beta)=0\}\cap \Omega=Z\cap \widetilde{e}_1.
	\end{equation}

	We conclude that $\widetilde{e}_1$ is a smooth $1$-dimensional regular submanifold of
	$\{(\Re{\alpha},\Im{\alpha},\Re{\beta},\Im{\beta})\in\mathbb{R}^4\}$. Furthermore $F_2|_T$
	maps $T$ diffeomorphically onto $I$. As $\mathcal{S}|_T=F_2|_T$, it follows that
	$\mathcal{S}|_{\widetilde{e}_1}:\widetilde{e}_1\rightarrow \widehat{e}_1$ is locally
	diffeomorphic at $(\alpha,\beta)=(\alpha^*,\beta^*)$ and in particular continuous.
	Combined with \eqref{eq:limitSn}, this shows that the extension of $\mathcal{S}$
	to $R\cup \widetilde{e}_1$ is indeed continuous. 
	
	Next, we will extend $\mathcal{S}$ to $R\cup \widetilde{e}_1\cup \widetilde{e}_2$ in much
	the same way. For $(\alpha,\beta)\in \widetilde{e}_2$, we define
	\begin{equation*}
	\mathcal{S}(\alpha,\beta)=(-is_1^{e_2}(\alpha,\beta),+\nu \pi),
	\end{equation*}
	where
	\begin{equation}\label{eq:defiS1}
	s_1^{e_2}(\alpha,\beta)=\int_{\gamma_1}\omega+\frac{i \pi (1-\nu)}{2}
	\end{equation}
	with $\gamma_1$ as defined in Figure \ref{fig:stokes_complexes} on the class $E_2$,
	and the branch chosen consistently with the definition of $s_1$. Analogously to the above,
	it is shown that $\widetilde{e}_2$ is a smooth $1$-dimensional regular submanifold of
	$\{(\Re{\alpha},\Im{\alpha},\Re{\beta},\Im{\beta})\in\mathbb{R}^4\}$, and that
	$\mathcal{S}|_{\widetilde{e}_2}\widetilde{e}_2\rightarrow \widehat{e}_2$ is locally diffeomorphic.
	
	Continuity of the extension $\mathcal{S}$ to $R\cup \widetilde{e}_1\cup \widetilde{e}_2$ on
	$\widetilde{e}_2$ is proven as it was on $\widetilde{e}_1$, noting that
	$\overline{\widetilde{e}_1}\cap \widetilde{e}_2=\emptyset$ by equation \eqref{eq:edge_isolation}.
	
	Similarly, it is shown that $\widetilde{e}_3$ and $\widetilde{e}_4$ are smooth $1$-dimensional
	regular submanifolds, we have $\overline{\widetilde{e}_i}\cap \widetilde{e}_j=\emptyset$
	for $1\leq i,j\leq 4$ with $i\neq j$, and we extend $\mathcal{S}$ to
	$R\cup \widetilde{e}_1\cup \widetilde{e}_2\cup \widetilde{e}_3\cup \widetilde{e}_4$,
	such that $\mathcal{S}(\widetilde{e}_i)\subseteq \widehat{e}_i$ and
	$\mathcal{S}|_{\widetilde{e}_i}\widetilde{e}_i\rightarrow \widehat{e}_i$ is locally diffeomorphic for $1\leq i\leq 4$.
	
	Finally, we define $\mathcal{S}(\alpha,\beta)=\widehat{c}_i$ if $(\alpha,\beta)\in \widetilde{c}_i$ for
	$1\leq i\leq 4$, so that $\mathcal{S}:K\rightarrow Q$, and it remains to be shown that $\mathcal{S}$
	is continuous at such points.
	
	Let $(\alpha,\beta)\in c_1$ and let us denote the simple and triple turning point of $V(\lambda;\alpha,\beta)$
	by respectively $\eta_1$ and $\eta_2$, in accordance with Figure \ref{fig:stokes_complexes}.	
	Suppose $(\alpha_n,\beta_n)_{n\geq 1}$ is a sequence in $R$ such that
	$(\alpha_n,\beta_n)\rightarrow (\alpha,\beta)$ as $n\rightarrow \infty$
	and let us denote the turning points of $V(\lambda;\alpha_n,\beta_n)$ by
	$\lambda_j^n=\lambda_j(\alpha_n,\beta_n)$ for $1\leq j\leq 4$ and $n\geq 1$.
	It follows from Lemmas \ref{lem:deformationlines} and \ref{lem:triple} that
	\begin{equation*}
	\lambda_1^n\rightarrow \eta_1,\quad \lambda_j^n\rightarrow \eta_2\quad (2\leq j\leq 4)
	\end{equation*}
	and hence
	\begin{equation*}
	\mathcal{S}(\alpha_n,\beta_n)=(-is_1(\alpha_n,\beta_n),-is_2(\alpha_n,\beta_n))\rightarrow \widehat{c}_1
	\end{equation*}
	as $n\rightarrow \infty$.
	
	Similarly, if $(\alpha_n,\beta_n)_{n\geq 1}$ is a sequence in $\widetilde{e}_1$ which converges to $(\alpha,\beta)$,
	we denote the turning points of $V(\lambda;\alpha_n,\beta_n)$ by $\mu_j^n$, $1\leq j\leq 3$
	in accordance with Figure \ref{fig:stokes_complexes}, and Lemma \ref{lem:triple} yields
	\begin{equation*}
	\mu_1^n\rightarrow \eta_1,\quad \mu_j^n\rightarrow \eta_2\quad (2\leq j\leq 3)
	\end{equation*}
	and  thus
	\begin{equation*}
	\mathcal{S}(\alpha_n,\beta_n)=(+\tfrac{1}{2}(1-\nu)\pi,-is_2^{e_1}(\alpha_n,\beta_n))\rightarrow \widehat{c}_1
	\end{equation*}
	as $n\rightarrow \infty$.
	
	We handle sequences in $\widetilde{e}_2$ converging to $(\alpha^*,\beta^*)$ similarly and finally note that the sets $\widetilde{c}_i$,$1,\leq i\leq 4$ 
	are discrete and $(\alpha,\beta)\notin \overline{\widetilde{e}_i}$ for $3\leq i \leq 4$.
	It follows that $\mathcal{S}$ is continuous at $(\alpha,\beta)\in \widetilde{c}_1$. Continuity at points in
	$\widetilde{c}_j$ for $2\leq j\leq 4$ are proven analogously and it follows that $\mathcal{S}$ is globally continuous.
\end{proof}
\end{lem}
\subsubsection*{Step \ref{step3}}
We start with proving that $K$ is compact.
\begin{lem}\label{lem:compact}
	The region $K=\overline{R}$ is compact.
\end{lem}
\begin{proof}
	Suppose $K$ is not compact, then $R$ is unbounded. Take a sequence $(\alpha_n,\beta_n)_{n>1}$ in
	$R$ such that $|\alpha_n|+|\beta_n|\rightarrow \infty$ as $n\rightarrow \infty$.
	Let us write $\widetilde{\beta}_n=\frac{\beta_n}{\alpha_n^3}$ for $n\geq 1$.
	By replacing $(\alpha_n,\beta_n)_{n>1}$ by an appropriate subsequence if necessary, we may assume
	that we are in one of the following four scenarios:
	\begin{enumerate}[label=(\roman*)]
		\item $\alpha_n\rightarrow \alpha^*\in\mathbb{C}$ and $|\beta_n|\rightarrow \infty$,
		\item $|\alpha_n|\rightarrow \infty$ and $|\widetilde{\beta}_n|\rightarrow \infty$,
		\item $|\alpha_n|\rightarrow \infty$ and $\widetilde{\beta}_n\rightarrow {\widetilde{\beta}}^*\in\mathbb{C}^*$, or
		\item $|\alpha_n|\rightarrow \infty$ and $\widetilde{\beta}_n\rightarrow 0$
	\end{enumerate}
	as $n\rightarrow \infty$.
	
	Each of the four cases leads to a contradiction in a similar fashion and we therefore limit our discussion to one of them, case (iii).
	Setting $\lambda=\alpha_n\mu$, we have
	\begin{equation}\label{eq:differential}
	\sqrt{V(\lambda;\alpha_n,\beta_n)}d\lambda=\alpha_n^2\sqrt{\mu^2+2\mu+1-\widetilde{\beta}_n\mu^{-1}-
		\alpha_n^{-2}+\alpha_n^{-4}\frac{\nu^2}{4}\mu^{-2}}d\mu.
	\end{equation}
	Let us write the turning points of ${V(\lambda;\alpha_n,\beta_n)}$ by $\lambda_j^n=\lambda_j(\alpha_n,\beta_n)$
	and define $\mu_j^n:=\alpha_n^{-1}\lambda_j^n$ 
	for $1\leq j\leq 4$ and $n\geq 1$. Then, by replacing $(\alpha_n,\beta_n)_{n>1}$ by an appropriate subsequence
	if necessary, we may assume that there exists a permutation $\sigma\in S_4$ such that
	\begin{equation*}
	\mu_{\sigma(4)}^n\rightarrow 0,\quad \mu_{\sigma(j)}^n\rightarrow u_j\in\mathbb{C}^*\quad (1\leq j\leq 3) 
	\end{equation*}
	as $n\rightarrow \infty$, where $\{u_1,u_2,u_3\}$ are the roots of $\mu^3+2\mu^2+\mu-{\widetilde{\beta}}^*$.
	
	However, $\mu=0$ is a simple pole of the differential \eqref{eq:differential}. Recalling definition
	\eqref{eq:defir} of $r_4,r_6\in\mathbb{R}_+$ in the proof of Lemma \ref{lem:range}, we hence obtain,
	as $n\rightarrow \infty$ and thus $\mu_{\sigma(4)}^n\rightarrow 0$, that either $r_4\rightarrow +\infty$
	or $r_6\rightarrow +\infty$, contradicting equation \eqref{eq:r4r6}.
\end{proof}
Step \ref{step3} is accomplished by the
the following proposition.
\begin{pro}\label{prop:Sextension}
	The continuous extension $\mathcal{S}:K\rightarrow Q$, defined in Lemma
	\ref{lem:extension}, is a homeomorphism and has the following further properties:
	\begin{enumerate}
		\item $\mathcal{S}$ maps the interior $R$ of $K$
		diffeomorphically onto the open rectangle $Q^\circ$;
		\item for $1\leq k\leq 4$,	 $\mathcal{S}$ maps part $\widetilde{e}_k$ of the border $\delta R$ diffeomorphically\footnote{diffeomorphically with respect to the geometric structure on $\widetilde{e}_k$
			as a smooth regular submanifold of $\{(\Re{\alpha},\Im{\alpha},\Re{\beta},\Im{\beta})\in\mathbb{R}^4\}$,
			see Lemma \ref{lem:extension}.} onto the open edge $\widehat{e}_k$ of $Q$;
		\item for $1\leq k\leq 4$, the part $\widetilde{c}_k$ of the border $\delta R$ is a singleton and
		$\mathcal{S}$ maps $\widetilde{c}_k\equiv\{\widetilde{c}_k\}$ to the corner $\widehat{c}_k$ of $Q$.
	\end{enumerate}
\begin{proof}
		As a first step, we prove that $\mathcal{S}$ is injective when restricted to $R$. We just give a rough sketch of the proof and refer the interested reader to \cite{bertola2009,jenkins}
		for a complete treatment of similar arguments.
		Suppose $\mathcal{S}(\alpha,\beta)=\mathcal{S}(\alpha',\beta')$ and let respectively
		$\lambda_1,\ldots,\lambda_4$ and $\lambda_1',\ldots,\lambda_4'$ be corresponding turning points. Furthermore, let us denote corresponding differentials by $\omega$ and $\omega'$.
		
		Note that the Stokes complex of the potential $V(\lambda,\alpha,\beta)$ naturally cuts
		the complex $\lambda$-plane into five disjoint open connected regions $I,\ldots,V$ as in Figure \ref{fig:contours}.
		By choosing the sign correctly, the action integral
		\begin{equation*}
		S(\lambda)=\int_{\lambda_1}^{\lambda}\omega
		\end{equation*}
		defines a uniformisation of $I$ (resp. homeomorphism from $\overline{I}$) onto the open
		(resp. closed) left-half plane, mapping the turning points $\lambda_1,\lambda_2,\lambda_3,\lambda_4$
		to the respective marked points
		$0$, $r_6i$, $(r_5+r_6)i$ and $(r_4+r_5+r_6)i$, where $r_4,r_5,r_6$ as defined in Lemma \ref{lem:range}.
		
		Similarly, the Stokes complex of the potential $V(\lambda,\alpha',\beta')$ cuts the complex plane into five pieces
		$I',\ldots, V'$ and 
		\begin{equation*}
		S'(\lambda)=\int_{\lambda_1'}^{\lambda}\omega'
		\end{equation*}
		defines a uniformisation of $I'$ (resp. homeomorphism from $\overline{I'}$) onto the open
		(resp. closed) left-half plane, which maps the turning points $\lambda_1',\lambda_2',\lambda_3',\lambda_4'$
		to the same respective marked points
		$0$, $r_6i$, $(r_5+r_6)i$ and $(r_4+r_5+r_6)i$, as $\mathcal{S}(\alpha,\beta)=\mathcal{S}(\alpha',\beta')$.
		
		Thus $\chi:\overline{I}\rightarrow \overline{I'}$, defined by
		\begin{equation*}
		\chi:=S'|_{\overline{I'}}^{-1}\circ S|_{\overline{I}},
		\end{equation*}
		is a homeomorphism which maps $I$ conformally and bijectively onto $I'$, mapping turning point
		$\lambda_k$ to $\lambda_k'$ for $1\leq k\leq 4$.
		
		Without going into further detail, one may prove that $\chi$ can be uniquely extended to an
		automorphism of $\mathbb{P}^1$, mapping each region $J\in \{I,\ldots V\}$ to its corresponding
		region $J'$ and in particular
		\begin{equation*}
		\chi(0)=0,\quad\chi(\infty)=\infty,\quad\chi(\lambda_i)=\lambda_i'\quad (1\leq i\leq 4).
		\end{equation*}
		It follows that $\chi$ is a dilation. There are only two dilations which preserve equations
		\eqref{eq:constraint1} and \eqref{eq:constraint2}, namely $\chi(\lambda)\equiv +\lambda$ and
		$\chi(\lambda)\equiv -\lambda$. Since by construction $\chi$ leaves the asymptotic direction
		$e^{-\frac{3}{4}\pi i}\infty$ invariant, we must have $\chi(\lambda)\equiv +\lambda$.
		In particular $(\alpha',\beta')=(\alpha,\beta)$ by equations \eqref{eq:alphaoflambdas} and
		\eqref{eq:betaoflambdas}, and we conclude that $\mathcal{S}$ is indeed injective on $R$.
	
	Next, we prove that $\mathcal{S}$ maps $R$ diffeomorphically onto the open rectangle $Q^\circ$. We already know that the restriction of $\mathcal{S}$ to $R$ is an injective smooth mapping
	from $R$ to $Q^\circ$, and hence, by the domain invariance theorem,
	$\widehat{R}=\mathcal{S}(R)$ is an open subset of $Q^\circ$ and $\mathcal{S}$
	maps $R$ diffeomorphically onto $\widehat{R}$. By construction of $\mathcal{S}$
	in Lemma \ref{lem:extension}, we have $\mathcal{S}(\delta R)\subseteq \partial Q$.	
	Since, by Lemma \ref{lem:compact}, $K$ is compact, $\partial \widehat{R}\subseteq\mathcal{S}(\delta R)$
	and hence $\partial \widehat{R}\subseteq \partial Q$.
	
	There exists only one non-empty open subset of $Q^\circ$ whose boundary is
	contained in $\partial Q^\circ$, namely $Q^\circ$ itself. We conclude that
	$\widehat{R}=Q^\circ$ and in particular $\mathcal{S}$ maps $R$ diffeomorphically onto the open rectangle $Q^\circ$.
	
	It follows that $\mathcal{S}(\delta R)=\partial \widehat{R}=\partial Q$ and we have
	$\mathcal{S}(K)=Q$, i.e. $\mathcal{S}$ is surjective.
	
	Next we consider $\mathcal{S}$ on the border $\delta R$ in more detail.
	Let $1\leq k\leq 4$, then, by construction of $\mathcal{S}$ in Lemma \ref{lem:extension},
	we have $\mathcal{S}^{-1}(\widehat{e}_k)\subseteq \widetilde{e}_k$,
	so the restricted mapping $\mathcal{S}|_{\widetilde{e}_k}:\widetilde{e}_k\rightarrow \widehat{e}_i$
	is surjective. Analogously to the injectivity argument above for $\mathcal{S}|_R$, it is shown
	that the restricted mapping $\mathcal{S}|_{e_k}$ is injective.
	Since, by Lemma \ref{lem:extension}, the restricted mapping
	$\mathcal{S}|_{e_k}:e_k\rightarrow \widehat{e}_k$ is locally diffeomorphic it
	follows that $\mathcal{S}$ maps $\widetilde{e}_k$ diffeomorphically onto
	$\widehat{e}_k$. 
	
	We conclude that $\mathcal{S}$ maps $R\cup \widetilde{e}_1\cup
	\widetilde{e}_2\cup \widetilde{e}_3 \cup \widetilde{e}_4$ bijectively onto
	$Q^{\circ}\cup \widehat{e}_1\cup \widehat{e}_2\cup \widehat{e}_3 \cup \widehat{e}_4$.
	
	Let $1\leq k\leq 4$, then, by construction of $\mathcal{S}$ in Lemma \ref{lem:extension},
	we have $\mathcal{S}^{-1}(\widehat{c}_k)\subseteq \widetilde{c}_k$. Hence the
	surjectivity of $\mathcal{S}$ implies that $\widetilde{c}_i$ is non-empty.
	On the other hand, it is easy to see that $\widetilde{c}_i$ has at most one element,
	applying the same reasoning as in the injectivity argument for $\mathcal{S}|_R$ above.
	So $\widetilde{c}_k$ is a singleton and $\mathcal{S}$ maps $\widetilde{c}_k\equiv\{\widetilde{c}_k\}$
	to the corner $\widehat{c}_k$ of $Q$.
	
	All together we obtain that $\mathcal{S}$ is a continuous bijection. Since $K$ is compact,
	$\mathcal{S}$ sends closed sets to closed sets and therefore $\mathcal{S}$ is a homeomorphism.
\end{proof}
\end{pro}
From here on, we will call $\widetilde{e}_k$, $1\leq k\leq 4$, and $\widetilde{c}_k$, $1\leq k\leq 4$,
the edges and corners of the region $K$, as justified by Proposition \ref{prop:Sextension}.
\subsubsection*{Step \ref{step4}}
We proceed with discussing Step \ref{step4},
i.e. the implicit parametrisation of the edges and corners of $K$. Let us introduce the following notation
\begin{equation*}
e_k=\Pi_a(\widetilde{e}_k),\quad c_k=\Pi_a(\widetilde{c}_k)\quad (1\leq k\leq 4).
\end{equation*}
Since $K$ is compact and $R_a=\Pi_a(R)$ is open, due to Lemma \ref{lem:deformation}, the boundary of the elliptic region is a subset of $\Pi_a(\delta R)$, i.e.
\begin{equation*}
\partial R_a\subseteq \Pi_a(\delta R)=e_1\cup e_2\cup e_3\cup e_4 \cup\{c_1,c_2,c_3,c_4\}.
\end{equation*}
However, at this point, it is by no means clear that $\partial R_a=\Pi_a(\delta R)$.

Recall from Lemma \ref{lem:partition} that points on the border $\delta R$ correspond to singular Boutroux curves, that is
\begin{equation*}
\delta R\subseteq \Omega \cap \{\Delta(\alpha,\beta)=0\},
\end{equation*}
where $\Delta(\alpha,\beta)$ denotes the discriminant \eqref{eq:discriminant}.

Our strategy to obtain a parametrisation of the border $\delta R$ goes by studying the set $\Omega \cap \{\Delta(\alpha,\beta)=0\}$ and constructing local parametrisations of this set in general after which we construct one for the subset of interest, $\delta R$.
We first briefly discuss the affine variety $\{\Delta(\alpha,\beta)=0\}$ and then consider its intersection with $\Omega$.

Take $(\alpha,\beta)\in \mathbb{C}^2$ such that $\Delta(\alpha,\beta)=0$. A priori, this means that the potential $V(\lambda;\alpha,\beta,\nu)$ can have one of four possible configurations of its turning points:
\begin{enumerate}
	\item the potential has one double turning point and two simple turning points,
	\item the potential has one triple turning point and one simple turning points,
	\item the potential has two double turning points,
	\item the potential has one quadruple turning point.
\end{enumerate}
One may check by direct computation that cases (4) and (5) are an algebraic impossibility for the potential $V(\lambda;\alpha,\beta,\nu)$ with $0<\nu \leq \tfrac{1}{3}$. Indeed, both cases can only be realised if $\nu^2=1$.

Let us denote by $X$ the unique double or triple turning point of the potential $V$. Then it follows by direct computation that $X$ and $\alpha$ are related by the following algebraic relation,
\begin{equation}\label{eq:defiX}
3 X^4+4\alpha X^3+(\alpha^2-1)X^2-\frac{\nu^2}{4}=0.
\end{equation}

For generic elements of $\{\Delta(\alpha,\beta)=0\}$ the multiple turning point $X$ is a double turning point. A triple turning point occurs when one of the remaining two simple turning points coincides with $X$. These points are characterised by the coincidence of two different branches of the quartic \eqref{eq:defiX}, and thus the discriminant of the quartic with respect to $X$ vanishes there, i.e.
\begin{equation}\label{eq:defiC}
C(\alpha):=\alpha^8-6(3\nu^2+1)\alpha^4+8(1-9\nu^2)\alpha^2-3(9\nu^4+6\nu^2+1),
\end{equation}
is zero there. More concretely, we have
\begin{align*}
\Omega_0:&=\{(\alpha,\beta) \in\mathbb{C}^2:V(\lambda;\alpha,\beta)\text{ has a triple turning point}\}\\
&=\{(\alpha,\beta)\in\mathbb{C}^2:C(\alpha)=0,\Delta(\alpha,\beta)=0\}\\
&=\{(\alpha,\beta)\in\mathbb{C}^2:C(\alpha)=0,\beta=f(\alpha)\},
\end{align*}
where the last equality follows by direct computation, with
\begin{equation}\label{eq:defif}
f(\alpha):=-\frac{\alpha^2(\alpha^2-2)+3\nu^2+1}{6\alpha}.
\end{equation}

Note that $\Omega_0\subset \Omega$ and that, for $1\leq k\leq 4$, by definition of the Stokes complex class $C_k$, the corner $\widetilde{c}_k$ of $K$ is an element of $\Omega_0$, namely
\begin{equation}\label{eq:cornerequation}
\widetilde{c}_k=(c_k,f(c_k)),\quad C(c_k)=0.
\end{equation}
Thus, to obtain an explicit description of the corners of $K$, all that is left to be done is to understand which of the eight roots of $C(\alpha)$ equals $c_k$, for $1\leq k\leq 4$.
 To facilitate the labeling of the eight roots of $C(\alpha)$, we have the following
 \begin{lem}\label{lem:zerospolynomial}
 	The polynomial $C(\alpha)$ has precisely two real roots and two purely imaginary roots.
 \end{lem}
 \begin{proof}
 	Application of Descartes' rule of signs on $C(\pm i \alpha)$ immediately gives that $C(\alpha)$ has precisely two purely imaginary roots. Note that
 	\begin{equation*}
 	C'(\alpha)=8\alpha(\alpha^2+2)(\alpha^2+3\nu-1)(\alpha^2-3\nu-1)
 	\end{equation*}
 	and direct computation gives that $C(\alpha)<0$ for each of the five real roots of $C'(\alpha)$. Since $C(\alpha)\sim 8\alpha^8$ as $\alpha\rightarrow \infty$, it follows that $C(\alpha)$ has precisely two real roots.
 \end{proof}
 By Lemma \ref{lem:zerospolynomial} and the fact that $C(\alpha)$ is real and symmetric, we may define $u_k$ as the unique root of $C(\alpha)$ in the $k$-th quadrant of the complex plane $\{\alpha\in\mathbb{C}\}$ for $1\leq k\leq 4$, so that
 \begin{equation}\label{eq:defiu}
 u_2=-\overline{u}_1,\quad u_3=-u_1,\quad u_4=\overline{u}_1.
 \end{equation}
 We also define $v_k$ as the unique root of $C(\alpha)$ within the semi-axis $i^{k-1}\mathbb{R}_+$, so that
 \begin{equation*}
 v_3=-v_1,\quad v_4=-v_2,
 \end{equation*}
 see Figure \ref{fig:complexl1}.  We conclude that $\Omega_0$ consists of the eight points $(u_k,f(u_k))$, $(v_k,f(v_k))$, $1\leq k\leq 4$. The final identification $c_k=u_k$, $1\leq k\leq 4$,
 requires some additional analysis (which we postpone to Proposition \ref{pro:boundary}).

Next we consider the intersection of $\{\Delta(\alpha,\beta)=0\}$ with $\Omega$.
The variety $\{\Delta(\alpha,\beta)=0\}$ has one complex dimension and the intersection $\Omega \cap \{\Delta(\alpha,\beta)=0\}$ with $\Omega$ imposes one additional real constraint, since the homology of the underlying singular elliptic curve (generically) has dimension one, and thus $\Omega \cap \{\Delta(\alpha,\beta)=0\}$ is locally (real) one-dimensional: it consists of smooth lines which can only meet at branching points of $\{\Delta(\alpha,\beta)=0\}$, i.e. points in $\Omega_0$.

Let us be more explicit, we have the following partition,
\begin{align}
\{\Delta(\alpha,\beta)=0\}\cap\Omega=\Omega_0\sqcup \Omega_1,
\end{align}
where
\begin{align*}
\Omega_0&=\{(\alpha,\beta) \in\Omega:V(\lambda;\alpha,\beta)\text{ has one triple turning point}\},\\
\Omega_1&=\{(\alpha,\beta) \in\Omega:V(\lambda;\alpha,\beta)\text{ has one double turning point}\}.
\end{align*}
We have already discussed the set $\Omega_0$ above, its a discrete set of eight elements. We proceed with a heuristic discussion of how to obtain local parametrisations of $\Omega_1$ which in particular show that it is made up of smooth lines. The full details can be found in the proof of Proposition \ref{pro:boundary}.

Consider a generic point in $\{\Delta(\alpha,\beta)=0\}$, so that the potential $V$ has a double turning point, and denote this double turning point by $X$ as before. Recall that $X$ is related to $\alpha$ via the algebraic equation \eqref{eq:defiX}.

Now, imposing $(\alpha,\beta)\in \Omega$, and hence $(\alpha,\beta)\in \Omega_1$, is then equivalent to imposing the real constraint
\begin{equation}\label{eq:realityconstraint}
\Re \int_{x_0}^X\sqrt{V(\lambda)}d\lambda=0,
\end{equation}
where $x_0$ denotes any of the remaining simple turning points. We may explicitly integrate the integral on the left-hand side, yielding locally an analytic function in $\alpha$, and thus locally $\Omega_1$ is described as the zero set of a harmonic function, i.e. a smooth line segment. In particular this gives us a local implicit parametrisation of $\Omega_1$.

Note that aforementioned implicit parametrisation depends crucially on the branch $X$ of \eqref{eq:defiX}. We proceed in introducing the branch relevant to the parametrisation of the entire border $\delta R$ of $R$, and stating rigorously the results concerning this parametrisation. The relevant branch $x=x(\alpha)$ is the one which solves the quartic \eqref{eq:defiX} analytically in the complex $\alpha$-plane cut along the diagonals
$[u_1,u_3]$ and $[u_2,u_4]$ plus midlines $[v_1,v_3]$ and $[v_2,v_4]$,
uniquely characterised asymptotically by
\begin{equation}\label{eq:defix}
x(\alpha)\sim\frac{\nu}{2}\alpha^{-1}\quad  (\alpha\rightarrow \infty).
\end{equation}

\begin{figure}
	\begin{tikzpicture}
	\node[anchor=south west,inner sep=0] (image) at (0,0) {\includegraphics[width=0.5\textwidth]{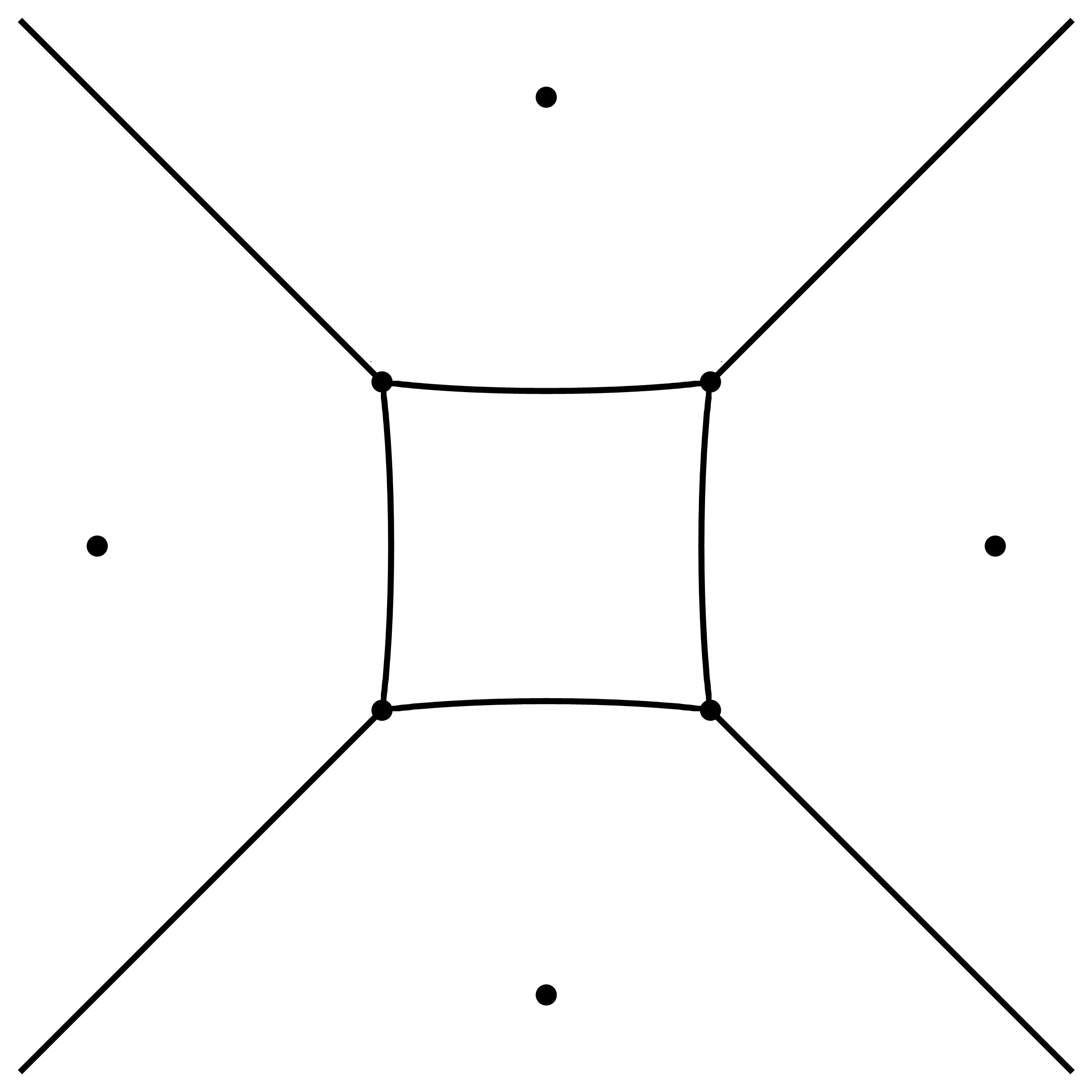}};
	\begin{scope}[x={(image.south east)},y={(image.north west)}]
	\node[above] at (0.5,0.5) {$0$};
	\node[below right] at (0.654,0.654) {$u_1$};
	\node[below left] at (0.346,0.654) {$u_2$};
	\node[above left] at (0.346,0.346) {$u_3$};
	\node[above right] at (0.654,0.346) {$u_4$};
	
	\node[above] at (0.915,0.5) {$v_1$};
	\node[above] at (0.5,0.915) {$v_2$};
	\node[above] at (0.085,0.5) {$v_3$};
	\node[above] at (0.5,0.085) {$v_4$};
	\draw[red,dashed] (0.36,0.36) -- (0.64,0.64);
	\draw[red,dashed] (0.36,0.64) -- (0.64,0.36);
	
	\node[right] at (0.64,0.5) {$\epsilon_1$};
	\node[above] at (0.5,0.64) {$\epsilon_2$};
	\node[left] at (0.36,0.5) {$\epsilon_3$};
	\node[below] at (0.5,0.36) {$\epsilon_4$};
	
	\node[above left] at (0.85,0.85) {$l_1$};
	\node[above right] at (0.15,0.85) {$l_2$};
	\node[below right] at (0.15,0.15) {$l_3$};
	\node[below left] at (0.85,0.15) {$l_4$};
	
	\end{scope}
	\end{tikzpicture}
	\caption{The level set $\{\psi(\alpha)=0\}$ and zeros of the polynomial $C(\alpha)$, defined in \eqref{eq:defiu}, with the branch cuts $[u_1,u_3]$ and $[u_2,u_4]$ in dashed red.}
	\label{fig:parametrisation2}
\end{figure}
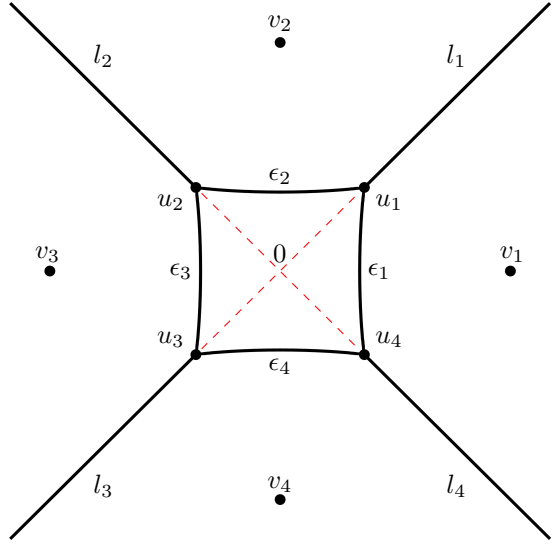

Since equation \eqref{eq:defiX} is of fourth order in $X$, we can compute all its branches explicitly. In particular, for $\alpha\in (v_1,\infty)$ all the branches of \eqref{eq:defiX} are real and only one is positive, namely $x(\alpha)$. Since $X=0$ is never a root of \eqref{eq:defix}, it follows that $x(\alpha)$ does not coincide with any other branch of \eqref{eq:defiX} at $\alpha =v_1$, so $x(\alpha)$ is single-valued at $\alpha=v_1$. Similarly it follows that $x(\alpha)$ is single-valued at $v_k$ for $2\leq k\leq 4$ and thus $x(\alpha)$ is analytic on the $\alpha$-plane merely cut along the diagonals $[u_1,u_3]$ and $[u_2,u_4]$, see Figure \ref{fig:parametrisation2}.

An analogous argument shows that there exists a unique algebraic function $y=y(\alpha)$ which solves
\begin{equation}\label{eq:ybranch}
y^2=\alpha^2+6 x\alpha+6x^2-1
\end{equation}
on the same cut $\alpha$-plane with $y(\alpha)\sim \alpha$ as $\alpha\rightarrow \infty$. Next, we reintroduce the function $\psi(\alpha)$ which allows us to parametrise the border $\delta R$,
\begin{equation}\label{eq:parametrisation2}
\psi(\alpha)=\tfrac{1}{2}\Re\left[\alpha y+\tfrac{1}{2}(1-\nu)\log(p_1)-\log(p_2)+\nu \log(p_3)\right],
\end{equation}
where
\begin{equation*}
p_1=1-2 x\alpha-2x^2,\quad p_2=2x+\alpha+y,\quad p_3=\frac{x(\alpha^2+5 x\alpha+4x^2-1) +\frac{1}{2}\nu  y}{x^2}.
\end{equation*}
This function constitutes the left-hand side of \eqref{eq:realityconstraint} for the branch $x$.
The precise results concerning the implicit parametrisation of the border $\delta R$ of $R$ in the following
proposition, which completes Step \ref{step4}.
\begin{pro}\label{pro:boundary}
	Considering the function $\psi(\alpha)$ defined in \eqref{eq:parametrisation2}, the following hold true.
	\begin{enumerate}
		\item The function $\psi(\alpha)$ is a harmonic function on the $\alpha$-plane cut along $[u_1,u_3]$ and $[u_2,u_4]$, and $\psi(\alpha)\rightarrow 0$ as $\alpha$ approaches any of the branch points $u_k$,$1\leq k \leq 4$ within the cut plane.
		\item The zero set $\{\psi(\alpha)=0\}$ takes the form depicted in Figure \ref{fig:parametrisation2}, namely it consists of the four roots $u_k$,$1\leq k \leq 4$ of $C(\alpha)$, and eight (mutually disjoint) level curves, $\epsilon_{k}$, $l_{k}$, $1\leq k\leq 4$.
		From each $u_k$ emanate three level curves, $\epsilon_{k},\epsilon_{k+1}$ and $l_k$, with $l_k$ going to infinity asymptotic to $e^{\frac{\pi i}{4}(2k-1)}\infty$, for $1\leq k\leq 4$.
		\item For $1\leq k\leq 4$, the internal radial angle between $\epsilon_{k}$ and $\epsilon_{k+1}$ at $u_k$ equals $\frac{2}{5}\pi$.	
		\item  Define
		\begin{equation}\label{eq:defiB}
		B_\Delta(\alpha)=x(-2+4 x^2+6 x\alpha+2\alpha^2),
		\end{equation}
		where $x=x(\alpha)$ is the algebraic function introduced in \eqref{eq:defix}, then $B_\Delta$ is an analytic function on the cut $\alpha$-plane with a well-defined limiting value $B_\Delta(\alpha)=f(\alpha)$ at the branching points $\alpha=u_k$, $1\leq k\leq 4$, where $f$ is defined in \eqref{eq:defif}, and	
		for $1\leq k\leq 4$, the corner $\widetilde{c}_k$ and edge $\widetilde{e}_k$ of $K$ are parametrised by
		\begin{subequations}
			\begin{align}
			\widetilde{c}_k&=(u_k,B_\Delta(u_k)),& c_k&=u_k,\label{eq:paracorner}\\
		\widetilde{e}_k&=\{(\alpha,B_\Delta(\alpha)):\alpha \in \epsilon_k\}, &	e_k&=\epsilon_k.\label{eq:paraedge}
			\end{align}
		\end{subequations}
	\end{enumerate}
\begin{proof}[Proof of Proposition \ref{pro:boundary}.]
	Let us first recall Lemma \ref{lem:partition}, which states in particular that
	\begin{equation*}
	\delta R\subseteq \Omega\cap \{\Delta(\alpha,\beta)=0\}.
	\end{equation*}
	Namely, on the border $\delta R$ two or three turning points have coalesced, so that the resulting Boutroux curve $\widehat{\Gamma}$ is singular. To prove the proposition, we first describe a method to compute $\Omega\cap \{\Delta(\alpha,\beta)=0\}$, and then we will specialise it to the subset of our interest, namely $\delta R$.
	
	Suppose $(\alpha,\beta)\in\mathbb{C}^2$ is such that $\Delta(\alpha,\beta)=0$, then there exist a unique turning point $\lambda=X$ of $V(\lambda;\alpha,\beta)$ which is not simple. Consider for the moment the generic case in which $X$ is double and thus $V(\lambda)$ has two remaining simple turning points, say $x_{1,2}$. Then 
	\begin{equation*}
	V(\lambda;\alpha,\beta)=\lambda^{-2}(\lambda-x_1)(\lambda-x_2)(\lambda-X)^2,
	\end{equation*}
	yielding the following helpful relations,
	\begin{subequations}\label{eq:helpful}
		\begin{align}
		&3 X^4+4\alpha X^3+(\alpha^2-1)X^2-\frac{\nu^2}{4}=0,\label{eq:sing1}\\
		&X(-2+4 X^2+6 X\alpha+2\alpha^2)=\beta,\label{eq:sing2}\\
		&x_1+x_2=-2(X+\alpha),\label{eq:sing3}\\
		&(x_2-x_1)^2=4(1-2 X\alpha-2 X^2).\label{eq:sing4}
		\end{align}
	\end{subequations}
	Here we recognise equations \eqref{eq:defix} and \eqref{eq:defiB} in (\ref{eq:sing1},\ref{eq:sing2}).
	
	Now clearly $\widehat{\Gamma}(\alpha,\beta)$ is a Boutroux curve if and only if
	\begin{equation*}
	\Re \int_{x_*}^X\frac{\sqrt{(\lambda-x_1)(\lambda-x_2)}(\lambda-X)}{\lambda}d\lambda=0,
	\end{equation*}
	for any choice of $x_*\in\{x_1,x_2\}$ and choice of contour. This integral can be evaluated explicitly, yielding
	\begin{equation}\label{eq:parametrisationgen}
	\Psi=\tfrac{1}{2}\Re\left[\alpha y+\tfrac{1}{2}(1-\nu)\log(P_1)-\log(P_2)+\nu \log(P_3)\right],
	\end{equation}
	with
	\begin{equation*}
	P_1=\tfrac{1}{4}(x_2-x_1)^2,\quad P_2=X-\tfrac{1}{2}(x_1+x_2)+Y,\quad P_3=\frac{X(x_1x_2-\tfrac{1}{2}(x_1+x_2)X) +\frac{1}{2}\nu  Y}{X^2},
	\end{equation*}
	where $Y$ is a branch of
	\begin{equation}\label{eq:Ybranch}
	Y^2=(X-x_1)(X-x_2)=\alpha^2+6 X\alpha+6X^2-1.
	\end{equation}
	Firstly note that $\Psi$ is invariant under interchanging of $x_1$ and $x_2$ and we may thus eliminate them from the formula using equations (\ref{eq:sing3},\ref{eq:sing4}), yielding equation \eqref{eq:parametrisation2} with $x\rightarrow X$ and $y\rightarrow Y$.
	
	In summary, we can compute $\Omega\cap \{\Delta(\alpha,\beta)=0\}$ as follows. Take a local branch $X=X(\alpha)$ of \eqref{eq:sing1} and choose the correct branch $Y=Y(\alpha)$ of \eqref{eq:Ybranch} realising equation \eqref{eq:parametrisationgen}, define the algebraic function $\beta_\Delta(\alpha)$ by the left-hand side of \eqref{eq:sing2}, 
	then $\Psi$ is (locally) a harmonic function in $\alpha$ and for any $\alpha_0$ in its zero set we have
	\begin{equation*}
	(\alpha_0,\beta_\Delta(\alpha_0))\in \Omega\cap \{\Delta(\alpha,\beta)=0\},
	\end{equation*} 
	yielding a local parametrisation of $\Omega\cap \{\Delta(\alpha,\beta)=0\}$. By considering each of the four branches $X$ of \eqref{eq:sing1} and following the above procedure we can in principle completely compute $\Omega\cap \{\Delta(\alpha,\beta)=0\}$. To prove the proposition this however will not be necessary.
	
	In the case of our interest the relevant branch of \eqref{eq:sing1} is given by $X=x$, introduced in equation \eqref{eq:defix}, and the corresponding correct branch of \eqref{eq:Ybranch} is given by $Y=y$ introduced in \eqref{eq:ybranch}. Let us first prove part (1) of the proposition. For $\alpha$ on the plane cut along $[u_1,u_3]$ and $[u_2,u_4]$ we have
	\begin{equation*}
	V(\lambda;\alpha,B_\Delta(\alpha))=\lambda^{-2}(\lambda-x_1)(\lambda-x_2)(\lambda-X)^2,
	\end{equation*}
	for some up to permutation unique $x_1$ and $x_2$. Recall that they cannot coincide on the cut plane, indeed, if $x_1=x_2$ equations \eqref{eq:helpful} imply $\nu^2=1$, which contradicts $0<\nu\leq \tfrac{1}{3}$.
	 We may thus specify each uniquely by choosing $x_{1,2}=x_{1,2}(\alpha)$ analytically on the cut plane with
	\begin{equation}\label{eq:x1x2asymptotic}
	x_1=-\alpha-\sqrt{1-\nu}+\mathcal{O}(\alpha^{-1}),\quad x_2=-\alpha+\sqrt{1-\nu}+\mathcal{O}(\alpha^{-1})
	\end{equation}
	as $\alpha\rightarrow \infty$, due to equations (\ref{eq:sing3},\ref{eq:sing4}).
	
	Since $x_1\neq x_2$ on the cut plane it follows from equation \eqref{eq:parametrisation2}, that 
	the term $p_1$ in formula \eqref{eq:parametrisationgen} does not vanish and thus $\Re\log(p_1)$ is finite and harmonic on the cut plane. Similarly, if $p_2=0$, then $Y^2=(X-\tfrac{1}{2}(x_1+x_2))^2$ which implies $x_1=x_2$. It follows that $p_2$ does not vanish and thus $\Re\log(p_2)$ is finite and harmonic on the cut plane. The same argument, using the identity $x_1x_2X^2=\frac{1}{4}\nu^2$, works for $p_3$ and we conclude that $\psi(\alpha)$ is a harmonic function on the entire cut plane.
	
	The branching points of $X(\alpha)$ are characterised by the coalescing of one of the two simple turning points $\{x_1,x_2\}$ with the double turning point $X(\alpha)$: at $\alpha=u_1,u_4$ the turning point $x_2$ coalesces with $X$ and at $\alpha=u_2,u_3$ the turning point $x_1$ coalesces with $X$. 
	
	The local behaviour of $\psi(\alpha)$ near the branching points is easily computed. Firstly, let us note that $\psi$ has the following symmetries,
	\begin{equation}\label{eq:psisym}
	\psi(-\alpha)=\psi(\alpha)=\psi(\overline{\alpha}),
	\end{equation}
	and we thus merely have to compute the behaviour near the branching point $\alpha =u_1$ where $x_2$ and $X$ coalesce. By direct computation we have the corresponding Puiseux series 
	\begin{align}
	x_1&=\eta_1+\mathcal{O}(\alpha-u_1),\\
	x_2&=\eta_2-2r(\alpha-u_1)^{\frac{1}{2}}+\mathcal{O}(\alpha-u_1),\\
	X&=\eta_2+r(\alpha-u_1)^{\frac{1}{2}}+\mathcal{O}(\alpha-u_1),
	\end{align}
	as $\alpha\rightarrow u_1$, where
	\begin{equation*}
	\eta_1=-3\eta_2-2u_1,\quad \eta_2=\frac{1+3\nu^2-u_1^4}{4u_1(2+u_1^2)},
	\end{equation*}
	for a unique root $r$ of $r^2=-\tfrac{1}{3}\eta_2$. Therefore
	\begin{align}
	\psi(\alpha)&=-\frac{4}{15}\Re\left[ \frac{(x_2-x_1)^{\frac{1}{2}}}{x_2}(X-x_2)^{\frac{5}{2}}\right](1+\mathcal{O}(\alpha-u_1)^{\frac{1}{2}})\nonumber\\
	&=
	-\frac{4}{15}3^{\frac{5}{4}}\Re\left[\kappa(\alpha-u_1)^{\frac{5}{4}}\right](1+\mathcal{O}(\alpha-u_1)^{\frac{1}{2}}) \label{eq:puiseux}
	\end{align}
	where the branch of $z^{\frac{5}{4}}$ is taken real and positive on $\mathbb{R}_+$ with branch cut $\{\arg z=-\frac{3\pi i}{4}\}$, for a unique $\kappa=\kappa(\nu)$ which satisfies
	\begin{equation}\label{eq:kappa}
	\kappa^4=-(\eta_2-\eta_1)^2\eta_2.
	\end{equation} 
	In particular $\psi$ vanishes near $u_1$ and due to symmetries \eqref{eq:psisym} we know that $\psi$ vanishes near the other branching points as well, finishing the proof of part (1) of the proposition.
	
	We now turn our attention to the zero set $\{\psi(\alpha)=0\}$. We first compute it locally near the branching point $u_1$.  
	To this end, note that $\kappa=\kappa(\nu)$ in equation \eqref{eq:puiseux} is an algebraic function in $\nu$. It is the unique branch of equation \eqref{eq:kappa} satisfying $\arg{\kappa(\frac{1}{3})}=\frac{3}{16}\pi$. By direct computation it can be shown that $\arg{\kappa}:(0,\frac{1}{3}]\rightarrow (\frac{1}{6}\pi,\frac{3}{16}\pi]$.
	
	The right-hand side of \eqref{eq:puiseux} can be expanded into a complete Puiseux series and since the argument of $\kappa$ is bounded by $\frac{1}{20}\pi<\arg{\kappa}<\frac{9}{20}\pi$, it follows that there are three level curves of $\{\psi(\alpha)=0\}$ emanating radially from $\alpha=u_1$ with angles
	\begin{equation}\label{eq:angles}
	\frac{4}{5}\left(\frac{1}{2}\pi-\arg{\kappa}+m \pi i\right)\quad (m\in\{-1,0,+1\}).
	\end{equation}   
	Due to the symmetries \eqref{eq:psisym}, it follows that each branching point $u_k$ has precisely three level curves emanating from it. We call them, going around $u_k$ in anti-clockwise direction starting from the branch-cut, $\epsilon_k,l_k^*,\epsilon_{k+1}^*$. Note in particular that the radial angle at $u_k$: between $\epsilon_{k}$ and $l_k^*$ equals $\frac{4}{5}\pi$, between $l_k^*$ and $\epsilon_{k+1}^*$ equals $\frac{4}{5}\pi$ and between $\epsilon_{k+1}^*$ and $\epsilon_{k}$ equals $\frac{2}{5}\pi$,
	for $1\leq k\leq 4$, due to equation \eqref{eq:angles}.   
	
	Clearly
	\begin{equation*}
	\psi(\alpha)=\tfrac{1}{2}\Re[\alpha^2]+\mathcal{O}(\log|\alpha|)\quad (\alpha\rightarrow \infty)
	\end{equation*}
	and it is relatively straightforward to show that there are precisely four level curves $\{\psi(\alpha)=0\}$ emanating from $\alpha=\infty$ along the asymptotic directions $e^{\frac{\pi i}{4}(2k-1)}\infty$, $1\leq k\leq 4$. We call these level curves $l_k$, $1\leq k\leq 4$, in accordance with Figure \ref{fig:parametrisation2}.
	
	The final piece of information to deduce the correctness of Figure \ref{fig:parametrisation2} is that $\psi(\alpha)$ is strictly monotonic along the branch cuts: for $k=1,2$ the lower (upper) branch of $\psi$ along $[0,u_{k}]$ is strictly increasing (decreasing) as $\alpha$ traverses along it in the direction of $u_{k}$ whereas for $k=3,4$ the lower (upper) branch of $\psi$ along $[0,u_{k}]$ is strictly decreasing (increasing) as $\alpha$ traverses along it in the direction of $u_{k}$. Indeed this implies that both the upper and lower branches of $\psi$ along the branch cuts are nonzero, except for at the branching points $u_k,1\leq k \leq 4$. Thus the $\epsilon_k,\epsilon_k^*,l_k,l_k^*$ make up all the level curves of $\{\psi(\alpha)=0\}$, and since $\psi$ is harmonic, it is easy to deduce that we must have $l_k=l_k^*$ and $\epsilon_k=\epsilon_k^*$ for $1\leq k\leq 4$, where $\epsilon_1^*:=\epsilon_5^*$, yielding part (2) of the proposition. Furthermore, since $\epsilon_k=\epsilon_k^*$, 
	it follows that the internal radial angle between $\epsilon_{k}$ and $\epsilon_{k+1}$ at $u_k$ equals $\frac{2}{5}\pi$, for $1\leq k\leq 4$, establishing part (3).
	
	Due to Lemma \ref{lem:deformboutroux}, we know that the isomorphism class of the Stokes complex $\mathcal{C}(\alpha,B_{\Delta}(\alpha))$ is constant along each of the level curves $\epsilon_k,l_k$, $1\leq k\leq 4$. To prove part (4), we have to determine the isomorphism class on each level curve. We first compute the isomorphism classes on the $l_k$, $1\leq k\leq 4$, after which the isomorphism classes at the branching points $u_k$ and curves $\epsilon_k$, $1\leq k\leq 4$, are easily deduced.
	
	We proceed in computing the isomorphism class on $l_1$. Firstly, due to Lemma \ref{lem:boutrouxclas}, we know that the inner Stokes complex is connected along $l_1$. By equations \eqref{eq:x1x2asymptotic}, setting $\lambda=-\alpha+t$ with $t=\mathcal{O}(1)$ and $x_k=-\alpha+\widetilde{x}_k$ for $k=1,2$, we have
	\begin{equation*}
	V(\lambda)=\widetilde{V}(t)(1-\mathcal{O}(\alpha^{-1})),\quad \widetilde{V}(t):=(t-\widetilde{x}_1)(t-\widetilde{x}_2),
	\end{equation*}
	as $\alpha\rightarrow e^{\frac{\pi}{4}i}\infty$ along $l_1$. Furthermore note that $\widetilde{X}=X+\alpha$ and $\widetilde{0}=0+\alpha$ are asymptotic to $e^{\frac{\pi}{4}i}\infty$ in the same limit. The Stokes complex of the leading order potential $\widetilde{V}(t)$ is depicted in Figure \ref{fig:complexl1} under $\widetilde{L}_1$.
	
	Similarly, setting $\lambda=|\alpha|^{-1}s$ with $s=\mathcal{O}(1)$ and writing $\widehat{X}=|\alpha|^{-1} X$, we have
	\begin{equation*}
	V(\lambda)=\widehat{V}(s)(1-\mathcal{O}(\alpha^{-1})),\quad \widehat{V}(s):=i\frac{(s-\widehat{X})^2}{s^2},
	\end{equation*}
	as $\alpha\rightarrow e^{\frac{\pi}{4}i}\infty$ along $l_1$. Furthermore $\widehat{x}_{1,2}=|\alpha|x_{1,2}$ are asymptotic to $e^{-\frac{3\pi}{4}i}\infty$ in the same limit. The Stokes complex of the leading order potential $\widehat{V}(s)$ is depicted in Figure \ref{fig:complexl1} under $\widehat{L}_1$.
	
	There is only one isomorphism class consistent with above two limiting behaviours, namely the class $L_1$ defined in Figure \ref{fig:complexl1}. Thus on $l_1$ the Stokes complex falls in the class $L_1$.

	\begin{figure}[htpb]
		\centering
		\begin{tabular}{| c | c | c |}
			\hline
			\begin{tikzpicture}[scale=0.5]
			
			\node     at (0,2.5) {$L_1$};
			
			\tikzstyle{star}  = [circle, minimum width=4pt, fill, inner sep=0pt];
			\node[blue,star]     at (1,0){};
			\node[right]     at (1,0) {$x_2$};

			\node[blue,star]     at (-1,0) {};
			\node[left]     at (-1,0) {$x_1$};

			\draw[blue,thick] (-1,0) -- (1,0);

			\draw[blue,thick] (1,0) -- (4,3);
			\draw[blue,thick] (1,0) -- (2,-1);
			\draw[blue,thick] (-1,0) -- (-2,1);
			\draw[blue,thick] (-1,0) -- (-2,-1);
			
			\node[blue,star]     at (3,2) {};
			\node[below right]     at (2.95,2.05) {$X$};
			
			\node[blue,star]     at (2.3,2.7) {};
			\node[above]     at (2.3,2.7) {$0$};
			
			\draw [blue,thick] plot [smooth,tension=1.0] coordinates 
			{(3,2) (2.85,3.25) (1.7, 3.3) (1.75, 2.15) (3,2) };
			
			\end{tikzpicture} 
			&
			
			\begin{tikzpicture}[scale=0.5]
			
			\node     at (0,2.6) {$\widetilde{L}_1$};
			
			\tikzstyle{star}  = [circle, minimum width=4pt, fill, inner sep=0pt];
			\node[blue,star]     at (1,0){};
			\node[right]     at (1,0) {$\hspace{0.1cm}\widetilde{x}_2$};

			\node[blue,star]     at (-1,0) {};
			\node[left]     at (-1,0) {$\widetilde{x}_1\hspace{0.1cm}$};

			\draw[blue,thick] (-1,0) -- (1,0);

			\draw[blue,thick] (1,0) -- (2,1);
			\draw[blue,thick] (1,0) -- (2,-1);
			\draw[blue,thick] (-1,0) -- (-2,1);
			\draw[blue,thick] (-1,0) -- (-2,-1);

			\end{tikzpicture} 
			&
			
			\begin{tikzpicture}[scale=0.5]
			
			\node     at (2.8,4.5) {$\widehat{L}_1$};
			
			\tikzstyle{star}  = [circle, minimum width=4pt, fill, inner sep=0pt];

			\draw[blue,thick] (2,1) -- (4,3);

			\node[blue,star]     at (3,2) {};
			\node[below right]     at (2.95,2.05) {$\widehat{X}$};
			
			\node[blue,star]     at (2.3,2.7) {};
			\node[above]     at (2.3,2.7) {$0$};
			
			\draw [blue,thick] plot [smooth,tension=1.0] coordinates 
			{(3,2) (2.85,3.25) (1.7, 3.3) (1.75, 2.15) (3,2) };
			
			\end{tikzpicture} \\
			\hline
		\end{tabular}
		\caption{Topological representation of isomorphism classes of Stokes complexes $L_1$, $\widetilde{L}_1$ and $\widehat{L}_1$ corresponding respectively to the potentials $V(\lambda;\alpha,B_\Delta(\alpha))$ with $\alpha \in l_1$, $\widetilde{V}(t)$ with $\alpha=e^{\frac{\pi i}{4}}\infty$  and $\widehat{V}(s)$ with $\alpha=e^{\frac{\pi i}{4}}\infty$.}
		\label{fig:complexl1}
	\end{figure}
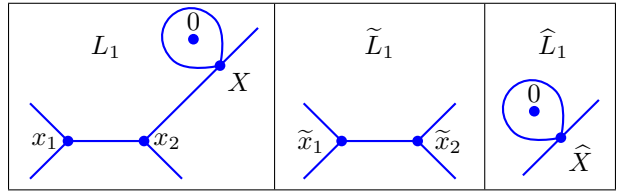

	As $\alpha\rightarrow u_1$ along $l_1$, we know that $x_2$ merges with $X$, and thus the resulting Stokes complex of $V(\lambda)$ at $(\alpha,\beta)=(u_1,B_\Delta(u_1))$ is given by $C_1$, defined in Figure \ref{fig:stokes_complexes}. We conclude that $c_1=u_1$ and hence $\widetilde{c}_1= (u_1,B_\Delta(u_1))$ by equation \eqref{eq:cornerequation}.
	Using the symmetries \eqref{eq:csym} we obtain equation \eqref{eq:paracorner} for all $1\leq k\leq 4$.
	
	Next we consider the Stokes complex of $V(\lambda;\alpha,B_\Delta(\alpha))$ along $\epsilon_2$. Again we know it's inner Stokes complex must be connected due to Lemma \ref{lem:boutrouxclas}. Furthermore, we have
	\begin{itemize}
		\item as $\alpha\rightarrow u_1=c_1$ along $\epsilon_2$, the turning point $x_2$ merges with $X$ and the resulting Stokes complex falls in the class $C_1$; 
		\item as $\alpha\rightarrow u_2=c_2$ along $\epsilon_2$, the turning point $x_1$ merges with $X$ and the resulting Stokes complex falls in the class $C_2$;
		\item The isomorphism class of the Stokes complex is invariant under reflection in the imaginary axes (and interchanging of marked points $\infty_1\leftrightarrow \infty_2$ and $\infty_3\leftrightarrow \infty_4$).
	\end{itemize}
	Indeed the third easily follows from the symmetries given in equation \eqref{eq:psisym}. Clearly, there is only one isomorphism class of Stokes complexes which satisfies these three conditions, namely $E_2$ defined in Figure \ref{fig:stokes_complexes}, with $\mu_{1,2}=x_{1,2}$ and $\mu_3=X$. Thus the Stokes complex along $\epsilon_2$ falls in the class $E_2$ and we have 
	\begin{equation}\label{eq:parasub}
	\epsilon_2\subseteq e_2,\quad\{(\alpha,B_\Delta(\alpha)):\alpha \in \epsilon_2\}\subseteq\widetilde{e}_2.
	\end{equation}
	Since $\widetilde{e}_2$ is a smooth non-self intersecting curve with end-points $\widetilde{c}_1$ and $\widetilde{c}_2$, and the same holds true for $\{(\alpha,B_\Delta(\alpha)):\alpha \in \epsilon_2\}$, equation \eqref{eq:parasub} implies $\{(\alpha,B_\Delta(\alpha)):\alpha \in \epsilon_2\}=\widetilde{e}_2$ and thus $\epsilon_2=e_2$. Analogously we prove equation \eqref{eq:paraedge} for the remaining $k\in\{1,3,4\}$ and with part (4) thus established, the proposition is proven.
\end{proof}

\end{pro}
\begin{rem}
	We note that $\psi(-\alpha)=\psi(\alpha)$. Furthermore, if $\nu=\frac{1}{3}$, then $\psi(e^{\frac{\pi i}{4}}\alpha)=-\psi(\alpha)$ and $l_k$ simply equals the straight diagonal line $\{u_k+te^{\frac{\pi i}{4}(2k-1)}:t\in \mathbb{R}_+\}$ for $1\leq k\leq 4$.
\end{rem}

\subsubsection*{Step \ref{step5}}
The following corollary constitutes step \ref{step5}, namely it justifies the description of the elliptic region given in Section \ref{sec:results}.
\begin{cor}\label{cor:ellipticregion}
 The elliptic region $K_a$ is the compact simply-connected quadrilateral domain whose boundary equals the Jordan curve
	\begin{equation}\label{eq:ellipticboundary}
	\partial K_a=e_1\sqcup e_2\sqcup e_3\sqcup e_4 \sqcup\{c_1,c_2,c_3,c_4\}.
	\end{equation}
	The boundary consists of four analytic edges $e_k$, $1\leq k\leq 4$, meeting at four corners $c_k$, $1\leq k\leq 4$ with internal angle $\frac{2\pi}{5}$, as in Figure \ref{fig:marked}. Furthermore, it is implicitly parametrised as part of the zero set $\{\psi(\alpha)=0\}$, as detailed in parts (2) and (4) of Proposition \ref{pro:boundary}.
\end{cor}
\begin{proof}
Recall that, by definition, $R_a=\Pi_a(R)$ and $K_a=\Pi_a(K)$. Since $\Pi_a$ is continuous, $K$ is compact and $R_a$ is open, by Lemma \ref{lem:deformation}, we know that
\begin{equation*}
\partial R_a\subseteq \Pi_a(\delta R).
\end{equation*}
 Due to Proposition \ref{pro:boundary}, we know that
 \begin{equation*}
 \Pi_a(\delta R)=e_1\sqcup e_2\sqcup e_3\sqcup e_4 \sqcup\{c_1,c_2,c_3,c_4\}
 \end{equation*} 
is a Jordan curve, and thus
\begin{equation*}
\partial R_a\subseteq  \Pi_a(\delta R)=e_1\sqcup e_2\sqcup e_3\sqcup e_4 \sqcup\{c_1,c_2,c_3,c_4\}.
\end{equation*}
Let us denote by $I$ the interior of the Jordan curve $\Pi_a(\delta R)$. We wish to show that $R_a=I$. To this end, let us note that $R_a\ni 0$ is non-empty, open and $\partial R_a\subseteq  \Pi_a(\delta R)$. But there exist only one such set, namely $R_a=I$. Therefore $\partial R_a=\partial I=\Pi_a(\delta R)$ which proves equation \eqref{eq:ellipticboundary}.

Furthermore, since $K$ is compact, 
\begin{equation*}
K_a=\overline{R_a}=\overline{I}=I\sqcup \Pi_a(\delta R),
\end{equation*}
and thus the elliptic region $K_a$ is indeed the compact simply-connected quadrilateral domain whose boundary equals the Jordan curve
$\Pi_a(\delta R)$. The remainder of the corollary is the content of (2) and (4) of Proposition \ref{pro:boundary}.
\end{proof}

\subsubsection*{Step \ref{step6}}
Step \ref{step6} is embodied by the following
\begin{pro}\label{pro:projection}
	The projection $\Pi_a:K\rightarrow K_a$ is a homeomorphism which maps the border $\delta R$ of $R$ homeomorphically onto the boundary $\partial K_a$ of the elliptic region and maps $R$ diffeomorphically onto $R_a$.
\end{pro}
\begin{proof}
	Firstly, note that Corollary \ref{cor:ellipticregion} and part (4) of Proposition \ref{pro:boundary} imply that
	\begin{equation*}
	\Pi_a|_{\delta R}:\delta R\rightarrow \partial R_a
	\end{equation*}
	is a homeomorphism.
	
	Consider now the mapping $g=\Pi_a\circ \mathcal{S}^{-1}:Q\rightarrow K_a$. It follows from Lemma \ref{lem:deformation} and Proposition \ref{prop:Sextension} that $g|_{Q^{\circ}}:Q^{\circ}\rightarrow R_a$ is a local diffeomorphism. Furthermore, $g$ maps the jordan curve $\partial Q$ homeomorphically onto the Jordan curve $\partial R_a$.
	
	 It follows from the latter and compactness of $K$ that $g|_{Q^{\circ}}:Q^{\circ}\rightarrow R_a$ is proper. Therefore $g^{-1}(\alpha)$ is finite for any $\alpha\in R_a$, and since $g|_{Q^{\circ}}$ is locally diffeomorphic, it follows that $g|_{Q^{\circ}}:Q^{\circ}\rightarrow R_a$ is a covering map. Since $Q^{\circ}$ and $R_a$ are simply connected, $g|_{Q^{\circ}}$ must be injective and thus $g$ is a continuous bijection which maps $Q^{\circ}$ diffeomorphically onto $R_a$. As $K$ is compact, $g$ maps closed sets to closed sets and thus $g$ is a homeomorphism.	The proposition now follows since $\Pi_a=g\circ \mathcal{S}$.
\end{proof}

\subsubsection*{Step \ref{step7}}
We can now put everything together to prove Theorem \ref{thm:mappingS}, and thus complete the final step, step \ref{step7}.

\begin{proof}[Proof of Theorem \ref{thm:mappingS}]	
Due to Lemma \ref{lem:extension}, the mapping $\mathcal{S}$ has a unique continuous extension to $K$ with range in $Q$ and, by Proposition \ref{prop:Sextension}, $\mathcal{S}$ maps $K$ homeomorphically to $Q$. Due to Proposition \ref{pro:projection}, the projection $\Pi_a:K\rightarrow K_a$ is a homeomorphism.

We may now consider the composition
\begin{equation*}
\mathcal{S}_a=\mathcal{S}\circ \Pi_a^{-1}:K_a\rightarrow Q,
\end{equation*}
which is a homeomorphism. Furthermore, due to Propositions \ref{prop:Sextension} and \ref{pro:projection}, it maps the interior $R_a$ of the elliptic region diffeomorphically onto the interior of the quadrilateral $Q$.

Let $1\leq k\leq 4$, then, due to part (4) of Proposition \ref{pro:boundary}, we know that $\Pi_a$ maps the edge $\widetilde{e}_k$ of $K$ diffeomorphically onto the edge $e_k$ of $K_a$ and similarly, due to Proposition \ref{prop:Sextension}, $\mathcal{S}$ maps the edge $\widetilde{e}_k$ diffeomorphically onto the open edge $\widehat{e}_k$ of $Q$. Thus $\mathcal{S}_a$ maps the edge $e_k$ of the elliptic region diffeomorphically onto the open edge $\widehat{e}_k$ of $Q$. Similarly, it follows that $\mathcal{S}_a$ maps the corner $c_k$ of the elliptic region to the corner $\widehat{c}_k$ of $Q$. The theorem follows.
\end{proof}

%

\subsection{Anti-Stokes Complexes}\label{subsec:anti}
In this subsection we discuss the concept of anti-Stokes complexes associated with the potential $V$ and characterise them topologically for all
$(\alpha,\beta) \in R$. Anti-Stokes complexes do not play an explicit role in the proof of
Theorem \ref{thm:mappingS}, but they do play a vital role in the WKB analysis in Section \ref{sec:WKB}.

Take any $\lambda_0$ which is not a critical point of the potential $V$, and consider the action integral \eqref{eq:actionintegral} on the universal covering of the complex plane minus critical points.
Consider the level curve $j_{\lambda_0}$ through $\lambda_0$, defined by $\Im S(\lambda_0,\lambda)=0$, on this universal covering. Due to equation \eqref{eq:actionregular}, this level curve is a line on the universal covering. Furthermore, since the action integral can differ by at most a sign on different sheets, it is clear the the projection of $\widetilde{j}_{\lambda_0}$ of $j_{\lambda_0}$ onto the complex plane minus critical points is diffeomorphic to either a line or a circle. Furthermore, since the residue of $\sqrt{V(\lambda)}$ near the only pole $\lambda=0$ is real, $\widetilde{j}_{\lambda_0}$ cannot be diffeomorphic to a circle. On the other hand, $\widetilde{j}_{\lambda_0}$ can emanate radially from $\lambda=0$. 

 Analogous to Lemma \ref{lem:stokes_curves}, we have the following trichotomy.
\begin{lem}\label{lem:stokes_curvesanti}
	Let $\lambda_0$ be a generic point of the potential \eqref{eq:potential} and $\widetilde{j}_{\lambda_0}$
	its corresponding projected level curve. Then $\widetilde{j}_{\lambda_0}$ is diffeomorphic to a line. Let $x\mapsto \widetilde{j}_{\lambda_0}(x)$
	be a diffeomorphism of $\mathbb{R}$ onto $\widetilde{j}_{\lambda_0}$, then, for $\epsilon\in \{\pm 1\}$,
	we have the following trichotomy: either
	\begin{enumerate}[label=(\roman*)]
		\item $\lim_{x\rightarrow \epsilon \infty}=\infty$ and the curve is asymptotic to one of the
		four semi-axes $i^k\mathbb{R}_+$, $k\in\mathbb{Z}_4$, in the $\lambda$-plane;
		\item $\lim_{x\rightarrow \epsilon \infty}=\lambda_*$ with $\lambda_*$ a turning point of the potential; or
		\item $\lim_{x\rightarrow \epsilon \infty}=0$.
	\end{enumerate}
\end{lem}
\begin{proof}
	See Strebel \cite{strebel}.
\end{proof}
In alignment with the above lemma, we make the following definition.
\begin{defi}\label{defi:dichotomyanti}
	Considering the trichotomy in Lemma \ref{lem:stokes_curvesanti} for a projected level curve $\widetilde{j}_{\lambda_0}$, we call,  in case (i), (ii) or (iii) respectively, the asymptotic direction $i^k\infty$, the turning point $\lambda_*$ or $0$ an endpoint of $\widetilde{j}_{\lambda_0}$.
	We call $\widetilde{j}_{\lambda_0}$ an anti-Stokes line if at least one endpoint is a turning point.
\end{defi}

To define the anti-Stokes complex associated with the potential $V$, we add four additional marked points to the space $\mathbb{C}_\infty$, see Definition \ref{def:stokescomplex},
\begin{equation*}
\widetilde{\mathbb{C}}_\infty=\mathbb{C}_\infty\sqcup\{\pm \infty, \pm i\infty\},
\end{equation*}
and endow it with the topology that makes the extension of the mapping $L$ defined in equation \eqref{eq:defiL}, by $L(\pm \infty)=\pm 1$ and $L(\pm i\infty)=\pm i$, a homeomorphism. Furthermore, we promote $\lambda=0$ to a marked point in $\widetilde{\mathbb{C}}_\infty$.

\begin{defi}\label{def:stokescomplexanti}
	The anti-Stokes complex $\widetilde{\mc{C}}=\widetilde{\mc{C}}(\alpha,\beta)\subseteq \widetilde{\mathbb{C}}_\infty$ of the potential
	\eqref{eq:potential} is the union of the marked points $\{0,\pm \infty,\pm i\infty\}$, the anti-Stokes lines and turning points
	of the potential. \\
	Topologically, the anti-Stokes complex $\widetilde{\mc{C}}$ is an embedded graph into $\widetilde{\mathbb{C}}_\infty$ with
	vertices equal to the turning points and the marked points $\{0,\pm \infty,\pm i\infty\}$, with edges given by the anti-Stokes lines.
\end{defi}

\begin{figure}[ht]
	\centering
	\includegraphics[width=10cm]{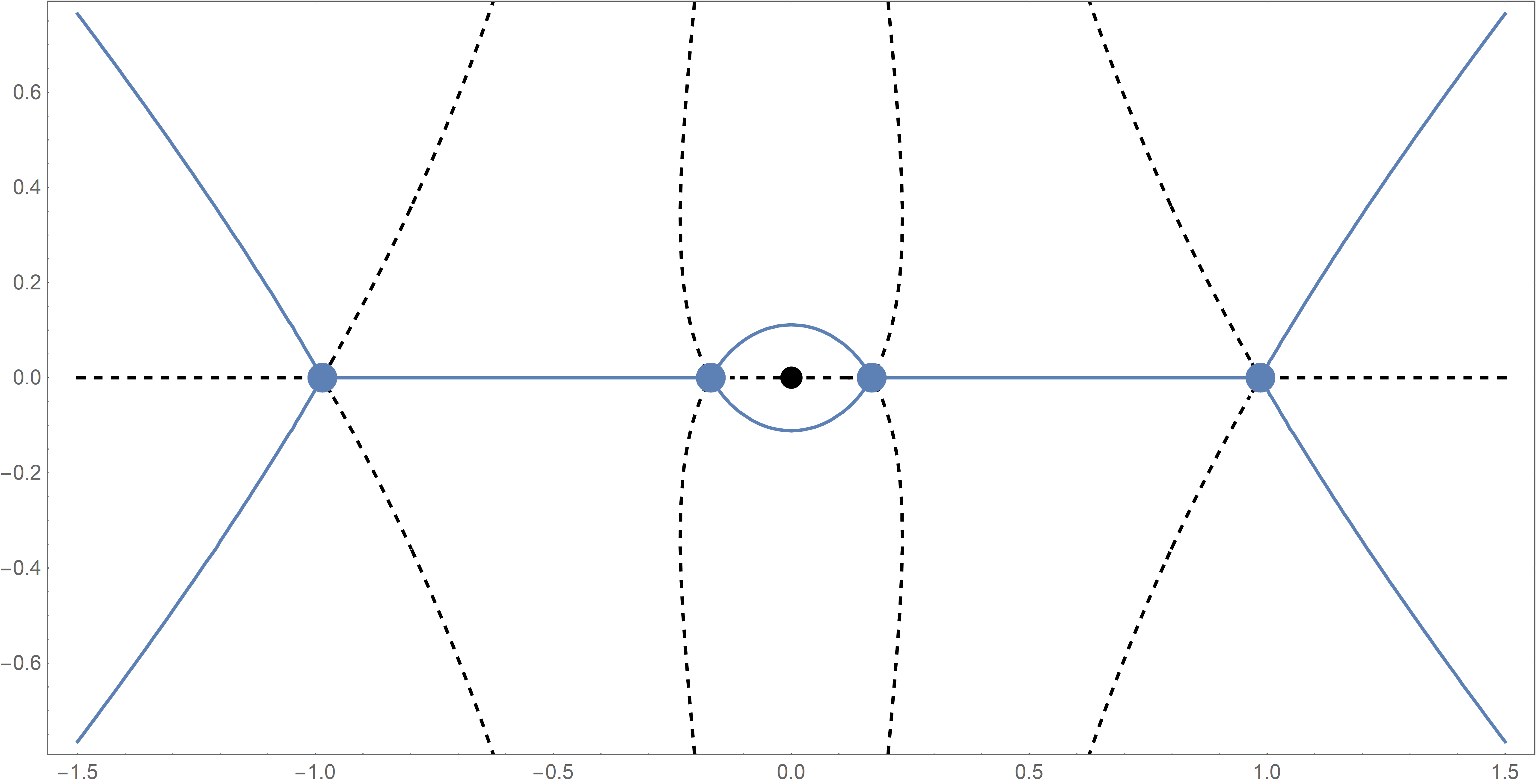} 
	\caption{Numerical plot of Stokes complex in blue and anti-Stokes complex in dashed purple of potential $V$ with $(\alpha,\beta)=(0,0)$ and $\nu=\tfrac{1}{3}$.
	} \label{fig:stokesrealplot}
\end{figure}

In Figure \ref{fig:stokesrealplot} a numerical plot of the Stokes complex and anti-Stokes complex is given of the potential $V$ with $(\alpha,\beta)=(0,0)$ and $\nu=\tfrac{1}{3}$. In the following proposition, we show that this numerical figure in fact represents, topologically speaking, the general picture when $(\alpha,\beta)\in R$.

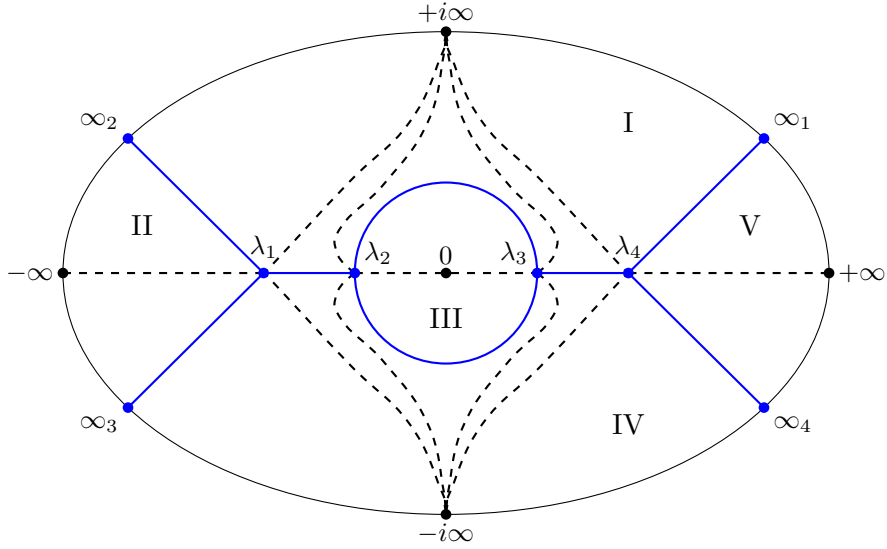
\begin{figure}[htpb]
	\centering
	\begin{tikzpicture}[scale=0.8]
	\tikzstyle{star}  = [circle, minimum width=4pt, fill, inner sep=0pt];
	
	\draw[black] (0,0) ellipse (6.3cm and 4cm);
	
	\draw[black,dashed,thick] (3,0) -- (6.3,0);
	\draw[black,dashed,thick] (-6.3,0) -- (-3,0);
	\draw[black,dashed,thick] (-1.5,0) -- (1.5,0);

	\draw [black,thick,dashed] plot [smooth,tension=0.8] coordinates {(1.5,0) (1.8,0.35) (1.7,0.8) (1,1.5) (0.3, 2.5) (0.04,3.6) (0.015,3.8) (0,4) };
	\draw [black,thick,dashed] plot [smooth,tension=0.8] coordinates {(3,0) (2.8,0.2) (1.6,1.5) (0.6, 2.5) (0.08,3.6) (0.03,3.8) (0,4) };
	
	\draw [black,thick,dashed] plot [smooth,tension=0.8] coordinates {(-1.5,0) (-1.8,0.35) (-1.7,0.8) (-1,1.5) (-0.3, 2.5) (-0.04,3.6) (-0.015,3.8) (0,4) };
	\draw [black,thick,dashed] plot [smooth,tension=0.8] coordinates {(-3,0) (-2.8,0.2) (-1.6,1.5) (-0.6, 2.5) (-0.08,3.6) (-0.03,3.8) (0,4) };
	
	\draw [black,thick,dashed] plot [smooth,tension=0.8] coordinates {(1.5,0) (1.8,-0.35) (1.7,-0.8) (1,-1.5) (0.3,-2.5) (0.04,-3.6) (0.015,-3.8) (0,-4) };
	\draw [black,thick,dashed] plot [smooth,tension=0.8] coordinates {(3,0) (2.8,-0.2) (1.6,-1.5) (0.6,-2.5) (0.08,-3.6) (0.03,-3.8) (0,-4) };
	
	\draw [black,thick,dashed] plot [smooth,tension=0.8] coordinates {(-1.5,0) (-1.8,-0.35) (-1.7,-0.8) (-1,-1.5) (-0.3,-2.5) (-0.04,-3.6) (-0.015,-3.8) (0,-4) };
	\draw [black,thick,dashed] plot [smooth,tension=0.8] coordinates {(-3,0) (-2.8,-0.2) (-1.6,-1.5) (-0.6,-2.5) (-0.08,-3.6) (-0.03,-3.8) (0,-4) };

	\node[blue,star]     at (-3,0){};
	\node[above]     at (-3,0.1) {$\lambda_1$};
	
	\node[blue,star]     at (-1.5,0) {};
	\node[above right]     at (-1.5,0) {$\lambda_2$};
	
	\node[black,star]     at (0,0) {};
	\node[above]     at (0,0) {0};
	
	\node[blue,star]     at (1.5,0) {};
	\node[above left]     at (1.5,0) {$\lambda_3$};
	
	\node[blue,star]     at (3,0) {};
	\node[above]     at (3,0.1) {$\lambda_4$};
	
	\draw[blue,thick] (-3,0) -- (-1.5,0);
	\draw[blue,thick] (1.5,0) -- (3,0);
	\draw[blue,thick] (0,0) circle (1.5cm);
	
	\draw[blue,thick] (-3,0) -- (-5.23,2.23);
	\node[blue,star]     at (-5.23,2.23) {};
	\node[above left]     at (-5.23,2.23) {$\infty_2$};
	
	\draw[blue,thick] (3,0) -- (5.23,2.23);
	\node[blue,star]     at (5.23,2.23) {};
	\node[above right]     at (5.23,2.23) {$\infty_1$};
	
	\draw[blue,thick] (-3,0) -- (-5.23,-2.23);
	\node[blue,star]     at (-5.23,-2.23) {};
	\node[below left]     at (-5.23,-2.23) {$\infty_3$};
	
	\draw[blue,thick] (3,0) -- (5.23,-2.23);
	\node[blue,star]     at (5.23,-2.23) {};
	\node[below right]     at (5.23,-2.23) {$\infty_4$};
	
	\node[right]     at (6.3,0) {$+\infty$};
	\node[black,star]     at (6.3,0) {};
	\node[left]     at (-6.3,0) {$-\infty$};
	\node[black,star]     at (-6.3,0) {};	
	\node[above]     at (0,4) {$+i\infty$};
	\node[black,star]     at (0,4) {};	
	\node[below]     at (0,-4) {$-i\infty$};	
	\node[black,star]     at (0,-4) {};
	
	\node     at (5,0.8) {{\large V}};
	\node     at (-5,0.8) {{\large II}};
	\node     at (3,2.5) {{\large I}};
	\node     at (3,-2.5) {{\large IV}};
	\node     at (0,-0.75) {{\large III}};

	\end{tikzpicture} 
	\caption{The Stokes complex $\mathcal{C}(\alpha,\beta)$ in blue and anti-Stokes complex $\widetilde{\mathcal{C}}(\alpha,\beta)$ in dashed black as embedded graphs in $\widetilde{\mathbb{C}}_\infty$, with $(\alpha,\beta)\in R$. Furthermore, I-V label the regions in the complex plane cut out by the Stokes complex, as in Figure \ref{fig:contours}.} \label{fig:complex00anti}
\end{figure}

\begin{pro}\label{pro:antiS}
	Let $(\alpha,\beta)\in R$ and label by I, II, III, IV and V the regions of the complex plane cut out by the Stokes complex $\mc{C}(\alpha,\beta)$ as in Figure \ref{fig:complex00}. Then, the anti-Stokes complex $\widetilde{\mc{C}}(\alpha,\beta)$ is topologically accurately described by Figure \ref{fig:complex00anti}, namely
	\begin{enumerate}[label=(\roman*)]
		\item The intersection of the Stokes and anti-Stokes complexes is given by the turning points of the potential;
		\item From each turning point emanate three Stokes and three anti-Stokes lines.
		\item From both the turning points $\lambda_2$ and $\lambda_3$ emanates precisely one anti-Stokes line that has $0$ as its other endpoint and lies	entirely in region III.
		\item From each turning point emanates exactly one anti-Stokes line which lies entirely in region I, each with other endpoint $+i\infty$.
		\item From each turning point emanates exactly one anti-Stokes line which lies entirely in region IV, each with other endpoint $-i\infty$.
		\item From $\lambda_2$ emanates precisely one anti-Stokes line which lies entirely in region II, with other endpoint $-\infty$.
		\item From $\lambda_4$ emanates precisely one anti-Stokes line which lies entirely in region V, with other endpoint $+\infty$.
	\end{enumerate}
\end{pro}
\begin{proof}
	Given any turning point, it is simple and a local analysis shows that it has three Stokes lines and three anti-Stokes lines emanating radially, with, going around the turning point, alternatively Stokes and anti-Stokes lines, separated radially by an angle $\frac{\pi}{3}$. In particular (ii) holds true.
	
	Now, let us take the turning point $\lambda_4$. Clearly it has one anti-Stokes line emanating into region I, one into region IV and one into region V, which we will denote respectively by $l_\text{I}^{(4)}$, $l_\text{IV}^{(4)}$ and $l_\text{V}^{(4)}$. Note that $l_\text{I}^{(4)}$ cannot escape region I and must therefore asymptotically connect with $+i\infty$. Similarly, $l_\text{IV}^{(4)}$ cannot escape region IV and must therefore asymptotically connect with $-i\infty$. Finally $l_\text{V}^{(4)}$ must lie entirely in region V and asymptotically connect with $+\infty$.
	
	Next, let us consider the turning point $\lambda_3$ and denote by $l_\text{I}^{(3)}$, $l_{\text{III}}^{(3)}$ and $l_\text{IV}^{(3)}$ the anti-Stokes lines emanating from it respectively into regions I, III and IV. Clearly  $l_\text{I}^{(3)}$ must lie entirely in region I and connect asymptotically with $+i\infty$. Similarly $l_\text{IV}^{(3)}$  lies entirely in region IV and connects asymptotically with $-i\infty$.
	
	Next we look at $l_{\text{III}}^{(3)}$. Let us denote by $W$ the closure of region III.
	  Then the action integral $S(\lambda_3,\lambda)$, see equation \eqref{eq:actionintegral}, is a well-defined analytic mapping on the universal covering of $W$ punctured at $\lambda=0$ to the complex plane. In fact, see Strebel \cite{strebel}, for a unique choice of sign, the action integral $S$ maps the universal covering of $W\setminus\{0\}$ homeomorphically onto the closed half-plane $\{\Re s\geq 0\}$, and the interior conformally onto the open half-plane $\{\Re s> 0\}$. Let $\pi$ be the projection of the universal covering onto $W\setminus\{0\}$. Then $S$ maps the points in the universal covering above $\lambda_3\in W$, i.e. $\pi^{-1}(\lambda_3)$, bijectively onto $\{k\nu \pi i:k\in\mathbb{Z}\}$, and $S$ maps $\pi^{-1}(\lambda_2)$ bijectively onto the set $\{r_2i+k\nu \pi i:k\in\mathbb{Z}\}$, where $r_2>0$ defined in the proof of Lemma \ref{lem:range}. Furthermore,
	  \begin{equation*}
	   \pi(S^{-1}(s_R+i s_I))\rightarrow 0\quad (s_R\rightarrow +\infty),
	  \end{equation*}
       for any fixed $s_I\in\mathbb{R}$.	   In particular, the anti-Stokes line $l_{\text{III}}^{(3)}$ is smoothly parametrised by 
	  \begin{equation*}
	  \mathbb{R}_+\rightarrow l_{\text{III}}^{(3)},s_R\mapsto \pi(S^{-1}(s_R))
	  \end{equation*}
	  and thus $l_{\text{III}}^{(3)}$ lies entirely in region III and has endpoints $\lambda_3$ and $0$.


	Similarly there are anti-Stokes lines $l_\text{I}^{(2)}$, $l_\text{III}^{(2)}$  and $l_\text{IV}^{(2)}$ emanating from $\lambda_2$, lying respectively in regions I, III and IV and connecting asymptotically with $+i\infty$, $0$ and $-i\infty$. Also, there are anti-Stokes lines $l_\text{I}^{(1)}$, $l_\text{II}^{(1)}$  and $l_\text{IV}^{(1)}$ emanating from $\lambda_1$, lying in regions I, II and IV respectively and connecting asymptotically with $+i\infty$, $-\infty$ and $-i\infty$.
	The topological description of the Stokes and anti-Stokes complexes in the proposition follows.	\end{proof}

\section{WKB asymptotics}\label{sec:WKB}
In the present section we study the inverse monodromy problem of equation \eqref{eq:scaled},
which characterises roots of generalised
Hermite polynomials in the rescaled variables $\alpha,\beta$, in the large $E$ limit.

We prove, by means of the complex WKB method, that the quantisation condition and the apparent singularity condition asymptotically read respectively
\begin{subequations}\label{eq:bohrsommerfeld}
 \begin{align}\label{eq:wkbcondition1}
& E\oint_{\gamma_1}\sqrt{V(\la;\alpha,\beta,\nu)}d\la = (2k+1) \pi, k \in \bb{Z} \\ \label{eq:wkbcondition2}
& E\oint_{\gamma_2}\sqrt{V(\la;\alpha,\beta,\nu)}d\la = 2l \pi , l \in \bb{Z} \mbox{ if } n \mbox { odd }, \;  l \in \bb{Z}+\frac12\mbox{ if } n \mbox { even }
 \end{align}
\end{subequations}
 where $V(\la;\alpha,\beta,\nu)$ is the potential \eqref{eq:scaledpotential}. Equations (\ref{eq:wkbcondition1},\ref{eq:wkbcondition2})
 are a pair of Bohr-Sommerfeld quantisation conditions, which form a complete system for $(\alpha,\beta)\in R$.
The same system, for the meromorphic differential $y \, dx$ over the elliptic curve $\{y^2=x^3+a x+b\}$,
was named the Bohr-Sommerfeld-Boutroux system in \cite{piwkb}, and it was shown to describe the poles of the Tritronqu\'ee solution of
Painleve-I equation in a neighbourhood of infinity.

In principle, the parameters $E$ and $\nu$ of the potential have the form $E=2m+n$ and $\nu=\frac{n}{E}$ with 
$n\leq m$ are positive integers. However, in the present section and unless otherwise stated, the parameters
$E,\nu$ are arbitrary positive real numbers, since most results hold true irrespective of their integer or rational character.

The precise formulation of equations (\ref{eq:wkbcondition1},\ref{eq:wkbcondition2}), together with an explicit error bound,
are stated and proven in theorems \ref{thm:wkbquantisation} and \ref{thm:wkbmonodromy}.

\subsection{Subdominant solutions}
Before we proceed with the WKB analysis to obtain quantitative estimates, we make a simple, albeit non-trivial,
qualitative study of \eqref{eq:scaled}.
Such a qualitative study is based on a well-established and well-known theory;
general references are \cite{piwkb,maro18,sibuya75,wasow12,eremenko}.

We consider here the following family of anharmonic potentials, which generalises \eqref{eq:scaled} slightly,
\begin{align}\label{eq:genharmonic}
 &y''(\la)=Q(\la)y(\la),\\
  &Q(\la)=a_0 \la^2+a_1 \la + a_2 + \frac{a_3}{\la}+ \frac{n^2-1}{4 \la^2} \qquad (a_0,n \in \bb{R}^+, a_1,a_2,a_3 \in \C).\nonumber
\end{align}
The above equation is regular outside $0$ and $\infty$: $0$ is a regular (or Fuchsian) singularity
and $\infty$ is an irregular singularity. Therefore solutions are in general
multi-valued. They are however single-valued if we restrict them to a simply-connected open subset of $\C^*$.
\begin{defi}
For $\varphi \in \bb{R}$, we
denote by $r_{\varphi}$  the ray of angle $\varphi$, $r_{\varphi}:=\lbrace e^{i\varphi} t,
t \in (0,\infty) \rbrace \subset \C^*$.
Moreover, we denote by $D_{\varphi}$ the plane cut by any line connecting $0$ and $\infty$, which for large $|\la|$ is asymptotic
to the ray $r_{\varphi}$.
\end{defi}

\subsubsection*{The singularity $\la=0$}
At $0$ a distinguished basis of solutions is provided by the Frobenius series
\begin{equation}\label{eq:Frobenius}
 \chi_+=\la^{\frac{n+1}{2}} \big( 1 + f(\la)\big),\quad \chi_-= \la^{\frac{-n+1}{2}} \big( 1 +g(\la)\big)+ C (\log \la) \chi_+,
\end{equation}
where $f,g$ are holomorphic and vanish at $\la=0$, and $C \in \C$ is a constant which is zero if $n\notin \mathbb{Z}$.
We notice that the solution $\chi_+$ is subdominant
at $\la=0$, irrespective of the direction along which the limit is taken, and this property characterises it up to a normalising constant.

Let $M$ be the monodromy operator: $M \psi(\la)=\psi(e^{2\pi i}\la)$.
If $n\notin \mathbb{Z}$, it is straightforward to verify that $\chi_+,\chi_-$ are eigenvectors of $M$:
$M \chi_+=-e^{ i \pi n }\chi_+$,
and $M \chi_-=-e^{- i \pi n }\chi_-$.
If $n$ is an integer, the monodromy operator $M$ is represented by the following matrix in the basis
$\lbrace \chi_+,\chi_-\rbrace$:
$$
M= (-1)^{n+1} \begin{pmatrix}1 & 0\\
2\pi i \, C & 1
\end{pmatrix}
$$
The above matrix has a single eigenvalue, $(-1)^{n+1}$, and it is not diagonalizable unless
$C=0$, in which case $\la=0$ is called an apparent singularity.
\subsubsection*{The singularity $\la=\infty$}
Infinity is an irregular singularity and since $Q(\la)=a_0 \la^2+ O(\la)$ as $\la \to \infty$ and $a_0 >0$, solutions look locally like
the solutions of the harmonic oscillator. We therefore define the four Stokes sectors
$\Sigma_k$, $k \in \mathbb{Z}_4$,
of  opening $\frac{\pi}{2}$ and centred at $\arg \la=\frac{\pi}{2}k$:
\begin{equation}\label{eq:stokessectors}
 \Sigma_k=\left\lbrace \la \in \bb{C}, \left|\arg \la -\frac{\pi}{2}k\right|<\frac{\pi}{4} \right\rbrace, \quad k \in \bb{Z}_4.
\end{equation}

\begin{pro}\label{pro:subdominants}\cite{sibuya75}
Choose a $k\in \Zq$ and a $\theta \in \bb{R}$ such that $r_\theta \subset \Sigma_k$. 
There exists a unique (up to a normalising constant) non-trivial solution $\psi:r_{\theta} \to \bb{C}$ of
equation \eqref{eq:genharmonic} such that
$\lim_{\la \to \infty}\psi(\la)=0$. 
Let $D=\Sigma_{k-1}\cup \overline{\Sigma}_k \cup \Sigma_{k+1}$, where $\overline{\Sigma}_k$ is the closure of $\Sigma_k$.
The analytic extension of $\psi(\la)$ to $D$
has the following asymptotics
\begin{equation*}
\lim_{\la \to \infty, \la \in r_{\theta'}} \psi(\lambda) =
\begin{cases}
0&  \text{if } r_{\theta'} \subset \Sigma_k,\\
\infty &  \text{if }  r_{\theta'} \subset \Sigma_{k\pm1}.
\end{cases} 
\end{equation*}
\end{pro}
A solution $\psi$ which vanishes at $\infty$ along a ray is called subdominant on this ray, while a solution which diverges along a ray is called dominant.
The previous proposition suggests that a subdominant solution is naturally defined in a Stokes sector $\Sigma_k$. However,
since the domain of the solutions is not $\C$ but the cut-plane $D_{\varphi}$, we need to distinguish and analyse
separately two cases, the case where
the cut $r_{\varphi}$ does not belong to the sector $\Sigma_k$, and the case where the cut $r_{\varphi}$ does belong to $\Sigma_k$.

In the first case, there exists a unique -- up to a normalising constant --
subdominant solution in the whole $\Sigma_k$ (which vanishes exponentially as $|\la| \to \infty$ in every closed subsector of $\Sigma_k$),
which we denote by $\psi_{k}(\la)$.

In the second case the ray $r_{\varphi}$ splits the Stokes sector $\Sigma_k$ into two non-empty open sub-sectors:
$\Sigma_k \setminus  r_{\varphi}= \Sigma_k^L \bigsqcup \Sigma_k^R$ (we define $\Sigma_k^L$ to be the subsector
containing the rays of argument $\varphi+\e$
for $\e$ small and positive.)
Due to Proposition \ref{pro:subdominants}, in both sub-sectors there exists a unique subdominant solution,
up to a normalising constant, which we denote by
$\psi^{L,R}_{k}(\la) $. The solutions $\psi^{L}_{k}(\la) $ and
$\psi^{R}_{k}(\la) $ may happen to be linearly dependent but in general they are not.
In fact, they are linearly dependent if and only if they are eigenvectors of the monodromy operator.
\begin{lem}\label{lem:monodromysubdominant}
Let $r_\varphi \subset \Sigma_k$.
The solutions $\psi^{L}_{k},\psi_k^R: D_{\varphi} \to \C $ subdominant in the sub-sectors $\Sigma_k^L,\Sigma_k^R$
are linearly dependent if and only if 
$\psi^{L}_{k} $ (and hence $\psi^{R}_k$) is an eigenvector of the monodromy operator $M$.
\end{lem}
\begin{proof}
The solution $\psi^{L}_k$ admits the decomposition $a \psi_k^{R}+
b \widehat{\psi}_k^{R}$ where $\widehat{\psi}_k^{R}$ is any dominant solution
in the right subsector (e.g. $\psi_{k-1}(\la)$), and $a,b$ are constants. Let us now compute the limit of the solution $\psi_k^L$
on the cut (the ray $r_{\varphi})$) when we approach the cut
from the left subsector or from the right subsector.
Due to Proposition \ref{pro:subdominants}, $\psi_k^{L},\psi_k^R, \widehat{\psi}_k^R $ can be analytically continued
in the whole sector $\Sigma_k$, without altering the property of being dominant or subdominant. Therefore
the from-the-left limit yields a subdominant solution on the cut, while the from the-right-limit yields
a dominant solution unless $b=0$. If $b=0$, then the from-the-right limit is proportional to the from-the-left limit,
since on any ray belonging to the Stokes sector the subdominant solution is uniquely defined up to a normalising constant.
By definition, the from-the-right and from-the-left limits are proportional if and only if $\psi^L$ is
an eigenvector of the monodromy matrix.
We can conclude that $\psi_k^R,\psi_k^L$ are linearly dependent 
if and only if $\psi_k^L$ (and hence $\psi_k^R$) is an eigenvector of the monodromy matrix. 
\end{proof}

We are now able to express the two conditions of Theorem \ref{thm:inversemonodromy} into a --
more manageable -- relations among subdominant solutions.
\begin{lem}\label{lem:conditionssub} We have the following equivalences for the quantisation and no-logarithm conditions
in Theorem \ref{thm:inversemonodromy}:
\begin{enumerate}
	\item Let $\chi_+,\psi_0: D_{\pi} \to \C$ be solutions of \eqref{eq:genharmonic}, $\chi_+$ subdominant at $0$,
	and $\psi_0$ subdominant in the Sector $\Sigma_0$. Then there exists a non-zero solution $\psi: D_{\pi} \to \C$ 
	of \eqref{eq:genharmonic} which solves
	the boundary value problem
	\begin{equation*}
	\lim_{\la \to +\infty} \psi(\la)=\lim_{\la \to 0^+} \psi(\la)=0,
	\end{equation*}
	if and only if
	$\chi_+$ and $\psi_{0}$ are linearly dependent.
	\item Let $n \in \bb{N}^*$ and $\chi_+,\psi^{L,R}_{1}:D_{\frac\pi2}\to \C$ be solutions of\eqref{eq:genharmonic} with $\chi_+$ subdominant at $0$
	and $\psi_1^{L,R}$ subdominant in the subsector $\Sigma_1^{R,L}$. Assume that $\chi_+$ and $\psi_1^L$ are linearly independent.
	Then the singularity $\la=0$ is apparent if and only if $\psi^L_{1}$ and $\psi^R_{1}$
	are linearly dependent.
\end{enumerate}	 
\end{lem}
\begin{proof}
(1) A solution of the boundary value problem must necessarily be proportional to $\chi_+$, since it vanishes at $0$, and proportional to
$\psi_0$, since it vanishes on a ray inside $\Sigma_0$.

(2) The singularity $\lambda=0$ is apparent if and only if every solution is an eigenvector of the monodromy operator $M$. Since $\chi_+$ is always an eigenvector of $M$ and, by hypothesis, $\psi_1^L$
is linearly independent from $\chi_+$, the latter holds if and only if $\psi_1^L$ is an eigenvector of $M$. By Lemma \ref{lem:monodromysubdominant}, $\psi_1^L$ is an eigenvector of the monodromy operator if and only if
$\psi_1^L$ and $\psi_1^R$ are linearly dependent.
\end{proof}

Next we transform the conditions obtained in Lemma \ref{lem:conditionssub} into equations among objects that can be 
explicitly computed in terms of the potential $Q$. To this aim we begin with the following definition.
\begin{defi}
Fix a $k\in \Zq$ and a $\varphi \in \bb{R}$ such that $r_{\varphi} \nsubseteq \Sigma_k$.
Let $f$ be a meromorphic function on $D_{\varphi}$.
We define the following quantities in $\bb{P}^1$, provided they exist:
\begin{equation}\label{eq:w0wk}
  w_{\underline{0}}(f)=\lim_{\la \to 0}f(\la) , \quad w_{k}(f):=\lim_{|\la| \to \infty,\la \in \Sigma_k} f(\la) \; .
\end{equation}
In the second case the limit is taken on any ray belonging to $\Sigma_k$ and the limit is defined to exist
if it is independent on the choice of the ray. 
\end{defi}
It is a basic fact that the above limits exist if $f$ is the ratio of a pair of solutions of the linear ODE \eqref{eq:genharmonic}
\cite{eremenko}. It is indeed a straightforward consequence of the fact that both at $0$ and in $\Sigma_k$
every solution admits a unique decomposition in a basis composed of a subdominant and a dominant solution.

We now prove Proposition \ref{prop:wronskiantovalues}, which states that the quantisation and apparent singularity conditions
can be expressed as an equality among pairs of values $w_{\underline{0}}(f),w_k(f)$ for some $f$ which is the ratio of
a specific pair of solutions of equation \eqref{eq:genharmonic}.
Such an idea is deeply rooted in the Nevanlinna theory of covering maps, which was explored by the first author in \cite{piwkb,pibelyi},
and by the first and second authors in \cite{maro18}.
\begin{pro}\label{prop:wronskiantovalues} We have the following equivalences for the quantisation and no-logarithm conditions in Theorem \ref{thm:inversemonodromy}:
	\begin{enumerate}
\item  Let $\chi_+,\psi_0,\psi_{\pm1}: D_{\pi} \to \C$ be solutions of \eqref{eq:genharmonic} with $\chi_+$ subdominant at $0$,
  $\psi_0$ subdominant in the sector $\Sigma_0$ and $\psi_{\pm1}$ subdominant in the sector $\Sigma_{\pm1}$.
  Assume that $\psi_{+1},\psi_{-1}$ are linearly independent. Then there exists a non-zero solution of \eqref{eq:genharmonic} which solves
		the boundary value problem
		\begin{equation*}
		 \lim_{\la \to +\infty} \psi(\la)=\lim_{\la \to 0^+} \psi(\la)=0,
		\end{equation*}
		if and only if
		\begin{equation}\label{eq:quantAB}
   w_{-1}\big( \frac{\chi_+}{\psi_{0}} \big)=
   w_{1}\big(\frac{\chi_+}{\psi_{0}}\big).
 \end{equation}
\item Let $n \in \bb{N}^*$ and $\chi_+,\psi_{-1},\psi^{L,R}_{1}:D_{\frac\pi2}\to \C$ be solutions of \eqref{eq:genharmonic} with
$\chi_+$ subdominant at $0$, $\psi_1^{L,R}$ subdominant in the subsector $\Sigma_1^{R,L}$ and $\psi_{-1}$
subdominant in $\Sigma_{-1}$. Assume that $\chi_+$ and $\psi_1^L$ are linearly independent, and that
  $\chi_{+}$ and $\psi_{-1}$ are linearly independent. Then the singularity $\la=0$ is apparent if and only if
 \begin{equation}\label{eq:ABRL}
 w_{-1}\big( \frac{\psi_{1}^{R}}{\psi_{1}^{L}}\big)=
  w_{\underline{0}}\big( \frac{\psi_{1}^{R}}{\psi_{1}^{L}} \big).
 \end{equation}
\end{enumerate}
In formulae \eqref{eq:quantAB},\eqref{eq:ABRL}, the quantities $w_{-1},w_1,w_{\underline{0}}$ are defined as per
 \eqref{eq:w0wk}.
\end{pro}
 \begin{proof}
(1) By hypothesis $\psi_{-1},\psi_1$ are linearly independent. It follows that, if
$\varphi,\psi$ is a pair of solutions, the equality $w_{-1}(\frac{\phi}{\psi})=w_{1}(\frac{\phi}{\psi})$
holds if and only if $\phi$ and $\psi$ are linearly dependent. Indeed, the if part of this claim is obvious. We proceed with proving the only if part,
by contradiction. Since $\psi_{-1}$ and $\psi_{1}$ are linearly independent,
$\psi_{-1}$ is dominant in $\Sigma_1$ and $\psi_{1}$ is dominant is $\Sigma_{-1}$, so that $w_{-1}(\frac{\psi_{-1}}{\psi_1})=0$
and $w_1(\frac{\psi_{-1}}{\psi_1})=\infty$. Assume that there is a basis of solutions $\phi,\psi$ such that
$w_{-1}(\frac{\phi}{\psi})=w_{1}(\frac{\phi}{\psi})$, then there exists $a,b,c,d \in \C$ with $ad-bc \neq 0$ such that
$\psi_{-1}=a \phi+b \psi,\psi_1=c \phi+ d \psi$. We deduce the contradicting
$$w_{1}(\frac{\psi_{-1}}{\psi_1})=\frac{a w_{1}(\frac{\phi}{\psi})+b}{c w_{1}(\frac{\phi}{\psi})+d}
=\frac{a w_{-1}(\frac{\phi}{\psi})+b}{c w_{-1}(\frac{\phi}{\psi})+d}=w_{-1}(\frac{\psi_{-1}}{\psi_1}).$$

The thesis now follows from Lemma \ref{lem:conditionssub}(1), which states that the quantisation condition is satisfied if 
and only if $\chi_+$ and $\psi_0$ are linearly dependent.

(2) Due to Lemma \ref{lem:conditionssub}(2), under the present hypotheses,
$\la=0$ is an apparent singularity if and only if $\psi_1^R$ and $\psi_1^L$ are linearly dependent. The proof follows from this fact
using the very same argument used to prove (1).
\end{proof}

We are interested in solving equations (\ref{eq:quantAB},\ref{eq:ABRL})
for the parameters $a_1,a_2,a_3$ of the potential $Q$ when $a_0,n$ are fixed to a certain value.
It is convenient to define a pair of meromorphic functions, of the parameters $a_1,a_2,a_3$,
which vanish when  (\ref{eq:quantAB},\ref{eq:ABRL}) are satisfied.
\begin{lem}\label{lem:analiticity}
Let the parameters $a_1,a_2,a_3$ of the potential $Q$ vary in $\C^3$ while $a_0,n$ are fixed.
The functions $U_1,U_2: \C^3 \to \Cb$ defined by the formulae
\begin{align*}
&U_1=1-\frac{w_{-1}\big( \frac{\chi_+}{\psi_{0}} \big)}
   {w_{1}\big(\frac{\chi_+}{\psi_{0}}\big)}\\
  & U_2=1-\frac{w_{-1}\big( \frac{\psi_{1}^{R}}{\psi_{1}^{L}}\big)}
  {w_{\underline{0}}\big( \frac{\psi_{1}^{R}}{\psi_{1}^{L}} \big)}
\end{align*}
 are well-defined holomorphic functions.
 In other words, $U_1,U_2$ are meromorphic functions on $\C^3$.
 \end{lem}
 \begin{proof}
We sketch the proof here, details can be found in the cited literature.
Let us analyse $U_1$, for example. It is a well-known fact that one can choose a normalisation of
$\psi_0$ in such a way it is an entire function of
$a_1,a_2,a_3,n$ of the potential $Q$ \cite{sibuya75}. Similarly, one can choose a normalisation of
$\chi_+$ such that $\chi_+$ is an entire function of $a_0,a_1,a_2,a_3$.
Under such a normalisation, one can show that the expression
$w_{\pm1}\big( \frac{\chi_+}{\psi_{0}}\big)$ is a holomorphic function,
with values in $\Cb$, 
of the coefficients $a_1,a_2,a_3$ of the potential; see \cite{eremenko,sibuya75}.
Since $U_1$ is invariant under re-scaling of $\chi_+,\psi_0$, the thesis follows.
 \end{proof}

\begin{rem}
We could have expressed conditions of linearly dependence of pairs of subdominant solutions
as conditions of the vanishing of their Wronskians. However, the direct asymptotic
analysis of a Wronskian is slightly more involved than the asymptotic analysis of a ratio of $w$ values,
because the former requires asymptotic estimates of solutions and of their derivatives.
\end{rem}

\subsection{WKB functions}
We consider the complex WKB method for an equation of the kind
\begin{equation}\label{eq:schrlan}
 \psi''(\lambda)=(E^2 V(\la)+r(\la))\psi(\la) \;.
\end{equation}
Specifying to $V=V(\la;\alpha,\beta,\nu)$ and $r=-\frac{1}{4 \la^2}$, we obtain the main object of our analysis,
namely equation (\ref{eq:scaled}).

The WKB asymptotic is based on the approximation of solutions by means of the (multivalued) WKB function
\begin{align}\label{eq:wkb}
& \Psi(\la;\la_0)= \exp\left\{\int^{\la}_{\la_0}E \sqrt{V(\mu)} - \frac{V'(\mu)}{4 V(\mu)} d \mu \right\}.
\end{align}

In order to measure the difference between a WKB function $\Psi(\la;\la_0)$ and a putative solution $\psi(\la)$ of
\eqref{eq:schrlan}, we introduce the
ratio
\begin{equation}\label{eq:z}
 z(\la)=\frac{\psi(\la)}{\Psi(\la;\la_0)}
\end{equation}
If we fix a point $\la' \in \Cb$, a path $\gamma$ starting at $\la'$ and passing through $\la_0$ and
the boundary conditions $$\lim_{\la \to \la'}z(\la)=1 \, , \quad
\lim_{\la \to \la'}z'(\la)=0 ,$$ 
$z(\la)$ is shown (after a few algebraic manipulations, see \cite[\S 4]{erdelyi10} for the details) to satisfy the following Volterra
integral equation 
\begin{align}\nonumber
 & z(\la)=1-\frac{1}{E}\int^{\la}_{\la',\gamma} B(\la,\mu) F(\mu) z(\mu) d\mu \, , \quad B(\la,\mu)= \frac{\exp\left\{ -
 2 E \int_{\mu,\gamma}^\la \sqrt{V(\nu)}d\nu \right\}-1}{2} \\
 \label{eq:volterra}
 &   F(\la)= V^{-\frac12}(\la) \left(-r(\la)+\frac{ -4 V''(\la) V(\la)+5V'^2(\la)}{16 V(\la)^2}\right),
\end{align}
for each $\la$ in the support of $\gamma$ such that $\int^{\la}_{\la',\gamma} B(\la,\mu) F(\mu)d\mu$ converges absolutely.

\begin{rem}
We notice that formula \eqref{eq:wkb} is the so-called \textbf{Langer-modified WKB approximation}. In fact,
 the standard WKB approximation is obtained by studying approximate solutions of the form
 $\widetilde{\Psi}(\la;\la_0)= \exp\big(\int^{\la}_{\la_0}\sqrt{Q(\mu)} - \frac{Q'(\mu)}{4 Q(\mu)} d \mu \big)$, with $Q=E^2 V+r$.
 This choice leads to an integral equation with forcing term $$\widetilde{F}(\la) d\la=\frac{ -4 Q''(\la) Q(\la)+5Q'^2(\la)}
 {16 Q^{\frac52}(\la)} d\la \; .$$
 The latter 1-form coincides asymptotically with $E^{-1}F(\la) d\la$ in all of $\hspace{0.1cm}\overline{\mathbb{C}}$ but for a neighbourhood of the Fuchsian
 singularity $\la=0$, where $F(\la)d\la$ is integrable while $\widetilde{F}(\la) d\la$ is not. For such a reason,
 the standard WKB approximation fails at that point, while the Langer modified WKB approximation provides the correct result.
\end{rem}

We begin the analysis of the integral equation by studying the singularities of the forcing form $F(\la) d\la$, in the case
$V=V(\la;\alpha,\beta,\nu)$ and $r=-\frac{1}{4 \la^2}$.
\begin{lem}\label{lem:integrability}
Let $\gamma:[0,1] \to \overline{\bb{C}}$ be a rectifiable curve with respect to the standard metric on the Riemann sphere, namely
\begin{equation}\label{eq:rectificable}
\int_0^1 \frac{|\dot{\gamma}(t)|}{1+|\gamma(t)|^2}dt < \infty.
\end{equation}
If $ V(\gamma(t))\neq 0$ for $t \in [0,1]$, then 
 \begin{equation}\label{eq:rhog}
  \rho_{\gamma}:= \int_0^1|F\big(\gamma(s)\big) \dot{\gamma}(s)|ds <	\infty  .
 \end{equation}
 \begin{proof}
  The form $F(\la) d\la$ is manifestly regular at all point in $\Cb$ where $V(\la)\neq0$ and $V(\la) \neq \infty$.
  The poles of $V(\lambda)$ are the points $\lambda=0$ and $\lambda=\infty$, and a simple computation shows that the form is integrable there.
  In fact $F(\la)=O(\lambda^{-3})$ as $\la \to \infty$ and
  $F(\la)=O(1)$ as $\la \to 0$.
 \end{proof}

\end{lem}

\begin{defi}\label{defi:admissible}
 Let $\la'$ be a regular point of $F(\la) d\la$. We say that $\gamma:[0,1] \to \overline{\bb{C}}, \gamma(0)=\la'$ is an admissible curve
 if
 \begin{enumerate}
  \item $ V(\gamma(t)) \neq \infty $, for all $t\in(0,1)$.
  \item $ V(\gamma(t))\neq 0$, for all $t \in [0,1] $.
  \item $\Re \int_{t'}^{t} \sqrt{V(\gamma(s))}\dot{\gamma}(s)ds\geq0 $ for all $0< t'\leq t<1 $,
 for one of the two branches of $\sqrt{V}$ definable on the support of $\gamma$.
 \item  $\gamma$ is rectifiable, i.e. equation \eqref{eq:rectificable} holds.
 \end{enumerate}
\end{defi}
For every admissible curve, there is an actual solution of \eqref{eq:schrlan} which converges uniformly to the WKB approximation
as $E$ grows large. More precisely, we have the following proposition.
\begin{pro}[Analytic continuation and WKB approximation]\label{pro:continuation}
Let $D \subset \C^*$ be an open simply-connected subset, and
$\gamma:[0,1] \to \Cb$ an admissible path such that $\gamma((0,1)) \subset D$.

Fix a point $t_0 \in (0,1)$ and
let $\sqrt{V}$ denote the branch such that $\Re \int_{t_0}^{t} \sqrt{V(\gamma(s))}\dot{\gamma}(s)ds\geq0$
for all $t\geq t_0$.

For all $E>0$, there exists a unique solution $\psi:D \to \C$ of \eqref{eq:schrlan} such that
 \begin{equation}\label{eq:wkbasymptotic}
  \norm{\frac{\psi\big(\gamma(t)\big)}{\Psi\big(\gamma(t);\gamma(t_0)\big)}-1} \leq \exp\big(\frac{\rho_\gamma(t)}{E}\big)-1
 , \; \forall t \in [0,1],
 \end{equation}
 where $\rho_{\gamma}(t)= \int_0^t|F\big(\gamma(s)\big) \dot{\gamma}(s)|ds $ (so that $\rho_{\gamma}(1)=\rho_{\gamma}$ as per \eqref{eq:rhog}),
 and
 $$\Psi\big(\gamma(t);\gamma(t_0)\big)= \exp\left\{ \int_{t_0}^t \big(E \sqrt{V(\gamma(s))} - 
 \frac{V'(\gamma(s))}{4 V(\gamma(s))} \big) \dot{\gamma}(s) ds \right\} ,$$
 as per formula \eqref{eq:wkb}.
 
\begin{proof}
Details can be found in \cite{piwkb}.
The integral equation \eqref{eq:volterra} is of the form $z=1+E^{-1}B[z]$, where
$$B[z](\gamma(t))=-\int^{t}_{0} B(\gamma(t),\gamma(s)) F(\gamma(s)) z(\gamma(s)) \dot{\gamma}(s) ds.$$
By hypothesis $\gamma$ is admissible path. It follows that
\begin{itemize}
 \item $ |B(\gamma(t),\gamma(s))|\leq1,\forall s \leq t$, by Definition \ref{defi:admissible}(3).
 \item $\rho_{\gamma}<\infty$, since, by Definition \ref{defi:admissible}(1,2,4), $\gamma$ satisfies all hypotheses
 of Lemma \ref{lem:integrability}.
\end{itemize}
These estimates imply that $B$ is a continuous linear operator on the Banach space of continuous functions
on $\gamma([0,t])$, equipped with the $L^{\infty}$ norm. In fact, by using H\"older's inequality, the operator norm
$\|{B}\|_{\infty}$ is readily seen to satisfy the inequality $\|{B}\|_{\infty}\leq \rho_{\gamma}(t)$. Moreover, we have that
$$\|{B^l}\|_{\infty}\leq \frac{\|{B}\|_{\infty}^l}{l!}\leq\frac{(\rho_{\gamma}(t))^l}{l!}.$$ Indeed
the $l$-th iteration of $B$ can be written explicitly as an integral on the standard simplex of dimension $l$,
which has volume $\frac{1}{l!}$, see e.g. \cite{piwkb} for more details.

We can therefore construct the unique solution of $z=1+E^{-1}B[z]$ by means of the Neumann series
$z=\sum_{l=0}^{\infty}E^{-l}B^l[1]$. The thesis immediately follows from the estimate on $\|{B^l}\|_{\infty}$.
\end{proof}
\end{pro}

Suppose that we have an admissible path $\gamma$ such that $\lim_{t \to 0}\gamma(t)$ equal to $0$ or $\infty$, and
in the latter case suppose that $\gamma(t)$ belongs, for $t$ small, to a closed sub-sector of a Stokes sector $\Sigma_k$ for some $k\in \Zq$.
A simple computation shows that the WKB function has limit
$\lim_{t \to 0}\Psi(\gamma(t);\gamma(t_0))=0$, therefore the corresponding solution $\psi(\la)$ constructed in the above Proposition
is subdominant at $0$ or in $\Sigma_k$.

Hence we are led to study admissible paths which start at $0$ or at $\infty$ (in a closed subsector of a Stokes sector).
For our specific purpose, we need to take into consideration a few more properties.
\begin{enumerate}
 \item We need to consider admissible curves that connect $0$ with a Stokes sector, or 
 a Stokes sector with another Stokes sector. This is because, after Proposition
 \ref{prop:wronskiantovalues}, the relevant objects are the asymptotic values of the ratio of two subdominant solutions.
 \item Since the WKB approximation is written in terms of the integral of $\sqrt{V}$, which is singular at poles and zeros of $V$,
 we need to keep track of the homotopy class of admissible curves in
 $\C^* \setminus \lbrace \la \in \C, V(\la)=0 \rbrace$.
\item We must restrict our study to those admissible paths
which are stable under small complex perturbation of the parameters $(\alpha,\beta) \in R$, in order to have enough room to apply
the implicit function theorem that we will develop in Section \ref{section:bulkasymptotic} of this paper.
\end{enumerate}

The admissidle curves that satisfy the above properties are the
level curves of the function $\Im \int^{\la}\sqrt{V(\mu)}d\mu$, which were introduced and studied in
Subsection \ref{subsec:anti}. These are also known horizontal trajectories of the quadratic differential $V(\la) d\la^2$,
and we stick to this second name for simplicity.
\subsection{Topology of horizontal trajectories}
Here we study the topology of the horizontal trajectories, namely of the level curves of the function $\Im \int^{\la}\sqrt{V(\mu)}d\mu$.
Recall from Subsection \ref{subsec:anti} that any horizontal trajectory, but for a finite
set of them known as anti-Stokes lines, can be indefinitely prolonged so that
its endpoints are poles of $V$, namely $0$ or $\infty$. Moreover,
if one endpoint is $\infty$, then the horizontal trajectory asymptotically lie in a closed
subsector of $\Sigma_k$ for some $k \in \Zq$; following the nomenclature introduced in Subsection \ref{subsec:anti},
we say that its endpoint is $e^{i\frac{\pi k}{2}}\infty$.

The following definition is quite natural.
\begin{defi}\label{def:Hab}
Let $a,b$ take value in the set of symbols $\lbrace \underline{0} \rbrace \cup \bb{Z}_4$, endowed with the function
$\pi:\lbrace \underline{0} \rbrace \cup \bb{Z}_4 \to \lbrace  0,+\infty,-\infty,+i \infty,-i\infty \rbrace$, defined by
$\pi(\uo)=0$, and $\pi(k)=e^{i\frac{\pi k}{2}}\infty$ for $k \in \bb{Z}_4$.

 We denote by $\overrightarrow{A}_{a,b}(\alpha,\beta)$
 the set of oriented horizontal trajectories of the potential $V(\la;\alpha,\beta,\nu)$ which starts at $\pi(a)$ and ends at $\pi(b)$.
We denote by $A_{a,b}(\alpha,\beta)$ the set of un-oriented trajectories with endpoints $a,b$.

Finally, we assume any oriented horizontal trajectory is given a smooth parametrisation $\gamma:[0,1] \to \Cb$,
such that $\lim_{t\to0}\gamma=\pi(a)$ and $\lim_{t\to1}\gamma=\pi(b)$.

By an abuse of notation, we use the same symbol $\gamma$ for a parametrised horizontal trajectory
and for its support. Moreover, we write $\overrightarrow{A}_{a,b},A_{a,b}$ for
$\overrightarrow{A}_{a,b}(\alpha,\beta),A_{a,b}(\alpha,\beta)$ when no confusion can arise.
\end{defi}\label{lem:antiadmissible}
\begin{lem}
If $\gamma \in \overrightarrow{A}_{a,b}$ for some $a,b \in \lbrace \underline{0} \rbrace \cup \bb{Z}_4$ then $\gamma$ is an admissible curve,
for the choice of $\sqrt{V}$ such that $\lim_{t\to 0}\Re \int_{t}^{t_0} \sqrt{V(\gamma(t))}\dot{\gamma}(t) dt =\infty $.
\begin{proof}
By construction $\gamma$ satisfies the first 3 conditions of Definition \ref{defi:admissible} (of an admissible curve).
We need to verify the last condition, namely that $\gamma$ is rectifiable.
To this aim we only need to verify that $\gamma$ is rectifiable near the poles of $V$,
that is near $0$ and $\infty$. Both properties follow from the local description of horizontal trajectories,
from Strebel's book:
\begin{itemize}
 \item After
 \cite[Theorem 7.2]{strebel} the horizontal trajectories can be rectified in a neighbourhood of $0$, since $\lim_{\la \to 0} \la^2V(\la)$
 is a positive constant.
 \item After
 \cite[Theorem 7.4]{strebel} the horizontal trajectories can be rectified in a neighbourhood of $\infty$, since
 $\infty$ is a pole of order greater or equal to $4$ (namely of order $6$) of the quadratic differential.\qedhere
\end{itemize}
\end{proof}
\end{lem}
After the above lemma and Proposition \ref{pro:continuation}, whenever we have a path $\gamma \in \overrightarrow{A}^{a,b}$,
we can estimate how the solution subdominant at $\pi(a)$
behaves at $\pi(b)$. In particular we know that the solution subdominant at $\pi(a)$ is dominant at $\pi(b)$ - along $\gamma$ - if $E$ is
large enough.
\begin{lem}\label{lem:domsub}
Let $D\subseteq \C^*$ be a simply connected open domain.
Suppose there exists an admissible path $\gamma \in \overrightarrow{A}^{a,b}$ such that $\gamma((0,1)) \subseteq D$.
For every $E>0$, there exists a solution $\psi: D \to \C$ of \eqref{eq:scaled} such that
 $$\lim_{t \to 0}\psi(\gamma(t))=0 .$$
Moreover, for all $E>\frac{\rho_{\gamma}}{\ln 2}$ -- where $\rho_{\gamma}$ is as per \eqref{eq:rhog} -- the same solution satisfies 
$$\lim_{t \to 1}\left| \psi(\gamma(t)) \right|=\infty .$$

\end{lem}
\begin{proof}
Let $\Psi(\gamma(t);\gamma(t_0))$ be the WKB function as per \eqref{eq:wkb}. A simple computation shows that
$\lim_{t\to 0}\Re \int_{t_0}^{t} \sqrt{V(\gamma(t))}\dot{\gamma}(t) dt =-\infty $
and $\lim_{t\to 1}\Re \int_{t_0}^{t} \sqrt{V(\gamma(t))}\dot{\gamma}(t) dt =\infty $. Therefore
$\lim_{t \to 0}\Psi(\gamma(t);\gamma(t_0))=0$ and $\lim_{t \to 1}\Psi(\gamma(t);\gamma(t_0))=\infty$.
Lemma \ref{lem:antiadmissible} states that $\gamma$ is admissible. Hence we can apply Proposition \ref{pro:continuation}
to deduce the following estimates which imply the thesis: there exists a $\psi:D \to \C$ such that
 \begin{itemize}
  \item $\lim_{t \to 0}\frac{\psi(\gamma(t))}{\Psi(\gamma(t);\gamma(t_0))}=1$ for all $E$.
  \item  $\lim_{t \to 1}\big|\frac{\psi(\gamma(t))}{\Psi(\gamma(t);\gamma(t_0))}\big|>0$
  for all $E>\frac{\rho_{\gamma}}{\ln 2}$.\qedhere
 \end{itemize}
\end{proof}

\begin{rem}\label{rem:paths} The following properties of $\overrightarrow{A}_{a,b} $ are either straightforward consequences of the definition or
	are well-known facts from Strebel's book \cite{strebel}.
\begin{enumerate}
\item The sets $\overrightarrow{A}_{a,b}$ and $\overrightarrow{A}_{b,a}$ are in natural bijection, and the set
$A_{a,b}$ coincides with $A_{b,a} $.
 \item
 The Anti-Stokes complex, defined in Subsection \ref{subsec:anti}, coincides, as a subset of $\mathbb{C}$, with the
 embedded graph defined by the set 
 \begin{equation*}
  \mathbb{C} \setminus \bigcup_{a \neq b}\lbrace x \in \gamma:
 \gamma \in A_{a,b} \rbrace.
 \end{equation*}
 \item
 $\forall k \in \Zq$ $A_{k,k\pm1}\neq \emptyset$, and
 $\inf_{\gamma \in A_{k,k\pm1}} \rho_{\gamma}=0$. Indeed,
 for every $\e>0$, there exists a path $\gamma \in A_{k,k\pm1}$ which has length less than $\e$ with respect to the
 standard metric on the sphere \cite[Section 7.4]{strebel}.
 \item For any ray $\theta \in [0,2 \pi) $ there exists a unique horizontal trajectory $\gamma_{\theta}$,
 which start at $\la=0$, tangent to the ray $r_{\theta}$. \cite[Theorem 7.2]{strebel}
\end{enumerate}
\end{rem}
The horizontal trajectories of the quadratic differentials are very useful tools in the WKB analysis, but they are slightly too rigid,
since distinct trajectories can only meet at $0$ and $\infty$.
We can however deform them quite freely in a neighbourhood of $0$ and $\infty$.
\begin{lem}\label{lem:pathdef}
(1) Let $\gamma:[0,1] \to \Cb$ be an admissible curve such that $\gamma(0)=e^{\frac{i \pi k}{2}}$ for some $k \in \bb{Z}_4$,
and $\varphi$ be an angle
$-\pi<\varphi<\pi$.
For all $\e>0$, there exists another admissible curve $\widetilde{\gamma}$ satisfying the following properties:
\begin{itemize}
 \item $\gamma$ and $\widetilde{\gamma}$ are homotopic in $\C^*\setminus \lbrace V(\la)=0 \rbrace $;
 \item $\arg \widetilde{\gamma}(t)=\frac{\pi k}{2}+\frac{\varphi}{4}$ as $t\to 0$;
 \item $|\rho_{\gamma}-\rho_{\widetilde{\gamma}}|\leq \e$, where $\rho_{\gamma}$ is as per \eqref{eq:rhog}.
 \end{itemize}
 
(2)Let $\gamma:[0,1] \to \Cb$ be an admissible curve such that $\gamma(1)=0$ and $\varphi \in \bb{R}$ be an arbitrary angle.
For all $\e>0$, there exists an admissible curve $\widetilde{\gamma}$ satisfying the following properties:
\begin{itemize}
 \item $\gamma$ and $\widetilde{\gamma}$ are homotopic in $\C^*\setminus \lbrace V(\la)=0 \rbrace $;
 \item $\arg \widetilde{\gamma}(t)=\varphi$ as $t \to 1$.
 \item $|\rho_{\gamma}-\rho_{\widetilde{\gamma}}|\leq \e$,  where $\rho_{\gamma}$ is as per \eqref{eq:rhog}.
 \end{itemize}
 \begin{proof}
 We prove (2) and leave (1) to the reader since the proof is almost identical.
 According to \cite[Theorem 7.2]{strebel}, there exists a holomorphic map in a neighbourhood $U$ of $0$,
 $\zeta:U \to \bb{C}, \zeta(0)=0,\zeta'(0)=1$
 such that  $V(\la) d\la^2=\frac{\nu^2}{4 \zeta^2} d\zeta^2$ so that  $\int \sqrt{V(\la)} d \la= \frac \nu2 \log \zeta$.
 Therefore,
 \begin{itemize}
  \item for any $\theta \in (0,\frac\pi{\nu})$
 or $\theta \in(-\frac\pi{\nu},0)$, any  finite part of the logarithmic spiral $\arg \log \zeta = \theta $ - in the direction
 converging to $\zeta=0$ -
 is admissible;
 \item if $\la$ is small enough the segment $[0,\la]$ is admissible.
 \end{itemize}

 Now we construct $\widetilde{\gamma}$ as follows. We choose a $\delta$ small and let $\widetilde{\gamma}$ coincide with $\gamma$
 on $\gamma([0,1-\delta])$.
 Starting from this point, $\widetilde{\gamma}$ coincides with the logarithmic spiral of 
 angle $\theta$ for a $\theta \in(0,\frac\pi{n})$
 or $\theta \in(-\frac\pi{n},0)$, until the point $\zeta^*$ such that
 $\arg \la(\zeta^*)=\varphi$. Finally, $\widetilde{\gamma}$
 coincide with segment connecting  $\la(\zeta^*)$ with $0$.
 As $\delta \to 0$, $\lim_{\delta \to 0}\int_{1-\delta}^1 |\dot{\widetilde{\gamma}}(t)|dt \to 0$ which implies the thesis since
 the form $F(\la) d\la$ is regular at $\la=0$.
 \end{proof}

\end{lem}

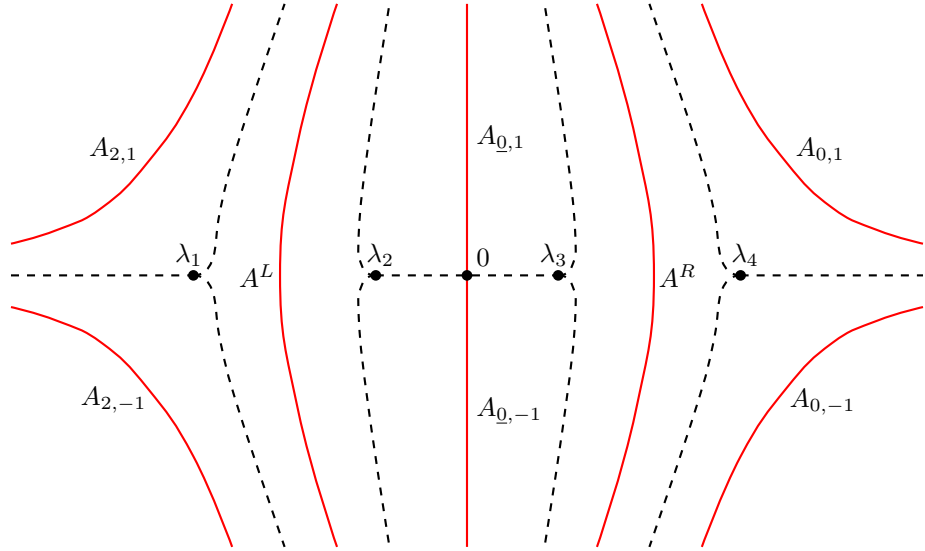
\begin{figure}
	\centering
		\begin{tikzpicture}[scale=0.6]

		\tikzstyle{star}  = [circle, minimum width=4pt, fill, inner sep=0pt];
	
		\node[above]     at (-6.1,0) {$\lambda_1$};

		\node[above]     at (-1.9,0) {$\lambda_2$};

		\node[above right]     at (0,0) {0};

		\node[above]     at (1.9,0) {$\lambda_3$};

		\node[above]     at (6.1,0) {$\lambda_4$};

		\draw[black,dashed,thick] (-2,0) -- (2,0);
		\draw[black,dashed,thick] (6,0) -- (10,0);
		\draw[black,dashed,thick] (-10,0) -- (-6,0);

		\draw [black,dashed,thick] plot [smooth,tension=0.8] coordinates {(2,0) (2.3,0.2)  (2.37,1) (2.15,3) (1.7,6)   };
		\draw [black,dashed,thick] plot [smooth,tension=0.8] coordinates {(-2,0) (-2.3,0.2)  (-2.37,1) (-2.15,3) (-1.7,6)   };
		\draw [black,dashed,thick] plot [smooth,tension=0.8] coordinates {(-2,0) (-2.3,-0.2)  (-2.37,-1) (-2.15,-3) (-1.7,-6)   };
		\draw [black,dashed,thick] plot [smooth,tension=0.8] coordinates {(2,-0) (2.3,-0.2)  (2.37,-1) (2.15,-3) (1.7,-6)   };

		\draw [black,dashed,thick] plot [smooth,tension=0.8] coordinates {(6,0) (5.6,0.3)  (5.4,1.5) (5,3) (4,6)   };
		\draw [black,dashed,thick] plot [smooth,tension=0.8] coordinates {(-6,0) (-5.6,0.3)  (-5.4,1.5) (-5,3) (-4,6)   };
		\draw [black,dashed,thick] plot [smooth,tension=0.8] coordinates {(6,0) (5.6,-0.3)  (5.4,-1.5) (5,-3) (4,-6)   };
		\draw [black,dashed,thick] plot [smooth,tension=0.8] coordinates {(-6,0) (-5.6,-0.3)  (-5.4,-1.5) (-5,-3) (-4,-6)   };

		\draw[red,thick] (0,-6) -- (0,6);
		
		\draw [red,thick] plot [smooth,tension=0.8] coordinates { (2.85,-6)  (3.7,-3)    (4.1,0) (3.7,3) (2.85,6)   };
		\draw [red,thick] plot [smooth,tension=0.8] coordinates { (-2.85,-6)  (-3.7,-3)    (-4.1,0) (-3.7,3) (-2.85,6)   };
		
		\draw [red,thick] plot [smooth,tension=0.8] coordinates { (10,0.7) (9,1)  (8,1.5)    (7,2.5) (6,4) (5.15,6)   };
		\draw [red,thick] plot [smooth,tension=0.8] coordinates { (-10,0.7) (-9,1)  (-8,1.5)    (-7,2.5) (-6,4) (-5.15,6)   };
		\draw [red,thick] plot [smooth,tension=0.8] coordinates { (10,-0.7) (9,-1)  (8,-1.5)    (7,-2.5) (6,-4) (5.15,-6)   };
		\draw [red,thick] plot [smooth,tension=0.8] coordinates { (-10,-0.7) (-9,-1)  (-8,-1.5)    (-7,-2.5) (-6,-4) (-5.15,-6)   };

			\node[black,star]     at (-6,0){};
			
			\node[black,star]     at (-2,0) {};
			
			\node[black,star]     at (0,0) {};

			\node[black,star]     at (2,0) {};

			\node[black,star]     at (6,0) {};

			\node    at (7.75,2.75) {$A_{0,1}$};
			\node    at (-7.75,2.75) {$A_{2,1}$};
			\node    at (-7.75,-2.75) {$A_{2,-1}$};
			\node    at (7.8,-2.8) {$A_{0,-1}$};
			
			\node[right]  at (0,3) {$A_{\underline{0},1}$};
			\node[right]  at (0,-3) {$A_{\underline{0},-1}$};

			\node[right] at (4,0)  {$A^R$};
			\node[left] at (-4,0)  {$A^L$};

	\end{tikzpicture} 
	\caption{Topological representation of the anti-Stokes complex and of the homotopy classes of horizontal trajectories.
	Dashed black lines: anti-Stokes lines. Red lines: representatives of homotopy classes of horizontal trajectories.} \label{fig:antiS}
\end{figure}

Proposition \ref{pro:antiS} of Section \ref{section:elliptic} states that,
for all $(\alpha,\beta) \in R$, the anti-Stokes complex $\widetilde{\mathcal{C}}$ is isomorphic to the the embedded graph
represented in Figure \ref{fig:antiS}. We notice that the totality of the elements in $A_{a,b}$ for all possible $a,b$ define
a foliation of $\C$ minus the anti-Stokes complex. By construction
every chamber of the foliation is foliated by curves which are homotopic:
chambers are in bijection with homotopy classes of horizontal trajectories which do not end at a turning point.
As a consequence the classification, for $(\alpha,\beta) \in R$, of the sets $A_{a,b}$ which are non-empty, and the partition
of these into homotopy classes can be read off from Figure \ref{fig:antiS}.
We thus have the following lemma.
\begin{lem}\label{lem:nonemptysets}
Suppose $(\alpha,\beta) \in R$.
\begin{itemize}
 \item The set $A_{a,b}$ is non-zero if and only if either $(a,b)$ or $(b,a)$ belongs to the following list:
$(\uo,1)$,$(\uo,-1)$, $(1,-1)$,$(0,1)$, $(1,2)$, $(2,-1)$,
$(-1,0)$ .
\item The set $A_{1,-1}$ is partitioned in two equivalence classes $A_{-1,1}^R,A_{-1,1}^L$ of paths which are mutually homotopic
in $\C^*\setminus \lbrace V(\la)=0 \rbrace $. 
For $(a,b) \neq (1,-1), (-1,1)$, $A_{a,b}$ is either empty or its elements belong to a unique homotopy equivalence class.
\end{itemize}
\end{lem}
The above lemma settles the topological classification of horizontal trajectories for the potential
$V(\la;\alpha,\beta,\nu)$ for all $(\alpha,\beta) \in R$: topologically we cannot distinguish the anti-Stokes complex or the
horizontal trajectories of two potentials if their parameters belong to $R$.
If we let $(\alpha,\beta)$ vary outside $R$, this may no longer be true. However,
for compact subset $K'$ of $R$, there exist an $\e>0$ such that on the $\epsilon$-neighbourhood $K'_{\e}$
of $K'$ the anti-Stokes complex $\mathcal{C}(\alpha,\beta)$ -- hence the homotopy classes of horizontal trajectories --
is topologically constant.

\begin{lem}[Horizontal trajectories in an $\e$ neighbourhood of $R$]\label{lem:defsets}
Let $K'$ be any compact subset of $R$, then there exist an $\e>0$ small enough such that on the $\epsilon$-neighbourhood $K'_{\e}$
of $K'$ the anti-Stokes complex $\widetilde{\mathcal{C}}(\alpha,\beta)$ is topologically constant.
Furthermore, we may assure that
\begin{enumerate}
   \item $\sup_{(\alpha,\beta) \in K'_\e}\inf_{\gamma \in  A_{0,\pm1}(\alpha,\beta) }\rho_{\gamma}<\infty$,
    \item $\sup_{(\alpha,\beta) \in K'_\e}\inf_{\gamma \in A^{R,L}(\alpha,\beta) }\rho_{\gamma}<\infty$,
 \end{enumerate}
 where $A^{R,L}(\alpha,\beta)$ is defined as in Figure \ref{fig:antiS}, through the anti-Stokes complex topology on $K'_{\e}$.
 \begin{proof}
Firstly, we prove that the anti-Stokes complex is topologically invariant under small perturbations around points in $R$. To this end, let us take any point $(\alpha^*,\beta^*)\in R$. To prove that the anti-Stokes complex is topologically invariant under small perturbations around $(\alpha^*,\beta^*)\in R$, it is enough to show that the same holds true for all of the individual anti-Stokes lines in $\tilde{\mathcal{C}}(\alpha^*,\beta^*)$. Let us consider for example the anti-Stokes line $\gamma^*$ in $\tilde{\mathcal{C}}(\alpha^*,\beta^*)$ which starts at $\lambda_4^*$ and connects to $+\infty$.
	
 	Firstly, it is easy to show that we can find an open ball $U_1$ around $(\alpha^*,\beta^*)$, an $R>0$ and a $0<\theta<\pi/4$, such that, for all $(\alpha,\beta)\in U_1$ and any horizontal trajectory $\gamma^\circ$, if $$\gamma^\circ \cap \{\lambda\in\mathbb{C}:|\lambda|>R, |\arg \lambda|<\theta\}\neq \emptyset,$$ then $\gamma^\circ$ must have $+\infty$ as an endpoint.
 		 	
 	Secondly, take a small open ball $U_2'$ around $(\alpha^*,\beta^*)$ such that $\lambda_4=\lambda_4(\alpha,\beta)$ does not coincide with any other turning point on $U_2'$. Let $\gamma(\alpha,\beta)$ denote the unique anti-Stokes line emanating from $\lambda_4$, with $\gamma(\alpha^*,\beta^*)=\gamma^*$, specified uniquely be the fact that the angle by which $\gamma(\alpha,\beta)$ emanates from $\lambda_4(\alpha,\beta)$ is a continuous function on $U_2'$.	
 	Then, given a point $p$ on the anti-Stokes line $\gamma^*$, and any small open ball $P$ around $p$ in the $\lambda$-plane, there exists an open ball $U_2\subseteq U_2'$ around $(\alpha^*,\beta^*)$ such that $\gamma(\alpha,\beta)\cap P\neq \emptyset$ for all $(\alpha,\beta)\in U_2$.
 	
 	Therefore, by applying the above with $p\in \gamma^*$ and $P \ni p$ an open ball with
 	 $$P\subseteq \{\lambda\in\mathbb{C}:|\lambda|>R, |\arg \lambda|<\theta\},$$
 	  we find that, setting $U:=U_1\cap U_2$, the anti-Stokes line $\gamma(\alpha,\beta)$ connects $\lambda_4$ with $+\infty$ for all $(\alpha,\beta)\in U$. In short, the anti-Stokes line $\gamma^*$ is topologically invariant under small perturbations around $(\alpha^*,\beta^*)$.
 	   	
 	The same argument works for the other anti-Stokes lines in $\tilde{\mathcal{C}}(\alpha^*,\beta^*)$, including those with endpoint $\lambda=0$, as we can find an open ball $V$ around $(\alpha^*,\beta^*)$ and a $\delta>0$, such that, for any horizontal trajectory $\gamma^\circ$, if 
 	$$\gamma^\circ \cap \{\lambda\in\mathbb{C}:|\lambda|<\delta\}\neq \emptyset,$$ then $\gamma^\circ$ must connect with $0$.
 It follows that the anti-Stokes complex is topologically invariant under small perturbations around points in $R$, that is, we can find an open disc around $(\alpha^*,\beta^*)$ such that on this disc the anti-Stokes complex $\tilde{\mathcal{C}}(\alpha,\beta)$ is topologically identical to the one at $(\alpha,\beta)=(\alpha^*,\beta^*)$.
 	
 	 By a straightforward compactness argument we may find an $\tilde{\epsilon}$-neighbourhood of $K'$ on which the anti-Stokes complex is topologically constant. We now set $\epsilon=\frac{1}{2}\tilde{\epsilon}$ as this will assure that the suprema in the lemma are also finite, since  $\overline{K_\epsilon'}\subseteq K_{\tilde{\epsilon}}'$ is compact. 		
 \end{proof}
\end{lem}

\subsection{Lifting WKB approximations to the elliptic curve}\label{sub:lift}
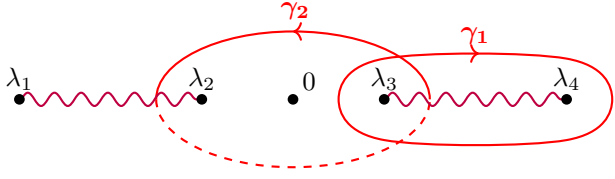
\begin{figure}
	\centering
	\begin{tikzpicture}[scale=0.6]

	\tikzstyle{star}  = [circle, minimum width=4pt, fill, inner sep=0pt];

	\node[above]     at (-6,0) {$\lambda_1$};

	\node[above]     at (-2,0) {$\lambda_2$};

	\node[above right]     at (0,0) {0};

	\node[above]     at (2,0) {$\lambda_3$};

	\node[above]     at (6,0) {$\lambda_4$};

	\draw[thick, red,decoration={markings, mark=at position 0.5 with {\arrow{>}}},postaction={decorate}] (0,0) [partial ellipse=0:180:3cm and 1.5cm];
	\draw[thick, red,dashed] (0,0) [partial ellipse=-180:0:3cm and 1.5cm];

	\draw [red,thick,decoration={markings, mark=at position 0.15 with {\arrow{>}}},
	postaction={decorate}] plot [smooth cycle,tension=0.6] coordinates {(1,0) (2,0.9)  (6, 0.9) (7,0) (6,-0.9)  (2,-0.9)  };	
	
	\node[red]     at (4,1.4) {$\boldsymbol{\gamma_1}$};
	\node[red]     at (0.05,1.95) {$\boldsymbol{\gamma_2}$};

	\draw[purple,thick,snake it] (-6,0) -- (-2,0);
	\draw[purple,thick,snake it] (2,0) -- (6,0);

	\node[black,star]     at (-6,0){};
	
	\node[black,star]     at (-2,0) {};
	
	\node[black,star]     at (0,0) {};

	\node[black,star]     at (2,0) {};

	\node[black,star]     at (6,0) {};
	
	\end{tikzpicture}
	\caption{Topological representation of cycles $\gamma_{1,2}$ in
	red on elliptic curve $\widehat{\Gamma}(\alpha,\beta)$, for $(\alpha,\beta)\in R$, 
	with	branch-cuts the Stokes lines with endpoints $\{\lambda_1,\lambda_2\}$ and $\{\lambda_3,\lambda_4\}$, depicted in wiggly purple and where
	$\lambda_k$, $1\leq k\leq 4$, denote the branching points of the elliptic curve.} \label{fig:cutplane_and_cycles}
\end{figure}

Recall from \eqref{eq:defigamma} the definition of the elliptic curve $\Gamma=\{(\la,y) \in \C^2, \, y^2=\la^2 V(\la)\}$
and of its compactification
$\widehat{\Gamma}$. It is convenient, for our computations, to realise the elliptic curve as a doubly
sheeted cover of the plane, see Figure \ref{fig:cutplane_and_cycles}: we choose two cuts
connecting two distinct pairs of zeros of $V$ so that the function $\sqrt{V}$
is single-valued on the cut plane. We name the lower sheet the one fixed by the requirement $\lim_{\la \to +\infty }\Re \sqrt{V(\la)}=+\infty$
(equivalently, $\res_{\la=0}\sqrt{V(\la)}=\frac{\nu}{2}$), and
to represent a curve in $\widehat{\Gamma}$ as a curve in the two-sheeted covering, we draw a solid line when 
the curve belong to the upper sheet, and a dashed line otherwise.
In case $(\alpha,\beta) \in R$, we choose the two cuts to be the Stokes lines connecting $\la_1$ with $\la_2$ and
$\la_3$ with $\la_4$ and the basis
$\gamma_1,\gamma_2$ of $H_1(\widehat{\Gamma},\bb{Z})$ as in Figure \ref{fig:cutplane_and_cycles}.

If $K'$ is a compact subset of $R$ and
$K'_\e$ an $\e$ neighbourhood of $K'$, small enough so that the roots of $V$ are all distinct,
this choice of the cuts and of the cycles $\gamma_1,\gamma_2$ extends continuously to $K'_\e$; as a consequence,
the elliptic integral $\int_{\gamma_i} \sqrt{V(\la)}d\la ,i=1,2$ is an analytic function on $K'_\e$.

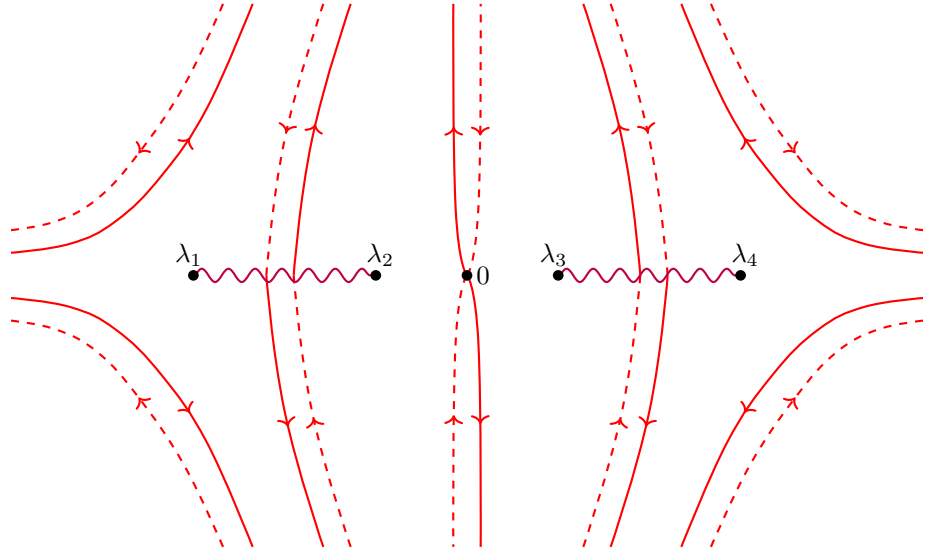
\begin{figure}
	\centering
	\begin{tikzpicture}[scale=0.6]
	
	\tikzstyle{star}  = [circle, minimum width=4pt, fill, inner sep=0pt];
	
	\node[above]     at (-6.1,0) {$\lambda_1$};

	\node[above]     at (-1.9,0) {$\lambda_2$};

	\node[right]     at (0,0) {0};

	\node[above]     at (1.9,0) {$\lambda_3$};

	\node[above]     at (6.1,0) {$\lambda_4$};
	
	%
	%
	%
	
	\draw [red,thick,dashed,decoration={markings, mark=at position 0.55 with {\arrow{<}}},
	postaction={decorate}] plot [smooth,tension=0.8] coordinates {  (0,0) (0.25,1.5) (0.3,6)   };
	\draw [red,thick,decoration={markings, mark=at position 0.55 with {\arrow{>}}},
	postaction={decorate}] plot [smooth,tension=0.8] coordinates {  (0,0) (-0.25,1.5) (-0.3,6)   };
	\draw [red,thick,dashed,decoration={markings, mark=at position 0.55 with {\arrow{<}}},
	postaction={decorate}] plot [smooth,tension=0.8] coordinates {  (0,0) (-0.25,-1.5) (-0.3,-6)   };
	\draw [red,thick,decoration={markings, mark=at position 0.55 with {\arrow{>}}},
	postaction={decorate}] plot [smooth,tension=0.8] coordinates {  (0,0) (0.25,-1.5) (0.3,-6)   };

	\draw [red,thick,dashed,decoration={markings, mark=at position 0.55 with {\arrow{<}}},
	postaction={decorate}] plot [smooth,tension=0.8] coordinates {  (4.4,0) (4.0,3) (3.15,6)   };
	\draw [red,thick,decoration={markings, mark=at position 0.55 with {\arrow{>}}},
	postaction={decorate}] plot [smooth,tension=0.8] coordinates { (4.4,0)  (4.0,-3)  (3.15,-6)   };
	
	\draw [red,thick,decoration={markings, mark=at position 0.55 with {\arrow{>}}},
	postaction={decorate}] plot [smooth,tension=0.8] coordinates {  (3.8,0) (3.4,3) (2.55,6)   };
	\draw [red,thick,dashed,decoration={markings, mark=at position 0.55 with {\arrow{<}}},
	postaction={decorate}] plot [smooth,tension=0.8] coordinates {(3.8,0)   (3.4,-3)  (2.55,-6)   };
	

	
	\draw [red,thick,dashed,decoration={markings, mark=at position 0.55 with {\arrow{<}}},
	postaction={decorate}] plot [smooth,tension=0.8] coordinates {  (-4.4,0) (-4.0,3) (-3.15,6)   };
	\draw [red,thick,decoration={markings, mark=at position 0.55 with {\arrow{>}}},
	postaction={decorate}] plot [smooth,tension=0.8] coordinates { (-4.4,0) (-4.0,-3) (-3.15,-6)    };
	
	\draw [red,thick,decoration={markings, mark=at position 0.55 with {\arrow{>}}},
	postaction={decorate}] plot [smooth,tension=0.8] coordinates {  (-3.8,0) (-3.4,3) (-2.55,6)   };
	\draw [red,thick,dashed,decoration={markings, mark=at position 0.55 with {\arrow{<}}},
	postaction={decorate}] plot [smooth,tension=0.8] coordinates {(-3.8,0)  (-3.4,-3)  (-2.55,-6)   };

	\draw [red,thick,dashed,decoration={markings, mark=at position 0.5 with {\arrow{<}}},
	postaction={decorate}] plot [smooth,tension=0.7] coordinates { (10,1) (9.0,1.2)  (8.1,1.7)    (7.2,2.7) (6.2,4.2) (5.35,6)   };
	\draw [red,thick,decoration={markings, mark=at position 0.6 with {\arrow{>}}},
	postaction={decorate}] plot [smooth,tension=0.7] coordinates { (10,0.5) (8.7,0.7)  (7.7,1.2)    (6.7,2.2) (5.7,3.7) (4.7,6)   };
	
	\draw [red,thick,dashed,decoration={markings, mark=at position 0.5 with {\arrow{<}}},
	postaction={decorate}] plot [smooth,tension=0.7] coordinates { (-10,1) (-9.0,1.2)  (-8.1,1.7)    (-7.2,2.7) (-6.2,4.2) (-5.35,6)   };
	\draw [red,thick,decoration={markings, mark=at position 0.6 with {\arrow{>}}},
	postaction={decorate}] plot [smooth,tension=0.7] coordinates { (-10,0.5) (-8.7,0.7)  (-7.7,1.2)    (-6.7,2.2) (-5.7,3.7) (-4.7,6)   };	
	
	\draw [red,thick,dashed,decoration={markings, mark=at position 0.5 with {\arrow{<}}},
	postaction={decorate}] plot [smooth,tension=0.7] coordinates { (-10,-1) (-9.0,-1.2)  (-8.1,-1.7)    (-7.2,-2.7) (-6.2,-4.2) (-5.35,-6)   };
	\draw [red,thick,decoration={markings, mark=at position 0.6 with {\arrow{>}}},
	postaction={decorate}] plot [smooth,tension=0.7] coordinates { (-10,-0.5) (-8.7,-0.7)  (-7.7,-1.2)    (-6.7,-2.2) (-5.7,-3.7) (-4.7,-6)   };
	
	\draw [red,thick,dashed,decoration={markings, mark=at position 0.5 with {\arrow{<}}},
	postaction={decorate}] plot [smooth,tension=0.7] coordinates { (10,-1) (9.0,-1.2)  (8.1,-1.7)    (7.2,-2.7) (6.2,-4.2) (5.35,-6)   };
	\draw [red,thick,decoration={markings, mark=at position 0.6 with {\arrow{>}}},
	postaction={decorate}] plot [smooth,tension=0.7] coordinates { (10,-0.5) (8.7,-0.7)  (7.7,-1.2)    (6.7,-2.2) (5.7,-3.7) (4.7,-6)   };


	\draw[purple,thick,snake it] (-6,0) -- (-2,0);
	\draw[purple,thick,snake it] (2,0) -- (6,0);

	\node[black,star]     at (-6,0){};
	
	\node[black,star]     at (-2,0) {};
	
	\node[black,star]     at (0,0) {};

	\node[black,star]     at (2,0) {};

	\node[black,star]     at (6,0) {};

	%
	%
	%
	%

	\end{tikzpicture} 
	\caption{In red, for both orientations and for every homotopy class,
		the lift on the elliptic curve
		$\widehat{\Gamma}$ of a horizontal trajectory of the potential $V(\la;\alpha,\beta,\nu)$,
		with $(\alpha,\beta) \in R$.
		The elliptic curve is realised as a two sheeted-cover of the cut plane, with branch-cuts depicted in wiggly purple.
		A solid line represents a curve on the upper sheet; A dashed line represents a curve on the lower sheet.} \label{fig:lifted}
\end{figure}

Now recall, from Lemma \ref{lem:domsub}, that a solution subdominant at $0$ or in a sub-sector $\Sigma_k$
is well-approximated on a curve $\gamma \in \overrightarrow{A}_{a,b}$ with $a=\uo$ or $a=k$ by the WKB function
\eqref{eq:wkb_functions}
$$
\Psi(\la;\la_0)= \exp\left\{\int^{\la}_{\la_0,\gamma}E \sqrt{V(\mu)} - \frac{V'(\mu)}{4 V(\mu)} d \mu \right\},
$$
provided the branch of $\sqrt{V}$ is chosen in such a way that
\begin{equation}\label{eq:sheetcondition}
 \lim_{t\rightarrow {0^+}} \Re \int_{t}^{t_0} \sqrt{V(\gamma(t))}\dot{\gamma}(t) dt =+\infty.
\end{equation}
Since the WKB function is written in terms of the integral of $\sqrt{V}$, it is naturally defined
on the lift of $\gamma$ to $\widehat{\Gamma}$, which we denote by $\widehat{\gamma}$.
This is not only natural, but also very convenient since there is a unique way of lifting $\gamma$ that enforces
condition \eqref{eq:sheetcondition}. In fact,
taking into consideration our choice of the branch-cuts of $\sqrt{V}$, 
for any $(a,b)$ and any path $\gamma \in \overrightarrow{A}_{a,b}$, the lift
$\widehat{\gamma}$ satisfies \eqref{eq:sheetcondition} provided
provided:
\begin{itemize}
	\item $\widehat{\gamma}(t)$ lies on the upper sheet as  $t \to 0$ if $a\in\{\underline{0},0,2\}$.
	\item $\widehat{\gamma}(t)$ lies on the lower sheet as  $t \to 0$ if $a\in\{1,-1\}$.
\end{itemize}
This situation is illustrated in Figure \ref{fig:lifted}.

\subsection{Bohr-Sommerfeld quantisation conditions.}
Recall from Lemma \ref{lem:analiticity} that there are two functions $U_1,U_2$, depending holomorphically on
the parameters $\alpha,\beta$ of the potential, which vanish whenever
the quantisation and apparent singularity conditions are met.
Here we prove asymptotic estimates for these two functions if $E$ is large and $(\alpha,\beta)$ belong to a neighbourhood of $R$.
As it turns out, the potentials $V(\la;\alpha,\beta,\nu)$, with $\nu$ fixed,
which satisfy the quantisation and apparent singularity conditions,
are approximate solutions of the Bohr-Sommerfeld-Boutroux system \eqref{eq:bohrsommerfeld}.

In the proofs below we use the following notation. Let $K' \subset R$ be a compact subset of $R$, and $\e$ a positive number.
We denote by $K'_\e$ the $\e$-neighbourhood of $K'$.

We begin with the analysis of the quantisation condition.
\begin{thm}\label{thm:wkbquantisation}
Fix a $\nu \in (0,\frac 13]$ (not necessarily rational), and a compact subset $K' \subset R$.
There exist positive constants $\e,E_1,C_1$ -- depending on $K'$ and $\nu$ -- and a function
$$U_1: [E_1,\infty) \times K'_\e \to \C, \; (E,\alpha,\beta) \mapsto U_1(E,\alpha,\beta),$$
such that 
\begin{itemize}
\item For all $E \in [E_1,\infty) $, $U_1(E,\cdot,\cdot):  K'_\e \to \C$ is analytic.
 \item If $(E,\alpha,\beta) \in [E_1,\infty) \times K'_\e$, then equation
 \eqref{eq:scaled} admits a non-trivial solution of the boundary value problem
		\begin{equation*}
		 \lim_{\la \to +\infty} \psi(\la)=\lim_{\la \to 0^+} \psi(\la)=0,
		\end{equation*}
		if and only if $U_1(E,\alpha,\beta)=0$.
 \item $\forall (E,\alpha,\beta) \in [E_1,\infty) \times K'_\e$, $U_1$ satisfies the estimate
 \begin{equation}\label{eq:quantisationwkb}
 \left| U_1(E,\alpha,\beta) - 1 - \exp\big(E\oint_{\gamma_1}\sqrt{V(\mu)}d\mu \big) \right|\leq \frac{C_1}{E}
 \exp\big(E \Re \oint_{\gamma_1}\sqrt{V(\mu)}d\mu \big) .
\end{equation}
\end{itemize}
\end{thm}
\begin{proof}
We consider solutions of \eqref{eq:scaled} on the cut-plane $D_{\varphi}$ with $\varphi=\pi$
(the cut coincides with $\bb{R}_-$ for large $\la$). We denote by $\chi_+,\psi_{\pm1},\psi_0:D_{\pi} \to \mathbb{C}$,
a solution subdominant at $0,\Sigma_{\pm1},\Sigma_0$ (which is uniquely defined up to a normalising constant).

According to Proposition \ref{prop:wronskiantovalues}(1),  under the hypothesis
   that $\psi_{+1},\psi_{-1}$ are linearly independent, the quantisation condition
   is equivalent to the equality $w_{-1}\big( \frac{\chi_+}{\psi_{0}} \big)=w_{1}\big(\frac{\chi_+}{\psi_{0}}\big)$.

The function $U_1$, which has been introduced in Lemma \ref{lem:analiticity},
\begin{equation}\label{eq:U1}
 U_1=1-\frac{w_{1}\big(\frac{\chi_+}{\psi_{0}}\big)}{w_{-1}\big( \frac{\chi_+}{\psi_{0}} \big)},
\end{equation}
vanishes when $w_{-1}\big( \frac{\chi_+}{\psi_{0}} \big)=w_{1}\big(\frac{\chi_+}{\psi_{0}}\big)$.

To prove the theorem, it is sufficient to show that there exists an $\e>0$ such that
   \begin{enumerate}
    \item If $E$ is large enough, $\psi_{+1},\psi_{-1}$ are linearly independent
    for all $(\alpha,\beta) \in K'_{\e}$ (so that $U_1$ vanishes if and only if the quantisation condition is satisfied).
    \item Fixed $E$, $U_1$ is an analytic function on $K'_{\e}$.
    \item Estimate \eqref{eq:quantisationwkb} holds.
   \end{enumerate}
   
In order to prove these three statements we make use of the horizontal trajectories, which are classified in Lemmas \ref{lem:nonemptysets}
and \ref{lem:defsets}. These lemmas in particular states that if $\e$ is small enough
the sets of horizontal trajectories $A_{0,\pm1}(\alpha,\beta),A_{\uo,\pm1}(\alpha,\beta),A^R(\alpha,\beta)$,
as depicted in Figure \ref{fig:antiS} are non empty.
Moreover, $$\sup_{\alpha,\beta \in K'_\e} \inf_{\gamma \in A(\alpha,\beta)} \rho_{\gamma} \leq C ,$$
where $A$ is one of the five sets of horizontal trajectories listed above, and $C$ a positive constant.

(1) For any $\gamma \in A^R$, we can choose the cut $l$, connecting $0$ to $\infty$ and asymptotic to the negative real semi-axis,
in such a way that $\gamma(]0,1[) \subset D_{\pi}:=\C \setminus l$. 
It follows from Lemma \ref{lem:domsub} that
 $\psi_{1}$ and $\psi_{-1}$ are linearly independent provided $E>\frac{C}{\ln2}$.

(2) After Lemma \ref{lem:analiticity}, for fixed $E,\nu$, the function $U_1$ is
a well-defined meromorphic function of the parameters $(\alpha,\beta)$ of the potential.
We are left to show that $U_1$ has no poles in $K'_{\e}$.
To this aim, it is sufficient to check that $w_{-1}\big( \frac{\chi_+}{\psi_{0}}\big) \neq 0$
and $w_{1}\big( \frac{\chi_+}{\psi_{0}}\big) \neq \infty$. We notice that
$w_{-1}\big( \frac{\chi_+}{\psi_{0}}\big)=0$ implies that $\chi_+$ and $\psi_{-1}$ are linearly dependent, while 
$w_{1}\big( \frac{\chi_+}{\psi_{0}}\big)= \infty$ implies that $\psi_0$ and $\psi_1$ are linearly dependent.
Since $A_{\uo,-1}$ and $A_{0,1}$ are non-empty,
neither equality can be satisfied for $E$ large enough (by the same reasoning as in (1) above).

(3) In order to derive estimate \eqref{eq:quantisationwkb}, we approximate the subdominant solutions $\chi_+$ and $\psi_0$
by the WKB function $$\Psi(\la;\la')=\exp\lbrace \int_{\la'}^\la E \sqrt{V(\mu)}-\frac{V'(\mu)}{4V(\mu)} d\mu \rbrace.$$
To this aim we select four admissible curves, as in the Figure \ref{fig:y1}:
$\alpha_{\pm1}$ connects $0$ with $\pm i \infty$,
and $\beta_{\pm1}$ connects $+\infty$ with $\pm i \infty$. We take some arbitrary paths
 $\widetilde{\alpha}_{\pm1} \in
\overrightarrow{A}_{\uo,\pm1}$, and $\widetilde{\beta}_{\pm1} \in
\overrightarrow{A}_{0,\pm1}$. According to Lemma \ref{lem:pathdef}, we can deform these four paths,
without increasing the corresponding the error
$\rho$'s, in such a way any two paths, which share the same endpoint, coincide in a neighbourhood of it. 
The paths $\alpha_{\pm1},\beta_{\pm1}$ are the result of such a deformation:
for some  $\delta>0$ we have that
\begin{itemize}
\item $\alpha_{1}(t)=\alpha_{-1}(t)$ for $t\in [0,\delta]$;
\item $\beta_{1}(t)=\beta_{-1}(t)$ for $t\in [0,\delta]$;
\item $\alpha_{1}(t)=\beta_{1}(t)$ for $t\in [1-\delta,1]$;
\item $\alpha_{-1}(t)=\beta_{-1}(t)$ for $t\in [1-\delta,1]$.
\end{itemize}
As explained in Subsection \ref{sub:lift}, it is convenient to lift admissible curves on the Riemann surface $\widehat{\Gamma}$,
in order to enforce
the correct choice of branch of $\sqrt{V}$. Since the paths $\alpha_{\pm1},\beta_{\pm1}$
starts at $\Sigma_0,0$, then the lift is defined by imposing that they belong to the upper sheet of the curve
as $t\to 0$.
 If we denote by
$\widehat{\alpha}_{\pm1},\widehat{\beta}_{\pm1}$ the lift of
$\alpha_{\pm1},\beta_{\pm1}$, the situation is as depicted in Figure \ref{fig:y2}.
We notice that, since the four lifted paths coincide pairwise near (the lift of) $0$ and $\infty$, their concatenation
$\widehat{\alpha}_{-1}-\widehat{\beta}_{-1}+\widehat{\beta}_1-
\widehat{\alpha}_1$ can be considered as a
closed loop in $\widehat{\Gamma}_{\alpha,\beta}$. In fact, we have
\begin{equation}\label{eq:homologygamma1}
\widehat{\alpha}_{1}-\widehat{\beta}_{1}+\widehat{\beta}_{-1}-
\widehat{\alpha}_{-1}= \gamma_1 \in H_1(\widehat{\Gamma}_{\alpha,\beta},\bb{Z})
\end{equation}
where $\gamma_1$ is the basis element of $H_1(\widehat{\Gamma}_{\alpha,\beta},\bb{Z})$, defined in Figure \ref{fig:cutplane_and_cycles}.

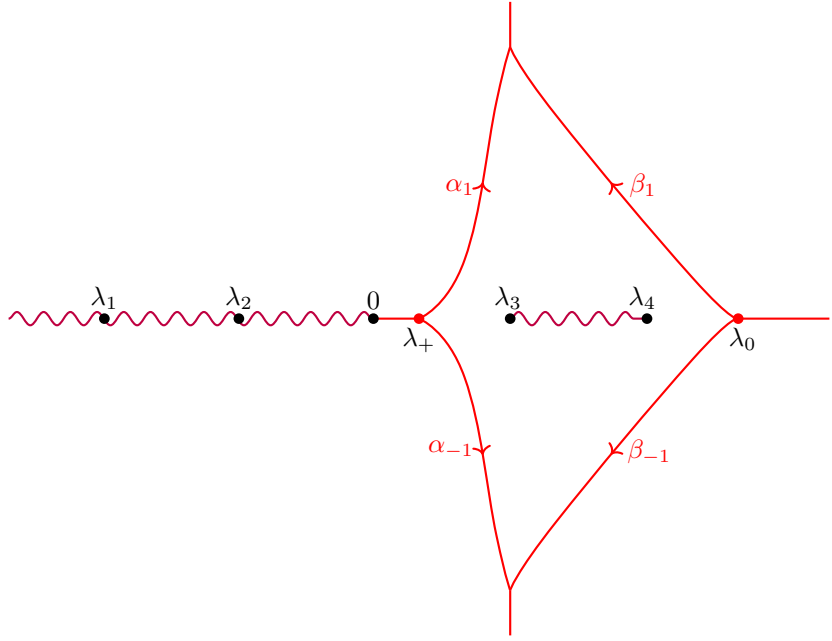
\begin{figure}
	\centering
	\begin{tikzpicture}[scale=0.6]

	\tikzstyle{star}  = [circle, minimum width=4pt, fill, inner sep=0pt];
	
	\node[above]     at (-5.9,0) {$\lambda_1$};

	\node[above]     at (-2.95,0) {$\lambda_2$};

	\node[above]     at (0,0) {0};

	\node[above]     at (2.95,0) {$\lambda_3$};

	\node[above]     at (5.9,0) {$\lambda_4$};
	
	\node[below]     at (1,0) {$\lambda_+$};
	\node[below]     at (8.1,0) {$\lambda_0$};
	
	\draw[red,thick] (0,0) -- (1,0);
	\draw[red,thick] (3,6) -- (3,7);
	\draw[red,thick] (3,-6) -- (3,-7);
	\draw[red,thick] (8,0) -- (10,0);

	\draw [red,thick,decoration={markings, mark=at position 0.5 with {\arrow{<}}},
	postaction={decorate}] plot [smooth,tension=0.8] coordinates { (3,6) (2.7,4.8)    (2,1.3)   (1,0) };
	\draw [red,thick,decoration={markings, mark=at position 0.5 with {\arrow{<}}},
	postaction={decorate}] plot [smooth,tension=0.8] coordinates { (3,-6) (2.7,-4.8)    (2,-1.3)   (1,0) };

	\draw [red,thick,decoration={markings, mark=at position 0.5 with {\arrow{<}}},
	postaction={decorate}] plot [smooth,tension=0.8] coordinates { (3,6) (4.0,4.5) (7,0.9) (8,0)  };
	\draw [red,thick,decoration={markings, mark=at position 0.5 with {\arrow{<}}},
	postaction={decorate}] plot [smooth,tension=0.8] coordinates { (3,-6) (4.0,-4.5) (7,-0.9) (8,0)  };

	\node[red]     at (1.9,2.9) {$\alpha_1$};
	\node[red]     at (5.91,2.95) {$\beta_1$};
	\node[red]     at (1.7,-2.9) {$\alpha_{-1}$};
	\node[red]     at (6.05,-2.95) {$\beta_{-1}$};

	\draw[purple,thick,snake it] (-8,0) -- (0,0);
	\draw[purple,thick,snake it] (3,0) -- (6,0);
	
	\node[red,star]     at (1,0){};
	\node[red,star]     at (8,0){};

	\node[black,star]     at (-5.9,0){};
	
	\node[black,star]     at (-2.95,0) {};
	
	\node[black,star]     at (0,0) {};

	\node[black,star]     at (3,0) {};

	\node[black,star]     at (6,0) {};

	\end{tikzpicture} 
	\caption{The integration paths used in the proof of Theorem \ref{thm:wkbquantisation}.} \label{fig:y1}
\end{figure}

\begin{figure}
	\centering
	\begin{tikzpicture}[scale=0.6]

	\tikzstyle{star}  = [circle, minimum width=4pt, fill, inner sep=0pt];
	
	\node[above]     at (-5.9,0) {$\lambda_1$};

	\node[above]     at (-2.95,0) {$\lambda_2$};

	\node[above]     at (0,0) {0};

	\node[above]     at (3,0) {$\lambda_3$};

	\node[above]     at (6,0) {$\lambda_4$};
	
	\node[below]     at (1.2,0) {$\widehat{\lambda}_+$};
	\node[below]     at (8.1,0) {$\widehat{\lambda}_0$};
	
	\draw[red,thick] (0,0) -- (1,0);
	\draw[red,thick] (3,6) -- (3,7);
	\draw[red,thick] (3,-6) -- (3,-7);
	\draw[red,thick] (8,0) -- (10,0);

	\draw [red,thick,decoration={markings, mark=at position 0.5 with {\arrow{<}}},
	postaction={decorate}] plot [smooth,tension=0.8] coordinates { (3,6) (2.7,4.8)    (2,1.3)   (1,0) };
	\draw [red,thick,decoration={markings, mark=at position 0.5 with {\arrow{<}}},
	postaction={decorate}] plot [smooth,tension=0.8] coordinates { (3,-6) (2.7,-4.8)    (2,-1.3)   (1,0) };

	\draw [red,thick,decoration={markings, mark=at position 0.5 with {\arrow{<}}},
	postaction={decorate}] plot [smooth,tension=0.8] coordinates { (3,6) (4.0,4.5) (7,0.9) (8,0)  };
	\draw [red,thick,decoration={markings, mark=at position 0.5 with {\arrow{<}}},
	postaction={decorate}] plot [smooth,tension=0.8] coordinates { (3,-6) (4.0,-4.5) (7,-0.9) (8,0)  };

	\draw [red,thick,decoration={markings, mark=at position 0.5 with {\arrow{>}}},
	postaction={decorate}] plot [smooth,tension=0.8] coordinates { (1,0)  (4.5,1.2)   (8,0) };
	\draw [red,thick,decoration={markings, mark=at position 0.5 with {\arrow{>}}},
	postaction={decorate}] plot [smooth,tension=0.8] coordinates { (1,0)  (4.5,-1.2)   (8,0) };	
	
	\node[red]     at (4.55,1.7) {$\gamma_1^+$};
	\node[red]     at (4.55,-1.65) {$\gamma_1^-$};	
		
	\node[red]     at (1.9,2.9) {$\widehat{\alpha}_1$};
	\node[red]     at (5.91,2.95) {$\widehat{\beta}_1$};
	\node[red]     at (1.7,-2.9) {$\widehat{\alpha}_{-1}$};
	\node[red]     at (6.05,-2.95) {$\widehat{\beta}_{-1}$};

	\draw[purple,thick,snake it] (-6,0) -- (-3,0);
	\draw[purple,thick,snake it] (3,0) -- (6,0);
	
	\node[red,star]     at (1,0){};
	\node[red,star]     at (8,0){};

	\node[black,star]     at (-5.9,0){};
	
	\node[black,star]     at (-2.95,0) {};
	
	\node[black,star]     at (0,0) {};

	\node[black,star]     at (3,0) {};

	\node[black,star]     at (6,0) {};

	\end{tikzpicture} 
	\caption{The integration paths used in the proof of Theorem \ref{thm:wkbquantisation}, lifted on the elliptic curve.
	The elliptic curve is realised as a two sheeted-cover of the
cut plane, with branch-cuts depicted in wiggly purple. A solid line represents a
curve on the upper sheet; A dashed line represents a curve on the lower sheet.} \label{fig:y2}
\end{figure}
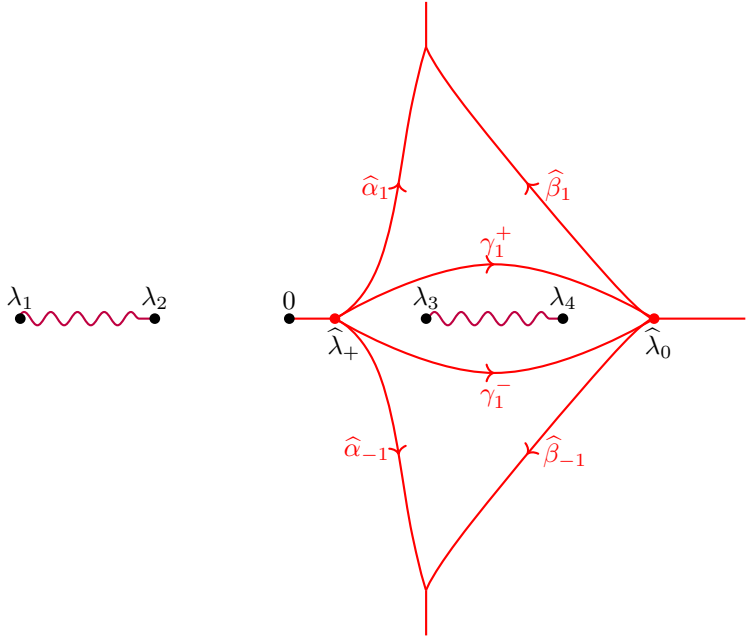

Having fixed the four admissible curves and their lift we can complete the proof of the theorem by proving \eqref{eq:quantisationwkb}.
We choose a $\la_+ \in \widehat{\alpha}_{1}\cap \widehat{\alpha}_{-1}$.
After Proposition \ref{pro:continuation}, there exists a solution 
$\chi_+(\la)$, subdominant at $0$, which admits the following estimate:
There exist $C_1, E_1$ such that for all $E\geq E_{1}$ and all $(\alpha,\beta) \in K'_{\e}$,
\begin{equation}\label{eq:chi+estimate}
 \left| \chi_+\big(\la \big)/\exp\left\lbrace\int_{\la_+,\widehat{\alpha}_{\pm1}}^{\widehat{\la}}
 E \sqrt{V(\mu)}-\frac{V'(\mu)}{4V(\mu)}d\mu \right\rbrace-1
 \right|\leq \frac{C_{1}}{E}, \; \forall \la \in \alpha_{\pm1},
\end{equation}
where $\widehat{\la}$ is the lift of $\la$ along $\widehat{\alpha}_{\pm1}$ .

Similarly, we choose a $\la_0 \in \widehat{\beta}_{1}\cap \widehat{\beta}_{-1}$.
There exists a solution 
$\psi_0(\la)$, subdominant in $\Sigma_0$, which admits the following estimate.
There exist $C_2, E_2$ such that for all $E\geq E_{2}$ and all $(\alpha,\beta) \in K'_{\e}$,
\begin{equation}\label{eq:psi0+estimate}
 \left| \psi_0\big(\la \big)/\exp\left\lbrace\int_{\la_0,\widehat{\beta}_{\pm1}}^{\widehat{\la}}
  E \sqrt{V(\mu)}-\frac{V'(\mu)}{4V(\mu)}d\mu\right\rbrace-1
 \right|\leq \frac{C_{2}}{E}, \; \forall \la \in \beta_{\pm1},
\end{equation}
where $\widehat{\la}$ is the lift of $\la$ along $\widehat{\beta}_{\pm1}$ .

We use (\ref{eq:chi+estimate},\ref{eq:psi0+estimate}) to derive the following estimate for $w_{1}\big( \frac{\chi_+}{\psi_{0}} \big)$.
There exist $C_+,E_+$ such that for all $(\alpha,\beta) \in K'_\e$,
\begin{equation}\label{eq:w1estimate}
 w_{1}\big( \frac{\chi_+}{\psi_{0}} \big)=\exp\left\lbrace\int_{\la_+,\gamma_1^+}^{\la_0}
 E \sqrt{V(\mu)}-\frac{V'(\mu)}{4V(\mu)} d\mu \right \rbrace \left( 1+ \e_+(\alpha,\beta) \right)
\end{equation}
where $\gamma_1^+=\widehat{\alpha}_1-\widehat{\beta}_1$, as depicted in Figure \ref{fig:y2}, and  $\e_+(\alpha,\beta)$  satisfies
the estimate
\begin{equation}\label{eq:estimatee+}
 |\e_+(\alpha,\beta)|\leq \frac{C_+}{E}\,, \quad \forall E\geq E_+.
\end{equation}
In order to derive (\ref{eq:w1estimate},\ref{eq:estimatee+}), we start by applying the definition of $w_1$:
$$
w_{1}\big( \frac{\chi_+}{\psi_{0}} \big)=\lim_{t \to 1}\frac{\chi_+(\alpha_1(t))}{\psi_0(\beta_1(t))}.$$
Using the inequalities (\ref{eq:chi+estimate},\ref{eq:psi0+estimate}), we immediately obtain that
there exist $C_+,E_+ >0$ such that for all $(\alpha,\beta) \in K'_\e$,
$$
\lim_{t \to 1}\frac{\chi_+(\alpha_1(t))}{\psi_0(\beta_1(t))}=
\lim_{t \to 1}\frac{\exp\left\lbrace\int_{\la_+,\widehat{\alpha}_{1}}^{\widehat{\alpha}_1(t)}
 E \sqrt{V(\mu)}-\frac{V'(\mu)}{4V(\mu)}d\mu \right\rbrace }
 {\exp\left\lbrace\int_{\la_0,\widehat{\beta}_{1}}^{\widehat{\beta}_1(t)}
 E \sqrt{V(\mu)}-\frac{V'(\mu)}{4V(\mu)}d\mu \right\rbrace }
  \left( 1+ \e_{+}(\alpha,\beta)\right)
$$
where $|\e_+|$ satisfies \eqref{eq:estimatee+}.

Finally we utilise the fact that $\widehat{\alpha}_1(t)=\widehat{\beta}_1(t)$ for $t \to 1$. Because of this,
$\widehat{\alpha}_1-\widehat{\beta}_1$ is to all effect a curve in $\widehat{\Gamma}$
homotopic to the curve $\gamma_1^+$, depicted in Figure \ref{fig:y2}. Therefore
$$
\lim_{t \to 1}\frac{\exp\left\lbrace\int_{\la_+,\widehat{\alpha}_{1}}^{\widehat{\alpha}_1(t)}
E \sqrt{V(\mu)}-\frac{V'(\mu)}{4V(\mu)}d\mu\right\rbrace}
{\exp\left\lbrace\int_{\la_0,\widehat{\beta}_{1}}^{\widehat{\beta}_1(t)}
 E \sqrt{V(\mu)}-\frac{V'(\mu)}{4V(\mu)}d\mu \right\rbrace }=
\exp\left\lbrace\int_{\la_+,\gamma_1^+}^{\la_0}
E \sqrt{V(\mu)}-\frac{V'(\mu)}{4V(\mu)}d\mu \right \rbrace,$$
from which (\ref{eq:w1estimate},\ref{eq:estimatee+}) follow.

By the same method, we obtain a similar estimate for  $w_{-1}\big( \frac{\chi_+}{\psi_{0}} \big)$.
There exists $C_-,E_->0$ such that for all $(\alpha,\beta)\in K'_\e$
\begin{equation}\label{eq:w-1estimate}
  w_{-1}\big( \frac{\chi_+}{\psi_{0}} \big)= \exp\left\lbrace 
 \int_{\la_+,\gamma_1^-}^{\la_0}E\sqrt{V(\mu)}-\frac{V'(\mu)}{4V(\mu)}d\mu \right\rbrace \left(1 +\e_-(\alpha,\beta)\right),
\end{equation}
where $\gamma_1^-=\widehat{\alpha}_{-1}-\widehat{\beta}_{-1}$, as depicted in Figure \ref{fig:y2}, and
\begin{equation}\label{eq:estimatee-}
 |\e_-(\alpha,\beta)|\leq \frac{C_-}{E}\,, \quad \forall E\geq E_-.
\end{equation}
 Finally we show that (\ref{eq:w1estimate},\ref{eq:estimatee+},\ref{eq:w-1estimate},\ref{eq:estimatee-}) implies \eqref{eq:quantisationwkb}.
 Since $\gamma_+-\gamma_-= \widehat{\alpha}_{1}-\widehat{\beta}_{1}+\widehat{\beta}_{-1}-
\widehat{\alpha}_{-1} $, then from \eqref{eq:homologygamma1} we obtain that $\gamma_+-\gamma_-= \gamma_1$.

Hence
from (\ref{eq:U1},\ref{eq:w1estimate},\ref{eq:w-1estimate}) we obtain
 $$
 U_1-1=-\exp\left\lbrace 
 \int_{\gamma_1}E\sqrt{V(\mu)}-\frac{V'(\mu)}{4V(\mu)}d\mu \right\rbrace \frac{1+\e_-(\alpha,\beta)}{1+\e_+(\alpha,\beta)}.
 $$
 By the residue theorem, $\exp\big\lbrace -\int_{\gamma_1}\frac{V'(\mu)}{4V(\mu)}d\mu\rbrace=-1$. Inserting this identity in the former
 display, we obtain that 
 $$
 U_1-1 - \exp\left\lbrace 
 \int_{\gamma_1}E\sqrt{V(\mu)}d\mu \right\rbrace=\exp\left\lbrace 
 \int_{\gamma_1}E\sqrt{V(\mu)}d\mu \right\rbrace \frac{\e_-(\alpha,\beta)-\e_+(\alpha,\beta)}{1+\e_+(\alpha,\beta)}.
 $$
 Using
 (\ref{eq:estimatee+},\ref{eq:estimatee-}), we obtain \eqref{eq:quantisationwkb}.
\end{proof}

\begin{thm}\label{thm:wkbmonodromy}
Fix a $\nu \in (0,\frac 13]$, and a compact subset $K' \subset R$.
There exist positive constants $\e,E_2,C_2$ -- depending on $K'$ and $\nu$ -- and a function
$$U_2: [E_2,\infty) \times K'_\e \to \C, \; (E,\alpha,\beta) \mapsto U_2(E,\alpha,\beta),$$
such that 
\begin{itemize}
\item For all $E \in [E_2,\infty) $, $U_2(E,\cdot,\cdot):  K'_\e \to \C$ is analytic.
 \item If $(E,\alpha,\beta) \in [E_2,\infty) \times K'_\e$ and $n:= E \nu$ is integer,
 the singularity $\la=0$ of equation \eqref{eq:scaled} is apparent
 if and only if $U_2(E,\alpha,\beta)=0$.
 \item $\forall(E,\alpha,\beta) \in [E_2,\infty) \times K'_\e$, the following estimate holds
 \begin{equation}\label{eq:monodromywkb}
 \left| U_2(E,\alpha,\beta)- 1-\exp{\left(E\oint_{\gamma_2}\sqrt{V(\la)}d\la+i\pi n\right)}\right|\leq \frac{C_2}{E} \exp\left(E
 \Re \oint_{\gamma_2}\sqrt{V(\la)}d\la\right) .
\end{equation}
\end{itemize}
\end{thm}
\begin{proof}
The proof follows the same steps of the proof of Theorem \ref{thm:wkbquantisation}. We therefore leave some details to the reader.
We consider solutions of \eqref{eq:scaled} on the cut-plane $D_{\varphi}$ with $\varphi=\frac\pi2$
(the cut coincides with the positive imaginary semi-axis for large $\la$).
 With such a choice of the domain,
the solutions $\chi_+,\psi_{1}^{R,L},\psi_{-1}$, subdominant at $0,\Sigma_{1}^{R,L},\Sigma_{-1}$ are uniquely defined up to a normalising
constant.

According to Proposition \ref{prop:wronskiantovalues}(2), if $n$ is integer and
   that $\chi_{+1},\psi_{-1}$ are linearly independent, $\la=0$ is an apparent singularity
   if and only if $w_{-1}\big( \frac{\psi_1^R}{\psi_{1}^L} \big)=w_{\uo}\big(\frac{\psi_1^R}{\psi_{1}^L}\big)$.

The function $U_2$, which has been introduced in Lemma \ref{lem:analiticity},
\begin{equation}\label{eq:U2}
 U_2=1-\frac{w_{\uo}\big(\frac{\psi_1^R}{\psi_{1}^L}\big)}{w_{-1}\big( \frac{\psi_1^R}{\psi_{1}^L} \big)}
\end{equation}
vanishes if and only if $w_{-1}\big( \frac{\psi_1^R}{\psi_{1}^L} \big)=w_{\uo}\big(\frac{\psi_1^R}{\psi_{1}^L}\big)$.

To prove the theorem, it is sufficient to show that we can find an $\e>0$ such that
   \begin{enumerate}
    \item If $E$ is large enough $\chi_{+1},\psi_{-1}$ are linearly independent for all $(\alpha,\beta) \in K'_{\e}$, so that 
    $U_2$ vanishes if and only if the singularity is apparent.
    \item $U_2$ is a well-defined analytic function on $K'_{\e}$.
    \item Estimate \eqref{eq:monodromywkb} holds.
   \end{enumerate}
   In order to prove these three statements we make use of the horizontal trajectories, which are classified in Lemmas \ref{lem:nonemptysets}
and \ref{lem:defsets}. These lemmas in particular states that if $\e$ is small enough
the sets of horizontal trajectories $A_{\uo,1}(\alpha,\beta),A^{R,L}(\alpha,\beta)$,
as depicted in Figure \ref{fig:antiS} are non empty.
Moreover, $$\sup_{\alpha,\beta \in K'_\e} \inf_{\gamma \in A(\alpha,\beta)} \rho_{\gamma} \leq C ,$$
where $A$ is one of the three sets of horizontal trajectories listed above, and $C$ a positive constant.

The proofs of (1) and (2) are carbon copies of the proofs of the analogous statements proven in Theorem \ref{thm:wkbquantisation}.
They are therefore
omitted.

(3)The proof of (3) is also analogous to the proof of the similar statement in Theorem \ref{thm:wkbquantisation},
but with a different choice of admissible curves.
The main idea is that we approximate $\psi_{1}^{R,L}$ 
by the WKB function $$\Psi(\la;\la')=\exp\lbrace \int_{\la'}^\la E \sqrt{V(\mu)}-\frac{V'(\mu)}{4V(\mu)} d\mu \rbrace,$$
along admissible curves which starts in $\Sigma_{1}^{R,L}$.
In fact, we choose admissible curves  $\alpha^{R,L}$ and $\beta^{R,L}$ as depicted in
Figure \ref{fig:y3}. We prove that this choice is possible.
We take two distinct horizontal trajectories $\widetilde{\alpha}^{L,R} \in \overrightarrow{A}_{1,\uo}$.
Since the trajectories $\widetilde{\alpha}_{R,L}$ do not intersect, but are asymptotic
as $t\to 0$ to positive imaginary axis,
we can define a cut $l$, asymptotic to the positive imaginary axis,
in such a way that $\alpha_{L,R}(]0,1[) \in D_{\frac{\pi}{2}}:= \C\setminus l$.
According to Lemma \ref{lem:pathdef}, we can deform $\widetilde{\alpha}_{L,R}$ to admissible curves
$\alpha_L,\alpha_R$, without increasing their error, such that
\begin{itemize}
 \item $\alpha_L(t)=\alpha_R(t)$ as $t \to 1$, i.e. in a neighbourhood of $0$;
 \item $\alpha_L(t)$ belongs to a closed subsector of $\Sigma_1^L$ as $t\to 0$;
 \item $\alpha_R(t)$ belongs to a closed subsector of $\Sigma_1^R$ as $t\to 0$;
\end{itemize}
Here we assume to have named $\widetilde{\alpha}_{L,R}$ so that $\widetilde{\alpha}_L$
lies on the left of the cut and $\widetilde{\alpha}_R$ lies on the right of the cut.

Similarly, the paths $\beta_{R,L}$ are appropriate deformations of horizontal trajectories belonging to the subsets $A^R,A^L$
(as defined in Lemma \ref{lem:defsets}) so that
\begin{itemize}
 \item $\beta_{L}(t)=\alpha_{L}(t)$ as $t \to 0$;
 \item $\beta_{R}(t)=\alpha_{R}(t)$ as $t \to 0$;
 \item $\beta_{L}(t)=\beta_{R}(t)$ as $t \to 1$.
\end{itemize}

\begin{figure}
	\centering
	\begin{tikzpicture}[scale=0.6]

	\tikzstyle{star}  = [circle, minimum width=4pt, fill, inner sep=0pt];
	
	\node[above]     at (-6,0) {$\lambda_1$};

	\node[above]     at (-3,0) {$\lambda_2$};

	\node[right]     at (0,0) {0};

	\node[above]     at (3,0) {$\lambda_3$};

	\node[above]     at (6,0) {$\lambda_4$};
	
	\draw[purple,thick,snake it] (0,0) -- (0,7);
	
	\draw[red,thick] (3,7) -- (3,6);
	\draw[red,thick] (-3,7) -- (-3,6);
	\draw[red,thick] (0,-7) -- (0,-6);
	\draw[red,thick] (0,0) -- (0,-1);

	\draw [red,thick,decoration={markings, mark=at position 0.32 with {\arrow{>}}},
	postaction={decorate}] plot [smooth,tension=0.8] coordinates { (3,6) (3.3,4.8)  (4.5,0)   (1,-4.5)   (0,-6) };
	
	\draw [red,thick,decoration={markings, mark=at position 0.49 with {\arrow{>}}},
	postaction={decorate}] plot [smooth,tension=0.8] coordinates { (3,6) (2.7,4.8)  (1.5,0)   (0.6,-1.5)   (0,-1) };
	
	\draw [red,thick,decoration={markings, mark=at position 0.32 with {\arrow{>}}},
	postaction={decorate}] plot [smooth,tension=0.8] coordinates { (-3,6) (-3.3,4.8)  (-4.5,0)   (-1,-4.5)   (0,-6) };

	\draw [red,thick,decoration={markings, mark=at position 0.49 with {\arrow{>}}},
	postaction={decorate}] plot [smooth,tension=0.8] coordinates { (-3,6) (-2.7,4.8)  (-1.5,0)   (-0.6,-1.5)   (0,-1) };

	\node[left]     at (-3,6.5) {$\lambda^L$};
	\node[right]     at (3,6.5) {$\lambda^R$};

	\node[red]     at (5.1,2) {$\beta_R$};
	\node[red]     at (-5.1,2) {$\beta_L$};
	\node[red]     at (1.35,1.9) {$\alpha_R$};
	\node[red]     at (-1.35,1.9) {$\alpha_L$};

	\node[red,star]     at (-3,6.5) {};

	\node[red,star]     at (3,6.5) {};
		
	\node[black,star]     at (-6,0){};
	
	\node[black,star]     at (-3,0) {};
	
	\node[black,star]     at (0,0) {};

	\node[black,star]     at (3,0) {};

	\node[black,star]     at (6,0) {};

	\end{tikzpicture} 
	\caption{The integration paths used in the proof of Theorem \ref{thm:wkbmonodromy}.} \label{fig:y3}
\end{figure}
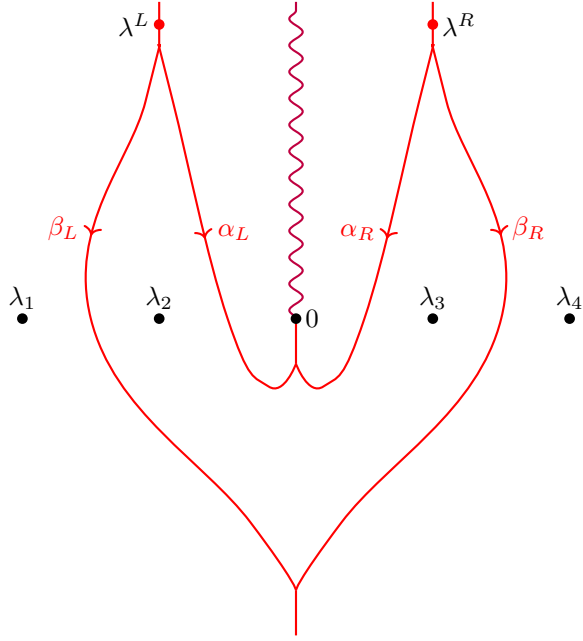

Having defined the four admissible curves, we lift them to the Riemann surface $\widehat{\Gamma}$ to enforce
the correct choice of the branch of $\sqrt{V}$. The paths $\widehat{\alpha}^{R,L},\widehat{\beta}^{R,L}$ are the lift of
$\alpha^{R,L},\beta^{R,L}$ defined by requiring that they belong to the lower sheet as $t\to 0$.
We are therefore in the situation depicted in Figure \ref{fig:y4}.
We notice that $\widehat{\alpha}_{R}-\widehat{\alpha}_{L}+\widehat{\beta}_L-
\widehat{\beta}_R$ makes a closed path and
\begin{equation}\label{eq:homologygamma2}
 \widehat{\alpha}_{R}-\widehat{\alpha}_{L}+\widehat{\beta}_L-
\widehat{\beta}_R= \gamma_2-\sigma \in H_1(\widehat{\Gamma}_{\alpha,\beta}\setminus\{ \la=0 \},\bb{Z})
\end{equation}
where $\gamma_2$ is the basis element of $H_1(\widehat{\Gamma}_{\alpha,\beta},\bb{Z})$, defined in Figure \ref{fig:cutplane_and_cycles},
and $\sigma$ is a small
positively-oriented loop about the lift of $\la=0$ on the lower sheet of $\widehat{\Gamma}$; see Figure \ref{fig:y4}.

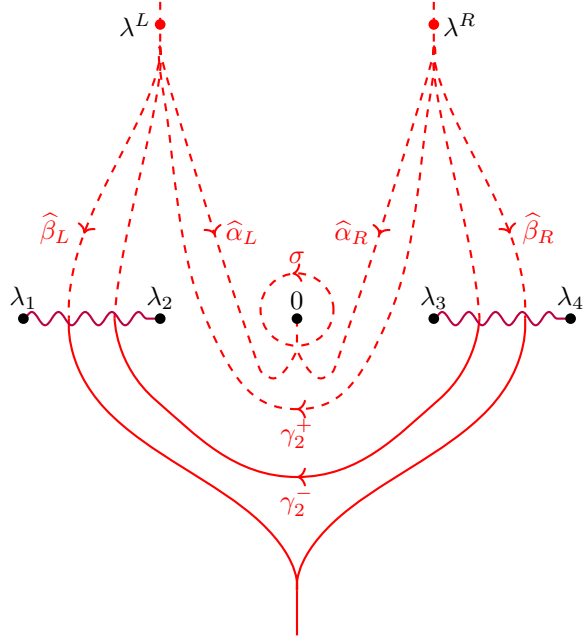
\begin{figure}
	\centering
	\begin{tikzpicture}[scale=0.6]

	\tikzstyle{star}  = [circle, minimum width=4pt, fill, inner sep=0pt];
	
	\node[above]     at (-6,0) {$\lambda_1$};

	\node[above]     at (-3,0) {$\lambda_2$};

	\node[above]     at (0,0) {0};

	\node[above]     at (3,0) {$\lambda_3$};

	\node[above]     at (6,0) {$\lambda_4$};

	\draw[red,thick,dashed] (3,7) -- (3,6);
	\draw[red,thick,dashed] (-3,7) -- (-3,6);
	\draw[red,thick] (0,-7) -- (0,-6);
	\draw[red,thick,dashed] (0,0) -- (0,-0.7);

	\draw [red,thick,dashed,decoration={markings, mark=at position 0.5 with {\arrow{>}}},
	postaction={decorate}] plot [smooth,tension=0.8] coordinates { (3,6) (2.7,4.8)  (1.2,0)   (0.6,-1.3)   (0,-0.7) };
	\draw [red,thick,dashed,decoration={markings, mark=at position 0.5 with {\arrow{>}}},
	postaction={decorate}] plot [smooth,tension=0.8] coordinates { (-3,6) (-2.7,4.8)  (-1.2,0)   (-0.6,-1.3)   (0,-0.7) };
	
	\draw [red,thick,dashed,decoration={markings, mark=at position 0.5 with {\arrow{>}}},
	postaction={decorate}] plot [smooth,tension=0.8] coordinates { (3,6) (2.9,4.8)  (2,0)   (0,-2)   (-2,0) (-2.9,4.8) (-3,6) };

	\draw [red,thick,dashed,decoration={markings, mark=at position 0.7 with {\arrow{>}}},
	postaction={decorate}] plot [smooth,tension=0.8] coordinates { (3,6) (3.3,4.8) (4.7,1.8) (5,0)  };
	\draw [red,thick,dashed,decoration={markings, mark=at position 0.7 with {\arrow{>}}},
	postaction={decorate}] plot [smooth,tension=0.8] coordinates { (-3,6) (-3.3,4.8) (-4.7,1.8) (-5,0)  };
	\draw [red,thick,dashed] plot [smooth,tension=0.8] coordinates { (3,6) (3.1,4.8) (3.8,1.5) (4,0)  };
	\draw [red,thick,dashed] plot [smooth,tension=0.8] coordinates { (-3,6) (-3.1,4.8) (-3.8,1.5) (-4,0)  };
	
	\draw [red,thick] plot [smooth,tension=0.8] coordinates { (5,0) (4.2,-2) (0.8,-4.5) (0,-6)  };
	\draw [red,thick] plot [smooth,tension=0.8] coordinates { (-5,0) (-4.2,-2) (-0.8,-4.5) (0,-6)  };
	
	\draw [red,thick,decoration={markings, mark=at position 0.5 with {\arrow{>}}},
	postaction={decorate}] plot [smooth,tension=0.8] coordinates { (4,0) (3.1,-1.8) (0,-3.5) (-3.1,-1.8) (-4,0)  };

	\draw[red,thick,dashed,decoration={markings, mark=at position 0.27 with {\arrow{>}}},
	postaction={decorate}] (0,0.2) ellipse (0.8cm and 0.8cm);
	
	\node[left]     at (-3,6.5) {$\lambda^L$};
	\node[right]     at (3,6.5) {$\lambda^R$};
	
	\node[red,below]     at (0,-3.5) {$\gamma_2^-$};
	\node[red,below]     at (0,-2) {$\gamma_2^+$};
	
	\node[red]     at (5.3,2) {$\widehat{\beta}_R$};
	\node[red]     at (-5.3,2) {$\widehat{\beta}_L$};
	\node[red]     at (1.2,1.9) {$\widehat{\alpha}_R$};
	\node[red]     at (-1.2,1.9) {$\widehat{\alpha}_L$};
	\node[red]     at (0,1.35) {$\sigma$};

	\draw[purple,thick,snake it] (-6,0) -- (-3,0);
	\draw[purple,thick,snake it] (3,0) -- (6,0);

	\node[red,star]     at (-3,6.5) {};

	\node[red,star]     at (3,6.5) {};
	
	\node[black,star]     at (-6,0){};
	
	\node[black,star]     at (-3,0) {};
	
	\node[black,star]     at (0,0) {};

	\node[black,star]     at (3,0) {};

	\node[black,star]     at (6,0) {};

	\end{tikzpicture} 
	\caption{The integration paths used in the proof of Theorem \ref{thm:wkbmonodromy}, lifted on the elliptic curve.
	The elliptic curve is realised as a two sheeted-cover of the
cut plane, with branch-cuts depicted in wiggly purple. A solid line represents a
curve on the upper sheet; A dashed line represents a curve on the lower sheet.} \label{fig:y4}
\end{figure}

We can now repeat the same steps of the proof of Theorem \ref{thm:wkbquantisation}.
We choose a $\la^R \in \widehat{\alpha}^R\cap \widehat{\beta}^R$ and 
$\la^L \in \widehat{\alpha}^L\cap \widehat{\beta}^L$, and
we obtain the following estimates for $w_{\uo}\big( \frac{\psi_1^R}{\psi_1^L} \big)$ and
$w_{-1}\big( \frac{\psi_1^R}{\psi_1^L} \big) $.

There exist a $C_+,E_+>0$ such that for all $(\alpha,\beta) \in K'_\e$,
\begin{equation}\label{eq:wuoestimate}
 w_{\uo}\big( \frac{\psi_1^R}{\psi_1^L} \big)=\exp\left\lbrace\int_{\la^R,\gamma_2^+}^{\la^L}
 E \sqrt{V(\mu)}-\frac{V'(\mu)}{4V(\mu)}d\mu \right \rbrace \left( 1+ \e_+(\alpha,\beta) \right)
\end{equation}
where $\gamma_2^+=\widehat{\alpha}^R-\widehat{\alpha}^L$, as depicted in Figure \ref{fig:y4}, and  $\e_+(\alpha,\beta)$  satisfies
the estimate
\begin{equation}\label{eq:estimatee+app}
 |\e_+(\alpha,\beta)|\leq \frac{C_+}{E}\,, \quad \forall E\geq E_+.
\end{equation}
There exists $C_-,E_->0$ such that for all $(\alpha,\beta)\in K'_\e$
\begin{equation}\label{eq:w-1estimateapp}
 w_{-1}\big( \frac{\psi_1^R}{\psi_1^L} \big)= \exp\left\lbrace 
 \int_{\la^L,\gamma_1^-}^{\la^R}E \sqrt{V(\mu)}-\frac{V'(\mu)}{4V(\mu)}d\mu \right\rbrace \left(1 +\e_-(\alpha,\beta)\right) 
\end{equation}
where $\gamma_2^-=\widehat{\beta}^R-\widehat{\beta}^L$, as depicted in Figure \ref{fig:y4}, and
\begin{equation}\label{eq:estimatee-app}
 |\e_-(\alpha,\beta)|\leq \frac{C_-}{E}\,, \quad \forall E\geq E_-.
\end{equation}
 Finally we show that (\ref{eq:wuoestimate},\ref{eq:estimatee+app},\ref{eq:w-1estimateapp},\ref{eq:estimatee-app})
 imply \eqref{eq:monodromywkb}.
 Since $$\gamma_2^+-\gamma_2^-=\widehat{\alpha}_{R}-\widehat{\alpha}_{L}+\widehat{\beta}_L-
\widehat{\beta}_R, $$ from \eqref{eq:homologygamma2} we get that
$\gamma_2^+-\gamma_2^-=\gamma_2-\sigma$.

 From equations (\ref{eq:U2},\ref{eq:wuoestimate},\ref{eq:w-1estimateapp}) we obtain
 $$
 U_2-1=-\exp\left\lbrace 
 \int_{\gamma_2-\sigma}E\sqrt{V(\mu)}-\frac{V'(\mu)}{4V(\mu)}d\mu \right\rbrace \frac{1+\e_-(\alpha,\beta)}{1+\e_+(\alpha,\beta)}.
 $$
	By the residue theorem, $$\int_{\gamma_2-\sigma}\frac{V'(\mu)}{4V(\mu)}d\mu=\pi i \mbox{ and }
	E \int_{-\sigma}\sqrt{V(\mu)}d\mu=i \pi n . $$
	Inserting these identities in the former
	display, we obtain
	$$
	U_2-1 - \exp\left\lbrace 
	\int_{\gamma_2}E\sqrt{V(\mu)}d\mu +i \pi n  \right\rbrace=\exp\left\lbrace 
	\int_{\gamma_2}E\sqrt{V(\mu)}d\mu  +i \pi n\right\rbrace \frac{\e_-(\alpha,\beta)-\e_+(\alpha,\beta)}{1+\e_+(\alpha,\beta)}.
	$$
	Using
	(\ref{eq:estimatee+app},\ref{eq:estimatee-app}), we obtain \eqref{eq:monodromywkb}.
\end{proof}

\section{Bulk Asymptotic}\label{section:bulkasymptotic}
In the present section we prove Theorem \ref{thm:bulkbeta}. In other words
we use the WKB analysis of equation \eqref{eq:scaled}
to provide an asymptotic solution for the inverse monodromy problem characterising the roots
of the generalised Hermite polynomials. 

In order to state and prove our results, for sake of simplicity we choose $p\geq q$, $p,q$ either equal
or co-prime, and fix a ratio $\frac{m}{n}=\frac{p}{q}$.
Hence, in the whole section, the numbers $m,n$ take values in the sequences $m=xq,n=xp, x \in \bb{N}$. 
Correspondingly $\nu=\frac{p}{2q+p}\in(0,\frac13]$ is fixed and the large parameter $E$ belongs to the sequence $(2q+p)x$.
However, all results are essentially unchanged if we let $\nu$ vary on $[\nu_0,\frac13]$, for some fixed $\nu_0 >0$.

For benefit of the reader we copy here the statement of Theorem \ref{thm:bulkbeta}.
The reader
should recall, from Section \ref{sec:results}, the following definitions:
\begin{itemize}
 \item The domains
 $K$, $K_a$ (i.e. the projection
of $K$ on the $\alpha$-plane, $K_a=\Pi_a(K)$); see Definition \ref{defi:RK}.
\item The function $\mc{S}$ and
the rectangle $Q=\left[-\tfrac{1}{2}(1-\nu)\pi,\tfrac{1}{2}(1-\nu)\pi\right]\times \left[-\nu\pi,\nu\pi\right]$,
which is the image of
$K$ under $\mc{S}$; see equation \eqref{def:S} and Theorem \ref{thm:mappingS}.
\item The finite sequence of integers $I_m=\lbrace -m+1,-m+3,\dots,m-1\rbrace$ and 
the point $(\alpha_{j,k},\beta_{j,k}) \in K$, with $j \in I_m, k \in I_n$, 
which is the unique solution of $S(\alpha,\beta)=(\frac{\pi j}{E},\frac{\pi k}{E})$; see Definition \ref{def:Im}.
\item For each filling fraction $\sigma \in (0,1)$, the subsets $Q^{\sigma} \subset Q$, which for large $E$ converges to $\sigma\cdot Q$,
$K^{\sigma}=\mathcal{S}^{-1}(Q^{\sigma}) \subset K$, and $K_a^{\sigma}=\Pi_a(K^{\sigma})\subset K_a$; see Definition \ref{def:Qsigma}.
\end{itemize}
Theorem \ref{thm:bulkbeta} is then stated as follows.
\begin{thm*}
 Fix $\sigma \in (0,1)$. Then there exists an $R_{\sigma}>0$ such that for $E$ large enough the following hold true:
 \begin{enumerate}
  \item In each
 ball of centre $(\alpha_{j,k},\beta_{j,k})$, $(j,k) \in I_m^{\sigma} \times I_n^{\sigma}$ and radius $R_{\sigma}E^{-2}$ there exists
 a unique point $(\alpha,\beta)$ such that the anharmonic oscillator \eqref{eq:scaled} satisfies
 the inverse monodromy problem characterising the roots of generalised Hermite polynomials.
 \item In the $\epsilon$ neighbourhood of $K^{\sigma}$ with radius $R_{\sigma}E^{-2}$, there are exactly
  $\lfloor \sigma m \rfloor \times \lfloor \sigma n \rfloor$ points $(\alpha,\beta)$ such that the anharmonic oscillator
 (\ref{eq:scaled}) satisfies the inverse monodromy problem characterising the roots of generalised Hermite polynomials.
 \end{enumerate}
\end{thm*}

In proving Theorem \ref{thm:bulkbeta}, we will use the multidimensional Rouch\'e theorem, which we state below.
\begin{thm}[\cite{aizenberg}]
Let $D,E$ be bounded domains in $\bb{C}^n $, $ \overline{D} \subset E$ , and let $f(z) , g(z)$ be holomorphic maps $ E \to \bb{C}^n$ such that
$\av{g(z)} < \av{f(z)}, \;\forall z \in \partial D$, 
then $u(z)=f(z)+g(z)$ and $f(z)$ have the same number (counting with multiplicities)
of zeroes inside $D$. Here $\av{\cdot}$ is any
norm on $\bb{C}^n$.
\end{thm}
\begin{proof}[Proof of Theorem \ref{thm:bulkbeta}] The proof is essentially an improved version of the proof of the analogous
result for poles of the Tritronqu\'ee solution of Painleve I \cite{piwkb2}.

For sake of definiteness, we use on $\bb{C}^2$ the euclidean norm and we denote it by $\| \,\|$. 
We denote by $\widetilde{K}^{\sigma}$ the epsilon neighbourhood of $K^{\sigma}$ where the radius epsilon is $R_{\sigma}E^{-2}$,
for some $R_{\sigma}>0$ to be determined.
By Theorem \ref{thm:wkbquantisation} and
Theorem \ref{thm:wkbmonodromy}, there exists an $E_{\sigma}>0$, such that
the solutions of the inverse monodromy problem, restricted to
$[E_{\sigma},\infty) \times \widetilde{K}^\sigma$, coincide with the simultaneous zeroes of the functions $U_1,U_2$
which satisfy the following estimate:
there exists a $C_\sigma$
such that for all $(E,\alpha,\beta) \in [E_{\sigma},\infty) \times \widetilde{K}^\sigma$,
\begin{align}\label{eq:quantproof}
 \left| U_1(E,\alpha,\beta) - 1 - \exp\left\{E\oint_{\gamma_1}\sqrt{V(\la)}d\la\right\}\right|&\leq \frac{C_{\sigma}}{E},  \\
 \label{eq:monodproof}
 \left| U_2(E,\alpha,\beta)-1 - \exp\left\{E\oint_{\gamma_2}\sqrt{V(\la)}d\la-i\pi n\right\}\right|&\leq \frac{C_{\sigma}}{E} .
\end{align}
The above inequalities stem from (\ref{eq:quantisationwkb},\ref{eq:monodromywkb}) and the fact that
$E \Re \oint_{\gamma_{1,2}}\sqrt{V(\la)}d\la =O(E^{-1})$ on $[E_{\sigma},\infty) \times \widetilde{K}^{\sigma}$.

For fixed $E \in [E_{\sigma},\infty)$, we define the function
$$u:\widetilde{K}^{\sigma} \to \bb{C}^2\, , \; (\alpha,\beta) \mapsto \left( U_1(E,\alpha,\beta),U_2(E,\alpha,\beta)\right) ,$$
whose zeroes characterise the solutions of the inverse monodromy problem for equation \ref{eq:scaled},
when $(\alpha,\beta)$ are restricted to $\widetilde{K}^{\sigma}$.

Similarly,for fixed $E \in [E_{\sigma},\infty)$, we define the function 
\begin{align*}
f:\widetilde{K}^{\sigma} \to \bb{C}^2 , \, 
f & = (\exp\big\{E\oint_{\gamma_1}\sqrt{V(\la)}d\la\big\}+1,
\exp\big\{E\oint_{\gamma_2}\sqrt{V(\la)}d\la-i\pi n\big\}+1) \\
& = \left(\exp\big\{i E s_1(\alpha,\beta)- i \pi m\big\}+1,\exp\big\{i E s_2(\alpha,\beta)-i\pi n\big\}+1\right),
\end{align*}
whose zeros coincide with
the points $(\alpha_{j,k},\beta_{j,k}), (j,k) \in I_m^{\sigma} \times I_n^{\sigma}$.

In order to prove the thesis by means of the Rouch\'e theorem,
it is sufficient to find an $R_{\sigma}$ such that
\begin{enumerate}
 \item $\forall (j,k) \in I_m^{\sigma} \times I_n^{\sigma}$,
 $\|f\|> \sqrt{2} E^{-1}C_{\sigma}$ on $\Sigma^{(j,k)}_{R_{\sigma}}$, where $\Sigma^{(j,k)}_{R_{\sigma}}$ is the sphere of centre
 $(\alpha_{j,k},\beta_{j,k})$, $j \in I_m^{\sigma}, k \in I_n^{\sigma}$ and
radius $R_{\sigma}E^{-2}$;
 \item $\|f\|> \sqrt{2} E^{-1}C_{\sigma}$ on the boundary of $\widetilde{K}^{\sigma}$.
\end{enumerate}
We do this by proving the corresponding estimates on the images of these sets under $\mathcal{S}$.

We choose the co-ordinates $(\hat{x},\hat{y})$ on $Q^{\sigma}$ and we embed $Q^{\sigma}$ in $\bb{C}^2$ by allowing $\hat{x},\hat{y}$
to be complex numbers.
We define the function $\hat{f}=(e^{i E \hat{x} - i\pi m},e^{iE \hat{y}-i \pi n})$ and
notice that $\hat{f}=f \circ \mc{S}$.

We begin by proving $(1)$. The following statement follows from a simple computation:
if $Y>0$ is big enough then for each $Y'>Y$ and each $(j,k)\in I_m \times I_n$, on the
spherical shell $\widehat{B}(j,k)_{Y,Y'} \subset \bb{C}^2$ of centre $(\frac{\pi j}{E},\frac{\pi k}{E})$, with
internal
radius $YE^{-2}$ and external radius $Y'E^{-2}$, the estimate $\|\hat{f}\|>\sqrt{2} E^{-1}C_{\sigma}$ is satisfied for $E$ big enough.
Let us then choose a $Y>0$ for which the above property holds and consider $S^{-1}(\widehat{B}(j,k)_{Y,Y'})$ for $(j,k) \in I_m^{\sigma} \times
I_n^{\sigma}$ and some $Y'>Y$. Denoting $J(j,k)$ the Jacobian of $\mc{S}$ at $(\alpha_{j,k},\beta_{j,k})$ and
$A(j,k)$ its inverse (recall that $S$ restricted to
$K^{\sigma}$ is a diffeomorphism),
we have that -up to a negligible $O(E^{-4})$
contribution- the counter-image of a sphere of radius $E^{-2} Y$ is just the ellipsoid with semi-axis $\sqrt{\la_i(j,k)} Y E^{-2},i=1,2$  where
$\la_{1,2}(j,k)>0$ are eigenvalues of the matrix $A^{\dagger}(j,k)A(j,k)$.
Since $K^{\sigma}$ is a compact subset of the open domain $R$, on wich $\mc{S}$ is a diffeomorphism,
the eigenvalues $\la_{1,2}$ are uniformly bounded. It means that there is a $R_{\sigma}$ and a $Y'>Y$ such that
$ \Sigma^{(j,k)}_{R_{\sigma}}\subset \mc{S}^{-1}(\widehat{B}(j,k)_{Y,Y'})$
for each $(j,k) \in I_{m}^{\sigma} \times I_n^{\sigma}$. Hence $\|f\|>\sqrt{2} E^{-1}C_{\sigma}$ for each of such
$\Sigma^{(j,k)}_{R_{\sigma}}$, and thus
the thesis is proven.

Part $(2)$ can be proven similarly. For $Y,E$ big enough, in the boundary of the $\e$ neighbourhood of $Q^{\sigma}$ with radius $E^{-2}Y$
the estimate $\|\hat{f}\|>\sqrt{2} E^{-1}C_{\sigma}$ is satisfied. By the same reasoning as above, we can find a $R_{\sigma}$
such that the estimate $\|f\|>\sqrt{2} E^{-1}C_{\sigma}$ is satisfied on the boundary of $\widetilde{K}^{\sigma}$. 
\end{proof}

\begin{rem}
Theorems \ref{thm:bulkbeta},\ref{thm:bulkintro}, and Corollary \ref{cor:smallalphaintro} hold true
unchanged if we let $\nu$ vary on the interval $[\nu_0,\frac13]$,
for any fixed $0<\nu_0<\frac13$. By this we mean that the constant $R_{\sigma}$ in the theorems above can be chosen to be independent on $\nu$.
\end{rem}

%
%

\appendix
\section{Elliptic Integrals}
\label{append:computation}
In this appendix we collect a number of explicit formulae for the elliptic integrals under consideration in the paper, which can be derived using standard elliptic function theory. The formulae are explicitly given in terms of the zeros of $V(\lambda;\alpha,\beta)$.
 For $(\alpha,\beta)$ close to $(0,0)$ the zeros $\lambda_k=\lambda_k(\alpha,\beta)$ of $V(\lambda;\alpha,\beta)$ do not coalesce and are analytic in $(\alpha,\beta)$. They are determined by equations \eqref{eq:constraint}, up to permutation, and we fix them unambiguously by the initial conditions \eqref{eq:turning} at $(\alpha,\beta)=(0,0)$.

Firstly, we have the following explicit formulas for $s_1$ and $s_2$,
\begin{align*}
s_1&=+\frac{2 i}{\sqrt{(\lambda_4-\lambda_1)(\lambda_3-\lambda_2)}}F(\lambda_1,\lambda_2,\lambda_3,\lambda_4)+\frac{1}{2} i\pi(1-\nu),\\
s_2&=-\frac{2 }{\sqrt{(\lambda_4-\lambda_3)(\lambda_2-\lambda_1)}}F(\lambda_4,\lambda_1,\lambda_2,\lambda_3)+i\pi\nu,
\end{align*}
with
\begin{align*}
F(\lambda_1,\lambda_2,\lambda_3,\lambda_4)=&-\tfrac{1}{4}(\lambda_4-\lambda_2)(\lambda_3-\lambda_2)(3\lambda_1-\lambda_2+\lambda_3+\lambda_4)\mathcal{K}(m)\\
&+\tfrac{1}{4}(\lambda_4-\lambda_1)(\lambda_3-\lambda_2)(\lambda_1+\lambda_2+\lambda_3+\lambda_4)\mathcal{E}(m)\\
&+(\lambda_4-\lambda_2)\Pi(n_1,m)+2\lambda_1\lambda_3(\lambda_4-\lambda_2) \Pi(n_2,m),
\end{align*}
where $\mathcal{K}(m)$, $\mathcal{E}(m)$ and $\Pi(n,m)$ denote the standard complete elliptic integrals of the respective first, second and third kind with parameter $m=k^2$ and elliptic characteristic $n$, in the above formula equal to
\begin{equation*}
m=\frac{(\lambda_2-\lambda_1)(\lambda_4-\lambda_3)}{(\lambda_3-\lambda_1)(\lambda_4-\lambda_2)},\quad n_1=-\frac{\lambda_4-\lambda_3}{\lambda_3-\lambda_2},\quad n_2=-\frac{(\lambda_4-\lambda_3)\lambda_2}{(\lambda_3-\lambda_2)\lambda_4}.
\end{equation*}
The formula for $s_1$ holds for $(\alpha,\beta)$ close to $(0,0)$ with all branches chosen principal, and is globally correct (on an open neighbourhood of $R$) via appropriate analytic continuation.
The formula for $s_2$ holds for $(\alpha,\beta)$ close to $(0,0)$, with all branches chosen principal except the one for $\Pi(n_2,m)$, namely
\begin{equation*}
\Pi(n_2,m)=\begin{cases}
\Pi^{(p)}(n_2,m) &\text{if $\Im{n_2}> 0$,}\\
\Pi^{(p)}(n_2,m)+\frac{i\pi}{\sqrt{(n_2-1)(1-m/n_2)}} &\text{if $\Im{n_2}\leq 0$,}\\
\end{cases}
\end{equation*}
where $\Pi^{(p)}(n_2,m)$ denotes the principal branch, so that $\Pi(n_2,m)$ is analytic in $n_2$ on an open neighbourhood of $(1,\infty)$ for $m\in \mathbb{C}\setminus [1,\infty)$.
It is globally correct (on an open neighbourhood of $R$) via appropriate analytic continuation.

Considering the Jacobian
\begin{equation}\label{eq:jacobianexplicit}
J_{(s_1,s_2)}(\alpha,\beta)=\begin{pmatrix}
\frac{\partial s_1}{\partial \alpha}& \frac{\partial s_1}{\partial \beta}\\
\frac{\partial s_2}{\partial \alpha}& \frac{\partial s_2}{\partial \beta}
\end{pmatrix},
\end{equation}
we have the following explicit expressions,
\begin{align*}
\frac{\partial s_1}{\partial \alpha}&=+\frac{2i}{\sqrt{(\lambda_4-\lambda_1)(\lambda_3-\lambda_2)}}\left[(\lambda_3-\lambda_2)(\lambda_4-\lambda_1)\mathcal{E}(m_{11})-(\lambda_1\lambda_2+\lambda_3\lambda_4)\mathcal{K}(m_{11})\right],\\
\frac{\partial s_2}{\partial \alpha}&=-
\frac{2}{\sqrt{(\lambda_4-\lambda_3)(\lambda_2-\lambda_1)}}\left[(\lambda_2-\lambda_1)(\lambda_4-\lambda_3)\mathcal{E}(m_{21})+(\lambda_2\lambda_3+\lambda_1\lambda_4)\mathcal{K}(m_{21})\right],\\
	\frac{\partial s_1}{\partial \beta}&=-\frac{2i}{\sqrt{(\lambda_3-\lambda_1)(\lambda_4-\lambda_2)}}\mathcal{K}(m_{12}),\\
	\frac{\partial s_2}{\partial \beta}&=+\frac{2}{\sqrt{(\lambda_3-\lambda_1)(\lambda_4-\lambda_2)}}\mathcal{K}(m_{22}),
\end{align*}
with parameters
\begin{equation*}
m_{12}=\frac{(\lambda_2-\lambda_1)(\lambda_4-\lambda_3)}{(\lambda_3-\lambda_1)(\lambda_4-\lambda_2)},\quad m_{22}+m_{12}=1,\quad m_{11}+m_{22}^{-1}=1,\quad m_{11}m_{21}=1,
\end{equation*}
which holds near $(\alpha,\beta)=(0,0)$ with all branches taken principally, and holds globally via appropriate analytic continuation.
In particular
\begin{equation}\label{eq:jacobiandet}
|J_{(s_1,s_2)}(\alpha,\beta)|\equiv 2\pi i
\end{equation}
due to Legendre's identity.


\end{document}